\documentclass[12pt,a4paper,reqno]{article} 

\pdfoutput=1

\usepackage{amssymb, amsmath, amsthm, a4} 
\usepackage{mathrsfs} 
\usepackage{graphicx}
\usepackage{enumitem}
\usepackage[titletoc,page]{appendix} 
\usepackage{setspace}

\usepackage[utf8]{inputenc}
\usepackage[english]{babel}
\usepackage{xpatch}
\usepackage[style=english]{csquotes} 
\usepackage[backend=biber,style=numeric]{biblatex}
\DeclarePrintbibliographyDefaults{heading=subbibliography} 

\numberwithin{equation}{section}

\theoremstyle{plain} 
\newtheorem{theorem}{Theorem}[section] 
\newtheorem{lemma}[theorem]{Lemma}  
\newtheorem{corollary}[theorem]{Corollary}  
\newtheorem{proposition}[theorem]{Proposition}
\newtheorem{conjecture}[theorem]{Conjecture}
\theoremstyle{definition} 
\newtheorem{definition}[theorem]{Definition}  
\theoremstyle{remark}
\newtheorem{remarks}[theorem]{Remarks}

\begin{document}

\title{On eigenvalues of the kernel  
$\frac{1}{2} + \lfloor \frac{1}{xy}\rfloor - \frac{1}{xy}$, II}

\author{Nigel Watt} 

\date{} 

\maketitle 

\begin{abstract} 
We study the eigenvalues $\lambda_1,\lambda_2,\lambda_3,\ldots$ (ordered by modulus) 
of the integral kernel $K\in L^2 ([0,1]\times [0,1])$ 
given by: $K(x,y) = \frac{1}{2} + \lfloor \frac{1}{x y}\rfloor - \frac{1}{x y}\,$ ($0<x,y\leq 1$). 
This kernel is of interest in connection 
with an identity of F. Mertens involving the 
M\"obius function. 
We establish that $\sum_{m=1}^{\infty} |\lambda_m|^{-1} = \infty$, 
and prove that $|\lambda_m| > m\log^{-3/2} m$ for all but finitely many $m\in{\mathbb N}$. 
The first of these results is an application of the theory of Hankel operators; 
the proof of the second result utilises a family of degenerate kernels 
$k_3,\allowbreak k_4,\allowbreak k_5,\ldots\,$ that are step-function approximations to $K$. 
Through separate computational work on eigenvalues of $k_N\,$ ($N=2^{21}$)  
we obtain numerical bounds, both upper and lower, for 
specific eigenvalues of $K$. Additional computational work, on eigenvalues of $k_N\,$ 
($N\in\{ 2^{10},\allowbreak 2^{11},\allowbreak \ldots ,2^{21}\}$), leads us to formulate a quite precise 
conjecture concerning where  on the real line the eigenvalues 
$\lambda_1,\lambda_2,\ldots ,\lambda_{767}$ are located:  
we discuss how this conjecture could (if it is correct) be viewed as supportive of 
certain interesting general conjectures concerning the eigenvalues of $K$. 
\end{abstract} 

\tableofcontents 

\begin{refsection}[refs0.bib]

\section{Introduction} 

This paper is a sequel to our recent work \cite{Wa2019}, in which it was shown that the integral kernel 
$K : [0,1]\times [0,1] \rightarrow {\mathbb R}$ given by 
\begin{equation}\label{DefK} 
K(x,y) := 
\begin{cases} 
\frac{1}{2} + \left\lfloor \frac{1}{xy}\right\rfloor - \frac{1}{xy} & \text{if $0<x,y\leq 1$} , \\ 
0 & \text{otherwise} , 
\end{cases}
\end{equation}
has infinitely many positive eigenvalues and infinitely many negative eigenvalues 
(note that, since $K$ is real and symmetric, its eigenvalues are non-zero real numbers). 
Its purpose is to add to what is known about the eigenvalues of $K$, and also 
to present empirical (or `experimental') evidence in support of several conjectures about them. 
\par 
Since $K$ is square-integrable, each eigenvalue $\lambda$ of $K$ 
has a finite `index', $\iota_K (\lambda)$ (say), which is the dimension of 
the complex vector space of all $\phi\in L^2 [0,1]$ satisfying $\phi(x)=\lambda\int_0^1 K(x,y)\phi(y)dy$ almost everywhere in $[0,1]$.  
We follow \cite[Section~3.8]{Tr1957}  in listing the eigenvalues of $K$ in a sequence 
$\lambda_1,\lambda_2,\ldots\ $ in which each eigenvalue $\lambda$ occurs 
exactly $\iota_K (\lambda)$ times, while the absolute value of $\lambda_m$ increases with $m$, so that one has 
\begin{equation}\label{lambdaOrder} 
0<\left| \lambda_m\right| \leq \left| \lambda_{m+1}\right| \quad \text{for $m\in{\mathbb N}$} .   
\end{equation} 
The sequence $|\lambda_1|,|\lambda_2|,\ldots\ $ is then independent of our 
choice of $\lambda_1,\lambda_2,\ldots\ $ (which is uniquely determined if and only if  
$\lambda_m\neq -\lambda_n$ for all $m,n\in{\mathbb N}$).  
\par 
The set  $\{ 0\}\cup\{ 1/\lambda_n : n\in{\mathbb N}\}$ is, of course, the spectrum of the operator 
\begin{equation}\label{K_as_HSoperator} 
f(x)\mapsto \int_0^1 K(x,y)f(y)dy\qquad \text{($f\in L^2 [0,1]$)} .   
\end{equation} 
This is a Hilbert-Schmidt integral operator. Its Hilbert-Schmidt norm is equal to the relevant $L^2$-norm of $K$,  which is  
\begin{equation*}
\| K\| := \left(\int_0^1 \int_0^1 |K(x,y)|^2 dx dy\right)^{\!\!1/2} < \infty\;.  
\end{equation*} 
\par 
Our interest in the kernel $K$ is connected with an elementary identity 
involving the M\"obius function  $\mu(n)$ that was proved in 1897 by F. Mertens \cite[Section~3]{Me1897}:  
\begin{equation*}
\sum_{k\leq n} \mu(k) = 2\sum_{k\leq\sqrt{n}} \mu(k) - \,\sum\sum_{\!\!\!\!\!\!\!\!\!\!\!\!\!\!\!\!{r,s\leq\sqrt{n}}} \mu(r)\mu(s) \left\lfloor\frac{n}{rs}\right\rfloor 
\quad\ \text{($n\in{\mathbb N}$)} .   
\end{equation*}
See the penultimate paragraph of \cite[Section~1]{Wa2019} for some discussion of how this result of Mertens motivates us: 
our unpublished result \cite[Theorem~3]{Wa2018} is also relevant. 

\subsection{Applications of the theory of Hankel operators}  
We observe in Section~4 that the operator \eqref{K_as_HSoperator} is unitarily equivalent to 
a specific Hankel operator on the space $L^2 (0,\infty)$,  
and is thereby also unitarily equivalent to Hankel operators on several other Hilbert spaces,  
including certain Hardy spaces and the space of `square summable' complex sequences. 
Thus, by applying relevant known results from the general theory of Hankel operators  
one can deduce several interesting facts about the eigenvalues of $K$. 
We find, in particular, that 
\begin{equation}\label{Intro-hankel1} 
\sum_{m=1}^{\infty} \frac{1}{|\lambda_m|} = \infty\;, 
\end{equation} 
that 
\begin{equation}\label{Intro-hankel2} 
\sum_{m=1}^{\infty} \frac{1}{\lambda_m^2} = \| K\|^2 = {\textstyle \frac14} - 2\zeta'(0)  + \zeta''(0)
\end{equation} 
(where $\zeta(s)$ is Riemann's zeta-function), and that if $\lambda$ is an eigenvalue of $K$ with $\iota_K (\lambda) \geq 2$ 
then $-\lambda$ is also an eigenvalue of $K$ and one has  
\begin{equation}\label{Intro-hankel3} 
\left| \iota_K(-\lambda) - \iota_K(\lambda)\right| \leq 1 \;. 
\end{equation} 
The first and last of these three results are, respectively, our Corollaries~\mbox{4.2} and~\mbox{4.5}.  
Regarding the proof of \eqref{Intro-hankel2}, see our Remarks~\mbox{6.11}~\mbox{(2)} and Section~\mbox{A.1}. 
\subsection{A general lower bound for $|\lambda_m|$} 
Sections~2 and~3 of this paper contain a proof of the following new result. 
\begin{theorem} 
There exists a (computable) $m_0\in{\mathbb N}$ such that one has 
\begin{equation}\label{MainBound} 
\left| \lambda_m\right| > m \log^{-\frac32} m 
\end{equation} 
for each integer $m\geq m_0$. 
\end{theorem}  
\par 
Theorem~\mbox{1.1} is  deduced from the following proposition, which we prove (in Section~3)  
by explicit construction of a kernel $k_N$ having the required properties. 
\begin{proposition} 
Let $N\in{\mathbb N}$ be sufficiently large. Then, 
for some symmetric bilinear combination 
\begin{equation}\label{Def-k_N}
k_N (x,y) := \sum_{i=1}^N \sum_{j=1}^N h_{i,j} \psi_i (x) \psi_j (y) 
\quad\text{($0\leq x,y\leq 1$)}  , 
\end{equation} 
in which $\psi_1,\ldots ,\psi_N$ are square-integrable and 
${\mathbb R}\ni h_{i,j}=h_{j,i}\,$ ($1\leq i,j\leq N$),  one has: 
\begin{equation}\label{Propn1Bound} 
\int_0^1 \int_0^1 \left( K(x,y) - k_N (x,y)\right)^2 dx dy < \frac{\log^3 N}{4N} \,. 
\end{equation}
\end{proposition}  
\subsection{Numerical bounds for $768$ eigenvalues of $K$} 
The specific kernel $k_N = k_N(x,y)$ constructed for our proof of Proposition~\mbox{1.2} 
also has a fundamental part to play in the quite separate work described in Sections~5, 6 and~7.  
This work, which is mostly of a numerical (or `experimental') character, 
has as its initial objective the determination of useful upper and lower bounds for some $768$ of the eigenvalues of $K$: specifically the 
first $384$ positive terms, and first $384$ negative terms, of the sequence $\{\lambda_m\}_{m=1}^{\infty}$. 
We achieve this in two stages, which are described in Sections~5 and~6, respectively. 
The first of these (Stage~1) involves the computation, for $N=2^{21}$, of upper and lower bounds for the relevant eigenvalues of
a kernel $k_N'$ that is a very close approximation to $k_N$, in the sense that $\| k_N' - k_N\|$ is very small. 
In Stage~2 we use the results of Stage~1 and sharp numerical upper bounds for $\| k_N' - k_N\|$ and $\| k_N' - K\|$  
to get our bounds on eigenvalues of $K$. The computations of Stage~2 are based on two theoretical results: Lemma~\mbox{6.7} 
(with its corollaries Lemma~\mbox{6.8} and \eqref{Weyl-4}) and Lemma~\mbox{6.12}. 
They include (ultimately) an application of the elementary principles described in Remarks~\mbox{6.14}~\mbox{(3)}. 
\par 
Our definition of $k_N$, in \eqref{epsilonDef}--\eqref{h_ijDef} and \eqref{Def-k_N}, 
involves a certain $N\times N$ real symmetric  matrix $H(N) = (h_{i,j})$ that (see Section~\mbox{5.2})  happens 
also to be a Hankel matrix: there are numbers $H_1,\ldots ,H_{2N-1}$ such that $h_{i,j} =  H_{i+j-1}$ for $1\leq i,j\leq N$. 
Our specification of $k_N'$ is similar (compare \eqref{k_N'Def} with \eqref{Def-k_N}), though  
the relevant Hankel matrix, $H'(N) = (h_{i,j}')$, differs slightly from $H(N)$. One has, in fact, 
$h_{i,j}' = H_{i+j-1}'\,$ ($1\leq i,j\leq N$), where $H_1',\ldots ,H_{2N-1}'\in{\mathbb R}$ are 
certain double precision approximations to the numbers  $H_1,\ldots ,H_{2N-1}$: the computation of these 
approximations is actually the first thing we do in Stage~1 (for how we do it, see Appendix~B). 
\par 
By a standard elementary argument (see \cite[Section~1.2]{Tr1957}, for example),  
the non-zero eigenvalues of $H'(N)\,$ (resp. $H(N)$) are the reciprocals of the eigenvalues of $k_N'\,$ (resp. $k_N$). 
Thus all that we need, in order to get good approximations to eigenvalues of $k_N'$, are some good approximations to  
the corresponding eigenvalues of $H'(N)$. This fact is key to our work in Stage~1, since  
good approximations to eigenvalues of real symmetric matrices such as $H'(N)$ can quite easily be obtained, 
with just the aid of a desktop computer and publicly available software. 
Also very helpful to us is the fact that $H'(N)$ is a Hankel matrix, since this allows calculations involving 
multiplication by $H'(N)$ to be greatly speeded up (through the use of a known algorithm that we detail in Appendix~C). 
Section~5 only summarises the Stage~1 computations: 
they are more fully described (and discussed) in Appendices~B, C, D and~E. 
\par 
In Section~\mbox{6.5} we give further details of the Stage~2 computations: the resulting bounds for eigenvalues of $K$  
are described and discussed in Remarks~\mbox{6.14}~\mbox{(1)--(4)}, at the end of the subsection. 
Some of these numerical bounds are shown in Tables~1 and~2 there. 
The principal achievements that we can point to, amongst these results, are as follows. 
For $m\leq 128$ we obtain, for the modulus of $\lambda_m$, upper and lower bounds whose quotient is less than $\frac65$.
For some small values of $m\in{\mathbb N}$ we get quite sharp upper and lower bounds for $\lambda_m$ itself, 
such as: 
\begin{align*} 
12.46<\lambda_1&<12.55\;,\quad 17.84<\lambda_5 <18.11\;,\quad 21.33<\lambda_7<21.79\;, \\
-17.65&<\lambda_4<-17.40\quad\text{and}\quad -19.76<\lambda_6<-19.41\;.   
\end{align*} 
In particular, with the help of \eqref{Intro-hankel3}, we can deduce that 
\begin{equation*} 
\iota_K(\lambda_m) = 1\quad\text{for the integers $m\in [1,14]\cup[25,26]$} . 
\end{equation*} 
There are disappointments, however. For example, although we can establish that $\lambda_2\lambda_3 < 0$, 
and that one has 
\begin{equation*} 
14.49>|\lambda_3|\geq |\lambda_2|>14.35\quad\text{and}\quad |\lambda_3|>14.42\;, 
\end{equation*}  
we nevertheless fail to determine which of the pair $\lambda_2,\lambda_3$ is the positive one: in particular, we 
fail to establish that $\lambda_3\neq -\lambda_2$.  Furthermore, when $m>268$ our upper and lower bounds 
for $|\lambda_m|$ differ by a factor that is greater than $2$, so that one of these bounds is rather weak (or both are). 
In fact the lower bounds that we obtain for each of $|\lambda_{269}|, |\lambda_{270}|, |\lambda_{271}|,\ldots$ have no  
intrinsic value, since these bounds are all equal to the lower bound obtained for $|\lambda_{268}|$.   
\par
We are confident that the above mentioned shortcomings are mainly due to a weakness in the method used    
to determine upper bounds for the moduli of eigenvalues of $H'(N)$: we point, specifically, to the use of $\| A^2\|^{1/2}$ 
as an upper bound for the spectral norm $\| A\|_2$ of the $N\times N$ real symmetric matrix $A$ defined in 
\eqref{Def-H''} and \eqref{Def-A(N)}. If resources of time and computing power had been available, then we 
might have tried some other (more effective) method of bounding this spectral norm. 
See Sections~\mbox{5.3}, \mbox{D.3} and~\mbox{D.4} for relevant details and discussion (the comments in Remarks~\mbox{D.3}~\mbox{(1)} and~\mbox{(3)}, in particular). 
\subsection{Probable bounds for eigenvalues} 
Our work in Section~7 assumes a certain hypothetical numerical upper bound for $\| A\|_2$. 
In fact it assumes such a bound not just in the case $N=2^{21}$, but also for each 
$N\in\{ 2^{10},2^{11},\ldots , 2^{20}\}$, with $A = A(N)$ being, in each case, a certain $N\times N$ real symmetric matrix 
that is defined in terms of the relevant $N\times N$ Hankel matrix $H'(N)$ and associated data (as detailed in Sections~\mbox{5.2} and~\mbox{5.3}). 
These hypothetical bounds are computed using a statistics-based algorithm (set out in Section~\mbox{D.5}). 
The first step of this algorithm involves the formation of an $N\times S$ matrix whose elements are chosen independently, and at random,  
from the set $\{ 1/\sqrt{N}, - 1/\sqrt{N}\}$.  
We deduce from Lemma~\mbox{D.5} that our usage of this algorithm should, most probably, yield valid bounds.   
Indeed, we find that if one were to repeatedly make independent applications of the algorithm, for a given $N$, and with $S\geq 30$ every time, 
then the rate at which invalid bounds would be produced would tend not to exceed one in a billion.  
For this reason we refer to our hypothetical numerical bounds for $\| A(N)\|_2\,$ ($N\in\{ 2^{10},2^{11},\ldots , 2^{21}\}$) 
as `probable upper bounds'. 
\par 
In resorting to the use of these probable upper bounds, we are motivated by 
a desire to transcend (however imperfectly) 
limits that  the practical considerations mentioned at the end of the previous subsection would otherwise 
force us to accept. 
\par 
Via Lemma~\mbox{D.4}, or its analogue for negative eigenvalues, 
our probable upper bound for $\| A(N)\|_2$ implies    
corresponding upper bounds  for the moduli of  eigenvalues of $H'(N)\,$ (henceforth also referred to as `probable upper bounds'). 
Whereas our use of the bound $\| A\|_2\leq \| A^2\|^{1/2}$ produces  
useful upper bounds for the moduli of just the first $100$ (or so) of the eigenvalues of $H'(2^{21})$, 
we find that for $N\in\{ 2^{10}, 2^{11}, \ldots , 2^{21}\}$ 
we get probable upper bounds,  sharp enough to be useful (for our purposes),  
for the $384$ least, and the $384$ greatest, of the eigenvalues of $H'(N)$ (these being, in fact, all of the eigenvalues of 
$H'(N)$ that we set out to estimate): Remarks~\mbox{D.7}~\mbox{(1)} and~\mbox{(3)} provide more detail on this. 
\par
Ultimately, since $\| k_N - k_N'\|$ is (in each case) very small, our lower bounds and probable upper bounds for the moduli of eigenvalues of $H'(N)$   
enable us to compute useful hypothetical (or `conjectural') estimates for the  $384$ least, and the $384$ greatest, of the eigenvalues of $H(N)$.  
That is, for $10\leq n\leq 21$, $N=2^n$ and $1\leq m\leq 384$, we determine a pair ${\mathcal X}^{+}_m (n), {\mathcal X}^{-}_m (n)$ of short real intervals 
that, if our probable upper bound for $\| A(N)\|_2$ is valid, 
must contain, respectively, the $m$-th greatest and $m$-th least of the eigenvalues of $H(N)$. It turns out 
that in all these cases $\inf{\mathcal X}^{+}_m (n) >0 > \sup{\mathcal X}^{-}_m (n)$. 
We base essentially the whole of our discussion in Section~7  on the premise (or `working hypothesis') 
that each of these intervals ${\mathcal X}^{\pm}_m (n)$ does indeed contain the relevant eigenvalue of $H(2^n)$. 
Subject to this working hypothesis, a new and appropriately revised application of Lemma~\mbox{6.12} (for $N=2^{21}$, as before) 
yields improved lower bounds for the moduli of eigenvalues of $K$: for some details of this see Remarks~\mbox{6.14}~\mbox{(5)}. 
However, what we do in Section~7 itself (and in the work discussed there) is to work at extracting even more from our working hypothesis by instead taking an empirical approach,  
in which the intervals ${\mathcal X}^{\pm}_m (n)\,$ ($10\leq n\leq 21$, $1\leq m\leq 384$, $\pm\in\{ +,-\}$) are treated as 
experimental data, from which we can extrapolate. This empirical approach is of course incapable of providing any definite conclusions, 
but it does allow us to come up with a number of interesting conjectures that our analysis of the data 
causes us to think are plausible. 
\subsection{Conjectures concerning $768$ eigenvalues of $K$}  
The discussion in Section~7 initially focusses on a study of the dependence on $n$ of the difference 
between the $m$-th greatest (or the $m$-th least) eigenvalue of $H(2^n)$ and the corresponding eigenvalue of $H(2^{n+1})$. 
There is (see our tentative conjecture \eqref{ForayHypothesis}) some indication that, for each fixed choice of $m$, these differences 
might be of size just $O(2^{-2n}n^A)$ as $n\rightarrow\infty$, 
where $A$ is a constant independent of $m$: 
this is also supported by observations noted in our Remarks~\mbox{C.3}~\mbox{(2)} and~\mbox{(3)}, though those 
observations are based on data different in kind to that considered in Section~7. 
However, the evidence we have for \eqref{ForayHypothesis} is a little weak, and this conjecture is also 
not explicit enough for our needs. The more explicit (but ultimately weaker) conjecture 
\eqref{SumCorr} is easier to justify; by adding it to our working hypothesis, 
we  arrive at \eqref{Extrapolation}, from which lower bounds for the relevant $767$ 
eigenvalues of $K$ follow, with the help of the data 
${\mathcal X}^{\pm}_m (20),{\mathcal X}^{\pm}_m (21)$ ($1\leq m\leq 384$, $\pm\in\{ +,-\}$). 
The missing ($768$-th) eigenvalue here (which is $384$-th positive term in the sequence $\lambda_1,\lambda_2,\ldots\ $)
does cause us more trouble: nevertheless, by an elaboration of the method used to get to \eqref{Extrapolation} 
(requiring a further extension of our working hypotheses), 
we arrive at a plausible lower bound for it. 
By Lemma~\mbox{6.7} the numbers $1/\inf{\mathcal X}^{+}_m (21),-1/\sup{\mathcal X}^{-}_m (21)$ ($1\leq m\leq 384$) 
are upper bounds for the moduli of the same $768$ eigenvalues of $K$: 
unlike the lower bounds that were just mentioned, these upper bounds are 
unconditionally valid (forming, as they do, one part of the set of unconditional results discussed in Section~6). 
\par
Combining these upper bounds and lower bounds, we obtain (as noted in \eqref{Work-Horse}, towards the end of Section~\mbox{7.1}) 
short real intervals, ${\mathcal B}^{+}_1,\ldots ,{\mathcal B}^{+}_{384}$ and 
${\mathcal B}^{-}_1,\ldots ,{\mathcal B}^{-}_{384}$ such that ${\mathcal B}_m^{+}$ (resp. ${\mathcal B}_m^{-}$) contains 
the $m$-th least (resp. $m$-th greatest) of the positive (resp. negative) eigenvalues of $K$, provided (of course) that 
our working hypotheses are correct. 
We find that for $1\leq m\leq 384$ both $(\inf{\mathcal B}_m^{+})/(\sup{\mathcal B}_m^{+})$ and 
$(\sup{\mathcal B}_m^{-})/(\inf{\mathcal B}_m^{-})$ are greater than $1-10^{-4}$, 
and that in many cases the interval ${\mathcal B}^{\pm}_m$ is even shorter than this would indicate: see \eqref{High-Precision}. 
The intervals ${\mathcal B}_m^{\pm}$ are, in particular, short enough for \eqref{Work-Horse} to serve   
as an interesting conjecture as to the locations of the relevant $768$ eigenvalues of $K\,$ (it specifies those locations 
with enough precision to reveal certain interesting patterns and trends, some of which are mentioned below). 
\par 
See Table~2 in Section~\mbox{6.5} for a selection of the numbers $\sup{\mathcal B}_m^{+}=1/\inf{\mathcal X}^{+}_m (21)$ 
and $|\inf{\mathcal B}_m^{-}|=-\inf{\mathcal B}_m^{-}=-1/\sup{\mathcal X}^{-}_m (21)$: these are, respectively, the numbers $1/{\mathcal L}^{+}_m$ and $1/{\mathcal L}^{-}_m$ tabulated there. 
In view of how very short the intervals ${\mathcal B}_m^{\pm}$  are, we decided it would not serve any useful purpose to   
tabulate the corresponding hypothetical lower bounds $\inf{\mathcal B}_m^{+}$ or $|\sup{\mathcal B}_m^{-}|=-\sup{\mathcal B}_m^{-}$.
\par 
Assuming that \eqref{Work-Horse} is correct, there are (see Section~\mbox{7.2}) a number of interesting findings concerning 
just the eigenvalues $\lambda_1,\ldots ,\lambda_{768}$ that can  
be extracted without much effort from the data ${\mathcal B}^{\pm}_m\,$ ($1\leq m\leq 384$, $\pm\in\{ +,-\}$). 
We mention just 3 of these conditional findings here. The first is that 
\begin{equation*} 
\iota_K(\lambda_m) = 1\quad\text{for $1\leq m\leq 767$} 
\end{equation*}
(this leads us to conjecture that all eigenvalues of $K$ are simple). 
The second conditional finding is (see the discussion around \eqref{Permed}--\eqref{TheFound}) 
a determination of ${\rm sgn}(\lambda_m):=\lambda_m / |\lambda_m|$ for each positive integer $m\leq 767$ with 
$\lfloor (m+1)/2\rfloor \not\in \{ 212,242,283,351,360,361,376,381\}$. 
One corollary of our being able to make such a determination is that   
\begin{equation*} 
\lambda_m\neq -\lambda_{\ell}\quad\text{when $\min\{ \ell,m\} < 423$} .
\end{equation*} 
The second conditional finding leads to the third, which (see \eqref{Shadowing}) is that 
\begin{equation}\label{SignSumBound-1}
\left| \sum_{\ell = 1}^m {\rm sgn}(\lambda_{\ell})\right|\leq 2\quad\text{for $1\leq m\leq 768$} . 
\end{equation}  
\subsection{General conjectures concerning the eigenvalues of $K$}  
In Sections~\mbox{7.3} and~\mbox{7.4} we formulate a number of wider conjectures, based on some analysis of the data  
${\mathcal B}^{\pm}_m\,$ ($1\leq m\leq 384$, $\pm\in\{ +,-\}$) and, of course, assuming the validity of \eqref{Work-Horse}.
In particular, we describe in Section~\mbox{7.3} how we are led to conjecture that  
for some $c_0\in (9.5, 10.2)$ one has 
\begin{equation*} 
\left|\lambda_m\right| \sim c_0 m \log^{-\frac32} m\quad\text{as $m\rightarrow\infty$} . 
\end{equation*} 
More speculatively (but still with some reason), 
our Remarks~\mbox{7.2} and~\mbox{7.4}~\mbox{(1)} embrace the conjectures that the above constant $c_0$ satisfies 
$9.85 < c_0 < 9.91$,  and that for any constant $\eta >0$ one has 
\begin{equation*} 
\left|\lambda_m\right| = c_0 m \log^{-\frac32} m + O\left( m^{\eta}\right) \quad\ \text{($m\in{\mathbb N}$)} . 
\end{equation*} 
\par 
Another of our conjectures, Conjecture~\mbox{7.3}, is simply a speculative extrapolation of our 
conditional finding \eqref{SignSumBound-1}: an equivalent 
conjecture is that, 
for any constant $\varepsilon >0$, one has 
\begin{equation}\label{Rails}
\left| \sum_{\ell = 1}^m {\rm sgn}(\lambda_{\ell})\right| + \iota_K\left( \lambda_m\right) = O\left( m^{\varepsilon}\right) \quad\  
\text{($m\in{\mathbb N}$)} . 
\end{equation}
We show in Section~\mbox{7.4} that Conjecture~\mbox{7.3}  implies (independently of any other hypotheses) 
that the series $\sum_{m=1}^{\infty}\lambda_m^{-1}$  is convergent: the proof is easy, although it does use Theorem~\mbox{1.1}. 
Some further (experimental) investigations relating to this series, in which \eqref{Work-Horse} is assumed,  
yield empirical evidence that appears quite supportive of the `naieve' conjecture that 
\begin{equation}\label{Mercer-INTRO}  
\sum_{m=1}^{\infty}\frac{1}{\lambda_m} = \int_0^1 K(x,x) dx\;. 
\end{equation} 
These investigations are detailed in Remarks~\mbox{7.4}~\mbox{(2)}, where there is also some discussion of 
plausible refinements of the conjecture \eqref{Mercer-INTRO} and supporting empirical evidence. 

\section*{Conventions and notation} 
\addcontentsline{toc}{section}{\protect\numberline{}Conventions and notation}
Subsections occur only within sections. Thus, for example, Section~\mbox{7.4} is the last subsection of Section~7, while Section~\mbox{A.1} is the first section of Appendix~A. There are no subsections in our  appendices. 
The sections, equations, remarks, tables and figures of each appendix are labelled with the 
appropriate capital letter (so `Table~\mbox{C-1b}', for example, means: `Table~\mbox{C-1b} of Appendix~C'). 
Each of Appendices A, C, D and E has its own independent list of references. 
\par 
A `double precision' number is a number that can be represented exactly in binary64 floating-point format. 
In Appendices~B, C and~D, we have ${\tt u} := 2^{-53}\,$ 
(the `machine epsilon', or `unit roundoff', for binary64 floating-point arithmetic). 
We use the same typewriter-like typeface for features of GNU Octave, and for objects defined in our Octave scripts and functions. 
Thus  `{\tt sort()}' denotes Octave's sort function, while `{\tt C}' is the name of a variable in a script. 
\par 
We sometimes write `$:=$' (or `$=:$'), instead of just `$=$' when defining a new constant or variable (i.e. 
`$A := B$' and `$B =: A$' both mean that $A$ is defined to be equal to $B$). The Bachmann-Landau (or `asymptotic')  notations `$A = O(B)$' and `$A\sim B$' have their usual meanings. 
We sometimes write `$A\ll B$' (resp. `$A\gg B$') 
to signify that $A=O(B)\,$  (resp. $B = O(A)$). 
The notation `$A\asymp B$' signifies that one has both 
$A\ll B$ and $B\ll A$. 
\par 
For $d\in{\mathbb N}$ and vectors ${\bf x}\in{\mathbb C}^d$, we have 
\begin{equation*} 
\| {\bf x}\|_p := 
\begin{cases}   (|x_1|^p+\ldots +|x_d|^p)^{1/p} & \text{when $1\leq p <\infty$} ,  \\  
\max\{|x_1|,\ldots ,|x_d|\} & \text{when $p=\infty$} . \end{cases} 
\end{equation*} 
We write  `$\| {\bf x}\|$' for the Euclidean norm $\| {\bf x}\|_2\,$ 
of a vector ${\bf x}\in{\mathbb C}^d$. 
For  $d,h\in{\mathbb N}$ and any $h\times d$ matrix $A = (a_{i,j})$, we have 
\begin{equation*} 
\| A\| := \sqrt{\sum_{i=1}^{h}\sum_{j=1}^d \left| a_{i,j}\right|^2}\qquad\text{(the Frobenius norm of $A$)}
\end{equation*} 
and $\| A\|_p := \max\{ \| A{\bf x}\|_p / \| {\bf x}\|_p : {\bf 0}\neq {\bf x}\in {\mathbb C}^d\}\,$ ($1\leq p\leq \infty$), so $\| A\|_2$ is the `spectral norm' of $A$.  
\par 
We sometimes write `$\{ x\}$' for the fractional part of $x\,$ (the number $x - \lfloor x\rfloor\in [0,1)$), though at other times 
`$\{ x\}$' may denote a singleton set (while `$\{ a_n\}_{n=0}^{\infty}$' would denote an infinite sequence). 
The `signum function', ${\rm sgn} : {\mathbb R}\backslash\{ 0\} \rightarrow \{ -1,1\}$, is given by ${\rm sgn}(x):=x/|x|$. 
We write `$I_d$' for the  $d\times d$ identity matrix, and write `$\zeta(s)$' for Riemann's zeta-function (evaluated at $s$). 
\par 
Additional notation is introduced later (closer to where it is needed). 

\section{Deducing Theorem~\mbox{1.1} from Proposition~\mbox{1.2}} 

We now establish the validity of Theorem~\mbox{1.1}, by showing that it is implied by Proposition~\mbox{1.2} 
(proof of which can be found in Section~3). 
\par 
Choose $m_0\in{\mathbb N}$ with $m_0 \geq 2$. Let $m\in{\mathbb N}$ satisfy $m\geq m_0$ and put $N=\lceil m/2\rceil$. 
Assuming that $m_0$ is chosen sufficiently large, we may apply Proposition~\mbox{1.2} to obtain the 
bound \eqref{Propn1Bound} for some kernel $k_N$ of the form \eqref{Def-k_N}. 
By \cite[Satz~VIII]{We1911}, the double integral on the left-hand side of \eqref{Propn1Bound} 
is greater than or equal to the sum 
\[ 
\Sigma'_N := \sum_{n=N+1}^{\infty} \frac{1}{\lambda_n^2} \,. 
\] 
Therefore, given that \eqref{lambdaOrder} implies 
\[ 
\Sigma'_N \geq \sum_{N<n\leq m} \frac{1}{\lambda_m^2} = \frac{m-N}{\lambda_m^2}  \,,
\] 
it follows from \eqref{Propn1Bound} that we have 
\[ 
\lambda_m^2 > \frac{4 (m-N)N}{\log^3 N} \geq \frac{m^2 - 1}{\log^3 \left(\frac{m+1}{2}\right)} \qquad\text{($N=\lceil m/2\rceil$)} . 
\] 
We have here $m\geq 2$ and 
$\log(\frac{2m}{m+1}) /\log(m) \geq (1-\frac{m+1}{2m}) / (m-1)=\frac{1}{2m}\geq m^{-2}$, 
and so 
$0 < \log(\frac{m+1}{2})  = \log(m) - \log(\frac{2m}{m+1}) \leq (1 - m^{-2})\log(m)$. 
Thus \eqref{MainBound} follows. \hfill $\square$ 

\section{Proving Proposition~\mbox{1.2}} 
\subsection{Defining $k_N$ and other preliminaries} 

Let $N\geq 3$ be an integer.
We construct a kernel $k_N : [0,1]\times[0,1]\rightarrow{\mathbb R}$ as follows. 
\par 
Firstly, we take $\varepsilon = \varepsilon(N)$ to be the unique real solution of 
the equation 
\begin{equation}\label{epsilonDef} 
2 \sinh (\varepsilon) = e^{-N\varepsilon} .
\end{equation} 
This implies $\varepsilon > 0$. Next, we put 
\begin{equation}\label{deltaDef} 
\delta = e^{-\varepsilon} , 
\end{equation}
and 
\begin{equation}\label{x_iDef} 
x_i = \delta^{i-1} \quad\text{for $i=1,\ldots ,N+1$.}
\end{equation}
Then, putting $x_{N+2} = 0$, we have: 
\begin{equation}\label{x_iOrder} 
1 = x_1 > x_2 > \text{ $\cdots $ } > x_{N+1} > x_{N+2}=0 .
\end{equation} 
When $(x,y)\in [0,1]^2$ is such that the sets 
$\{ x , y\}$ and $\{ x_1 , x_2, \ldots , x_{N+2}\}$ 
are not disjoint, we put $k_N(x,y) = 0$. 
In all other cases (where, still, $(x,y)\in [0,1]^2$), there will exist 
a unique ordered pair $(i,j)\in\{ 1,2,\ldots , N+1\}^2$ such that 
the open rectangle 
\begin{equation}\label{DefR_ij} 
{\mathcal R}_{i,j} := \left( x_{i+1} , x_i\right) \times \left( x_{j+1} , x_j\right) 
\end{equation} 
contains the point $(x,y)$. If, here, $\max\{ i,j\} > N$, then we put 
$k_N(x,y) =0$ (as before).
If instead $i,j\leq N$, then we put
$k_N(x,y) = \mu_{i,j}$, where 
\begin{equation}\label{mu_ijDef} 
\mu_{i,j} = \mu_{i,j}(N) := 
\frac{1}{\operatorname{vol}\left( {\mathcal R}_{i,j}\right)} \iint\limits_{{\mathcal R}_{i,j}} K(x,y) dx dy \,, 
\end{equation}
where, of course, 
$\operatorname{vol}( {\mathcal R}_{i,j}) = (x_i - x_{i+1})(x_j - x_{j+1})$.
Note that, since $K(x,y)$ is symmetric, the substitution of $x$ for $y$ (and vice versa) in \eqref{mu_ijDef} 
shows (given \eqref{DefR_ij}) that we have 
\begin{equation}\label{muSymmetry} 
\mu_{i,j} = \mu_{j,i} 
\quad\text{($i,j\in{\mathbb N}$, $i,j\leq N$).} 	 
\end{equation}
\par
The kernel $k_N$ (just defined) is of the form \eqref{Def-k_N},  with 
\begin{equation}\label{psi_iDef} 
\psi_i (x) := 
\begin{cases} 
\frac{1}{\sqrt{x_i - x_{i+1}}} & \text{if $x_{i+1}< x < x_i$} , \\ 0 & \text{otherwise} , 
\end{cases} 
\end{equation}
and 
\begin{equation}\label{h_ijDef} 
h_{i,j} := \mu_{i,j} \cdot \sqrt{(x_i - x_{i+1})(x_j - x_{j+1})}\in{\mathbb R} \,, 
\end{equation}
for $0\leq x\leq 1$ and $1\leq i,j\leq N$. 
Note that we have here $\int_0^1 \psi_i^2 (x) dx = 1<\infty\,$ ($1\leq i\leq N$), and that, 
by virtue of \eqref{muSymmetry}, the real numbers $h_{i,j}$ defined in \eqref{h_ijDef} 
are such that $h_{i,j}=h_{j,i}\,$ ($1\leq i,j\leq N$). Therefore all we need do, 
in order to complete our proof of Proposition~\mbox{1.2}, is show that 
the kernel $k_N$ that we have defined in the last paragraph is such that 
the inequality \eqref{Propn1Bound} will hold if $N$ is large enough (in absolute terms). 
The next three lemmas (in which $N$, $\varepsilon$ and $\delta$ are the numbers 
appearing in the definition of $k_N$) will enable us to achieve this. 

\begin{lemma} 
We have 
\[ 
\int_0^1 \int_0^1 \left( K(x,y) - k_N (x,y)\right)^2 dx dy 
< {\textstyle\frac12}\delta^N + \sum_{n=1}^{2N-1} \delta^{n-1} I_n\cdot\min\{ n,2N-n\} \,, 
\] 
with: 
\[ 
I_n := \min_{\kappa\in{\mathbb R}} \int_{\delta}^1 \int_{\delta}^1  
\left( K\left( \delta^{n-1} uv , 1\right) - \kappa\right)^2 du dv \,. 
\] 
\end{lemma} 

\begin{proof} 
Let ${\mathcal D} = [0,1]^2 \backslash [\delta^N , 1]^2 \subset {\mathbb R}^2$ and put 
$K_N (x,y) = K(x,y) - k_N (x,y)\,$ ($0\leq x,y\leq 1$). 
By the definitions of $K$ and $k_N$, \eqref{DefK} and \eqref{epsilonDef}--\eqref{mu_ijDef}, one obtains:  
\[
\int_0^1 \int_0^1 K_N^2 (x,y) dx dy  
=  \int_{x_{N+1}}^1 \int_{x_{N+1}}^1 K_N^2 (x,y) dx dy + \iint\limits_{\mathcal D} K^2 (x,y) dx dy \,, 
\] 
\begin{align*} 
\int_{x_{N+1}}^1 \int_{x_{N+1}}^1 K_N^2 (x,y) dx dy  
 &= \sum_{i=1}^N \sum_{j=1}^N \iint\limits_{{\mathcal R}_{i,j}} \left( K(x,y) - \mu_{i,j}\right)^2 dx dy  \\ 
 &= \sum_{i=1}^N \sum_{j=1}^N \min_{\kappa\in{\mathbb R}} \int_{\delta^i}^{\delta^{i-1}}\!\!\!\int_{\delta^j}^{\delta^{j-1}}  
\left( K(x,y) - \kappa\right)^2 dx dy  \\ 
 &= \sum_{i=1}^N \sum_{j=1}^N  \delta^{i+j-2} I_{i+j-1} \,,
\end{align*} 
and 
\begin{equation*} 
\iint\limits_{\mathcal D} K^2 (x,y) dx dy \leq {\textstyle\frac14} \operatorname{vol}({\mathcal D}) <  {\textstyle\frac12} \delta^N \,. 
\end{equation*} 
We deduce the required result by observing that, 
if $n$ is an integer, then one has 
$|\{ (i,j)\in {\mathbb N}^2 : i,j\leq N\ \text{and}\ i+j-1=n\}| = \max\{ 0,\min\{ n,2N-n\}\}$.   
\end{proof} 

\begin{lemma} 
Let $n\leq 2N-1$ be a positive integer. Let $I_n$ be as defined in Lemma~\mbox{3.1} and put  
\[ 
B_n := \left( \frac{1}{\delta^{n-1}} \,,\,  \frac{1}{\delta^{n+1}}\right]\subset{\mathbb R}\,. 
\]
Then 
\[ 
I_n \leq  \delta^{-2(n+1)} (1-\delta)^4 \quad\text{if $B_n \cap {\mathbb Z} = \emptyset$} . 
\] 
\end{lemma} 

\begin{proof} 
Suppose that there are no integers lying in the interval $B_n$. Then, by virtue of the definition \eqref{DefK}, the function 
$(u,v) \mapsto K(\delta^{n-1} uv , 1)$ is continuous on $[\delta , 1]\times [\delta , 1]\subset {\mathbb R}^2$, 
and one has, in particular, 
\[ 
K\left(\delta^{n-1} uv , 1\right) = c_n - \frac{1}{\delta^{n-1} uv} \qquad\text{($\delta\leq u,v\leq 1$)} , 
\] 
where $c_n := \frac12 + \lfloor \delta^{1-n} \rfloor$ is independent of the point $(u,v)$. 
Therefore 
\begin{align*} 
I_n &= \min_{\kappa\in{\mathbb R}} \int_{\delta}^1 \int_{\delta}^1  
\left( \left( c_n - \kappa\right) - \frac{1}{\delta^{n-1} uv}\right)^2 du dv \\ 
 &\leq \int_{\delta}^1 \int_{\delta}^1  
\left(\frac{1}{\delta^n} - \frac{1}{\delta^{n-1} uv}\right)^2 du dv  
=\frac{1}{\delta^{2n-2}}  \int_{\delta}^1 \int_{\delta}^1  
\left( \frac{1}{uv} - \frac{1}{\delta}\right)^2 du dv  \,. 
\end{align*} 
Since the last of the integrands here is bounded above by $(\delta^{-2} - \delta^{-1})^2 = \delta^{-4} (1-\delta)^2\,$ (for $\delta\leq u,v\leq 1$), 
the required result follows. 
\end{proof} 

\begin{remarks} 
In light of the definition \eqref{DefK} one has $-\frac12 < K(x,y)\leq \frac12\,$ ($0\leq x,y\leq 1$), 
and so it is trivially the case that $I_n$ (defined as in Lemma~\mbox{3.1}) must 
satisfy: 
\begin{equation}\label{I_n_TrivBound} 
I_n \leq {\textstyle\frac14} (1-\delta)^2 \,. 
\end{equation} 
By the definitions \eqref{epsilonDef} and \eqref{deltaDef}, we have  
\begin{equation}\label{deltaAlg} 
1-\delta^2 = \delta^{N+1}\, . 
\end{equation} 
Using this we find that Lemma~\mbox{3.2} conditionally yields a bound for $I_n$ that is $4(1+\delta)^{-2} \delta^{2(N-n)}$ times 
the trivial bound given in \eqref{I_n_TrivBound}. Thus, given that $1+\delta > 2\delta >0$, we may say 
that Lemma~\mbox{3.2} is `non-trivial' in respect of cases where $n<N$ and $B_n \cap {\mathbb Z}=\emptyset$. 
\end{remarks}

\begin{lemma} 
We have 
\[ 
N\varepsilon = \log(N) -\log\log(N) + O(1)\,. 
\]
\end{lemma} 

\begin{proof} 
Recall that we assume $N\geq 3$. Thus, since \eqref{epsilonDef} implies $\varepsilon>0$, 
we have $N\varepsilon = \vartheta \log N$ for some $\vartheta\in (0,\infty)$. 
The equation \eqref{epsilonDef} may therefore be reformulated as: 
\begin{equation}\label{BWa1} 
2\sinh\left( \frac{\vartheta\log N}{N}\right) = N^{-\vartheta}\,. 
\end{equation} 
\par 
Since $\sinh(x)>x$ for $x>0$, we deduce from \eqref{BWa1}  that 
$\vartheta N^{\vartheta - 1} \log N < \frac12$. It follows (since $\log N > 1 > 0$) that we must have $\vartheta < 1$. 
If we had also $\vartheta \leq \frac12$ then, since $N\geq 3$ and $\sinh(x) < x\cdot \cosh(x)$ for $x>0$, the equation \eqref{BWa1} would 
imply $N^{-1/2}\leq N^{-\vartheta} < N^{-1} \log(N)\cdot\cosh(\frac12\cdot 3^{-1}\log 3) < \frac54 N^{-1}\log N$. 
The last three inequalities would imply $N^{-1/2}\log(N^{1/2}) > \frac25$, which is absurd (since $x^{-1} \log x \leq e^{-1} < \frac25$ for $x>0$). 
Therefore $\vartheta$ cannot satisfy $\vartheta \leq \frac12$. Thus we have $\frac12 < \vartheta < 1$, and so 
\begin{equation}\label{BWa2} 
\frac{\log N}{2N} < \varepsilon < \frac{\log N}{N} \leq\frac{1}{e} \,. 
\end{equation} 
\par 
Since $\sinh(\varepsilon) = \sum_{r=1}^{\infty} \varepsilon^{2r-1} / (2r-1)!$, the bounds in \eqref{BWa2} imply 
that we have $\sinh(\varepsilon) = (1 + O(\varepsilon^2))\varepsilon = (1+O(N^{-1}))\varepsilon$. 
Thus, assuming that $N$ is sufficiently large, it will follow from the equation \eqref{BWa1} (where 
we have $\vartheta := \varepsilon N/\log N$)  that  
\[ 
\log(2\vartheta) + \log\log(N) - \log(N) = -\vartheta \log N + O\left( N^{-1}\right) \,. 
\] 
Since we have here $1<2\vartheta <2$ and $\vartheta \log N = N\varepsilon$, the required result follows. 
\par 
In the remaining case, one has both $(0, \log N)\supset\{ N\varepsilon , \log\log N\}$ and $N=O(1)$, 
and it follows (trivially) that  $N\varepsilon + \log\log(N) - \log(N) = O(1)$. 
\end{proof} 

\begin{remarks} 
In the above proof it was established that the inequalities in \eqref{BWa2} are valid whenever $N\geq 3$. 
Those inequalities have a couple of useful corollaries, namely the bounds  $N\varepsilon > \frac12$ 
and $N\varepsilon^2 < 1$. 
There is some implicit use of the latter bounds, and of \eqref{BWa2}, in the next subsection. 
\end{remarks}  

\subsection{Completing the proof of Proposition~\mbox{1.2}} 

We apply Lemma~\mbox{3.1}, and then bound (individually) the terms of the sum over $n$ that occurs there,  
by means of Lemma~\mbox{3.2} or the estimate \eqref{I_n_TrivBound}.     
We use the latter (trivial) estimate only when $n\in\{ 1,2,\ldots ,2N-1\}$ is such that  
the interval $B_n := ((1/\delta)^{n-1} , (1/\delta)^{n+1}]$ contains an integer: 
in all remaining cases the bound from Lemma~\mbox{3.2} is used. 
Note that, by \eqref{deltaAlg}, the interval $B_n$ has length $\ell_n = \delta^{N-n}$, 
and so (since $B_n\cap{\mathbb Z}$ is not empty if $\ell_n\geq 1$) 
the bound \eqref{I_n_TrivBound} is employed in those cases where $n\geq N$. 
Note also that, if $n\in\{ 1,2,\ldots , N-1\}$ and $k\in B_n\cap {\mathbb Z}$, then one has both 
\[ 
2\leq k\leq \frac{1}{\delta^{n+1}}\leq \delta^{-N} = e^{N\varepsilon} 
\qquad\text{and}\qquad 
\frac{\log k}{\varepsilon} - 1 \leq n < \frac{\log k}{\varepsilon} + 1 \,.
\] 
In particular, for $k\in{\mathbb Z}$, one has 
$|\{ n\in{\mathbb N} : n\leq N-1\ \text{and}\ B_n\ni k\}| \leq 2$. 
We therefore obtain the bound  
\begin{equation}\label{BWb1} 
\int_0^1 \int_0^1 \left( K(x,y) -k_N(x,y)\right)^2 dx dy < \frac12 \delta^N + \Sigma' + \Sigma_1 + \Sigma_2 \,,
\end{equation} 
where: 
\[ 
\frac{\Sigma_1}{{\textstyle\frac14} (1-\delta)^2} 
= \sum_{n=N}^{2N-1} (2N - n) \delta^{n-1} 
\leq \frac{N \delta^{N-1}}{1-\delta} = \frac{(1+\delta) N}{ \delta^2} \,,  
\] 
\[ 
\frac{\Sigma_2}{{\textstyle\frac14} (1-\delta)^2}  
= 2 \sum_{1<k\leq\exp(N\varepsilon)} \frac{1 +\varepsilon^{-1} \log k}{\delta^2 k} 
\leq \frac{2}{\delta^4 \varepsilon} \sum_{1<k\leq\exp(N\varepsilon)} \frac{\log k}{k} \,, 
\] 
and 
\[ 
\Sigma' =  \frac{(1-\delta)^4}{\delta^3} \sum_{n=1}^{N-1} \delta^{-n} n 
\leq\frac{(1-\delta)^4}{\delta^3} \cdot \frac{N \delta^{1-N}}{(1-\delta)} 
=\frac{(1-\delta)^2 N}{\delta (1+\delta)}  \,. 
\] 

\par

Since we have $1>\delta = e^{-\varepsilon} > 1 - \varepsilon > 1 - e^{-1} > 0$, we may deduce that
\[ 
\Sigma_1 , \Sigma' \leq \frac{N\varepsilon^2}{2\delta^2} \ll N\varepsilon^2 < \frac{2 N^2 \varepsilon^3}{\log N}
\qquad\text{and}\qquad 
{\textstyle\frac12}\delta^N =\frac{1-\delta^2}{2\delta} \ll \varepsilon < \frac{4 N^2\varepsilon^3}{\log^2 N}  \,, 
\] 
and (similarly) that 
\begin{align*} 
\Sigma_2 \leq \frac{{\textstyle\frac12} (1-\delta)^2}{\delta^4 \varepsilon} \sum_{1<k\leq \exp(N\varepsilon)} \frac{\log k}{k} 
 &\leq {\textstyle\frac12} e^{4\varepsilon} \varepsilon \cdot \left( 
\int_1^{\exp(N\varepsilon)} (\log x) x^{-1} dx + O(1)\right) \\ 
 &=  {\textstyle\frac14} e^{4\varepsilon} N^2 \varepsilon^3  \cdot \left( 1 + O\left( (\log N)^{-2}\right)\right) \,. 
\end{align*} 
It follows by these estimates, \eqref{BWb1} and Lemma~\mbox{3.4}, that one has  
\begin{align*} 
\int_0^1 \int_0^1 \left( K(x,y) -k_N(x,y)\right)^2 dx dy 
 &\leq {\textstyle\frac14} N^2 \varepsilon^3  \cdot \left( 1 + O\left( (\log N)^{-1}\right)\right)  \\ 
 &\leq {\textstyle\frac14} N^{-1} (N\varepsilon)^2  \cdot \left( N\varepsilon + O(1)\right) \,. 
\end{align*} 
These bounds yield the required result, for Lemma~\mbox{3.4} implies that if $N$ is sufficiently large (in absolute terms), 
then the last of the bounds just obtained will be less than 
the upper bound, $\frac14 N^{-1} (\log N)^3$, that appears in \eqref{Propn1Bound}. 
\hfill $\square$ 

\begin{remarks}
By elaboration of the proof just completed, 
one can show that if  $N$ is sufficiently large (so that $\varepsilon = \varepsilon(N)$ is sufficiently small)  
then 
\begin{align*} 
\frac{\log^3 N}{4N} > \| K - k_N\|^2 &\gg \sum_{\varepsilon^{-2/3} < k \leq \frac14 \varepsilon^{-1}} 
\frac{\varepsilon \log k}{k}\cdot 
\left| {\mathbb Z}\cap\Bigl( \frac{\log k}{\varepsilon} - \frac14 , \frac{\log k}{\varepsilon} + \frac14\Bigr]  \right| \\ 
 &\gg \varepsilon\log^2 (1/\varepsilon)
\end{align*} 
(the last of these bounds following by \cite[Lemma~5.4.3, Corollary~2]{Hu1996}). It follows that one has 
$\| K - k_N\|^2 \asymp N^{-1}\log^3 N$ for $N\geq 3$ (this relation  
being trivially valid when $3\leq N \ll 1$, since one can never have 
$\| K - k_N\| = 0$). 
\par 
With more effort it can be shown that 
\begin{equation}\label{Uncond_Frob_Asymp}
\| K - k_N\|^2 =\left( 1 + O\bigl( \log^{-\frac23} N\bigr)\right) \cdot {\textstyle\frac{7}{60}} N^{-1} \log^3 N 
\quad\ \text{($N\geq 3$)} , 
\end{equation} 
and that if the Riemann Hypothesis is correct then one may substitute 
\begin{equation*}
- \left( c + 3\log\log N\right) (\log N)^{-1} + O\left( (\log\log N)^2 (\log N)^{-2} \right) , 
\end{equation*}
with $c := 2\log(2\pi) + \frac{23}{7}\log(2) - \frac{137}{30}\approx 1.38657106$, 
for the $O$-term in \eqref{Uncond_Frob_Asymp}. We have a proof of these assertions that 
utilises Lemmas~\mbox{B.3}, \mbox{B.5} and~\mbox{B.6} from our Appendix~B. Although this proof is relatively straightforward, 
it is rather long, and so best omitted: 
all our other results are independent of it. It is however worth noting  
the overlap between this proof and that of \cite[Theorem~4]{Wa2018} presented in \cite[Section~6]{Wa2018}:  
both proofs depend on certain bounds for the sum 
\begin{equation*} 
\sum_{\quad 1<k<e^{(N+1)\varepsilon}}   \left( \frac{\log k}{k}\right)\cdot B_4\left(\left\{ \frac{\log k}{\varepsilon}\right\}\right)\;, 
\end{equation*} 
where $\varepsilon = \varepsilon(N)$ and $B_4(x)$ is the Bernoulli polynomial $x^4 - 2x^3 + x^2 - \frac{1}{30}$. The required 
bounds are deduced, via \cite[Equation~(6.54) and Lemma~6.8]{Wa2018}, from upper bounds for the absolute value 
of $\zeta'(s)$ on the line ${\rm Re}(s) = 1$. 
For further details, and references to relevant results in the literature, see \cite[Lemma~6.9]{Wa2018} and 
the proof and `Remark' following it. 
\end{remarks}

\section{The connection with Hankel operators} 
\subsection{The linear integral operator with kernel $K$} 
We define $B_K $ to be the linear integral operator on the space 
$L^2 [0,1]$ that has kernel function $K$. 
Thus we have: 
\begin{equation*} 
\left( B_K  f\right)(x) := \int_0^1 K(x,y) f(y) dy\qquad \text{($f\in L^2 [0,1]$ and $0\leq x\leq 1$)} . 
\end{equation*} 
Since $K$ is symmetric and square-integrable on $[0,1]\times[0,1]$, the operator $B_K $ is 
Hermitian and compact (see for example \cite[(3.9-3)]{Tr1957}, which implies compactness). 
The spectrum of  $B_K $ is  
$\{ 0\}\cup\{ 1/\lambda : \lambda\ {\rm is\ an\ eigenvalue\ of}\ K\}$, 
and, since  $B_K $ is Hermitian, its singular values, 
$\sigma_1\geq \sigma_2\geq\sigma_3\geq\ \ldots $,   
satisfy $\sigma_m = 1/|\lambda_m|\,$ ($m\in{\mathbb N}$). 
Given  the definition \eqref{DefK}, it follows by \cite[(3.10-8), (2.5-17) and (2.1-3)]{Tr1957} 
that $B_K $ is a Hilbert-Schmidt operator, with Hilbert-Schmidt norm $\| B_K \|_{\rm HS}\geq 0$ 
satisfying: 
\begin{equation}\label{KopHS} 
\| B_K \|_{\rm HS}^2 
:= \sum_{m=1}^{\infty} \sigma_m^2 = \int_0^1 \int_0^1 K^2 (x,y) dx dy < \infty\;.
\end{equation} 
\subsection{An equivalent Hankel operator on $L^2 (0,\infty)$} 
Let $S$ be the linear operator given by 
\begin{equation}\label{Def_opU} 
\left( S f\right) (t) := e^{-\frac{1}{2} t} f\left( e^{-t}\right) \qquad \text{($f\in L^2 [0,1]$ and $0<t<\infty$)} . 
\end{equation} 
Then $S$ maps the space $L^2 [0,1]$ isometrically onto $L^2 (0,\infty)$. 
Therefore the operators $B_K  : L^2 [0,1] \rightarrow L^2 [0,1]$ and  
$S B_K  S^{-1} :  L^2 (0,\infty) \rightarrow L^2 (0,\infty)$ are unitarily equivalent. 
A calculation shows that the latter operator (necessarily compact and Hermitian) is the Hankel operator 
with kernel 
\begin{equation}\label{kernel-h} 
h(t) = \left(  S K(1, \cdot)\right) (t) =  e^{-\frac{1}{2} t} K\left( 1,  e^{-t}\right) \qquad \text{($0<t<\infty$)} . 
\end{equation} 
That is, we have $S B_K  S^{-1} = \Gamma_h$, where 
\begin{equation}\label{DefGamma_h} 
\left(\Gamma_h g\right) (t) := \int_0^{\infty} h(t+u) g(u) du 
\qquad \text{($g\in L^2 (0,\infty)$ and $0<t<\infty$)} . 
\end{equation} 
Since $\Gamma_h$ is unitarily equivalent to  
$B_K $, it has the same spectrum and singular values as that operator, 
and so is (like $B_K $) a Hilbert-Schmidt operator. We have  
\begin{equation}\label{GammaHSnorm} 
\infty > \| B_K \|_{\rm HS} = \| \Gamma_h\|_{\rm HS} 
= \left( \int_0^{\infty} h^2 (t) t dt\right)^{\frac12} \;, 
\end{equation} 
in which the final equality (a special case of a result  \cite[Theorem~7.3~(ii)]{Pa1988} 
in the theory of Hankel operators) may here be obtained  
as a corollary of \eqref{KopHS}, \eqref{kernel-h} and the first equality in 
\eqref{GammaHSnorm}. 

\subsection{Equivalent Hankel operators on other spaces} 
The space $L^2 (0,\infty)$ is isomorphic (as a Hilbert space) to several other 
spaces that have been much studied. 
These include: 
\begin{description}[leftmargin = 9.25em, rightmargin = 2.75em, 
labelwidth=4em, labelsep = 1.0em, itemindent=-2.5em, align=right]
\item[$\ell^2$] ---\quad the space  of `square summable' complex sequences $\{ a_n\}_{n=0}^{\infty}$ 
with norm $\| \{ a_n\}_{n=0}^{\infty}\|_{\ell^2} := (\sum_{n=0}^{\infty} |a_n|^2)^{\frac12} < \infty$;  
\item[$H^2  ({\mathbb D})$] ---\quad the Hardy space  of functions $g(z) = \sum_{n=0}^{\infty} a_n z^n$ analytic on the open disc 
${\mathbb D} := \{ z\in{\mathbb C} : |z|<1\}$ with norm $\| g\|_{H^2  ({\mathbb D})} := (\sum_{n=0}^{\infty} |a_n|^2)^{\frac12} < \infty$; 
\item[$H^2  ( {\mathbb C}_{+})$] ---\quad the Hardy space 
 of functions $G(s)$ analytic on the half plane ${\rm Re}(s) >0$ 
with norm $\| G\|_{H^2  ( {\mathbb C}_{+})} := \sup_{x>0} ( \int_{-\infty}^{\infty} |G(x+iy)|^2 dy)^{\frac12} < \infty$. 
\end{description} 
One has, in particular,  
the surjective isomorphisms $L : L^2 (0,\infty) \rightarrow H^2  ({\mathbb C}_{+})$, 
$V :  H^2  ({\mathbb D}) \rightarrow  H^2  ({\mathbb C}_{+})$ and 
$X : \ell^2\rightarrow H^2  ({\mathbb D})$ given by: 
\begin{align*} 
(L\varphi )(s) &:= (2\pi)^{-1/2} \mathscr{L}\left( \varphi (t) ; s\right) \qquad\quad\ \,\text{($\varphi \in L^2 (0,\infty)$)} ,\\ 
(Vg)(s)  &:=  \pi^{-1/2} (1+s)^{-1} g(Ms)\qquad\text{($g\in  H^2  ({\mathbb D})$)} , \\ 
( X\alpha ) (z) &:= {\textstyle\sum_{m=0}^{\infty} a_m z^m } 
\qquad\qquad\ \ \,\text{($\alpha = \left\{ a_m\right\}_{m=0}^{\infty}\in\ell^2$)} ,
\end{align*} 
where 
\begin{equation}\label{Moebius_Transform} 
Ms := \frac{1 - s}{1+s} 
\end{equation} 
and $ \mathscr{L}( \varphi (t) ; s)$ denotes the Laplace transform of $\varphi$: 
see \cite[Pages~15 and~23--26]{Pa1988} for relevant discussion and proofs. 
It follows that the operators $\widehat \Gamma_h := L \Gamma_h L^{-1} = \mathscr{L} \Gamma_h \mathscr{L}^{-1}$, 
$\widetilde \Gamma_h := V^{-1} \widehat \Gamma_h V$ and $\ddot \Gamma_h := X^{-1}\widetilde \Gamma_h X$ 
(acting on the spaces $H^2  ({\mathbb C}_{+})$, $H^2  ({\mathbb D})$ and $\ell^2$, repectively) 
are each unitarily equivalent to the operator $\Gamma_h$. 
\par 
Regarding $\widehat \Gamma_h$, it can be shown
that, with 
\begin{equation}\label{Def-H(s)} 
H(s) :=  \mathscr{L} (h(t) ; s)\qquad\text{(${\rm Re}(s) > 0$)}  , 
\end{equation}  
one has, in a certain sense,  
\begin{equation*} 
\big(\widehat \Gamma_h G\big)(s) = P_{+}  (H(s)G(-s))\qquad\text{($G\in H^2 ( {\mathbb C}_{+})$)} , 
\end{equation*} 
where $P_{+}$ denotes orthogonal projection from $L^2 (i{\mathbb R})$ onto a subspace identified with 
$H^2  ( {\mathbb C}_{+})$:  see \cite[Chapter~4]{Pa1988} for a proof, 
and \cite[Pages 13--18 and~23--28]{Pa1988} for relevant definitions and theory.
Since the definitions \eqref{kernel-h} and \eqref{DefK} imply $h\in L^2 (0,\infty)$, it follows from 
\eqref{Def-H(s)} that we have $H=\sqrt{2\pi} Lh\in H^2 ( {\mathbb C}_{+})$.  
\par 
One can express $H(s)$ explicitly, in terms of Riemann's zeta function, $\zeta(s)$. 
By the definitions \eqref{kernel-h} and \eqref{DefK}, 
we have $h(t) = - e^{-\frac{1}{2} t} \widetilde B_1 (e^t)$, 
where $\widetilde B_1 (x)$ is the periodic Bernoulli function 
$\{ x\} - \frac12 = (x -\lfloor x\rfloor) -\frac12$. By this and \eqref{Def-H(s)}, it follows that  
$H(s) = - \mathscr{L}(\widetilde B_1 (e^t) ; s+\frac12)$. 
Through the 
substitution $t=\log x\,$ ($t\geq 0$) we find that, for ${\rm Re}(s) > 0$, $s\neq 1$, 
one has:  
\begin{align}\label{LaplaceB1} 
\mathscr{L}\left( \widetilde B_1 (e^t) ; s\right) := 
\int_0^{\infty} e^{-st} \widetilde B_1 (e^t) dt 
 &=\int_1^{\infty} x^{-(s+1)} \left( \{ x\} - \textstyle{\frac12}\right)  dx \nonumber\\ 
 &= -\frac{\left(\zeta(s) - \frac{1}{s-1} - \frac12\right)}{s} 
\end{align} 
(see \cite[(2.1.4)]{Ti1986} regarding the final equality here). 
Thus, for ${\rm Re}(s) > -\frac12$, $s\neq\frac12$, we have: 
\begin{equation}\label{Laplace-h}
H(s) = \frac{\zeta(s+\frac12) - \frac{1}{s-\frac12} - \frac12}{s+\frac12} \; . 
\end{equation} 
Through another application of \eqref{LaplaceB1} we obtain:  
\begin{align}\label{Gamma_h-HSnormEval}  
\|\Gamma_h \|_{\rm HS} &= \sqrt{\textstyle{\frac14} - 2\zeta'(0) + \zeta''(0)} \nonumber\\ 
 &= \sqrt{0.0815206105007606323505594460\ldots\ } \nonumber\\ 
 &= 0.285518143908159770648629633\ldots\ .
\end{align} 
This result is used in some of of the numerical work that we describe in Section~6 and Appendix~B; and we find a different use for it in Section~\mbox{7.3}. 
\par  
We prove the first equality of \eqref{Gamma_h-HSnormEval} in Section~\mbox{A.1}.   
Regarding the remaining (numerical) parts of \eqref{Gamma_h-HSnormEval}, we note that \cite[(2.4.5)]{Ti1986} gives us  
$\zeta' (0) = - \frac12 \log(2\pi)$. 
We compute our estimates for $\frac14 - 2\zeta' (0) + \zeta'' (0)$ and 
its positive square root by using 
known estimates \cite[A061444, A257549]{Sl2022} for $\log(2\pi)$ and $-\zeta'' (0)$. 
Our computations were aided by use of the 
GNU Octave programming language, and  (in order to surpass GNU Octave's level of precision) were partly done by hand.  
\par 
We shall now describe, in concrete terms, how $\widetilde \Gamma_h$ and  
$\ddot \Gamma_h$ act on their respective domains, $H^2  ({\mathbb D})$ and $\ell^2$. 
We begin by defining, for $m\in{\mathbb Z}$, 
\begin{equation*} 
p_m (z) := z^m\qquad\text{($z\in{\mathbb D}$)} . 
\end{equation*} 
It is well-known (see, for example, \cite[Lemma~2.3]{Pa1988}) that  
the family $\{ p_m\}_{m=0}^{\infty}$ is an orthonormal basis for $H^2 ({\mathbb D})$. 
Therefore there is a  unique sequence $\{ c_n\}_{n=0}^{\infty}\in\ell^2$ such that 
$\widetilde \Gamma_h \, p_0 = \sum_{n=0}^{\infty} c_n p_n$.   
It can be shown that this $\{ c_n\}_{n=0}^{\infty}$ is a real sequence satisfying 
\begin{equation}\label{Parseval-A} 
\sum_{m=0}^{\infty} (m+1) \left| c_m\right|^2 =  \| \Gamma_h\|_{\rm HS}^2  < \infty \;, 
\end{equation}
and that one has  
\begin{equation}\label{H2(D)HankelAction} 
\widetilde \Gamma_h \sum_{m=0}^{\infty} a_m p_m 
= \sum_{n=0}^{\infty} \left( \sum_{m=0}^{\infty} a_m c_{m+n}\right) p_n 
\end{equation} 
for all sequences $\{ a_m\}_{m=0}^{\infty}\in\ell^2$. One has moreover: 
\begin{equation}\label{c_m-in-series} 
\frac{H(1)}{z} + \sum_{m=0}^{\infty} c_m z^m = \frac{H(Mz)}{z} \qquad\text{($z\in{\mathbb D}$)} , 
\end{equation} 
with $Mz$ and $H(s)$ as defined in \eqref{Moebius_Transform} and \eqref{Def-H(s)}. 
Proofs of \eqref{Parseval-A}, \eqref{H2(D)HankelAction}  and \eqref{c_m-in-series} are given in Section~\mbox{A.2}. 
\par 
We have, in the three results \eqref{H2(D)HankelAction},  \eqref{c_m-in-series} and  \eqref{Laplace-h},  
a fairly explicit description of the action of $\widetilde \Gamma_h$ on $H^2  ({\mathbb D})$. 
There are  several commonly used alternative descriptions: see \cite[Chapter~3]{Pa1988}. 
In particular, if one identifies 
$H^2 ({\mathbb D})$ with a closed subspace of $L^2 (\partial{\mathbb D})$, via the linear and isometric embedding that, 
for $m=0,1,2,\ldots\ $,  maps  
$p_m$ to the function $\tilde p_m : \partial{\mathbb D}\rightarrow{\mathbb C}$ 
satisfying $\tilde p_m (z) = z^m$ for $z=e^{i\theta}$, $0<\theta\leq 2\pi$,  
then, for certain functions $F\in L^2 (\partial{\mathbb D})$, known as `symbols' for  $\widetilde\Gamma_h$, one has 
\begin{equation*} 
\big(\widetilde\Gamma_h \varphi\big) (z) = P(F(z)\varphi(z^{-1})) \qquad
\text{($\varphi\in H^2 ({\mathbb D})$)} , 
\end{equation*} 
where $P$ is the orthogonal projection from $L^2 (\partial{\mathbb D})$ onto $H^2 ({\mathbb D})$. 
By \eqref{H2(D)HankelAction} and \eqref{c_m-in-series}  (or by  \cite[Theorem~4.6]{Pa1988}), 
it may be seen that the function $z\mapsto H(Mz)/z$ is a symbol  for $\widetilde\Gamma_h$. 
\par 
An immediate consequence of \eqref{H2(D)HankelAction} is that, for $\alpha = \{ a_m\}_{m=0}^{\infty}\in \ell^2$,  
one has 
\begin{equation*} 
\ddot\Gamma_h \alpha = X^{-1}\widetilde\Gamma_h X \alpha = \beta\;, 
\end{equation*} 
where $\beta = \{ b_m\}_{m=0}^{\infty}\in \ell^2$ is given by the infinite matrix equation: 
\begin{equation}\label{InfiniteHankelMatrix}
\begin{pmatrix}  b_0 \cr b_1 \cr b_2 \cr \vdots \cr\end{pmatrix} = 
\begin{pmatrix} 
c_0 &c_1 &c_2&\cdots \cr c_1 &c_2 &c_3&\cdots \cr c_2 &c_3 &c_4&\cdots \cr \vdots &\vdots &\vdots &\ddots \cr 
\end{pmatrix} 
\begin{pmatrix} a_0 \cr a_1 \cr a_2 \cr \vdots \cr\end{pmatrix} \;. 
\end{equation}

\subsection{Consequences for the eigenvalues of $K$} 
Thanks to the unitary equivalence of 
$B_K$, $\Gamma_h$  and $\ddot \Gamma_h$ there are a number of quite general results 
on Hankel operators (previously established by experts in the field) that have significant  
implications concerning the spectrum of $B_K$. We shall consider here just two such results. 
\par
We recall that a compact linear operator $\Gamma$ on a Hilbert space $H$ 
is `trace-class' if and only if one has $\sum_{m=1}^{\infty} \sigma_m (\Gamma) < \infty$ when  
$\sigma_1(\Gamma)\geq \sigma_2(\Gamma) \geq\ \ldots\ $
are the singular values of $\Gamma$. Since trace-class operators have especially nice properties 
(see Partington's book \cite[Chapter~1]{Pa1988} for details) it is very natural to ask if the operator $B_K$ is trace-class.
We can answer this question with the help of the following general result, taken from \cite{Pa1988}. 
\begin{theorem}  
Let $f\in L^1 (0,\infty)\cap L^2 (0,\infty)$ and suppose that the Hankel operator $\Gamma_f : L^2 (0,\infty)\rightarrow L^2 (0,\infty)$  
defined by 
\begin{equation*} 
(\Gamma_f g)(t) := \int_0^{\infty} f(t+u)g(u)du\qquad\text{($g\in L^2 (0,\infty)$, $0<t<\infty$)}  
\end{equation*}  
is trace-class. Then $f$ is equal almost everywhere to a function $f_0$ that is continuous on $(0,\infty)$. 
\end{theorem} 
\begin{proof} 
This theorem is the first part of \cite[Corollary~7.10]{Pa1988}, which is a corollary of  
results first obtained in  \cite{Pe1980} and \cite{CoRo1980}. 
\end{proof} 
\begin{corollary} 
The operators $\Gamma_h$ and $B_K $  are not trace-class. One has: 
\begin{equation}\label{NotTraceClass} 
\sum_{m=1}^{\infty} \frac{1}{\left|\lambda_m\right|} = \infty\;. 
\end{equation} 
\end{corollary} 
\begin{proof} 
We know (see Sections~\mbox{4.1} and~\mbox{4.2}) that $\Gamma_h$  is a compact linear operator on the Hilbert space $L^2 (0,\infty)$, and 
that its singular values are the numbers $1/|\lambda_1|\geq 1/|\lambda_2|\geq\ \ldots\ $. Therefore 
$\Gamma_h$ is trace-class if and only if the series $\sum_{m=1}^{\infty} \lambda_m^{-1}$ is absolutely convergent. 
\par 
Note  that $\Gamma_h$ is the Hankel operator defined on $L^2 (0,\infty)$ by \eqref{DefGamma_h}, 
and that, by \eqref{kernel-h} and \eqref{DefK}, the kernel function $h$ lies in both of the spaces $L^1 (0,\infty)$ and $L^2 (0,\infty)$. 
Theorem~\mbox{4.1} therefore implies that if the series $\sum_{m=1}^{\infty} \lambda_m^{-1}$ does converge absolutely 
(so that $\Gamma_h$ is trace-class) 
then $h$ is equal almost everywhere to a function $h_0$ that is continuous on $(0,\infty)$, and so one has 
$\lim_{a\rightarrow c-} (c-a)^{-1}\int_a^c h(t) dt = h_0 (c) = \lim_{b\rightarrow c+} (b-c)^{-1}\int_c^b h(t) dt$ 
for all $c\in (0,\infty)$. We observe, however, that it follows from 
\eqref{kernel-h} and \eqref{DefK} that, whenever $c\in\{ \log 2, \log 3, \log 4,\ \ldots\ \}$,   
the last two limits are not equal: the first equaling $-\frac12 e^{-\frac12 c} < 0$, 
while the second equals $\frac12 e^{-\frac12 c} > 0$. The series 
$\sum_{m=1}^{\infty} \lambda_m^{-1}$ must therefore not be absolutely convergent, and so 
neither $\Gamma_h$ nor the (unitarily equivalent) operator $B_K$ is trace-class. 
\end{proof} 

The authors of \cite{MPT1995} have succeeded in characterizing the spectra of bounded Hermitian Hankel operators. 
We state here the special case of \cite[Theorem~1]{MPT1995} applicable to compact operators, followed by a corollary 
concerning its implications for both the operator $B_K$ and the indices of eigenvalues of the kernel $K$. 

\begin{theorem}[Megretskii, Peller and Treil] 
Let $B$ be a compact Hermitian linear operator on a Hilbert space with identity operator $I$. 
Then $B$ is unitarily equivalent to a Hankel operator if and only if the following three conditions are satisfied: 
\begin{description} 
\item{(C1)} \quad $\dim\ker(B) \in \{ 0 , \infty\}$; 
\item{(C2)} \quad $B$ is non-invertible; 
\item{(C3)} \quad $|\dim\ker (B + \varkappa I) - \dim\ker (B - \varkappa I) | \leq 1$ for $ \varkappa\in{\mathbb R}\backslash\{ 0\}$. 
\end{description} 
\end{theorem} 
\begin{proof} 
The hypotheses imply that every non-zero complex number contained in the spectrum of $B$ 
is both real and an eigenvalue of $B$ 
(see, for example, \cite[Theorems~8.4-4 and~9.1-1]{Kr1978}). It follows (see \cite[Introduction, Condition~(3)]{MPT1995}) 
that \mbox{(C3)}, above, holds if and only if 
the corresponding condition \mbox{(C3)} occurring in the statement of \cite[Theorem~1]{MPT1995} is satisfied in the case  
where the operator $\Gamma$ there is equal to our $B$. Since the remaining conditions on $\Gamma$ 
in \cite[Theorem~1]{MPT1995} are exactly 
those imposed on $B$ in \mbox{(C1)} and~\mbox{(C2)} above,  that theorem therefore 
contains this one. 
\end{proof} 
\begin{remarks} 

\item{\it 1)}\quad The necessity of Condition~\mbox{(C1)}  follows from Beurling's theorem \cite[Theorem~IV]{Be1949} on 
`shift-invariant' closed subspaces of $H^2({\mathbb D})$: see \cite[Lemma~6.3 and Theorem~6.4]{Pa1988}.

\item{\it 2)}\quad  The necessity of \mbox{(C2)} may be understood by considering the case 
of the Hankel operator $\widetilde \Gamma_h$ defined in Section~\mbox{4.3}: observe that, by   
\eqref{H2(D)HankelAction}, one has  
\begin{equation*} 
\left\| \widetilde \Gamma_h \, p_m - 0 p_m\right\|_{H^2({\mathbb D})} =  \left\| \sum_{n=0}^{\infty} c_{m+n} p_n\right\|_{H^2({\mathbb D})} 
= \left(\sum_{\ell =m}^{\infty} |c_{\ell}|^2 \right)^{\!\!\!1/2} \rightarrow 0\quad\text{as $\,m\rightarrow\infty$} ,
\end{equation*} 
so that $0$ lies in the approximate 
point spectrum of $\widetilde\Gamma_h$.

\item{\it 3)}\quad The necessity of \mbox{(C3)}  was first proved (subject to $B$ being compact) in \cite{Pe1990} and, 
by a different method, in \cite{HW1993}. 
\end{remarks} 

\begin{corollary} The operator $B_K$ is non-invertible and the dimension of the space  
$\ker(B_K)$ is either $0$ or $\infty$. 
Each eigenvalue $\lambda$ of the kernel $K$ has index $\iota_K (\lambda)\in{\mathbb N}$ satisfying 
\begin{equation}\label{Balanced_Indices}
\iota_K (\lambda) \leq   
\begin{cases}  1 + \iota_K (-\lambda)\;, &\text{if  $-\lambda$ is an eigenvalue of $K$} ; \\ 1\;, &\text{otherwise} . \end{cases} 
\end{equation}
\end{corollary} 
\begin{proof} 
Put $B=B_K$. Then we know (see Section~\mbox{4.1}) that the hypotheses of Theorem~\mbox{4.3} are satisfied. 
Therefore, since we have seen (in Sections~\mbox{4.2} and~\mbox{4.3}) that $B$ is unitarily equivalent to 
the Hankel operator $\ddot\Gamma_h$, it follows that each of the conditions \mbox{(C1)}, \mbox{(C2)} and \mbox{(C3)}  occurring 
in the statement of Theorem~\mbox{4.3} is satisfied. By \mbox{(C1)} and \mbox{(C2)}, we have 
the first two assertions of the corollary. By \mbox{(C3)}, we have also 
$\dim\ker (B - \varkappa I) \leq 1 + \dim\ker (B + \varkappa I)$ for 
all non-zero eigenvalues $\varkappa$ of $B$ (these being real, since $B$ is Hermitian). 
By the relevant definitions, the non-zero eigenvalues of $B$ are the reciprocals of the eigenvalues of the kernel $K$, 
and each eigenvalue $\lambda$ of $K$ has index $\iota_K (\lambda) =  \dim\ker (B - \lambda^{-1} I)$. Thus, for 
each eigenvalue $\lambda$ of $K$, we have 
$1\leq \iota_K (\lambda) \leq 1 + \dim\ker (B + \lambda^{-1} I) = 1 + \dim\ker (B - (-\lambda)^{-1} I)$. 
The remaining part of the corollary follows from this, 
since $\dim\ker (B - (-\lambda)^{-1} I)$ equals $\iota_K (-\lambda)$ if $-\lambda$ is an eigenvalue of $K$, 
and is otherwise equal to $0$. 
\end{proof} 

\begin{remarks} 
One can construct a quite different proof of \eqref{Balanced_Indices} using  
results from \cite[Section~4]{WaPP} concerning derivatives of eigenfunctions of $K$. 
\end{remarks} 

\section{Some degenerate approximations to $K$}
In the next section we shall describe our work on obtaining numerical approximations to eigenvalues of~$K$. 
We prepare for that, in this section, 
by first revisiting the degenerate kernel $k_N$ defined in Section~\mbox{3.1}. 
We then define a related $N\times N$ matrix $H(N)$ such that the eigenvalues of $k_N$ 
are simply the reciprocals of the non-zero eigenvalues of $H(N)$. 
\par 
The kernel $k_N$ may be considered an approximation to $K$, 
if $N$ is large enough. 
This suggests that reciprocals of some of the larger eigenvalues of $H(N)$ might 
serve as useful  approximations to eigenvalues of $K$. In practice  
we must work instead with an $N\times N$ Hankel matrix 
$H'(N)$ whose elements are approximations to those of $H(N)$. 
We denote by $k_N'$ the corresponding (very close) approximation to the kernel $k_N$. 
\par
In Section~\mbox{5.3} we summarise our work on obtaining, for some specific values of $N$, 
accurate numerical approximations 
to eigenvalues of $H'(N)$. 
This is covered in much greater detail in 
Appendices~C, D and~E: the computation of the elements of $H'(N)$ is covered in Appendix~B. 
See Tables~\mbox{D-1} and~\mbox{D-3} for a selection of the numerical results obtained. 
\subsection{Notation} 
We introduce notation for the eigenvalues of an arbitrary real symmetric integral kernel $k\in L^2 ([0,1]\times [0,1])$ 
similar to that which we have used earlier, for $k=K$, but including (as an extra argument) the symbol denoting the relevant kernel. 
Thus, for example, the eigenvalues of the kernel $k_N$ are denoted by $\lambda_1 (k_N),\lambda_2(k_N),\ldots\ $. 
We shall in each case assume that $\lambda_1 (k),\lambda_2(k),\ldots\ $ are ordered by absolute value, similarly to how 
the eigenvalues $\lambda_1,\lambda_2,\ldots\ $ are ordered (via \eqref{lambdaOrder}). 
For each kernel $k$ encountered, we denote by $\lambda^{+}_n (k)\,$ 
(resp. $\lambda^{-}_n (k)\,$) the $n$-th positive (resp. negative) term in the 
sequence $\lambda_1 (k),\lambda_2(k),\ldots\ $ (though, when $k=K$ we omit the 
argument `$k$' from this notation). Thus the sequence of positive (resp. negative) 
eigenvalues $\lambda^{+}_1 (k),\lambda^{+}_2 (k),\ldots\ $ 
(resp. $\lambda^{-}_1 (k),\lambda^{-}_2 (k),\ldots\ $) 
will always be monotonic increasing (resp. decreasing). 
\par 
In certain cases the sequence of positive (resp. negative) eigenvalues might be finite, or even empty. 
Taking $\omega^{+}(k),\omega^{-}(k)\in {\mathbb N}\cup\{ 0, \infty\}$ to denote  
(respectively) the number of positive eigenvalues of $k$ and the number of 
negative eigenvalues of $k$, we follow Weyl \cite{We1911} in defining, 
for each positive integer $n$ and either (consistent) choice of sign ($\pm$),
the `reciprocal eigenvalue': 
\begin{equation}\label{Def-nu_n} 
\varkappa^{\pm}_n (k)  := 
\begin{cases} 
1/\lambda^{\pm}_n (k) & \text{if $n\leq\omega^{\pm}(k)$} , \\ 0 & \text{otherwise} .  
\end{cases} 
\end{equation}
Similarly, with  $\omega(k) := \omega^{+}(k) + \omega^{-}(k)$, we define 
reciprocal eigenvalues $\varkappa_1 (k),\varkappa_2(k),\ldots\ $ by putting:
\begin{equation}\label{Def-nu_n-unsigned} 
\varkappa_n (k)  := 
\begin{cases} 
1/\lambda_n (k) & \text{if $n\leq\omega(k)$} , \\ 0 & \text{otherwise} ,  
\end{cases} 
\end{equation} 
for $n\in{\mathbb N}$. 
The argument `$k$' in these notations will be omitted when $k=K$.  

\begin{remarks} 
The sequences 
$\{ \varkappa^{+}_n (k)\}_{n\in \mathbb N}$, $\{ \varkappa^{-}_n (k)\}_{n\in \mathbb N}$ and 
$\{ | \varkappa_n (k) | \}_{n\in \mathbb N}$ are each uniquely determined by the choice of $k$.  
One has $\infty > \varkappa^{+}_n(k)\geq \varkappa^{+}_{n+1}(k)\geq 0$, 
$-\infty < \varkappa^{-}_n(k)\leq \varkappa^{-}_{n+1}(k)\leq 0$, 
$\varkappa_n(k)\in{\mathbb R}$ and $|\varkappa_n(k)|\geq |\varkappa_{n+1}(k)|$, for $n\in{\mathbb N}$. 
\end{remarks} 

\subsection{The Hankel matrices $H(N)$ and $H'(N)$}  

Let $N\geq 3$ be an integer. As in Section~3, we take $k_N (x,y)$ to be the kernel given by  
\eqref{Def-k_N}, with $h_{i,j}$ and $\psi_i (x)\,$ ($1\leq i,j\leq N$) determined by  
\eqref{epsilonDef}--\eqref{mu_ijDef}, \eqref{psi_iDef} and \eqref{h_ijDef}. 
Note, in particular, the dependence of $h_{i,j}$ on $\mu_{i,j}$ that is evident in \eqref{h_ijDef}. 
By \eqref{DefK}, \eqref{x_iDef}--\eqref{DefR_ij} and a change of the variable of 
integration in \eqref{mu_ijDef}, we have that 
\begin{equation}\label{m_ijHankelformed}
\mu_{i,j} = M_n := \frac{1}{(1-\delta)^2} \int_{\delta}^1 \int_{\delta}^1 K\left( 1 , \delta^{n-1} uv\right) du dv 
\qquad\text{($n = i + j - 1$)} , 
\end{equation} 
whenever $i,j\in\{ 1, 2, \ldots , N\}$. 
By this, \eqref{h_ijDef} and \eqref{x_iDef}, we get:
\begin{equation}\label{H_nDef} 
h_{i,j} = H_n := (1 - \delta) \delta^{(n-1)/2} M_n \qquad\text{($n = i + j - 1$)} , 
\end{equation} 
whenever $i,j\in\{ 1, 2, \ldots , N\}$. 
Thus the numbers $h_{i,j}\,$ ($1\leq i,j\leq N$) are the elements of the unique 
finite ($N\times N$) Hankel matrix, 
$H(N)$ (say), that has first column $(H_1 , H_2 , \ldots , H_N)^{\rm T}$ and last row 
$(H_N , H_{N + 1} , \ldots , H_{2N - 1})$. 
\par 
By exploiting (similarly to \cite[Section~1.2]{Tr1957}) the fact that $k_N (x,y)$ is a step-function, one can show 
that the eigenvalues of the kernel $k_N$ are the reciprocals of the non-zero eigenvalues of 
the matrix $H(N)$: one finds, in particular, that any eigenfunction $\phi$ of $k_N$  must have the form 
$\phi = \sum_{n=1}^N v_n \psi_n$, where $(v_1 , \ldots , v_N)^{\rm T}\in {\mathbb C}^N$ is an eigenvector of $H(N)$  
and $\psi_1,\ldots ,\psi_N$ are the functions defined on $[0,1]$ by \eqref{psi_iDef}. 
Thus, with $P = P(N) := \omega^{+} (k_N)$ and  $Q = Q(N) := \omega^{-} (k_N)$, we have that 
$P + Q \leq N$, and that 
$\varkappa_1^{+} (k_N) , \ldots , \varkappa_P^{+} (k_N)\,$ (resp. $\varkappa_1^{-} (k_N) , \ldots , \varkappa_Q^{-} (k_N)$) 
are the positive (resp. negative) eigenvalues of the matrix $H(N)$. 
\par
In Section~3 we showed, in effect, that the sequence $k_3 , k_4, k_5 , \ldots\ $ converges to $K$ in the Hilbert space 
$L^2 ([0,1]\times [0,1])$. This (as we shall see) implies that the eigenvalues of $k_N$ will approximate those of $K$, 
provided that $N$ is sufficiently large. Thus it follows (in light of our observations in the previous paragraph) that by 
computing eigenvalues of the matrix $H(N)$, for some sufficiently large $N$, one will obtain good approximations 
to eigenvalues of $K$. This, however, is not quite as straightforward as it seems: for we can only 
obtain estimates for the numerical values of the elements of the matrix $H(N)\,$  (that is, the numbers $H_1,\ldots ,H_{2N - 1}$ 
given by \eqref{m_ijHankelformed} and \eqref{H_nDef}), and are therefore forced to consider how much the  
errors in these estimates may affect the final results. 
\par 
It is helpful to distinguish between the `theoretical' Hankel matrix $H(N)$, and our `best approximation' to it: another $N\times N$ 
Hankel matrix, $H'(N)$. 
The elements $h_{i,j}'$ of $H'(N)$ have the form $H_{i+j-1}'$, where (for $n=1,\ldots ,2N - 1$) the number $H_n'\in{\mathbb R}$ 
is our `best estimate' of the numerical value of $H_n$: 
see Appendix~B for details of how we compute $H_n'$.  
The non-zero eigenvalues of the matrix $H'(N)$ are the reciprocals of the eigenvalues 
of the symmetric integral kernel $k_N'$ that is given by 
\begin{equation}\label{k_N'Def} 
k_N' (x,y) := \sum_{i=1}^N \sum_{j=1}^N H_{i + j - 1}' \psi_i (x) \psi_j (y) 
\quad\text{($0\leq x,y\leq 1$)}  , 
\end{equation} 
with $\psi_1 , \ldots , \psi_N$ given by \eqref{psi_iDef} and \eqref{epsilonDef}--\eqref{x_iOrder}:   
the proof of this is the same, in essence, as that of the relationship between the eigenvalues of $k_N$ 
and the non-zero eigenvalues of $H(N)$. 

\subsection{Approximations to eigenvalues of $H'(N)$}

Since we are only interested in the larger sized eigenvalues of $H'(N)$,   
rather than those that are close to (or equal to) $0$,   
the notation introduced in Section~\mbox{5.1} is applicable:  recall that the non-zero eigenvalues of $H'(N)$ are the 
reciprocals of the eigenvalues of the kernel $k_N'$.  
We observe in particular that, with $P' = P'(N) := \omega^{+} (k_N')$ and $Q' = Q'(N) := \omega^{-} (k_N')$, 
the matrix $H'(N)$ has rank $R' = P' + Q' \leq N$ and 
non-zero eigenvalues $\varkappa_1 (k_N') , \ldots , \varkappa_{R'} (k_N')$, while 
its positive (resp. negative) eigenvalues are $\varkappa_1^{+} (k_N') , \ldots , \varkappa_{P'}^{+} (k_N')\,$ 
(resp. $\varkappa_1^{-} (k_N') , \ldots , \varkappa_{Q'}^{-} (k_N')$). 
\par 
For $n\in\{ 4, 5, \ldots , 21\}$, $N=2^n$ and $M:=\min\{ 384, N-2\}$ we have used the `{\tt eigs()}' 
function of GNU Octave to compute approximations to the $M$ greatest, 
and $M$ least, of the eigenvalues of $H'(N)$. 
In these computations we opted not 
to have the $N\times N$ matrix $H'(N)$ be passed directly (as an argument) to the function {\tt eigs()}. 
We instead pass to {\tt eigs()} a function handle referencing the function {\tt fast\_hmm()} 
that is discussed in Appendix~C: 
this gives the function {\tt eigs()} a fast and accurate means 
of estimating products of the form $H'(N)X$.
For $n = 21\,$ ($N = 2^{21}$) the computation of the relevant $768$ approximate eigenvalues 
took just under 4 hours. The prior computation of 
$H_1', H_2', \ldots , H_{2N-1}'$, with $N = 2^{21}$, took far longer (about 180 hours),  due in part to 
the use made there of Octave's {\tt interval} package. 
One constraint that influenced our choice of $M$ 
was the limited working memory of our computer (16 gigabytes of RAM), 
which was just enough to cope with the relevant calls to {\tt eigs()}  in the case $n=21$.
\par 
These computations yielded numerical approximations 
\begin{equation}\label{alphasOK}
\alpha^{+}_1 > \alpha^{+}_2 > \cdots > \alpha^{+}_M 
\end{equation} 
(resp. $\alpha^{-}_1 < \alpha^{-}_2 < \cdots < \alpha^{-}_M$) 
to the $M$ greatest (resp. $M$ least) eigenvalues of $H'(N)$, and specific  
vectors ${\bf v}^{+}_1,{\bf v}^{+}_2,\ldots ,{\bf v}^{+}_M\in{\mathbb R}^N\,$ 
(resp. ${\bf v}^{-}_1,{\bf v}^{-}_2,\ldots ,{\bf v}^{-}_M\in{\mathbb R}^N$) 
approximating an orthonormal system of eigenvectors for those $M$ eigenvalues. 
To best explain our use of this data (particularly in the cases where $n\leq 9$) some further notation 
is helpful. We define: 
\begin{equation*}
p' = p'(N) := \left| \left\{ j\leq M : \alpha^{+}_j > 0\right\}\right| 
\quad\text{and}\quad 
q' = q'(N) := \left| \left\{ j\leq M : \alpha^{-}_j < 0\right\}\right| \;.
\end{equation*}
Upon examination of the relevant data, it turns out that we have 
\begin{equation*}
| p' - q' |\leq 2 \quad \text{and}\quad p' + q' = \min\{ 2M , N\} = N' \quad \text{(say)} 
\end{equation*} 
in all cases (i.e. for $4\leq n\leq 21$). 
In particular,  we have $p' + q' = N' = N$ when $4\leq n\leq 9$, 
and $p' = q' = M = \frac12 N' \leq \frac38 N$ when $10\leq n\leq 21$. 
We now write 
\begin{equation*}
\left( \alpha^{+}_1 , \ldots , \alpha^{+}_{p'}, \alpha^{-}_{q'} , \ldots , \alpha^{-}_1\right) 
= \left( \alpha_1,\ldots , \alpha_{N'}\right)\quad \text{(say)} , 
\end{equation*}
\begin{equation*} 
\left( {\bf v}^{+}_1 , \ldots , {\bf v}^{+}_{p'}, {\bf v}^{-}_{q'} , \ldots , {\bf v}^{-}_1\right) 
= \left( {\bf v}_1,\ldots , {\bf v}_{N'}\right)\quad \text{(say)} , 
\end{equation*}
and put: 
\begin{equation} \label{Def-H''}
H''(N,N') := \sum_{j=1}^{N'} \alpha_j {\bf v}_j {\bf v}_j^{\rm T}
\end{equation} 
and 
\begin{equation} \label{Def-A(N)}
A = A(N) := H'(N) - H''(N,N')  
\end{equation} 
(so that both $H''(N,N')$ and $A(N)$ are real and symmetric $N\times N$ matrices). 
This notation will be of use in our subsequent discussion (and is used in Appendix~D): 
note that it ensures that we have $\alpha_1 > \alpha_2 > \ldots > \alpha_{N'}$.
\par
Ideally, we would like to be certain that   
$\alpha^{+}_1 , \ldots , \alpha^{+}_{p'}\,$ (resp. $\alpha^{-}_1 , \ldots , \alpha^{-}_{q'}$) 
are indeed accurate approximations to the eigenvalues 
$\varkappa^{+}_1 (k_N'), \ldots , \varkappa^{+}_{p'}(k_N')\,$ 
(resp. $ \varkappa^{-}_1(k_N') , \ldots ,  \varkappa^{-}_{q'}(k_N')$).  
We have been able to establish that this is the case when $4\leq n\leq 9$.     
For each such $n$ we can compute 
(using the Gershgorin disc theorem and an upper bound for $\| A\|_2$) 
a pairwise disjoint set of very short intervals, 
$[a_j, b_j]\ni\alpha_j\,$ ($1\leq j\leq N$), each containing exactly one of the eigenvalues of $H'(N)$: 
for further details see Section~\mbox{D.1}. 
\par 
For $n\geq  10$, we have $N' = 2M = 768\leq \frac34 N$, so that the Gershgorin disk theorem is not 
readily applicable. Nevertheless, we have managed to compute, for $10\leq n\leq 21$, satisfactory 
lower bounds for the moduli of the $M$ greatest (resp. $M$ least) eigenvalues of $H'(N)$.
These lower bounds are positive numbers $L^{+}_1 , \ldots , L^{+}_M > 0\,$ 
(resp. $L^{-}_1 , \ldots ,  L^{-}_M$), dependent on $N$,  such that one has 
\begin{equation}\label{LBs_for_eigenmoduli} 
\left| \varkappa^{\pm}_j (k_N')\right| \geq L^{\pm}_j\quad\text{($1\leq j\leq M$)} , 
\end{equation}
for either consistent choice of sign ($\pm$). 
Section~\mbox{D.2} contains further details. 
These bounds are considered satisfactory due 
to $L^{\pm}_j$ being, in each case, a number that is only very slightly smaller than $|\alpha^{\pm}_j|$. 
\par 
We have also computed, just for $n=21$, certain non-trivial 
upper bounds $U^{\pm}_1 , \ldots ,  U^{\pm}_M\,$ 
for the moduli of eigenvalues of $H'(N)$, complementary to the 
lower bounds \eqref{LBs_for_eigenmoduli}. Thus, when $N=2^{21}$ and $1\leq j\leq M=384$, one has  
\begin{equation}\label{n=21intervals}
\left[ L^{+}_j , U^{+}_j\right] \ni \varkappa^{+}_j (k_N') 
\quad\text{and}\quad 
\left[ -U^{-}_j , -L^{-}_j\right] \ni \varkappa^{-}_j (k_N') \;. 
\end{equation}
The bounds $U^{\pm}_j\,$ ($1\leq j\leq M$) result from 
an upper bound for $\| A\|_2$ that (in turn) is deduced from a sharp numerical upper bound for $\| A^2\|$. 
In computing the latter bound, for $n=21\,$ ($N=2^{21}$), we faced certain practical problems  
(due to the size of $A$, and the limitations of the hardware and software 
that we had at our disposal): these we overcame 
by adopting a somewhat elaborate method. 
We split the computation into two parts: the first part (computation of a sharp upper  
bound for $\| (H'(2^{21}))^2 \|^2$) is described in Section~\mbox{C.3}, 
while the second part (application of the bound for $\| (H'(2^{21}))^2 \|^2$) 
is covered in Section~\mbox{D.3}. See  \eqref{SpecNormA21}  for 
an explicit statement of  our numerical bound for $\| A(2^{21})\|_2$.  
For some discussion of the pros and cons of this bound see Remarks~\mbox{D.3}.  
\par 
In Section~\mbox{D.4} we describe how (with the help of our upper bound for $\| A(2^{21})\|_2$) 
we have computed the above mentioned bounds   $U^{\pm}_j\,$ ($1\leq j\leq M$). 
The algorithms employed ensure that  we have both \eqref{n=21intervals} 
and $[L^{+}_j , U^{+}_j]\times [-U^{-}_j , -L^{-}_j]\ni (\alpha^{+}_j , \alpha^{-}_j )$,  for $1\leq j\leq M$. 
Although  the intervals $[ L^{+}_j , U^{+}_j]$ and $[ -U^{-}_j , -L^{-}_j]$ turn out to be quite `long' 
when $50\leq j\leq M=384\,$ ($U^{\pm}_j / L^{\pm}_j > \exp(0.01)$ in these cases), 
we find (in contrast) that  $U^{\pm}_j / L^{\pm}_j < \exp(10^{-8})$ when $j\leq 49$: 
consequently we can (at least) say for certain that, for $N=2^{21}$, the number 
$\alpha^{+}_j\,$ (resp. $\alpha^{-}_j$) is a fairly accurate approximation to 
$\varkappa^{+}_j (k_N')\,$ (resp. $\varkappa^{-}_j (k_N')$) when $j\leq 49$. 
These results do not satisfy us completely, 
since they tell us little about the accuracy of the 
majority of the approximations $\alpha^{\pm}_1,\ldots ,\alpha^{\pm}_M$. 
\par 
In  Section~\mbox{D.5} we describe 
a $4$-step statistics-based algorithm that helps us to  
further probe the accuracy of $\alpha^{\pm}_1,\ldots ,\alpha^{\pm}_M$. 
This algorithm produces a result $\tilde R = \tilde R(N) > 0$, to be used as a hypothetical 
upper bound for $\| A(N)\|_2$: 
since the probability of it generating results $\tilde R(2^{10}), \ldots , \tilde R(2^{21})$ 
that satisfy $\tilde R(2^n) \geq \| A(2^n)\|_2\,$ ($10\leq n\leq 21$) is, arguably, very close to~$1$,  
we call $\tilde R\,$ (i.e. $\tilde R(N)$) a `probable upper bound' for $\| A(N)\|_2$.  
There are three caveats this. There is, firstly, the 
obvious point that even a very high probability of success does not make success certain.  
Secondly, the algorithm that we use requires an initial input of random data 
(some fairly long random sequence of elements of the set $\{ 0, 1\}$), 
and so the validity of our claim concerning the probability of having 
$\tilde R\geq \| A\|_2\,$ ($10\leq n\leq 21$) is dependent on our having 
convenient access to a valid source of randomness: 
we have therefore employed a random number generator that is 
the least obviously  defective amongst those that are known to us (see Appendix~E for details). 
Thirdly, we have taken some (highly expedient)  
short cuts in carrying out the theoretical rounding error analysis on which 
Step~3 of the algorithm depends, by making certain assumptions; 
we might, thereby, have come to false conclusions:  
for relevant details, see Remarks~\mbox{D.6} and~\mbox{C.2}~\mbox{(1)} and~\mbox{(2)}. 
\par 
As we describe (briefly) in Section~\mbox{D.6}, the methods of Section~\mbox{D.4}  
can be adapted so as to yield, when $10\leq n\leq 21$, certain hypothetical numerical 
bounds $\tilde U^{\pm}_j \geq |\varkappa_j^{\pm} (k_N')|\,$ ($1\leq m\leq M$) 
that are valid if $\tilde R(N) \geq \| A(N)\|_2$. Extending the terminology used for $\tilde R$, 
we call the number $\tilde U^{\pm}_j$ our `probable upper bound' for the 
modulus of the eigenvalue $\varkappa_j^{\pm} (k_N')$. 
\par
Our combined 
numerical results for $10\leq n\leq 21\,$ (i.e. the numbers $L^{\pm}_j$, $U^{\pm}_j$ and $\tilde U^{\pm}_j$ 
whose computation is described in Sections~\mbox{D.2--D.6}) 
are summarised in Section~\mbox{D.7}.
The relevant probable upper bounds $\tilde R(2^{10}), \ldots , \tilde R(2^{21})$ 
are stated (to an accuracy of $16$ significant digits) in Table~\mbox{D-2}, which is in Section~\mbox{D.5}. 
\par 
In Section~6, below, we discuss an application of the numerical data 
$L^{\pm}_j, U^{\pm}_j\,$ ($1\leq j\leq M=384$, $\pm\in\{ + , -\}$) computed in 
our work on the case $N=2^{21}$. We explain there how this data, when combined with 
sharp upper bounds for $\| k_N' - K\|$ and $\| k_N' - k_N\|$, leads to  
fairly accurate estimates for several of the smaller eigenvalues of the kernel $K$.
\par 
In Section~7 we work on the assumption that each probable upper bound  
$\tilde R(2^n)\,$ ($10\leq n\leq 21$) is in fact a valid upper bound for the spectral norm  
of the corresponding matrix $A(2^n)$, so that we have  
$\tilde U^{\pm}_j\geq |\varkappa^{\pm}_j (k_N')|\geq L^{\pm}_j$ 
when $1\leq j\leq M=384$, $\pm\in\{ + , -\}$, $N=2^n$ and $10\leq n\leq 21$. 
The numerical data $L^{\pm}_j=L^{\pm}_j(n)$, $\tilde U^{\pm}_j=\tilde U^{\pm}_j(n)\,$ 
($1\leq j\leq M=384$, $\pm\in\{ + , -\}$, $10\leq n\leq 21$) 
is used in support of conjectures about specific eigenvalues of the kernel $K$. 
We are led (ultimately) to formulate some wider conjectures concerning 
the distribution of the eigenvalues of $K$. 

\section{Bounds for eigenvalues of $K$} 

In this section we describe our work on obtaining bounds 
for specific eigenvalues of $K$. 
Relevant notation and results 
from existing literature are introduced in Sections~\mbox{6.1} and~\mbox{6.2}. 
In Sections~\mbox{6.3} and~\mbox{6.4} we show that, when $N\geq 3$, 
the reciprocal eigenvalues of $K$ 
are bounded (above and below) by certain expressions 
involving the Hilbert-Schmidt norms $\| k_N' - K\|$, $\| k_N' - k_N\|$ and eigenvalues of $H'(N)$.  
Our main results there, Lemmas~\mbox{6.8} and~\mbox{6.12}, 
improve upon what can be got by direct application of Weyl's inequalities. 
We have applied these results in computational work that utilises  
numerical bounds for $\| k_N' - K\|$, $\| k_N' - k_N\|$  and eigenvalues of $H'(N)$,  
obtained in respect of the case $N=2^{21}$. This computational work and its results 
are summarised in Section~\mbox{6.5}.

\subsection{Notation}
For any integral kernel $k\in L^2 ([0,1]\times [0,1])$ and 
any $N\times N$ matrix $C$ with elements $c_{i,j}\in{\mathbb C}$, the norms 
$\| k\| , \| C\| \geq 0$ are given by: 
\begin{equation*} 
\| k\|^2 = \int_0^1 \int_0^1 |k(x,y)|^2 dx dy\quad\text{and}\quad 
\| C\|^2 =  \sum_{i=1}^N\sum_{j=1}^N \left| c_{i,j}\right|^2 \;. 
\end{equation*} 
We use the notation $\langle \alpha , \beta\rangle$ for the inner product on $L^2 [0,1]$ given by: 
\begin{equation*} 
\langle \alpha , \beta\rangle := \int_0^1 \alpha(x) \overline{\beta(x)} dx\qquad\text{($\alpha , \beta\in L^2 [0,1]$)} . 
\end{equation*} 
Functions $\alpha , \beta\in L^2 [0,1]$ are `orthogonal' if and only if $\langle \alpha , \beta\rangle  = 0$. 
\par
Note also that $k_N$, $H(N)$, $k_N'$ and $H'(N)$ all have the same meaning as in Section~5, and that 
we continue using the notation defined in Section~\mbox{5.1}. 

\subsection{Statements of requisite known results} 

Let the $N\times N$ matrix $C$ and the integral kernels $k,k'\in  L^2 ([0,1]\times [0,1])$ be real and symmetric. 
Then the following three results hold. 

\begin{lemma} 
One has  
\begin{equation}\label{TraceLemma-k} 
\| k\|^2 = \sum_{n=1}^{\infty} \left( \varkappa_n (k)\right)^2
= \sum_{n=1}^{\infty} \left( \varkappa^{+}_n (k)\right)^2 + \sum_{n=1}^{\infty} \left( \varkappa^{-}_n (k)\right)^2 
\end{equation} 
and 
\begin{equation}\label{TraceLemma-A} 
\| C\|^2 = \sum_{n=1}^N \Lambda_n^2 \;, 
\end{equation} 
where $\Lambda_1,\ldots ,\Lambda_N$ are the eigenvalues of $C$. 
\end{lemma} 
\begin{proof} 
By the relevant definitions, $\| C\|^2 = {\rm Trace}( C^* C)$, with $C^*:= {\overline C}^{\rm T}=C\,$ (since $C$ is real and symmetric). 
Using the spectral decomposition of normal matrices  
(for which see \cite[Section~5.4]{Os1990}) one can establish that 
$|\Lambda_1|^2, \ldots , |\Lambda_N|^2$ are the eigenvalues of  $C^* C$.  
Thus we have $ {\rm Trace}( C^* C) = \sum_{n=1}^N  |\Lambda_n|^2$, and so obtain \eqref{TraceLemma-A} 
(since, by virtue of $C$ being real and symmetric, its eigenvalues $\Lambda_1,\ldots ,\Lambda_N$ are real).  
\par 
The result \eqref{TraceLemma-k} is \cite[Section~3.10 (8)]{Tr1957} (expressed in our notation), and is  
a corollary of the `bilinear formula' \cite[Section~3.9, (3)]{Tr1957}. 
\end{proof} 

\begin{lemma} 
There exists an orthonormal sequence $\phi_n\,$ ($\,{\mathbb N}\ni n\leq\omega(k)$) of eigenfunctions of~$k$,     
each satisfying $\phi_n (x) = \lambda_n (k) \int_0^1 k(x,y) \phi_n (y) dy$ almost everywhere in $[0,1]$.  
Given any such sequence $\phi_n\,$ ($\,{\mathbb N}\ni n\leq\omega(k)$), one has
\begin{equation}\label{HSlemma}
\int_0^1 \int_0^1 k(x,y)\psi(x)\overline{\xi (y)} dx dy   
= \sum_{{\mathbb N}\ni n\leq\omega(k)} \varkappa_n (k) 
\left\langle \psi , \phi_n\right\rangle \left\langle \phi_n , \xi \right\rangle
\end{equation} 
whenever $\psi,\xi\in L^2[0,1]$. 
\end{lemma} 
\begin{proof} 
Suppose, firstly, that $\omega(k)\neq 0\,$ (so that $\omega(k)\in{\mathbb N}\cup\{\infty\}$). 
For the first assertion of the lemma see \cite[Section~3.8]{Tr1957}. 
The result \eqref{HSlemma} is obtained in \cite[Section~3.11]{Tr1957}, as a corollary of a 
theorem of Hilbert and Schmidt (discussed in \cite[Section~3.10]{Tr1957}).
\par 
If $\omega(k) = 0$, then by \eqref{Def-nu_n-unsigned} and \eqref{TraceLemma-k} 
one must have $\| k\| = 0$, so that the double integral in \eqref{HSlemma} has to equal $0\,$ 
(trivially also the value of the sum over $n$ in \eqref{HSlemma}). 
\end{proof} 

\begin{lemma} 
For either consistent choice of sign ($\pm$), one has: 
\begin{equation}\label{Weyl-0} 
\left| \varkappa^{\pm}_{m + n - 1} \left( k + k'\right)\right| \leq 
\left| \varkappa^{\pm}_m (k) + \varkappa^{\pm}_n \left( k'\right)\right| 
\qquad\text{($m,n\in{\mathbb N}$)} . 
\end{equation} 
\end{lemma} 
\begin{proof} 
The case where the sign is `$+$' is \cite[Satz~III]{We1911}, and this case implies the other case (since one has 
$\varkappa^{-}_m (k) = - \varkappa^{+} (-k)\leq 0$, and similar relations involving the kernels $k'$ and $k+k'$). 
\end{proof} 

\begin{remarks}  

\item{\it 1)}\quad It follows immediately from Lemma~\mbox{6.1} that, for either choice of sign ($\pm$), one has 
\begin{equation}\label{Root-n-bound}
\left| \varkappa^{\pm}_m (k)\right| \leq m^{-1/2} \| k\|\qquad\text{($m\in{\mathbb N}$)} . 
\end{equation} 
Indeed, for $m\in{\mathbb N}$ one has: $\| k\|^2 \geq \sum_{n=1}^m  (\varkappa_n (k))^2 \geq m \cdot (\varkappa_m (k))^2$.   

\item{\it 2)}\quad Since  $\max\{ |\varkappa_n (k)|   :  n\in{\mathbb N}\} =  |\varkappa_1(k)| < \infty$,     
one can show (using Bessel's inequality) that for  $\psi,\xi\in L^2[0,1]$ 
the sum over $n$ in \eqref{HSlemma} is absolutely convergent. 
\end{remarks}

\subsection{Lower bounds for $|\varkappa^{\pm}_m|$} 

Let $m\geq 1$ and $N\geq 3$ be integers. Let `$\pm$' be a fixed choice of sign ($+$ or $-$), and take `$\mp$' to denote 
the opposite sign. By \cite[Theorem~1.1]{Wa2019} we know that one has $\omega^{\pm} = \infty$, 
and so $1/\lambda^{\pm}_m = \varkappa^{\pm}_m\neq 0$.  
The reciprocal of any strictly positive lower bound for $|\varkappa^{\pm}_m|$ is an upper bound for $|\lambda^{\pm}_m|$.
\par
By applying \eqref{Weyl-0} with $k=k_N'$, $k' = K - k_N'$ and $n=1$ one gets 
the bound $|\varkappa^{\pm}_m| \leq  |\varkappa^{\pm}_m (k_N')| + |\varkappa^{\pm}_1  (K - k_N')|$. 
By Weyl's inequality \eqref{Weyl-0}, with $k = K$, $k' = k_N' - K$ and $n=1$, one has also 
$|\varkappa^{\pm}_m (k_N')|\leq |\varkappa^{\pm}_m| + |\varkappa^{\pm}_1 (k_N' - K)|$. 
Thus, given that $\varkappa^{\pm}_1 (k_N' - K) =  -\varkappa^{\mp}_1  (K - k_N')$, 
we have:  
\begin{equation}\label{Weyl-1} 
|\varkappa^{\pm}_m (k_N')| + |\varkappa^{\pm}_1  (K - k_N')| \geq |\varkappa^{\pm}_m|  
\geq  |\varkappa^{\pm}_m (k_N')| - |\varkappa^{\mp}_1  (K - k_N')|\;. 
\end{equation} 
By \eqref{Root-n-bound}, one has $ \max\{ |\varkappa^{+}_1  (K - k_N')| ,  |\varkappa^{-}_1  (K - k_N')|\} \leq \| K - k_N'\|$. 
Therefore the bounds \eqref{Weyl-1} certainly imply: 
\begin{equation}\label{Weyl-2} 
\left| \varkappa^{\pm}_m  -   \varkappa^{\pm}_m (k_N')\right| \leq \| K - k_N'\|\;. 
\end{equation} 
One can also show (similarly) that 
\begin{equation}\label{Weyl-3} 
\left| \varkappa^{\pm}_m (k_N)  -   \varkappa^{\pm}_m (k_N')\right| \leq \| k_N - k_N'\|\;. 
\end{equation} 
\par 
In this subsection we shall obtain an improvement of the lower bound for  $|\varkappa^{\pm}_m|$ that is 
implicit in  \eqref{Weyl-2}. We state first a definition and two lemmas that yield the improved lower bound. 

\begin{definition} 
We put
\begin{equation}\label{J_xiDef}
J_{\xi} (k) := \int_0^1 \int_0^1 k(x,y) \xi (x) \overline{\xi(y)} dx dy \;, 
\end{equation} 
when $\xi\in L^2 [0,1]$ and  $k\in L^2 ([0,1]\times [0,1])$.  
\end{definition} 

\begin{lemma} 
Let $\xi ,\varrho \in L^2 [0,1]$. 
Suppose that $\xi (x) = \int_0^1 k_N(x,y) \varrho (y) dy$ almost everywhere in~$[0,1]$.  
Then $J_{\xi} (K) = J_{\xi}(k_N)$. 
\end{lemma} 

\begin{proof} 
By the definition \eqref{J_xiDef}, we have $J_{\xi} (K) - J_{\xi}(k_N) = J_{\xi}(K - k_N)$. 
Therefore the lemma follows if we can show that $J_{\xi}(K - k_N)$ equals $0$. 
\par 
Recall that the kernel $k_N (x,y)$ is given by \eqref{Def-k_N},  
with $h_{i,j}$ and $\psi_i (x)\,$ ($1\leq i,j\leq N$) determined by  
\eqref{epsilonDef}--\eqref{mu_ijDef}, \eqref{psi_iDef} and \eqref{h_ijDef}. 
By this and the hypotheses of the lemma, it follows that  $\xi = \sum_{i=1}^N b_i \psi_i$ 
(in $L^2 [0,1]$), where 
$b_i = \sum_{j=1}^N  \langle h_{i,j} \psi_j , \overline{\varrho}\rangle \in{\mathbb C}\,$ 
($1\leq i\leq N$). From this and \eqref{J_xiDef}, we deduce that one has 
$J_{\xi} \left( K - k_N\right)  = \sum_{i=1}^N \sum_{j=1}^N G_{i,j} b_i \overline{b_j}$, 
where 
\begin{equation*} 
G_{i,j} := \int_0^1 \int_0^1 (K(x,y) - k_N (x,y)) \psi_i (x) \overline{\psi_j (y)} dx dy\qquad 
\text{($1\leq i,j\leq N$)} . 
\end{equation*} 
Therefore the lemma follows if we show that $G_{i,j} = 0\,$  ($1\leq i,j\leq N$). 
\par
Let $i,j\in\{ 1,\ldots ,N\}$. By \eqref{psi_iDef} and \eqref{x_iOrder}, we find that 
\begin{equation}\label{BWe1}
G_{i,j} = \frac{1}{\sqrt{\operatorname{vol}\left( {\mathcal R}_{i,j}\right)}} 
\iint\limits_{{\mathcal R}_{i,j}} \left( K(x,y) - k_N (x,y)\right) dx dy \,, 
\end{equation}
where ${\mathcal R}_{i,j}\subset {\mathbb R}^2$ is the rectangle defined in \eqref{DefR_ij}. 
By \eqref{Def-k_N}, \eqref{psi_iDef}, \eqref{x_iOrder} and \eqref{h_ijDef}, we have here 
$k_N (x,y) = \mu_{i,j}$ for all $(x,y)\in {\mathcal R}_{i,j}$, and so, after recalling also 
the definition \eqref{mu_ijDef} of $\mu_{i,j}$, we find that one has 
\begin{equation*} 
\iint\limits_{{\mathcal R}_{i,j}} k_N (x,y) dx dy = \operatorname{vol}\left( {\mathcal R}_{i,j}\right) \cdot \mu_{i,j} 
= \iint\limits_{{\mathcal R}_{i,j}} K(x,y) dx dy \;. 
\end{equation*} 
By this and \eqref{BWe1}, we have $G_{i,j} = 0$. 
\end{proof} 

\begin{lemma} 
For either (consistent) choice of sign ($\pm$) one has: 
\begin{equation*} 
\left| \varkappa^{\pm}_m\right|\geq \left| \varkappa^{\pm}_m\left( k_N\right)\right|  
\qquad\text{($m,N\in{\mathbb N}$, $N\geq 3$)} . 
\end{equation*} 
\end{lemma} 

\begin{proof} 
We shall adapt methods used by Weyl in his proof of \cite[Satz~I]{We1911}. Let $m\geq 1$ and $N\geq 3$ be integers. 
We consider only the case in which the sign ($\pm$) is `$+$', since the proof of 
the other case is similar: in doing so we may assume that $\varkappa^{+}_m ( k_N)$  is strictly positive, 
since otherwise one would necessarily have $|\varkappa^{+}_m ( k_N)| = 0 \leq |\varkappa^{+}_m|\,$ 
(in view of the relevant definitions). Thus (recalling also relevant points noted in Section~\mbox{5.2}) 
we have $m\leq P := \omega^{+} (k_N) \leq N$. 
\par
Let $\xi : [0,1]\rightarrow{\mathbb C}$ be some finite linear combination of eigenfunctions of $k_N$ associated 
with eigenvalues that are positive (the choice of this linear combination will be refined later). 
It follows that $\xi\in L^2 [0,1]$, and that $\xi$ is orthogonal to every 
eigenfunction of $k_N$ that is associated with a negative eigenvalue ($k_N$ being a real and symmetric integral kernel). 
Therefore, given that $K(x,y)$ and $k_N(x,y)$ are real and symmetric kernels, 
it follows by Lemma~\mbox{6.2} and \cite[Theorem~1.1]{Wa2019} that the integrals $J_{\xi} (K)$ and $J_{\xi} (k_N)\,$  
(defined as in \eqref{J_xiDef}) must satisfy 
\begin{equation}\label{BWd1} 
J_{\xi} (K) \leq\sum_{n=1}^{\infty}  
\varkappa^{+}_n \cdot \left| \left\langle \xi , \phi^{+}_n \right\rangle\right|^2 
\quad\text{and}\quad 
J_{\xi} \left( k_N\right) = \sum_{n=1}^{P}  
\varkappa^{+}_n \left( k_N \right)\cdot\left| \left\langle \xi , \theta^{+}_n \right\rangle\right|^2 \;, 
\end{equation} 
where $\phi^{+}_1 , \phi^{+}_2 , \ldots\ $ (resp. $\theta^{+}_1 , \theta^{+}_2 , \ldots ,  \theta^{+}_{P}$) is some 
orthonormal sequence of eigenfunctions of $K$ (resp. $k_N$) in $L^2 [0,1]$. 
\par 
We now become more specific regarding our choice of $\xi$: since $1\leq m\leq P$, we may suppose that 
\begin{equation}\label{BWd2} 
\xi = \sum_{r = 1}^m c_r \theta^{+}_r\;, 
\end{equation} 
where $c_1 , \ldots , c_m$ are complex constants, and are not all zero. 
If $m > 1$ then the column vectors 
$( \langle \theta^{+}_r , \phi^{+}_1\rangle , \ldots , \langle \theta^{+}_r , \phi^{+}_{m-1}\rangle )^{\rm T}\,$ 
($1\leq r\leq m$) are elements of ${\mathbb C}^{m-1}$, and so these $m$ vectors cannot be linearly independent. 
Thus one can choose $c_1 , \ldots , c_m\in{\mathbb C}\,$ (not all zero) so as to have 
\begin{equation}\label{BWd3} 
\left\langle \xi , \phi^{+}_n\right\rangle = 0 \qquad\text{for all $n\in{\mathbb N}$ with $n\leq m - 1$} 
\end{equation} 
(when $\xi$ is given by \eqref{BWd2}). By \eqref{BWd2} and the orthonormality of $\theta^{+}_1 , \ldots ,\theta^{+}_m$, 
we have here $\|\xi\|^2 = \sum_{r=1}^m |c_r|^2 > 0$, and so (by replacing $\xi$ with $\| \xi\|^{-1} \xi$) 
may suppose that (in addition to \eqref{BWd2} and \eqref{BWd3}) one has: 
\begin{equation}\label{BWd4} 
\| \xi\|^2 = \sum_{r=1}^m |c_r|^2 = 1\;. 
\end{equation} 
Using now \eqref{BWd1}, \eqref{BWd3}, Bessel's inequality and \eqref{BWd4}, we get:  
\begin{equation*} 
J_{\xi} (K) \leq \sum_{n=m}^{\infty}  \varkappa^{+}_n\cdot \left| \left\langle \xi , \phi^{+}_n \right\rangle\right|^2 
\leq \varkappa^{+}_m\cdot \sum_{n=m}^{\infty}  \left| \left\langle \xi , \phi^{+}_n \right\rangle\right|^2 
\leq \varkappa^{+}_m\cdot \|\xi\|^2  = \varkappa^{+}_m\;. 
\end{equation*} 
By  \eqref{BWd1}, \eqref{BWd2}, the orthonormality of $\theta^{+}_1 , \ldots ,  \theta^{+}_P$, 
and  \eqref{BWd4},  one gets also: 
\begin{equation*} 
J_{\xi} \left( k_N\right) = \sum_{n=1}^m \varkappa^{+}_n \left( k_N \right)\cdot \left| c_n\right|^2 
\geq \varkappa^{+}_m \left( k_N \right)\cdot \sum_{n=1}^m \left| c_n\right|^2 
= \varkappa^{+}_m \left( k_N \right) > 0\;. 
\end{equation*} 
\par 
We recall (see Lemma~\mbox{6.2}) that  $\theta^{+}_1 , \ldots ,\theta^{+}_m$ are eigenfunctions of $k_N$ associated with the 
eigenvalues $\lambda^{+}_1 (k_N) , \ldots ,\lambda^{+}_m (k_N)\,$ (respectively).  It therefore follows from 
\eqref{BWd2} that, with $\varrho := \sum_{r = 1}^m c_r \lambda^{+}_r (k_N) \theta^{+}_r \in L^2 [0,1]$, 
one has $\int_0^1 k_N(x,y) \varrho (y) dy = \xi (x)$ almost everywhere in $[0,1]$. Therefore Lemma~\mbox{6.6} implies 
that we have $J_{\xi} (K) = J_{\xi} (k_N)$. By combining this with the bounds on $J_{\xi} (K)$ and $J_{\xi} (k_N)$ obtained 
at the end of the last paragraph, we may deduce that 
one has $0 < \varkappa^{+}_m (k_N) \leq J_{\xi} (K) \leq \varkappa^{+}_m$, and 
so $0 < \varkappa^{+}_m (k_N) \leq \varkappa^{+}_m$. 
\end{proof} 
 
By \eqref{Weyl-3}, we have 
$|\varkappa^{\pm}_m (k_N)| \geq  |\varkappa^{\pm}_m (k_N')| - \| k_N - k_N'\|$ for $m\in{\mathbb N}$.  
Lemma~\mbox{6.7} therefore has the following immediate corollary. 
 
\begin{lemma}
Let $m,N\in{\mathbb N}$, with $N\geq 3$. Then, for either consistent choice of sign ($\pm$), one has: 
\begin{equation*} 
\left|\varkappa^{\pm}_m\right| \geq 
|\varkappa^{\pm}_m (k_N')| - \| k_N - k_N'\| 
\;. 
\end{equation*} 
\end{lemma} 

\begin{remarks} 

\item{\it 1)}\quad When $\| k_N - k_N'\| <  \| K - k_N'\|$, Lemma~\mbox{6.8} gives us a better lower bound for $| \varkappa^{\pm}_m|$ than 
can be obtained directly from \eqref{Weyl-2} alone. 
In our numerical applications of Lemma~\mbox{6.8} (for which see Section~\mbox{6.5})  
the ratio $\| k_N - k_N'\| /  \| K - k_N'\|$ never exceeds $10^{-9}$. 

\item{\it 2)}\quad Let $m\geq 1$ and $N\geq 3$ be integers. Put $\alpha = |\varkappa_m (k_N)|$, so that $\alpha\in [0,\infty)$. 
By Lemma~\mbox{6.7} and our conventions concerning the ordering of eigenvalues and reciprocal eigenvalues 
(for which see \eqref{lambdaOrder}  and Section~\mbox{5.1}), one has  
\begin{align*} 
m &\leq \big| \{ n\in{\mathbb N} : |\varkappa_n (k_N)| \geq \alpha\}\big| \\ 
 &= \big| \{ n\in{\mathbb N} : |\varkappa_n^{+} (k_N)| \geq \alpha\}\big| 
+ \big| \{ n\in{\mathbb N} : |\varkappa_n^{-} (k_N) | \geq \alpha\}\big| \\ 
 &\leq \big| \{ n\in{\mathbb N} : |\varkappa_n^{+}| \geq \alpha\}\big| 
+ \big| \{ n\in{\mathbb N} : |\varkappa_n^{-}| \geq \alpha\}\big|  \\
 &=\big| \{ n\in{\mathbb N} : |\varkappa_n| \geq \alpha\}\big| \;, 
\end{align*}
and so must have $|\varkappa_m |\geq \alpha$. Thus Lemma~\mbox{6.7} implies the bounds:  
\begin{equation}\label{RELBlemma-v2} 
\left| \varkappa_m\right|\geq \left| \varkappa_m\left( k_N\right)\right|  
\qquad\text{($m,N\in{\mathbb N}$, $N\geq 3$)} . 
\end{equation} 

\item{\it 3)}\quad By a very slight elaboration of the reasoning used in our last remark above, 
it may be deduced from 
the result \eqref{Weyl-2} that one has 
\begin{equation}\label{Weyl-2-unsigned} 
\bigl| |\varkappa_m|  -   |\varkappa_m (k_N')|\bigr| \leq \| K - k_N'\| 
\end{equation} 
for all $m,N\in{\mathbb N}$ with $N\geq 3$. Similarly,  \eqref{Weyl-3} implies: 
\begin{equation}\label{Weyl-3-unsigned} 
\bigl| |\varkappa_m (k_N)|  -   |\varkappa_m (k_N')|\bigr| \leq \| k_N - k_N'\|\qquad \text{($m,N\in{\mathbb N}$, $N\geq 3$)} .  
\end{equation} 
It therefore follows from \eqref{RELBlemma-v2} that 
\begin{equation}\label{Weyl-4}
\left|\varkappa_m\right| \geq 
|\varkappa_m (k_N')| - \| k_N - k_N'\| 
\qquad \text{($m,N\in{\mathbb N}$, $N\geq 3$)} . 
\end{equation}

\end{remarks} 

\subsection{Upper bounds for $|\varkappa_m|$} 

The main result of this subsection is Lemma~\mbox{6.12}, which, in cases where both $\| k_N - K\| / |\varkappa_m|$ 
and $\| k_N' - k_N\| / \| k_N - K\|$ are sufficiently small (in absolute terms), gives us an    
upper bound for  $|\varkappa_m|$ stronger than that which is implicit in  \eqref{Weyl-2-unsigned}.

\begin{lemma} 
One has 
\begin{equation}\label{ParallelAxis-1} 
\left\| K - k_N'\right\|^2 = \left\| K - k_N\right\|^2 + \left\| k_N - k_N'\right\|^2 
\end{equation} 
and 
\begin{equation}\label{ParallelAxis-2} 
\| K\|^2 = \left\| K - k_N\right\|^2 + \left\| k_N\right\|^2 \;. 
\end{equation} 
\end{lemma} 

\begin{proof} 
Put $f(t) := \| K - k_N + t(k_N - k_N') \|^2\,$ ($t\in{\mathbb R}$). 
By the relevant definition in Section~\mbox{6.1}, this function $f(t)$ is a polynomial with leading 
term $ \| k_N - k_N' \|^2 t^2$. 
The specification of $\mu_{i,j}\in{\mathbb R}$ in \eqref{mu_ijDef} is such that each one of the integrals 
$\iint_{{\mathcal R}_{i,j}} (K(x,y) - \mu_{i,j})^2 dx dy\,$ ($1\leq i,j\leq N$) 
is minimised. Therefore, recalling (see Sections~\mbox{3.1} and~\mbox{5.2}) the definition of  
$k_N (x,y)$ and the (similar) form of the kernel $k_N' (x,y)$, one can deduce that $f(t)$ attains its global minimum at the point $t=0$, and that  
$f'(0)$ must therefore equal zero. It follows that $f(t) = f(0) +  \| k_N - k_N' \|^2 t^2\,$ ($t\in{\mathbb R}$), and so 
we have in particular the equation $f(1) = f(0) +  \| k_N - k_N' \|^2$, which is the result \eqref{ParallelAxis-1}. 
\par 
The result \eqref{ParallelAxis-2} is just the 
special case of \eqref{ParallelAxis-1} in which $k_N'(x,y)$ equals zero everywhere in $[0,1]\times [0,1]$: note 
that we are free to choose the coefficients $H_1',\ldots ,H_{2N - 1}'\in{\mathbb R}$ that occur in our specification 
\eqref{k_N'Def} of $k_N'$, and so could choose to have $H_n' = 0\,$ ($1\leq n\leq 2N - 1$). 
\end{proof} 

\begin{remarks} 

\item{\it 1)}\quad Recall the definitions \eqref{psi_iDef} and \eqref{epsilonDef}--\eqref{x_iDef}:  
since these ensure that $\psi_1,\ldots , \psi_N$  is an orthonormal system in $L^2 [0,1]$,   
one can deduce from the definitions \eqref{Def-k_N} and \eqref{k_N'Def} 
that $\| k_N\| = \| H(N) \|$, $\| k_N'\| = \| H'(N) \|$ and $\| k_N - k_N' \| = \| H(N) - H'(N) \|$, where $H(N)$ and $H'(N)$ 
are the matrices mentioned in Section~\mbox{5.2}. It follows (since $H(N)$ and $H'(N)$ are Hankel matrices) that one has  
$\| k_N\|^2 = \sum_{n=1}^{2N - 1} \min\{ n , 2N -n\} H_n^2$, 
$\| k_N'\|^2 = \sum_{n=1}^{2N - 1} \min\{ n , 2N -n\} ( H_n')^2$ and  
$\| k_N - k_N' \|^2 =  \sum_{n=1}^{2N - 1} \min\{ n , 2N -n\} ( H_n - H_n')^2$.  

\item{\it 2)}\quad In view of \eqref{KopHS} and \eqref{DefK}, the number $\| K\|$ equals 
the Hilbert-Schmidt norm of the operator $B_K$ defined in Section~\mbox{4.1}. 
Thus, by \eqref{GammaHSnorm} and \eqref{Gamma_h-HSnormEval}, one has 
$\| K \|^2 = \frac14 - 2\zeta'(0) + \zeta''(0) = 0.0815206105007606323505594460\ldots\ $. 

\item{\it 3)}\quad By \eqref{ParallelAxis-1} and \eqref{ParallelAxis-2}, one has 
\begin{equation*} 
\left\| K - k_N'\right\|^2 = \| K\|^2 - \left( \left\| k_N\right\|^2 - \left\| k_N - k_N'\right\|^2\right) .
\end{equation*}
Since 
\begin{align*} 
\left\| k_N\right\|^2 - \left\| k_N - k_N'\right\|^2 &= \sum_{1\leq n < 2N} \min\{ n , 2N -n\} \left( H_n^2 - ( H_n - H_n')^2\right) \\
 &=\sum_{1\leq n < 2N} \min\{ n , 2N -n\} \left( \left( H_n'\right)^2 + 2 H_n' \left( H_n - H_n'\right)\right) \;, 
\end{align*}
one can deduce that    
\begin{equation}\label{HSnormBound} 
\left\| K - k_N'\right\|^2 \leq \| K\|^2 - \left\| k_N'\right\|^2   
+ 2 \sum_{1\leq n < 2N} \min\{ n , 2N -n\} \left| H_n'\right|\cdot \left| H_n - H_n'\right| \;. 
\end{equation} 
These observations (combined with those in Points~\mbox{(1)}  and~\mbox{(2)}, above) enable  
a numerical upper bound for $\| K - k_N'\|$ to be computed, once 
$N$ and $k_N'$ have been specified and upper bounds for the numbers 
$|H_n - H_n'|\,$ ($1\leq n < 2N$) have been obtained: see Section~\mbox{B.2} for details.

\end{remarks} 

We now state and prove our main result in this subsection. 

\begin{lemma} 
Let $\ell$ and $m$ be positive integers satisfying $\ell\leq m$. For $p>0$, put 
$S_p = S_p (k_N';\ell,m)  := \sum_{n=\ell}^m \left| \varkappa_n (k_N')\right|^p$. 
Then, provided that $S_1 \neq 0$, one has:  
\begin{equation}\label{ELBlemma} 
\left| \varkappa_m \right| \leq \frac{S_2 + {\textstyle\frac12} \left\| K - k_N'\right\|^2 
+ {\textstyle\frac12} (m - \ell) \left\| k_N - k_N'\right\|^2}{S_1}  + \left\| k_N - k_N'\right\| \;.  
\end{equation} 
\end{lemma} 
    
\begin{proof} 
By \eqref{RELBlemma-v2}, one has 
\begin{equation}\label{BWf1} 
\left| \varkappa_n\right| = \left| \varkappa_n (k_N) \right| + \delta_n (k_N) 
\qquad\text{($n\in{\mathbb N}$)} , 
\end{equation} 
with real numbers $\delta_1 (k_N), \delta_2 (k_N), \ldots \ $ satisfying: 
\begin{equation}\label{BWf2} 
\delta_n (k_N) \geq 0\qquad \text{($n\in{\mathbb N}$)} . 
\end{equation} 
By this and \eqref{TraceLemma-k}  (applied firstly for $k=K$, and then for $k=k_N$) we get: 
\begin{align*} 
\| K \|^2 &= \sum_{n=1}^{\infty} \left( \left| \varkappa_n (k_N) \right| + \delta_n (k_N) \right)^2 \\ 
 &\geq \sum_{n=1}^{\infty} \left( \varkappa_n (k_N) \right)^2 
 + \sum_{n=\ell}^m \left( 2  \left| \varkappa_n (k_N) \right|  + \delta_n (k_N)\right)  \delta_n (k_N) \\     
 &= \left\| k_N \right\|^2  
 + \sum_{n=\ell}^m \left( 2  \left| \varkappa_n (k_N) \right|  + \delta_n (k_N)\right)  \delta_n (k_N) \;. 
\end{align*} 
It follows from this and \eqref{ParallelAxis-2} that one has 
\begin{equation}\label{BWf3} 
\sum_{n=\ell}^m \left( \left| \varkappa_n (k_N) \right|  + {\textstyle \frac12}\delta_n (k_N)\right)  \delta_n (k_N) 
\leq {\textstyle \frac12}\left\| K - k_N \right\|^2  \;. 
\end{equation} 
\par 
By \eqref{Weyl-3-unsigned} and \eqref{BWf2}, one has 
\begin{multline*} 
 \left( \left| \varkappa_n (k_N) \right|  + {\textstyle \frac12}\delta_n (k_N)\right)  \delta_n (k_N)  \\ 
\begin{aligned}
&\geq   \left| \varkappa_n (k_N') \right| \delta_n (k_N)   
+  {\textstyle \frac12}\left( \delta_n (k_N) - 2 \left\| k_N - k_N'\right\| \right)  \delta_n (k_N)  \\ 
&=  \left| \varkappa_n (k_N') \right| \delta_n (k_N)   
+ {\textstyle \frac12}\left( \delta_n (k_N) - \left\| k_N - k_N'\right\| \right)^2 
- {\textstyle \frac12}\left\| k_N - k_N'\right\|^2 
\end{aligned} 
\end{multline*} 
for $n\in{\mathbb N}$. 
By this and \eqref{BWf3}, it follows (given that $x^2\geq 0$ for $x\in{\mathbb R}$) that one must have 
\begin{equation} \label{BWf4} 
{\textstyle \frac12}\left\| K - k_N \right\|^2  
\geq \sum_{n=\ell}^m\left| \varkappa_n (k_N') \right| \delta_n (k_N)   
- {\textstyle \frac12}(m + 1  - \ell) \left\| k_N - k_N'\right\|^2 
\end{equation} 
\par
We observe that, by \eqref{BWf1}, \eqref{Def-nu_n-unsigned}  
and our conventions on the ordering of eigenvalues (for which see Section~\mbox{5.1}), one has 
$|\varkappa_n (k_N)| + \delta_n (k_N) = |\varkappa_n| \geq |\varkappa_m|$ when $1\leq n\leq m$. 
By this and \eqref{Weyl-3-unsigned}, it follows that we have  
\begin{equation*} 
\delta_n (k_N) \geq  |\varkappa_m| -  |\varkappa_n (k_N)| 
\geq   |\varkappa_m| -  |\varkappa_n (k_N')|   - \| k_N - k_N'\|  
\end{equation*} 
when $1\leq n\leq m$. Using this to bound the sum in \eqref{BWf4}, 
one gets:   
\begin{equation*} 
{\textstyle \frac12}\left\| K - k_N \right\|^2  
\geq \left( \left| \varkappa_m\right|  -  \| k_N - k_N'\| \right) S_1  - S_2 
- {\textstyle \frac12}(m + 1  - \ell) \left\| k_N - k_N'\right\|^2 \;. 
\end{equation*} 
The required result  \eqref{ELBlemma} follows directly from this and \eqref{ParallelAxis-1}. 
\end{proof} 

\begin{remarks} 
Let $m$ and $N$ be integers, with $N\geq 3$. 
From \eqref{Weyl-2-unsigned} and the case $\ell = m$ of \eqref{ELBlemma}, we obtain 
the bounds 
\begin{equation*} 
|\varkappa_m| \leq |\varkappa_m (k_N')| + \| K - k_N' \|
\end{equation*}
and 
\begin{equation*} 
|\varkappa_m| \leq  |\varkappa_m (k_N')| + \frac{\| K - k_N' \|^2}{2  |\varkappa_m (k_N')|} + \| k_N - k_N' \|\;,
\end{equation*} 
respectively. 
The latter bound is an improvement upon the former one if one has 
\begin{equation*} 
 |\varkappa_m (k_N')| > {\textstyle\frac12} \| K - k_N' \| \cdot \left( 1 - \frac{\| k_N - k_N' \|}{ \| K - k_N' \| }\right)^{-1} 
\end{equation*} 
(note that, by \eqref{ParallelAxis-1} and the inequality $\| K - k_N\| >0$, 
we have $\| K - k_N' \| > \| k_N - k_N' \|$ here). 
\end{remarks} 

\subsection{Computations and numerical Results}

The computations described in this subsection utilise only data associated with the matrix $H'(2^{21})$ 
and the corresponding kernel: $k_N'$ for $N=2^{21}$. 
Thus, in what follows, it is to be understood that $n=21$, $N=2^n$ and $M=384$ when there is no indication to the contrary.
\par
Recall (from Section~\mbox{5.3}) that we have computed, for $N=2^{21}$, $\pm\in\{ +, -\}$ and $1\leq m\leq M=384$,  
certain binary64 numbers 
$L^{\pm}_m, U^{\pm}_m >0$ with $L^{\pm}_m \leq |\varkappa^{\pm}_m (k_N')|\leq U^{\pm}_m\,$ 
(so that the relations in \eqref{n=21intervals} hold). Decimal approximations to some of these numbers 
can be seen in Table~\mbox{D-3}, which is in Section~\mbox{D.7}. We also have at our disposal 
certain numerical upper bounds, $E=E(N)$ and $F=F(N)$, for $\| k_N' - k_N\|$ and $\| k_N' - K\|$, respectively: these 
bounds were obtained via computations based on points noted in Remarks~\mbox{6.11} (for a little more detail 
see Remarks~\mbox{B.8}~\mbox{(4)} and Table~\mbox{B-1}, which lists decimal approximations to $E(N)$ and $F(N)$ in 
the cases where $N\in\{ 2^4, 2^5, \ldots , 2^{21}\}$). 
\par
The above mentioned bounds for the numbers $|\varkappa^{\pm}_m (k_N')|\,$ ($\pm\in\{ +, -\}$ and $1\leq m\leq M$)  
immediately imply similar bounds for $|\varkappa_1 (k_N')|,\ldots , |\varkappa_{2M} (k_N')|$.
Indeed, with $W_j := W_0 :=\max\{ U_M^{+} , U_M^{-}\}$ for $j\geq 0$, one has 
\begin{equation}\label{WIDE-1}
L_m'\leq |\varkappa_m (k_N')|\leq U_m' \quad\text{($1\leq m\leq 2M$)} , 
\end{equation}
where $L_m'\,$ (resp. $U_m'$) is the $m$-th largest term of the sequence 
$L^{+}_1,\ldots ,L^{+}_M,  L^{-}_1,\ldots L^{-}_M\,$ 
(resp. $U^{+}_1, \ldots ,U^{+}_M, U^{-}_1,\ldots ,U^{-}_M, W_0, W_1,\ldots , W_{2M-2}$). 
We used Octave's {\tt sort()} function to compute the relevant sequences 
$L_1',\ldots ,L_{2M}'$ and $U_1',\ldots ,U_{2M}'$. 
\par
With the help of the {\tt interval} package we computed 
sharp lower bounds ${\mathcal L}_1,\ldots ,{\mathcal L}_{2M}$ for the differences $L_1' - E,\ldots ,L_{2M}' - E$. 
By \eqref{Weyl-4} and \eqref{WIDE-1}, the numbers ${\mathcal L}_1,\ldots ,{\mathcal L}_{2M}$ 
are such that 
\begin{equation}\label{WIDE-2}
{\mathcal L}_m \leq |\varkappa_m|\quad\text{($1\leq m\leq 2M$)} . 
\end{equation}
Numbers ${\mathcal L}^{+}_1\ldots , {\mathcal L}^{+}_M$ and ${\mathcal L}^{-}_1\ldots , {\mathcal L}^{-}_M$ satisfying 
\begin{equation}\label{WIDE-3}
{\mathcal L}^{\pm}_m \leq |\varkappa^{\pm}_m|\quad\text{($\pm\in\{ + , -\}$ and $1\leq m\leq M$)}
\end{equation}
were computed in a similar fashion: ${\mathcal L}^{\pm}_m$ being, in each case, a sharp numerical lower bound for the difference 
$L^{\pm}_m - E$. Lemma~\mbox{6.8} and~\eqref{LBs_for_eigenmoduli} provide the necessary justification for this. 
We remark that ${\mathcal L}^{\pm}_m$ is, in all relevant cases, a positive number. 
Unsurprisingly, for $1\leq m\leq 2M$, the $m$-th largest term of the sequence 
${\mathcal L}^{+}_1\ldots , {\mathcal L}^{+}_M, {\mathcal L}^{-}_1\ldots , {\mathcal L}^{-}_M$ 
equals ${\mathcal L}_m$.
\par 
Finally, using an algorithm based on Lemma~\mbox{6.12}, we computed numbers 
${\mathcal U}_1,\ldots ,{\mathcal U}_{2M}$ satisfying 
\begin{equation}\label{WIDE-4}
|\varkappa_m|\leq {\mathcal U}_m\quad\text{($1\leq m\leq 2M$)} .
\end{equation}
As input data, for this computation, 
we used the above mentioned  
numerical bounds $E\geq \| k_N' - k_N\|$ and $F\geq \| k_N' - K\|$, and   
the numerical bounds \eqref{WIDE-1} 
for $|\varkappa_1 (k_N')|,\ldots , |\varkappa_{2M} (k_N')|$. 
Note that our algorithm initially generates $m$ upper bounds for $|\varkappa_m|$:  
one bound, ${\mathcal U}_{m,\ell}\in (0, +\infty]$, for each $\ell\in{\mathbb N}$ satisfying the condition $\ell\leq m$ of Lemma~\mbox{6.12}. 
It provisionally assigns to ${\mathcal U}_m$ the value $\min_{\ell\leq m} {\mathcal U}_{m,\ell}\,$ 
(see below regarding final adjustments made to the value of ${\mathcal U}_m$). 
The final step of the algorithm is to 
recursively assign  to ${\mathcal U}_m$ the value $\min\{ {\mathcal U}_m , {\mathcal U}_{m-1}\}\,$ 
(with $m$ here going from $2$ to $2M$, in increments of $1$). This is justified 
by virtue of \eqref{lambdaOrder} and \eqref{Def-nu_n-unsigned}, 
and ensures that we 
end up having ${\mathcal U}_1\geq {\mathcal U}_2\geq \ldots \geq {\mathcal U}_{2M}>0$. 
Where $\min_{\ell\leq m} {\mathcal U}_{m,\ell} < +\infty$, 
we use $\ell_m$ as our notation for the associated (optimal) value of $\ell\in\{ 1,\ldots ,m\}$.  
\par
In order to ensure the validity of the numerical bounds \eqref{WIDE-4}  
all of the relevant non-integer calculations were carried out with the aid of Octave's {\tt interval} package. 
An important part of these calculations is concerned with the computation of a reasonably sharp upper bound 
for the maximum of 
$\phi(x_{\ell},\ldots , x_m) := (\frac12 (m-\ell) E^2 + \frac12 F^2 + \sum_{n=\ell}^m x_n^2) / \sum_{n=\ell}^m x_n$ 
subject to the constraints $L_n'\leq x_n\leq U_n'\,$ ($\ell\leq n\leq m$). 
One can show that this maximum, $\Phi_{\ell, m}\,$ (say), must equal $\phi (U_{\ell}',\ldots ,U_m')$ if the latter number  
is less than $2\min\{ L_{\ell}',\ldots , L_m'\} = 2 L_m'$. 
Our algorithm exploits this fact whenever there is a clear opportunity to do so:   
but if it fails to verify that $\phi (U_{\ell}',\ldots ,U_m') < 2 L_m'$  
then it simply puts ${\mathcal U}_{m,\ell}=+\infty$. 
\par
Reformulating \eqref{WIDE-2}--\eqref{WIDE-4} in terms of eigenvalues of $K$, one obtains: 
\begin{equation}\label{WIDE-5} 
1/{\mathcal U}_m\leq |\lambda_m|\leq 1/{\mathcal L}_m\quad\text{($1\leq m\leq 2M$)} 
\end{equation}
and 
\begin{equation}\label{WIDE-6}
|\lambda^{\pm}_m|\leq 1/{\mathcal L}^{\pm}_m\quad\text{($\pm\in\{ + , -\}$ and $1\leq m\leq M$)} .
\end{equation}
A selection of these bounds, rounded to 6 significant digits, are shown in Tables~1 and~2 (below). 

\begin{remarks}

\item{1)}\quad 
The value of $\ell_m$ is defined only when $1\leq m\leq 268$, since  $\min_{\ell\leq m} {\mathcal U}_{m,\ell} = +\infty$ 
for $269\leq m\leq 2M$. 
The final step of our algorithm has no effect on the values of 
${\mathcal U}_1,\ldots , {\mathcal U}_{268}\,$ 
(it turns out that the function 
$m\mapsto \min_{\ell\leq m} {\mathcal U}_{m,\ell}$ is strictly decreasing for $m\leq 268$); 
its effect on ${\mathcal U}_{269},\ldots , {\mathcal U}_{2M}$ is that 
we end up having ${\mathcal U}_m = {\mathcal U}_{268}$ when $268\leq m\leq 2M=768$.

\item{2)}\quad As a rough guide to the strength of the bounds \eqref{WIDE-5}, 
we note that one has:
\begin{equation*}
\frac{{\mathcal U}_m}{{\mathcal L}_m} \leq 1 + 
\begin{cases} 
1 & \text{if}\ m\leq 268\,; \\
\frac18 & \text{if}\ m\leq 118\,; \\
\frac{1}{32} & \text{if}\ m\leq 10\,; \\
\frac{1}{128} & \text{if}\ m\in\{ 1, 3\} \,.
\end{cases}
\end{equation*} 

\item{3)}\quad Given our definition of the relevant notation (for which see Section~\mbox{5.1}), it follows immediately from \eqref{lambdaOrder} 
that $|\varkappa^{\pm}_m|\leq |\varkappa_m|$ for $\pm\in\{ +,-\}$ and $m\in{\mathbb N}$. 
Thus we may certainly deduce from \eqref{WIDE-4} that $|\varkappa^{\pm}_m|\leq {\mathcal U}_m$  
when $\pm\in\{ +,-\}$ and $m\leq M$. Observe, furthermore, that if $m,r\in\{ 1,\ldots ,M\}$ are such that ${\mathcal U}_{m+r} < {\mathcal L}^{-}_r$ 
then one must have $\varkappa^{+}_m \leq {\mathcal U}_{m+r}\,$ 
(else, by \eqref{WIDE-3}, one would have ${\mathcal U}_{m+r} < \min\{ |\varkappa^{-}_1|, \ldots , |\varkappa^{-}_r|, |\varkappa^{+}_1|, \ldots , |\varkappa^{+}_m|\}$, 
and so ${\mathcal U}_{m+r} < |\varkappa_{m+r}|$, which would contradict \eqref{WIDE-4}). Therefore, when $m\leq M$, we have 
\begin{equation*}
\varkappa^{+}_m \leq {\mathcal U}_{m+\rho(m,+)}\;,
\end{equation*}
where  $\rho(m,+)$ is defined to be the maximum of the set $\{ 0\}\cup\{ r\in{\mathbb N} : r\leq M\ \text{and}\ {\mathcal U}_{m+r} < {\mathcal L}^{-}_r\}$. 
Using this bound for $\varkappa^{+}_m\,$ ($m\leq M$), and the analogous bounds for $|\varkappa^{-}_1|,\ldots ,|\varkappa^{-}_M|$, 
we determine sets 
$\{{\mathcal U}^{+}_1\ldots , {\mathcal U}^{+}_M\}, \{{\mathcal U}^{-}_1\ldots , {\mathcal U}^{-}_M\}\subset\{ {\mathcal U}_1\ldots , {\mathcal U}_{2M}\}$ such that  
\begin{equation}\label{WIDE-7}
 |\varkappa^{\pm}_m|\leq {\mathcal U}^{\pm}_m \quad\text{($\pm\in\{ + , -\}$ and $1\leq m\leq M$)} .
\end{equation}
Note that it follows from this and \eqref{WIDE-3} that one has 
\begin{equation*}
1 / {\mathcal U}^{\pm}_m \leq |\lambda^{\pm}_m| \leq 1 / {\mathcal L}^{\pm}_m \quad\text{($\pm\in\{ + , -\}$ and $1\leq m\leq M$)} .
\end{equation*}
A selection of these bounds are included in Table~2.

\item{4)} \quad Let ${\mathcal J} (\varkappa) := \{ j\in{\mathbb N} : |\varkappa_j| = |\varkappa|\}\,$  ($\varkappa\in{\mathbb R}$). 
By computing sets of the form $\{  j\leq 2M :  [{\mathcal L}_j, {\mathcal U}_j]\cap [{\mathcal L}^{\pm}_m, {\mathcal U}^{\pm}_m]\neq \emptyset\}$ 
one can determine real intervals ${\mathcal I}^{+}_1,{\mathcal I}^{-}_1,\ldots ,{\mathcal I}^{+}_M,{\mathcal I}^{-}_M$ 
with ${\mathcal I}^{\pm}_m \supseteq{\mathcal J}(\varkappa^{\pm}_m)$ 
for $\pm\in\{ -1 , 1\}$ and $1\leq m\leq M$ (some are shown in Table~2). 
Since the sequence $|\varkappa_1|,|\varkappa_2|,\ldots$ is monotonic decreasing, 
one must have $|\varkappa| > |\varkappa'|$ whenever $\max{\mathcal J}(\varkappa) < \min{\mathcal J}(\varkappa')$. 
Thus we find (for example) that, based on the sets ${\mathcal I}^{+}_1,{\mathcal I}^{-}_1,\ldots ,{\mathcal I}^{+}_8,{\mathcal I}^{-}_8$ 
shown in Table~2, one certainly has: 
\begin{align*}
|\varkappa^{+}_m| &> |\varkappa^{-}_m|\quad\text{when $1\leq m\leq 7$ and $m\neq 6$} , \\
|\varkappa^{-}_m| &> |\varkappa^{+}_{m+1}|\quad\text{when $2\leq m\leq 7$ and $m\neq 4$} ,
\end{align*}
and 
\begin{equation}\label{WIDE-8}
|\varkappa^{\pm}_m|  > |\varkappa^{\pm}_{m+1}| \quad\text{when $\pm\in\{ +, -\}$ and $1\leq m\leq 7$} . 
\end{equation} 
Since we have ${\mathcal I}^{+}_1 = [1, 1] = \{ 1\}$, and so also ${\mathcal J}(\varkappa^{+}_1) = \{ 1\}$, 
it is clear that $\varkappa_1 = \varkappa^{+}_1$. We find, similarly, that $\varkappa_4 = \varkappa^{-}_2$, 
$\varkappa_5 = \varkappa^{+}_3$, $\varkappa_6 = \varkappa^{-}_3$, $\varkappa_7 = \varkappa^{+}_3$, $\varkappa_{10} = \varkappa^{-}_5$, 
$\varkappa_{13} = \varkappa^{+}_7$ and $\varkappa_{14} = \varkappa^{-}_7$. 
Since ${\mathcal I}^{\pm}_m\subset [1,14]$ for $\pm\in\{ +,-\}$ and $1\leq m\leq 7$, 
we have also: $\{ \varkappa_m : m\leq 14\} = \cup_{m=1}^7 \{ \varkappa^{+}_m , \varkappa^{-}_m\}$. 
By this and \eqref{WIDE-8}, each one of $\lambda_1,\ldots ,\lambda_{14}$ is a simple eigenvalue of $K$, 
and one has $|\lambda_{15}| > |\lambda_{14}|\,$ (as one can check by glancing at Table~1).
\par
By further inspection of the data in Table~2, we find that the sets ${\mathcal J}(\varkappa^{+}_m)\cap\{ 25, 26\}$ 
and ${\mathcal J}(\varkappa^{-}_m)\cap \{ 26\}$ are empty unless $m=13$, 
and that ${\mathcal J}(\varkappa^{-}_m)\ni 25$ only if $m\in\{ 12, 13\}$.   
It follows that the eigenvalue $\lambda_{26}$ is simple, and that either $\lambda_{26} = \lambda^{-}_{13}$ 
and $\lambda_{25}\in\{ \lambda^{+}_{13}, \lambda^{-}_{12}\}$, or else 
$\lambda_{26} = \lambda^{+}_{13}$ 
and $\lambda_{25}\in\{ \lambda^{-}_{13}, \lambda^{-}_{12}\}\,$ (in the former case $\lambda_{25}$ is simple; in the latter it is of multiplicity 
$1$ or $2$). 
\par 
More generally, our best upper bounds for multiplicities were obtained with the help of the relation \eqref{Balanced_Indices}. 
In cases where ${\mathcal L}_m > {\mathcal U}_{2M}\,$ (so that one has  
${\mathcal J}(\varkappa_m) \subseteq \{ 1,\ldots ,2M-1\}$), it follows from \eqref{Balanced_Indices}, \eqref{WIDE-2} and \eqref{WIDE-4} 
that if $t_m$ is the greatest odd integer $t$ such that there 
exists an integer $s$ with $1\leq s\leq m\leq s+t-1\leq 2M$ and $\cap_{j=s}^{s+t-1} [{\mathcal L}_j, {\mathcal U}_j] \neq \emptyset$ 
then the number $\mu_m := (t_m + 1)/2$ is an upper bound for the multiplicity of $\lambda_m$.   
We computed this bound, $\mu_m$, for each positive integer $m\leq 98\,$ 
(having found that ${\mathcal L}_{98} > {\mathcal U}_{2M} > {\mathcal L}_{99}$): for a selection of the results thereby obtained, 
see the rightmost column of Table~1. We find, in particular, that $\mu_m = 1$ when either $1\leq m\leq 14$ or $m\in\{ 25, 26\}\,$ 
(which mostly just confirms results already noted, though it does establish that $\lambda_{25}$ is a simple eigenvalue). 
We find also that $\mu_m\leq 3$ whenever $m\leq 46$, and that $\mu_m\leq 4$ whenever $m\leq 65$. 

\begin{table}[bh]
\centering 
\begin{minipage}{55mm}
\resizebox{50mm}{!}{ 
\begin{tabular}[c]{|r|r|c|c|c|r|} 
\hline 
$m$ & $\ell_m$  & $1/{\mathcal U}_m$ & $1/{\mathcal L}_m$  & $\mu_m$ \\ \hline  
1 & 1 &  12.4612 & 12.5488  & 1 \\ \hline  
2 & 2  & 14.3515 & 14.4860  & 1 \\ \hline  
3 & 2 & 14.4201 & 14.4893  & 1 \\ \hline  
4 & 4  & 17.4022 & 17.6441  & 1 \\ \hline  
5  & 5  & 17.8465 & 18.1078  & 1 \\ \hline  
6 & 6  & 19.4178 & 19.7562  & 1 \\ \hline  
7 & 7  & 21.3352 & 21.7875  & 1 \\ \hline  
8 & 8  & 22.1484 & 22.6561  & 1 \\ \hline  
9  & 8  & 22.4261 & 22.7117  & 1 \\ \hline  
10 & 10  & 23.3777 & 23.9780  & 1 \\ \hline  
11 & 11  & 25.6811 & 26.4856  & 1 \\ \hline  
12 & 11  & 26.3138 & 26.9915  & 1 \\ \hline  
13  & 13  & 28.4755 & 29.5888  & 1 \\ \hline  
14 & 14  & 30.2660 & 31.6172  & 1 \\ \hline  
15 & 15  & 32.9808 & 34.7604  & 2 \\ \hline  
16 & 15  & 34.1616 & 35.4565  & 3 \\ \hline  
24 & 22  & 36.3086 & 37.9903  & 2 \\ \hline  
25 & 24  & 37.8582 & 40.4782  & 1 \\ \hline  
26 & 25  & 39.6664 & 41.9071  & 1 \\ \hline  
27 & 27  & 42.2006 & 46.2279  & 3 \\ \hline  
32  & 28  & 46.5603 & 48.1481  & 3 \\ \hline  
48 & 44 & 60.5511 & 63.9656 & 4 \\ \hline 
64 & 59  & 71.1528 & 77.1839  & 4 \\ \hline  
98 & 89  & 93.6489 & 101.382  & 45 \\ \hline  
99 & 92  & 94.4381 & 103.037  & - \\ \hline  
128 & 95  & 100.448 & 118.324  & - \\ \hline  
268 & 99  & 101.637 & 202.768  & - \\ \hline 
269 & -  & 101.637 & 203.448  & - \\ \hline 
512 & -  & 101.637 & 323.956  & - \\ \hline 
768 & -  & 101.637 & 446.298  & - \\ 
\hline 
\end{tabular} 
}
\vskip 3mm 
\ Table~1  
\end{minipage}
\begin{minipage}{95mm}
\resizebox{90mm}{!}{ 
\begin{tabular}[c]{|r|c|c|c|c|c|c|} 
\hline 
$m$  &   $1/{\mathcal U}_m^{+}$  & $1/{\mathcal L}_m^{+}$ & ${\mathcal I}^{+}_m$ 
&  $1/{\mathcal U}_m^{-}$ & $1/{\mathcal L}_m^{-}$ & ${\mathcal I}^{-}_m$  \\ \hline 
1 & 12.4612 & 12.5488 & [1,1] & 14.3515 & 14.4893  & [2,3] \\ \hline
2 & 14.3515 & 14.4860 & [2,3] & 17.4022 & 17.6441  & [4,4] \\ \hline
3 & 17.8465 & 18.1078 & [5,5] & 19.4178 & 19.7562  & [6,6] \\ \hline
4 & 21.3352 & 21.7875 & [7,7] & 22.1484 & 22.6561  & [8,9] \\ \hline
5 & 22.1484 & 22.7117 & [8,9] & 23.3777 & 23.9780  & [10,10] \\ \hline
6 & 25.6811 & 26.4856 & [11,12] & 25.6811 & 26.9915  & [11,12] \\ \hline
7 & 28.4755 & 29.5888 & [13,13] & 30.2660 & 31.6172  & [14,14] \\ \hline
8 & 32.9808 & 34.7604 & [15,17] & 32.9808 & 35.4838  & [15,21] \\ \hline
9 & 34.1616 & 35.4565 & [15,20] & 34.1616 & 35.7961  & [15,22] \\ \hline
10 & 34.5881 & 35.7926 & [15,22] & 34.9184 & 36.1243  & [16,23] \\ \hline
11 & 34.9184 & 36.1987 & [16,23] & 35.1328 & 36.3612  & [16,24] \\ \hline
12 & 35.1328 & 36.8747 & [16,24] & 35.4676 & 37.9903  & [17,25] \\ \hline
13 & 35.3246 & 40.4782 & [16,26] & 37.8582 & 41.9071  & [24,26] \\ \hline
14 & 42.2006 & 46.2279 & [27,31] & 42.2006 & 46.7331  & [27,32] \\ \hline
15 & 44.3390 & 46.9532 & [27,33] & 44.3390 & 47.8324  & [27,34] \\ \hline
16 & 45.1717 & 47.8932 & [27,34] & 45.1717 & 48.1481  & [27,35] \\ \hline
24 & 57.3340 & 63.9656 & [41,54] & 57.3340 & 63.9162  & [41,54] \\ \hline
32 & 65.9454 & 76.7727 & [54,70] & 66.8455 & 77.1839  & [54,71] \\ \hline
48 & 87.6573 & 100.770 & [79,137] & 87.6573 & 101.081  & [79,153] \\ \hline
49 & 88.1822 & 101.106 & [79,155] & 88.1822 & 101.382 & [79,186] \\ \hline
50 & 88.6548 & 103.037 & $[79,\infty)$ & 88.6548 & 103.064  & $[79,\infty)$ \\ \hline
64 & 98.8998 & 118.264 & $[94,\infty)$ & 99.0621 & 118.324  & $[94,\infty)$ \\ \hline
96 & 100.895 & 156.912 & $[96,\infty)$ & 100.915 & 156.717  & $[96,\infty)$ \\ \hline
128 & 101.313 & 197.050 & $[98,\infty)$ & 101.321 & 196.603  & $[98,\infty)$ \\ \hline
192 & 101.585 & 264.853 & $[99,\infty)$ & 101.585 & 264.711  & $[99,\infty)$ \\ \hline
256 & 101.637 & 323.956 & $[99,\infty)$ & 101.637 & 320.259  & $[99,\infty)$ \\ \hline
384 & 101.637 & 446.298 & $[99,\infty)$ & 101.637 & 445.246  & $[99,\infty)$ \\ 
\hline 
\end{tabular} 
} 
\vskip 3mm 
\ Table~2 
\end{minipage}
\end{table}

\item{5)}\quad  The calculations outlined in this subsection can easily be adapted so as to 
make use of the $2M$ `probable upper bounds' $\tilde U^{\pm}_1,\ldots , \tilde U^{\pm}_m$ that are mentioned 
in Section~\mbox{5.3} (in connection with a certain number $\tilde R =\tilde R(N)$, 
which we associate with the $N\times N$ matrix $A$ given by  
\eqref{Def-H''} and \eqref{Def-A(N)}). 
This yields numbers  
$\tilde {\mathcal U}_1, \ldots , \tilde {\mathcal U}_{2M}$ and $\tilde {\mathcal U}^{+}_1,\tilde {\mathcal U}^{-}_1,\ldots , \tilde {\mathcal U}^{+}_M,\tilde {\mathcal U}^{-}_M$  
which, if $\tilde R\geq \| A\|_2$, 
are such that 
$\tilde {\mathcal U}_m \geq |\varkappa_m|\,$ ($1\leq m\leq 2M$) and 
$\tilde {\mathcal U}^{\pm}_m \geq |\varkappa^{\pm}_m|\,$ ($\pm\in\{ +,-\}$ and $1\leq m\leq M$).
Provided that $\tilde R\geq \| A\|_2$, one has both 
$[{\mathcal L}_m , \tilde {\mathcal U}_m] \ni |\varkappa_m|\,$ ($1\leq m\leq 2M$)
and $[{\mathcal L}^{\pm}_m , \tilde {\mathcal U}^{\pm}_m] \ni |\varkappa^{\pm}_m|\,$ ($\pm\in\{ +,-\}$ and $1\leq m\leq M$). 
The latter set of results are deduced from the former set, through calculations similar to those indicated in Point~\mbox{(3)}, above.  
This leads to some loss of precision; we find, nevertheless, that
\begin{equation*} 
\log\bigl( \tilde {\mathcal U}^{\pm}_m / {\mathcal L}^{\pm}_m \bigr) 
\leq {\textstyle \frac52} \log\bigl( \tilde {\mathcal U}_{2m} / {\mathcal L}_{2m} \bigr) 
\quad \text{for $\pm\in\{ + , -\}$ and $1\leq m\leq M$} .
\end{equation*}
\par
When $m\geq 120\,$   (say) the difference between the conditional bound $\tilde {\mathcal U}_m$ and the unconditional bound ${\mathcal U}_m$ 
is significant: we find, for example, that $\tilde {\mathcal U}_m / {\mathcal L}_m \leq 23/20$ for $1\leq m\leq 2M$, and that 
$\tilde {\mathcal U}_m / {\mathcal L}_m \leq 9/8$  whenever $m\leq 301\,$ (compare this with Point~\mbox{(2)}, above).
\par 
We find also that $\tilde {\mathcal U}_{2M} < {\mathcal L}_{650}$. 
Thus, using the numerical data ${\mathcal L}_m , \tilde {\mathcal U}_m\,$ ($1\leq m\leq 2M$), 
in conjunction with the method described in the final paragraph of Point~\mbox{(4)}, above, 
we get conditional upper bounds $\tilde \mu_1, \ldots ,\tilde \mu_{650}$ for the multiplicities 
of the first $650$ eigenvalues of $K$ (i.e. these bounds are valid if $\tilde R\geq \| A\|_2$). 
To give some idea of the strength of these conditional bounds, we mention here 
that 
\begin{equation*}
\left\lfloor m/15\right\rfloor \leq \tilde\mu_m\leq \lceil m/10\rceil \quad \text{for $21\leq m\leq 650$},
\end{equation*}
and that $\tilde\mu_m = \mu_m$ for $m\leq 93$:  Table~1 lists  some of the sequence $\mu_1,\ldots ,\mu_{98}$. 
\par
Our tables omit the data 
$\,\tilde {\mathcal U}_m$, $\tilde {\mathcal U}^{\pm}_m$ 
and $\tilde\mu_m$, which plays no further part in our investigations. 
\end{remarks}

\section{Conjectures about eigenvalues of $K$} 
In this section $M$ denotes just the number $384$, and $N$ is `shorthand' for $2^n$. 
We assume throughout that our `probable upper bounds' 
(discussed in Section~\mbox{5.3}, and, in greater detail, in Sections~\mbox{D.5--D.7}) are valid upper bounds. 
Thus we suppose, in particular, that 
$|\varkappa^{\pm}_j (k_N')|\leq \tilde U^{\pm}_j = \tilde U^{\pm}_j(n)$ 
when $1\leq j\leq M$, $\pm\in\{ + , -\}$ and $10\leq n\leq 21$. 
Complementing this, there are the lower bounds \eqref{LBs_for_eigenmoduli}, the validity of which 
is not in question; note that the relevant numerical data there, $L^{+}_1, L^{-}_1,\ldots ,L^{+}_M, L^{-}_M$, depends on $n$, 
and so (since we shall need to compare results obtained for different choices of $n$) we now switch to 
using the notation $L^{+}_j (n)$ (resp. $L^{-}_j (n)$), in place of just $L^{+}_j$ (resp. $L^{-}_j$):  
for the same reason, we shall write $\tilde U^{+}_j (n)$ (resp. $\tilde U^{-}_j (n)$) for $\tilde U^{+}_j$ (resp. $\tilde U^{-}_j$). 
Our working hypothesis is, then, that 
\begin{equation} \label{W-H}
\big[ L^{+}_j (n) , \tilde U^{+}_j (n)\big] \ni \varkappa^{+}_j (k_N') 
\quad\text{and}\quad 
\big[ -\tilde U^{-}_j (n) , -L^{-}_j (n)\big] \ni \varkappa^{-}_j (k_N') \;,
\end{equation} 
for $1\leq j\leq M$ and $n\in\{ 10,11,\ldots ,21\}$. 
As we shall explain, in Section~\mbox{7.1} (below), this working hypothesis leads us to make certain (quite specific) conjectures
regarding the behaviour, for $n\geq 21$, of the differences 
\begin{equation*} 
\Delta^{\pm}_m (n) := |\varkappa^{\pm}_m (k_{N})| - |\varkappa^{\pm}_m (k_{N/2})|\quad 
\text{($1\leq m\leq M$, $\pm\in\{ +,-\}$)} .
\end{equation*}
We shall also discuss there certain conditional (or `conjectural') lower bounds for the moduli of eigenvalues of $K$ that hold if the latter 
conjectures, and our working hypothesis, are correct. These bounds follow by virtue of the Weyl bound \eqref{Weyl-3} 
and the case $n=21$ of the relations 
\begin{equation}\label{SummedDifferences} 
\varkappa^{\pm}_m = \varkappa^{\pm}_m ( k_N ) \pm \sum_{r=1}^{\infty} \Delta^{\pm}_m (n+r)\quad 
\text{($\pm\in\{ + , -\}$, $m\geq 1$ and $n\geq 2$)} ,
\end{equation} 
which are themselves direct consequences of our having both $|\varkappa^{\pm}_m - \varkappa^{\pm}_m ( k_{2^R N})|\leq \| K - k_{2^R N}\|\,$ 
($\pm\in\{ +,-\}$ and $m,R\in{\mathbb N}$) and  
$\| K - k_{2^R N}\| \rightarrow 0$ as $R\rightarrow +\infty$: 
the former inequalities being obtained similarly to \eqref{Weyl-2}, while the latter result is a corollary of the  proof of Proposition~\mbox{1.2} that we give in Section~3.

\subsection{Conjectural lower bounds for $|\lambda^{\pm}_m|$} 

Recall, from Section~\mbox{6.5}, that  we  computed, for $n=4,5,\ldots ,21$,  
a numerical upper bound $E(N)$ for $\| k_N' - k_N\|$.  
It follows from the Weyl inequality \eqref{Weyl-3} that by increasing the lengths of  
the intervals appearing in \eqref{W-H} by $2E(N)$, while keeping the midpoints of the intervals fixed, 
one obtains intervals containing $\varkappa^{+}_j (k_N)$ and $\varkappa^{-}_j (k_N)$: 
using, in this way, the cases of  \eqref{W-H} allowed us by our working hypothesis  
we compute, for $\pm\in\{ + , -\}$, $1\leq m\leq M$ and $10\leq n\leq 21$, 
an explicit real interval ${\mathcal X}^{\pm}_m (n)$ containing the number $\varkappa^{\pm}_m (k_N)$. 
A further simple computation then yields, for $\pm\in\{ + , -\}$, $1\leq m\leq M$ and $11\leq n\leq 21$, 
an explicit real interval ${\mathcal D}^{\pm}_m (n)$ containing the difference $\Delta^{\pm}_m (n)$. 
We find that nearly all of the latter $22M = 8448$ intervals have a lower end point (i.e. infimum) that is positive:  
the only exceptions to this rule occur when $\pm = +$, $m=M$ and $18\leq n\leq 21$. 
An obvious corollary of this observation is that, when $\pm\in\{ + , -\}$ and $1\leq m\leq M$, 
one has:  
\begin{equation}\label{Happy=Positive}
\Delta^{\pm}_m (n) > 0
\quad\text{for}\ 11\leq n\leq\begin{cases} 17 & \text{if $\{\pm\} = \{ +\}$ and $m = M$} , \\ 21 & \text{otherwise} . \\ \end{cases} 
\end{equation}
In the $4$ exceptional cases we have ${\mathcal D}^{\pm}_m (n)\ni 0$, and so we fail to determine the sign of $\Delta^{\pm}_m (n)$. 
\par
In all cases where $\Delta^{\pm}_m (n) > 0$ we put: 
\begin{equation*} 
Y^{\pm}_m (n) := \log_2 \left( \Delta^{\pm}_m (n)\right) \in{\mathbb R}\;. 
\end{equation*}
When $n$, $m$ and the sign $\pm$ satisfy the conditions attached to \eqref{Happy=Positive}, 
we can use the data ${\mathcal D}^{\pm}_m (n)$ to compute a real interval ${\mathcal Y}^{\pm}_m (n)$ containing $Y^{\pm}_m (n)$. 
Thus, for $\pm\in\{ +, -\}$ and $1\leq m\leq M$, we have the pair of functions 
$n\mapsto \inf {\mathcal Y}^{\pm}_m (n)$ 
and $n\mapsto \sup {\mathcal Y}^{\pm}_m (n)$, whose common domain is the set of integers $n$ satisfying the conditions 
in \eqref{Happy=Positive}. 
Using Octave's {\tt plot()} function to view their graphs, we find in general that these functions are decreasing, in an approximately linear fashion,  
on a significant part of their domain (especially where the independent variable $n$ is near to the top end of its range). 
There are relatively few choices of the pair 
$(\pm, m)$ for which the graphs do not fit this description. Moreover, upon comparing graphs associated with 
distinct choices of the pair $(\pm , m)$, we observe a clear tendency for the graphs to run approximately parallel to one 
another (this is most obvious when the relevant values of $m$ are small, and is, again, more evident where $n$ is near the top of its range). 
\par 
Further examination of the data ${\mathcal Y}^{\pm}_m (n)$ leads us to 
conjecture, tentatively, that there exist absolute constants $\alpha\in (-2.02, -1.98)$ and $\beta\in (3,4)$ such that, 
for $\pm\in\{ +,-\}$ and $m\in{\mathbb N}$, one has   
\begin{equation}\label{ForayHypothesis}
Y^{\pm}_m (n) - Y^{\pm}_m (n-1) = \alpha  + \frac{\beta + o_m (1)}{n}\;, 
\end{equation}
in the limit as $n\rightarrow\infty$. 
We prefer to omit the evidence gathered in support this conjecture,  
since \eqref{ForayHypothesis} serves only to guide subsequent work 
(and is not taken for a working hypothesis). We instead argue that there is evidence that the difference 
\begin{equation*} 
\dot Y^{\pm}_m (n) := Y^{\pm}_m (n) - Y^{\pm}_m (n-1) = \log_2 \left( \Delta^{\pm}_m (n) / \Delta^{\pm}_m (n-1) \right) 
\end{equation*}
behaves (in certain respects) similarly to how it would if the conjecture \eqref{ForayHypothesis} were correct. 
\par 
If it was known that \eqref{ForayHypothesis} holds, with $\alpha\approx -2$ and $\beta\approx 3.5\,$ (say), 
then it would be reasonable to hope that, for some `sufficiently small' values of $m$, 
one would have 
$\max\{ \dot Y^{\pm}_m (n) : n\geq 22\} \leq -1.8$  for $\pm\in\{ + , -\}$: we remark here 
that our numerical data supports, to some extent, the hypothesis that this does hold whenever one has $m\leq 192\,$ (say). 
The upper bound $-1.8$ is (in practice) an overly ambitious target to set ourselves, 
given the limited data that we have to work with,  
especially if we are to make claims concerning what occurs for an arbitrary choice of $m\in\{ 1,\ldots , M\}$. 
We therefore select a set of  three easier targets: ${\mathcal T} := \{ 0, -1, -\omega\}$, 
with 
\begin{equation*}
\omega := \log_2(3) = -1.58496\ldots 
\end{equation*} 
(both here and below). 
For each $\tau\in{\mathcal T}$, we attempt to determine 
when it is that one may reasonably expect to have $\dot Y^{\pm}_m (n)\leq \tau$. 
In particular, 
for $\pm\in\{ + , -\}$, $12\leq n\leq 21$  and  
\begin{equation*}
1\leq m\leq M_0(\pm n) := \begin{cases} M-1 & \text{if $\pm n\geq 18$} , \\  M & \text{otherwise} , \end{cases} 
\end{equation*}
we use the data ${\mathcal Y}^{\pm}_m (n)$ and ${\mathcal Y}^{\pm}_m (n-1)$ to compute  
an interval $\dot{\mathcal Y}^{\pm}_m (n)$ containing $\dot Y^{\pm}_m (n)$. 
We then go on to compute, for $\tau\in{\mathcal T}$, $\pm\in\{ + , -\}$ and $12\leq n\leq 21$, 
an upper bound $F^{\pm}_{\tau} (n)$ for the cardinality of the set 
\begin{equation*} 
{\mathcal F}^{\pm}_{\tau} (n) := \bigl\{ m\leq M_0 (\pm n) : \dot Y^{\pm}_m (n) > \tau\bigr\}  
\end{equation*} 
(more specifically, the difference $M_0 (\pm n) - F^{\pm}_{\tau} (n)$ is here the number of 
positive integers $m\leq M_0 (\pm n)$ for which we find that    one has 
$\sup\dot{\mathcal Y}^{\pm}_m (n) < \tau$). At the same time (and in the same cases) we 
compute a lower bound $f^{\pm}_{\tau} (n)\in{\mathbb N}$ for  ${\mathcal F}^{\pm}_{\tau} (n)\cup \{ M_0 (\pm n) +1\}$: 
this lower bound is the least $m\in{\mathbb N}$ for which we fail to determine that $\sup\dot{\mathcal Y}^{\pm}_m (n) < \tau$. 
The results of these computations are shown in Tables~\mbox{3a--3c} (below). 

\begin{remarks}[On Tables~\mbox{3a--3c}]
\item{\it 1)}\quad Regarding these tables, we observe firstly that for each of the relevant choices of $\tau$ 
the functions $n\mapsto f^{+}_{\tau} (n)$  
and $n\mapsto f^{-}_{\tau} (n)$ show a clear tendency to be increasing on (at least) a significant part of their domain. 
In particular, $f_0^{+}(n)$ and $f_0^{-}(n)$ are increasing for $15\leq n\leq 21$; while, but for the fact that 
$f_{-1}^{+} (20)$ is less than $f_{-1}^{+} (19)$, it could be said that both  $f_{-1}^{+}(n)$ and $f_{-1}^{-}(n)$ are increasing 
for $12\leq n\leq 21$. 
\par 
In the case $\tau = -\omega$, we observe (against the trend) that $f_{-\omega}^{+}(21)$ 
is very much smaller than $f_{-\omega}^{+}(20)$: we argue that this outcome is misleading, since it occurs as a result of 
our having $\dot{\mathcal Y}^{+}_1  (21)\approx [-2.3229, -1.3992]$, an interval that contains the number 
$-\omega\approx -1.585$, and that is exceptionally long 
in comparison with the corresponding data $\dot{\mathcal Y}^{\pm}_m  (n)$ obtained in most other cases. 
This item of data $\dot{\mathcal Y}^{+}_1  (21)$  is also exceptional in one other respect, for in every 
other case we find either that $\dot{\mathcal Y}^{\pm}_m  (n) \subset (-\infty , -\omega)$, 
or else that $\dot{\mathcal Y}^{\pm}_m  (n) \subset (-\omega, \infty)$: we consequently have  both 
$f^{\pm}_{-\omega} (n) = \min {\mathcal F}^{\pm}_{-\omega} (n)$ 
and $F^{\pm}_{-\omega} (n) = |{\mathcal F}^{\pm}_{-\omega} (n)|$ in all of those cases.  
\par
In view of the above observations, we are tempted to conjecture that one has $\dot Y^{+}_1  (21) < -\omega$, and so 
$f^{+}_{-\omega}(21) = 1\not\in {\mathcal F}^{+}_{-\omega} (21)$. If we assume that this conjecture is valid then 
it follows (by re-examination of our data) that $\min{\mathcal F}^{+}_{-\omega} (21) = 211$. 
We therefore think it is plausible that, if we had been able to 
carry out the relevant numerical calculations with greater precision, 
then we might well have ended up having $f^{+}_{-\omega} (21)$ equal to $211\,$   (instead of $1$), 
so that both $f_{-\omega}^{+}(n)$ and $f_{-\omega}^{-}(n)$ would then have been 
found to be increasing for $12\leq n\leq 21$. This allows us to regard the 
results on $f^{\pm}_{-\omega} (n)$ in Table~\mbox{3c} as being generally supportive 
of the conjecture that, for $\pm\in\{ +, -\}$ and $m\leq 192\,$ (say), one has 
$\dot Y_m^{\pm} (n) < -\omega$ for all integers $n\geq 20$. 
\item{\it 2)}\quad Examining the results in Tables~\mbox{3a--3c}, 
we find that it is almost the case that the function $n\mapsto F^{\pm}_{\tau} (n+1) / F^{\pm}_{\tau} (n-1)$ is decreasing for 
\begin{equation*}
n\geq \begin{cases} 16 & \text{when $\tau = 0$} , \\ 15 & \text{when $\tau = -1$} , \\ 13 & \text{when $\tau = -\omega$} . \end{cases}
\end{equation*}
The only exception to this is that $F^{+}_{-\omega} (18) / F^{\pm}_{-\omega} (16) \not\leq F^{+}_{-\omega} (17) / F^{\pm}_{-\omega} (15)$, 
and this one violation of the general rule is quite marginal,   
since $F^{+}_{-\omega} (18) / F^{\pm}_{\tau} (16)$ is only 
about $1\%$ greater than $F^{+}_{-\omega} (17) / F^{\pm}_{\tau} (15)$.  
This (together with our earlier observation, in Point~\mbox{(1)} above, concerning 
cases where $F^{\pm}_{-\omega} (n) = |{\mathcal F}^{\pm}_{-\omega} (n)|$) 
leads us to conjecture that for $\pm\in\{ + , -\}$ and all integers $n\geq 20$ one has  
\begin{equation*} 
|{\mathcal F}^{\pm}_{-\omega} (n+1)|  \leq  c^{\pm} \cdot |{\mathcal F}^{\pm}_{-\omega} (n-1)|\;, 
\end{equation*} 
where  
\begin{equation*} 
c^{\pm} = \frac{F^{\pm}_{-\omega} (21)}{F^{\pm}_{-\omega} (19)}  
=\begin{cases} 1/9 & \text{if $\pm\in\{ +\}$} , \\ 4/31 & \text{if $\pm\in\{ -\}$} . \end{cases} 
\end{equation*}
Similar considerations lead us to conjecture also that 
$|{\mathcal F}^{\pm}_{-1} (n+1)|  \leq  |{\mathcal F}^{\pm}_{-1} (n-1)| / 11$ 
for $\pm\in\{ = , -\}$ and all integers $n\geq 20$. 
\par 
By results noted in Tables~\mbox{3b--3c}, we  have: 
$|{\mathcal F}^{\pm}_{-1} (20)|\leq 5$, $|{\mathcal F}^{\pm}_{-1} (21)|\leq 1$, 
$|{\mathcal F}^{\pm}_{-\omega} (20)|\leq 11$ and $|{\mathcal F}^{\pm}_{-\omega} (21)|\leq 4$, 
for $\pm\in\{ + , -\}$. Therefore, if both the conjectures just mentioned are correct, then it follows (by induction) 
that we have $|{\mathcal F}^{\pm}_{-1} (22)|\leq\frac{5}{11}$,  $|{\mathcal F}^{\pm}_{-\omega} (22)|\leq\frac{44}{31}$, 
and $|{\mathcal F}^{\pm}_{-1} (n)|\leq\frac{1}{11}$ and $|{\mathcal F}^{\pm}_{-\omega} (n)|\leq\frac{16}{31}$ 
for all integers $n\geq 23$: that is, we would have 
${\mathcal F}^{\pm}_{-1} (n) = {\mathcal F}^{\pm}_{-\omega} (n+1)=\emptyset$ for all integers $n\geq 22$, 
and $|{\mathcal F}^{\pm}_{-\omega} (22)|\leq 1$.  
In summary, we record here our opinion that the results in Tables~\mbox{3b--3c} provide a fair degree of support 
for the conjecture that,  when $\pm\in\{ +, -\}$, $1\leq m\leq M_0 (\pm 18)$ and $n\in{\mathbb N}$,  one has   
\begin{equation*} 
\dot Y^{\pm}_m (n) < \begin{cases}  -1 & \text{if $n=22$} , \\ -\omega &\text{if $n\geq 23$} . \end{cases} 
\end{equation*} 
This concludes our remarks on Tables~\mbox{3a-3c}. 
\end{remarks} 

\begin{table}[ht]
\centering 
\begin{minipage}{55mm}
\resizebox{!}{27mm}{ 
\begin{tabular}[c]{|r|r|r|r|r|} 
\hline 
$n$  & $F^{+}_0 (n)$  & $F^{-}_0 (n)$ & $f^{+}_0 (n)$  & $f^{-}_0 (n)$ \\ \hline 
12  & 7  & 10  & 18  & 19  \\ \hline
13  & 7    & 6  & 35     & 34 \\ \hline
14    & 13    & 11  & 60       & 60 \\ \hline
15     & 8    & 14  & 57     & 57 \\ \hline
16     & 11    & 14  & 107      & 95 \\ \hline
17     & 26    & 24  & 159     & 159  \\ \hline 
18     & 16    & 19  & 185     & 201 \\ \hline
19      & 5    & 8  & 303     & 262 \\ \hline
20      & 2    & 0  & 331  & 385 \\ \hline
21     & 0     & 0  & 384 & 385   \\ 
\hline 
\end{tabular} 
}
\vskip 3mm 
\ Table~3a ($\tau = 0$) 
\end{minipage}
\begin{minipage}{58mm}
\resizebox{!}{27mm}{ 
\begin{tabular}[c]{|r|r|r|r|r|} 
\hline 
$n$  & $F^{+}_{-1} (n)$  & $F^{-}_{-1} (n)$ & $f^{+}_{-1} (n)$  & $f^{-}_{-1} (n)$ \\ \hline 
12    & 268     & 271  & 14      & 19 \\ \hline 
13    & 211     & 211  & 25      & 19 \\ \hline 
14    & 100     & 99  & 32      & 45 \\ \hline  
15    & 79    & 78  & 57     & 45 \\ \hline 
16    & 75    & 70  & 91     & 81 \\ \hline 
17     & 54    & 53  & 140      & 96 \\ \hline 
18     & 36    & 36  & 159     & 159 \\ \hline 
19     & 11    & 11  & 226     & 262 \\ \hline 
20      & 5    & 4  & 211     & 303 \\ \hline 
21      & 1    & 1  & 331     & 374 \\ 
\hline 
\end{tabular} 
}
\vskip 3mm
\ Table~3b ($\tau = -1$) 
\end{minipage}
\begin{minipage}{57mm}
\resizebox{!}{27mm}{ 
\begin{tabular}[c]{|r|r|r|r|r|} 
\hline 
$n$  & $F^{+}_{-\omega} (n)$  & $F^{-}_{-\omega} (n)$ & $f^{+}_{-\omega} (n)$  & $f^{-}_{-\omega} (n)$ \\ \hline 
12    & 344      & 351  & 8       & 8 \\ \hline
13    & 302     & 300  & 14      & 15 \\ \hline
14    & 239     & 243  & 25      & 15 \\ \hline
15    & 187     & 195  & 26      & 33 \\ \hline
16    & 128     & 135  & 45      & 33 \\ \hline
17     & 84     & 98  & 45      & 67 \\ \hline
18     & 58    & 59  & 111      & 67 \\ \hline
19     & 27    & 31  & 192     & 169 \\ \hline
20     & 11    & 8  & 211     & 303 \\ \hline
21      & 3      & 4  & 1     & 331  \\ 
\hline
\end{tabular} 
}
\vskip 3mm
\ Table~3c ($\tau = - \omega$)
\end{minipage} 
\end{table} 

We now expand our working hypothesis (that \eqref{W-H} holds when $10\leq n\leq 21$ and $j\leq M$) 
by assuming also that the conjectures mentioned at the conclusions of Remarks~\mbox{7.3}~\mbox{(1)} and~\mbox{7.3}~\mbox{(2)} are correct. 
Recalling \eqref{Happy=Positive} and relevant definitions, we find (by induction) that our expanded working hypotheses 
imply that when $\pm\in\{ + , -\}$, $1\leq m\leq M_0 (\pm 18)$ and $n\geq 22$ one has 
\begin{equation}\label{SumCorr} 
0 < \Delta^{\pm}_m (n) \leq 3^{21 - n} \Delta^{\pm}_m (21) \cdot 
\begin{cases} 1 & \text{if $m\leq 192$} , \\ \frac32 & \text{otherwise} . \end{cases} 
\end{equation}
Using the case $n=21$ of \eqref{SummedDifferences}, we can 
deduce from the upper bound for $|\Delta^{\pm}_m (n)|$ implied by \eqref{SumCorr} 
that when $\pm\in\{ + , -\}$, $1\leq m\leq M_0 (\pm 18)$ and $N=2^{21}$ one has  
\begin{equation}\label{Extrapolation}
\left| \varkappa^{\pm}_m\right| \leq 
\begin{cases}  \frac32 |\varkappa^{\pm}_m (k_N)| - \frac12 |\varkappa^{\pm}_m (k_{N/2})|   &  \text{if $m\leq 192$} ,   \\  
\frac74 |\varkappa^{\pm}_m (k_N)| - \frac34 |\varkappa^{\pm}_m (k_{N/2})|    &  \text{otherwise} .   \end{cases} 
\end{equation}
We can get a useful numerical bound for $|\varkappa^{\pm}_m|$ from this, since we have (in the explicit real intervals 
${\mathcal X}^{\pm}_m (20)\ni \varkappa^{\pm}_m (k_{2^{20}})$ and  ${\mathcal X}^{\pm}_m (21)\ni \varkappa^{\pm}_m (k_{2^{21}})$)    
sufficient data to compute a sharp upper bound, $\hat{\mathcal U}^{\pm}_m$, 
for the relevant expression on the right-hand side of \eqref{Extrapolation}.  
Thus we obtain, for $\pm\in\{ + , -\}$ and $1\leq m\leq M_0(\pm 18)$, a conditional numerical bound  
\begin{equation}\label{Extrap-1}
\left| \varkappa^{\pm}_m\right| \leq \hat{\mathcal U}^{\pm}_m 
\end{equation}
that is valid if our assumptions  (i.e. the expanded working hypotheses mentioned at the start of this paragraph) are indeed correct. 
\par 
The condition $m\leq M_0(\pm 18)$ attached to \eqref{Extrap-1} is a minor nuisance, since we would like to replace it 
with the weaker (and  simpler) condition $m\leq M$. This might be achieved by putting $\hat{\mathcal U}^{+}_M = {\mathcal U}^{+}_M\,$ 
(the upper bound for $\varkappa^{+}_M$ mentioned in Remarks~\mbox{6.14}~\mbox{(3)}), 
or by instead putting $\hat{\mathcal U}^{+}_M = \tilde{\mathcal U}^{+}_M\,$ (the conditional upper bound for $\varkappa^{+}_M$ mentioned 
in Remarks~\mbox{6.14}~\mbox{(5)}). Neither of these upper bounds, however, is very close to the best lower bound 
that we have for $\varkappa^{+}_M$, which is  ${\mathcal L}^{+}_M\approx 2.2406558\times 10^{-3}$ 
(whereas $\tilde{\mathcal U}^{+}_M\approx 2.8408\times 10^{-3}$ and ${\mathcal U}^{+}_M\approx 9.84\times 10^{-3}$). 
We have therefore preferred to give $\hat{\mathcal U}^{+}_M$ a value (closer to ${\mathcal L}^{+}_M$) obtained 
through an adaptation of the approach which gave us \eqref{Extrapolation} (for $m\leq M_0(\pm 18)$). 
This  involves a study of the differences 
\begin{equation*}
\Delta^{*} (n) := \varkappa^{+}_M (k_N) - \varkappa^{+}_M (k_{N/16}) = \sum_{0\leq t\leq 3} \Delta^{+}_M (n-t)\quad\  \text{($14\leq n\leq 21$)} . 
\end{equation*} 
We can establish, firstly, that these $8$ differences are all positive numbers.   
An examinination of the behaviour of the function $n\mapsto \log_2(\Delta^{*} (n))$, 
yields results appearing to support of the conjecture that one has 
$\Delta^{*} (t+13) \leq 5^{-2}\cdot 2^{-(7/4)t}\,$ ($t\in{\mathbb N}$). 
Assuming that this conjecture is correct, it follows 
by the relevant case of \eqref{SummedDifferences} that one has   
$\varkappa^{+}_M - \varkappa^{+}_M (k_{2^{21}}) \leq (2^7 -1)^{-1} \cdot 2^{-14} \cdot 5^{-2} \approx 1.92 \times 10^{-8}$. 
We incorporate the last conjecture into our working hypothesis. 
Thus, using the numerical bound  $\varkappa^{+}_M (k_{2^{21}})\leq \sup{\mathcal X}^{+}_M (21)\approx 2.2406562\times 10^{-3}$, 
we get a conditional upper bound $\hat{\mathcal U}^{+}_M \approx 2.2406754\times 10^{-3}$ for the number $\varkappa^{+}_M$.
\par 
Recalling now the unconditional numerical bounds \eqref{WIDE-3}, we may conclude that if the working 
hypotheses that we have adopted (up to this point) are correct then, 
for $\pm\in\{ + , -\}$ and $1\leq m\leq M$,  one has 
\begin{equation}\label{Work-Horse}
{\mathcal B}^{\pm}_m\ni\lambda^{\pm}_m 
\end{equation} 
with  
\begin{equation*}
{\mathcal B}^{+}_m := [1 / \hat{\mathcal U}^{+}_m , 1 / {\mathcal L}^{+}_m]\quad \text{and} \quad 
{\mathcal B}^{-}_m := [-1 / {\mathcal L}^{-}_m , -1 / \hat{\mathcal U}^{-}_m]\;. 
\end{equation*}
If correct, these conditional bounds tell us far more precisely   
where the relevant eigenvalues of $K$ are located than do the unconditional 
bounds discussed in Section~6: we find, in particular, that 
for $\pm\in\{ +,-\}$ and $1\leq m\leq M=384$, one has 
\begin{equation}\label{High-Precision}
0<\frac{\hat{\mathcal U}^{\pm}_m}{{\mathcal L}^{\pm}_m} - 1 \leq (6.4\times 10^{-10}) m^2 < 10^{-4} 
\end{equation}
unless $m=1$ and $\pm\in\{ -\}\,$ (in which case the upper bound is not $6.4\times 10^{-10}$, but 
$1.1\times 10^{-9}$). The  data sets 
$\{ (\hat{\mathcal U}^{+}_m / {\mathcal L}^{+}_m - 1) / m^2 : 1\leq m\leq M\}$ and 
$\{ (\hat{\mathcal U}^{-}_m / {\mathcal L}^{-}_m - 1) / m^2 : 1\leq m\leq M\}$ have medians that are 
are less than $5.2\times 10^{-11}$.  
Note that \eqref{High-Precision} compares very favourably 
with what is observed in Remarks~\mbox{6.14}~\mbox{(2)} regarding the unconditional bounds 
\eqref{WIDE-5}.
\par 
Table~2 (in Section~\mbox{6.5}) contains some $6$-digit  approximations to terms of the sequence 
$1/{\mathcal L}^{+}_1,\allowbreak 1/{\mathcal L}^{-}_1,\allowbreak 1/{\mathcal L}^{+}_2,\allowbreak 1/{\mathcal L}^{-}_2, \ldots\ $; 
we have not felt it necessary to include in this paper a corresponding table of approximations to terms of the sequence $1/\hat{\mathcal U}^{+}_1,\allowbreak 1/\hat{\mathcal U}^{-}_1,\allowbreak 1/\hat{\mathcal U}^{+}_2,\allowbreak 1/\hat{\mathcal U}^{-}_2, \ldots\ $, 
since it is clear (see above) that in every case under discussion one has 
$0.9999 < {\mathcal L}^{\pm}_m/\hat{\mathcal U}^{\pm}_m < 1$.  
We shall instead describe (in the next subsection) some significant features of the data 
$\,{\mathcal B}^{+}_1,\allowbreak {\mathcal B}^{-}_1,\allowbreak {\mathcal B}^{+}_2,\allowbreak {\mathcal B}^{-}_2,\allowbreak\ldots ,{\mathcal B}^{-}_M$. 

\subsection{Preliminary examination of the data: some observations}

In this subsection, and the next, we assume the validity of the accumulated working hypotheses of the previous  subsection. 
On this assumption, the relation \eqref{Work-Horse} holds for $\pm\in\{ +,-\}$ and $1\leq m\leq M$; we shall 
discuss some consequences. 
\par
It is worth noting, firstly, that we observe nothing obvious in our data that prevents \eqref{Work-Horse} from being true. 
We have checked, in particular, that each interval ${\mathcal B}^{\pm}_m$ is non-empty,  
and that these intervals are consistent with our having $0<\lambda^{+}_1 \leq \lambda^{+}_2 \leq \ldots \ $ 
and $0>\lambda^{-}_1\geq \lambda^{-}_2 \geq \ldots \ $.  We find, in fact, that 
$\max{\mathcal B}^{-}_1 < 0 < \min{\mathcal B}^{+}_1$, 
and that one has both $\max{\mathcal B}^{+}_m < \min{\mathcal B}^{+}_{m+1}$ and 
$\min{\mathcal B}^{-}_m > \max{\mathcal B}^{-}_{m+1}$ for $1\leq m < M$. 
By this last observation and \eqref{Work-Horse}, it follows that each of  
$\lambda^{+}_1,\allowbreak\lambda^{-}_1,\allowbreak\lambda^{+}_2,\allowbreak\lambda^{-}_2,\allowbreak\ldots ,\allowbreak\lambda^{+}_{M-1},\allowbreak \lambda^{-}_{M-1}$ 
is a simple eigenvalue of $K$: one has, in particular, 
both $\lambda^{+}_{M-1} < \lambda^{+}_M$ and $\lambda^{-}_{M-1} > \lambda^{-}_M$. 
\par 
From the observations that ${\mathcal B}^{+}_{M-1} \subseteq  [445.004 , 445.011]$, 
${\mathcal B}^{-}_M\subseteq [-445.247 , -445.240]$ 
and ${\mathcal B}^{+}_M\subseteq [446.294 , 446.298]$  
we deduce that $\lambda^{+}_{M-1} < |\lambda^{-}_M| < \lambda^{+}_M$. 
It follows that $\lambda_{2M-1}=\lambda^{-}_M$, and that the sequence  
$\lambda_1,\allowbreak\lambda_2,\allowbreak\ldots,\allowbreak\lambda_{2M-1}$ must be some permutation of the $2M-1$ distinct eigenvalues 
$\lambda^{-}_M,\allowbreak\lambda^{+}_{M-1},\allowbreak\lambda^{-}_{M-1},\allowbreak\ldots ,\allowbreak\lambda^{+}_1,\allowbreak\lambda^{-}_1$.  
We have no data concerning $\lambda^{-}_{M+1}$. All we can say regarding this eigenvalue is 
that it less than $\lambda^{-}_M\,$ (if it were equal to $\lambda^{-}_M$ then, by \eqref{Balanced_Indices}, 
the number  $|\lambda^{-}_M|$ would be a positive eigenvalue of $K$, and so could not possibly lie strictly between 
$\lambda^{+}_{M-1}$ and $\lambda^{+}_M$, as it does). 
If $\lambda^{-}_{M+1}$ is less than or equal to $-\lambda^{+}_M$ then 
$|\lambda_{2M}| = \lambda^{+}_M$, but if instead it 
lies in the (non-empty) interval $(-\lambda^{+}_M, \lambda^{-}_M)$ then 
$\lambda_{2M} = \lambda^{-}_{M+1} \in (-\lambda^{+}_M, 0)$. Thus it remains 
an open question as to whether or not one might have $\lambda_{2M} = \lambda^{+}_M$. 
\par
It is also worth noting that 
$\min{\mathcal B}^{+}_{m+2} + \min{\mathcal B}^{-}_m > 0 > \max{\mathcal B}^{+}_{m} + \max{\mathcal B}^{-}_{m+2}$ 
for $1\leq m\leq M-2$. It follows from this (and the observations of the preceding paragraph) that 
\begin{equation}\label{LeapFrog}
|\lambda^{\pm}_{m'}| > |\lambda^{\mp}_m|\quad\text{when $\pm\in\{ -,+\}$, $m,m'\in{\mathbb N}$ and $M+1\geq m' > m+1$} . 
\end{equation}  
This (together with our previous observations) implies something quite interesting about the sum function  
\begin{equation*} 
S(m) := \sum_{\ell = 1}^m {\rm sgn}(\lambda_{\ell})\;. 
\end{equation*} 
Namely it implies that 
\begin{equation}\label{Shadowing} 
|S(m)|\leq 2\quad\text{for $1\leq m\leq 2M = 768$} . 
\end{equation} 
To see this we note firstly that, since $\lambda_{2M-1} = \lambda^{-}_M\in (\lambda^{-}_{M+1}, \lambda^{-}_{M-1})$,  
we must have $S(2M-1) = (M-1)\cdot (1) + M\cdot (-1) = -1$ and $S(2M-2) = 0$, and so certainly will have 
$|S(m)|\leq 2$ when either $1\leq m\leq 2$ or $2M\geq m\geq 2M-4\,$ (given that $|S(m+1) - S(m)|=1$ for $m\in{\mathbb N}$).  
In the remaining cases, where $3\leq m\leq 2M-5$, one can see (with the help of \eqref{LeapFrog}) 
that the ordered pair $(\min_{\ell\leq m} \lambda_{\ell}, \max_{\ell\leq m} \lambda_{\ell})$ 
must equal $(\lambda^{-}_{(m-r)/2} , \lambda^{+}_{(m+r)/2})\in{\mathbb R}^2$ for some choice of integer 
$r\equiv m\pmod{2}$ satisfying $|r|\leq 2$. From this one can deduce 
that $S(m) = \frac{m+r}{2}\cdot 1 + \frac{m-r}{2}\cdot (-1) = r\in\{ -2, -1, 0, 1, 2\}$ in those remaining cases. 
\par
A second consequence of \eqref{LeapFrog} (and observations preceding it) is that, for $\pm\in\{ +,-\}$ 
and $1\leq m\leq M-1$, the set $\{ \ell\in{\mathbb N} : \lambda_{\ell} = \lambda^{\pm}_m\}$ is 
a singleton $\{ \ell^{\pm}_m\}\subset {\mathbb N}$ with 
\begin{equation}\label{Split-Shuffle} 
-2\leq   \ell^{\pm}_m - 2m\leq 1
\end{equation} 
(note, however, that the integer $\ell^{\pm}_m$ here is defined relative to a sequence  
$\lambda_1,\lambda_2,\ldots\ $ whose definition is ambiguous if one has $\lambda_k = - \lambda_j$ 
for some  $j,k\in{\mathbb N}$). We omit the proof of this consequence of \eqref{LeapFrog}, leaving it as an easy exercise. 
\par
We now discuss where in the sequence  
$\lambda^{-}_M,\allowbreak\lambda^{+}_{M-1},\allowbreak\lambda^{-}_{M-1},\allowbreak\ldots ,\allowbreak\lambda^{+}_1,\allowbreak\lambda^{-}_1$ the eigenvalues $\lambda_m\,$ ($1\leq m < 2M$) might occur.  
It is useful to note, firstly, that the sequence of compact intervals 
${\mathcal B}^{-}_M,\allowbreak {\mathcal B}^{+}_{M-1},\allowbreak {\mathcal B}^{-}_{M-1},\allowbreak\ldots ,\allowbreak {\mathcal B}^{+}_1,\allowbreak {\mathcal B}^{-}_1$
happens to have a (necessarily unique) permutation 
${\mathcal B}_1,\allowbreak\ldots ,\allowbreak {\mathcal B}_{2M-1}$ satisfying both   
$\min\{ |x| : x\in {\mathcal B}_{m-1}\} < \min\{ |x| : x\in {\mathcal B}_m\}$ 
and $\max\{ |x| : x\in {\mathcal B}_{m-1}\} < \max\{ |x| : x\in {\mathcal B}_m\}$, for $2\leq m\leq 2M-1$. 
This permutation is (of course) easily computed, and so we are able to find  
the  permutation $\varpi$ of the set $\{ 1,\ldots ,2M-1\}$ that (uniquely) satisfies:
\begin{equation}\label{Permed}
{\mathcal B}_m = 
\begin{cases} 
{\mathcal B}^{+}_{\varpi (m) / 2} & \text{if $\varpi (m)$ is even} , \\ {\mathcal B}^{-}_{(\varpi(m) + 1)/2} & \text{otherwise} .
\end{cases} 
\end{equation}
From  what has been noted (up to this point) one can readily deduce that 
\begin{equation}\label{abs(lam)_intervals}
\{ |x| : x\in{\mathcal B}_m\}\ni |\lambda_m| \quad\text{($1\leq m\leq 2M-1$)} . 
\end{equation} 
\par 
It is  desirable that the relation in \eqref{abs(lam)_intervals} 
be sharpened, where possible, to ${\mathcal B}_m\ni \lambda_m$. 
It is worth noting here, for use later, that we certainly do have 
\begin{equation}\label{Catchall} 
{\mathcal B}_1 \cup {\mathcal B}_1 \cup \ \cdots \ \cup {\mathcal B}_{2M-1} \supseteq \left\{ \lambda_1, \lambda_2 , \ldots , \lambda_{2M-1}\right\} 
\end{equation}
(this follows from \eqref{Work-Horse} and points noted earlier in this subsection). 
It is also clear, in view of previous observations, 
that for $1\leq \ell <m<2M$ either $\max\{ |x| : x\in {\mathcal B}_{\ell}\}< \min\{ |x| : x\in {\mathcal B}_m\}$,  
or else $\{ |x| : x\in {\mathcal B}_{\ell}\}\cap\{ |x| : x\in {\mathcal B}_m\}\neq\emptyset$ and 
$\{ {\mathcal B}_{\ell}, {\mathcal B}_m\} = \{ {\mathcal B}^{-}_j, {\mathcal B}^{+}_k\}$ 
for some positive integers $j\leq M$, $k<M$. By examining our data we find that 
the only pairs $(j,k)\in{\mathbb N}^2$ with $\max\{ j,k\}\leq M$ and $\{ |x| : x\in{\mathcal B}^{-}_j\}\cap {\mathcal B}^{+}_k \neq\emptyset$ 
are those that have $j=k\in W$, where $W := \{ 212, 242, 283, 351, 360, 361, 376, 381 \}$.  
One can deduce from this (and some of the earlier observations) that 
$\{ {\mathcal B}^{-}_j , {\mathcal B}^{+}_j\} = \{ {\mathcal B}_{2j-1}, {\mathcal B}_{2j}\}$ for $j\in W$, 
and that $\max\{ |x| : x\in {\mathcal B}_{\ell}\}< \min\{ |x| : x\in {\mathcal B}_m\}$ 
whenever $\ell$ and $m$ are integers with  
$1\leq \ell <m<2M$ and either $\ell\not\in\{ 2j-1 : j\in W\} = W_{O}\,$ (say) or $m\not\in\{ 2j : j\in W\} = W_{E}\,$ (say). 
By this, combined with \eqref{Permed}--\eqref{Catchall}, 
it follows that one has both $\{ \varpi(2j-1), \varpi(2j)\} = \{ 2j-1, 2j\}$ 
and $\{ \lambda_{2j-1}, \lambda_{2j}\} = \{ \lambda^{+}_j, \lambda^{-}_j\} \subset {\mathcal B}^{+}_j\cup {\mathcal B}^{-}_j$
for all $j\in W$, and that for $m\in\{ 1,\ldots , 2M-1\}\backslash (W_{O} \cup W_{E})$ one has 
\begin{equation}\label{TheFound}
{\mathcal B}_m\ni \lambda_m = \begin{cases} 
\lambda^{+}_{\varpi (m) / 2} & \text{if $\varpi (m)$ is even} , \\ \lambda^{-}_{(\varpi(m) + 1)/2} & \text{otherwise} .
\end{cases} 
\end{equation}
Thus one has:  $\{\ell^{+}_m , \ell^{-}_m\} = \{ 2m-1, 2m\} = \{ \varpi^{-1} (2m-1), \varpi^{-1} (2m)\}$   
when $m\in W$; $\ell^{-}_m = \varpi^{-1} (2m-1)$ when $m\in\{ 1,\ldots , M\}\backslash W$; 
and $\ell^{+}_m = \varpi^{-1} (2m)$ when $m\in\{ 1,\ldots , M-1\}\backslash W$.
Also implicit in what has already been noted is the empirical fact that one has
\begin{equation*} 
-1 \leq (-1)^m \left( \varpi (m) - m\right) \leq 2\quad\ \text{($1\leq m\leq 2M-1$)} .
\end{equation*}

\subsection{A conjecture on the asymptotic behaviour of $|\lambda_m|$ as $m\rightarrow\infty$} 

The results \eqref{MainBound} and \eqref{NotTraceClass} (our Theorem~\mbox{1.1} and Corollary~\mbox{4.2}) 
imply that 
\begin{equation*}
\liminf_{m\rightarrow\infty} \frac{\log(|\lambda_m|/m)}{\log\log m} \in [-{\textstyle\frac32},1]\;. 
\end{equation*}
In this subsection we discuss our conjecture that one has, in fact, 
\begin{equation}\label{MC-1}
\left| \lambda_m\right| \sim c_0 m \log^{-\frac32} m\quad\text{as $\,m\rightarrow\infty$} , 
\end{equation}
where $c_0$ is some positive absolute constant.  We present (below) our arguments in support of this conjecture: 
these rest on an analysis of the data discussed in the previous subsection. 
\par
Let 
\begin{equation*}
C_m := \frac{\left| \lambda_m\right| \log^{\frac32} m}{m}\quad\ \text{($m\in{\mathbb N}$)} . 
\end{equation*}
Our conjecture \eqref{MC-1} is correct if and only if the sequence $C_1,C_2,\ldots\ $ converges to a finite limit:  
by \eqref{MainBound}, we can be sure that if this limit exists then it is some number $c_0\geq 1$. 
We may therefore reformulate \eqref{MC-1} as the conjecture that one has 
\begin{equation}\label{MC-2} 
\lim_{m\rightarrow\infty} C_m = c_0 \quad\text{for some $c_0\in [1,\infty)$} . 
\end{equation}
\par
\begin{figure}[h]
\begin{center}
\includegraphics*[scale=0.98, bb=48 200 549 593]{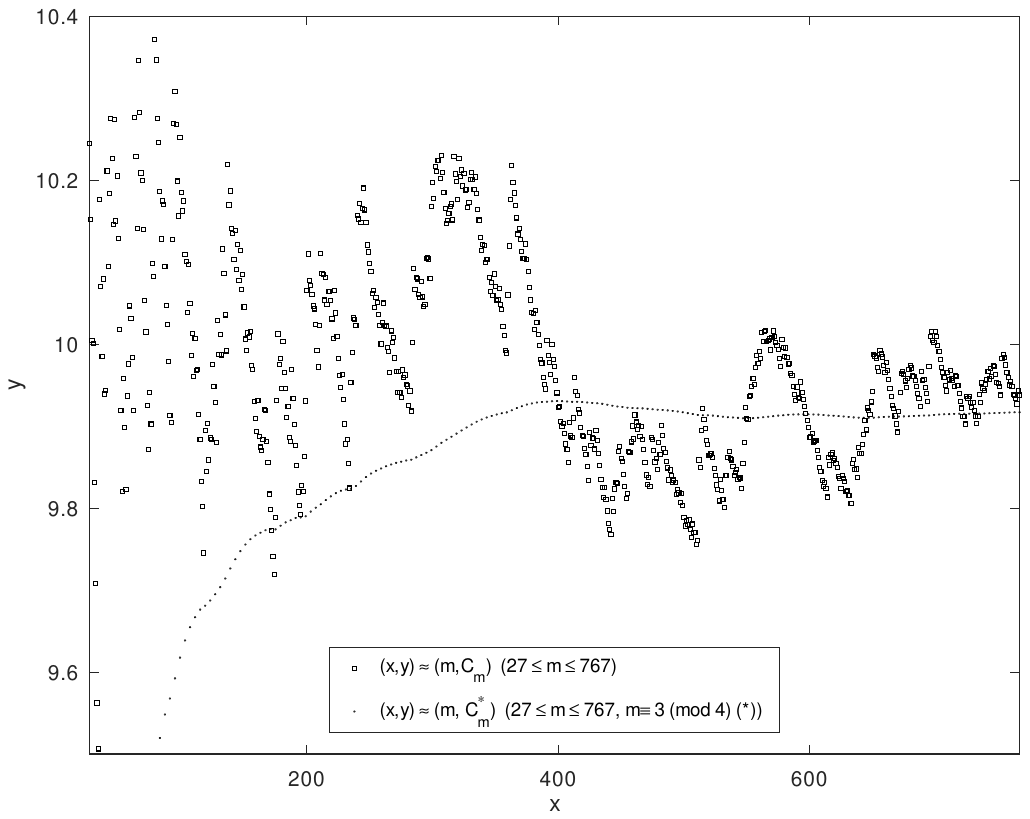} 	
\end{center}
\caption{\small Partial graphs of $m\mapsto C_m$ and $m\mapsto C^{*}_m\,$ [${}^{\star}$points $(m, C^{*}_m)$ with $C^{*}_m < 9.5$ omitted].} 
\label{fig:MC2}
\end{figure}

 
Since we have the relations \eqref{abs(lam)_intervals}, our data 
${\mathcal B}_m\,$ ($1\leq m\leq 2M-1$) enables us to compute short intervals 
${\mathcal C}_1,\ldots ,{\mathcal C}_{2M-1}$ satisfying ${\mathcal C}_m\ni C_m\,$ ($1\leq m\leq 2M-1$). 
For $m\in{\mathbb N}$ we define: 
\begin{equation*}  
C^{*}_{m} := \frac{1}{m}\sum_{\ell =1}^m C_{\ell} \;.  
\end{equation*}
We use the data ${\mathcal C}_m\,$ ($1\leq m\leq 2M-1$) in computing intervals 
${\mathcal C}^{*}_1,\ldots ,{\mathcal C}^{*}_{2M-1}$ with ${\mathcal C}^{*}_{m} \ni C^{*}_{m}$ for $1\leq m\leq 2M-1$. 
These computations yield satisfactorily precise information regarding where 
$C_1,\ldots ,C_{2M-1}$ and $C^{*}_1,\ldots ,C^{*}_{2M-1}$ are to 
be found on the real line. Indeed, we find that for $2\leq m\leq 2M-1$ we have both 
$\log\sup{\mathcal C}_m - \log\inf{\mathcal C}_m < 7.0\times 10^{-5}$ 
and 
$\log\sup{\mathcal C}^{*}_m - \log\inf{\mathcal C}^{*}_m < 4.6\times 10^{-6}$, 
while $(2,11)\supset {\mathcal C}_m, {\mathcal C}^{*}_m\,$ (of course, we have also 
${\mathcal C}_1 =  {\mathcal C}^{*}_1 = \{ 0\}$). Thus we can plot, with reasonable accuracy,  
parts of the graphs of the functions $m\mapsto C_m$ and $m\mapsto C^{*}_m$. 
One such plot is shown in Figure~\ref{fig:MC2} (below). 

Note the restriction to the domain $\{ m\in{\mathbb N} : 27\leq m\leq 2M-1\}$, in Figure~\ref{fig:MC2}. The purpose of the constraint $m\geq 27$ is to 
restrict the range of the $y$-coordinates, thus making it easier to spot small variations in those coordinates. 
In the case of points $(m, C^{*}_m)$ there is the further (quite artificial) condition $m\equiv 3\pmod{4}$: without this most of the 
graph of $m\mapsto C^{*}_m$ would appear (in print) as an unbroken line, rather than as a set of distinct points.  
The position of each omitted point, $(m, C^{*}_m)$ with $m\not\equiv 3\pmod{4}$, can be satisfactorily 
estimated by linear interpolation between $(\ell -1, C^{*}_{\ell -1})$ and $(\ell + 3, C^{*}_{\ell +3})$, for $\ell = 4\lfloor m/4\rfloor$. 
\par
\begin{figure}[ht]
\begin{center}
\includegraphics*[scale=1, bb=75 201 549 582]{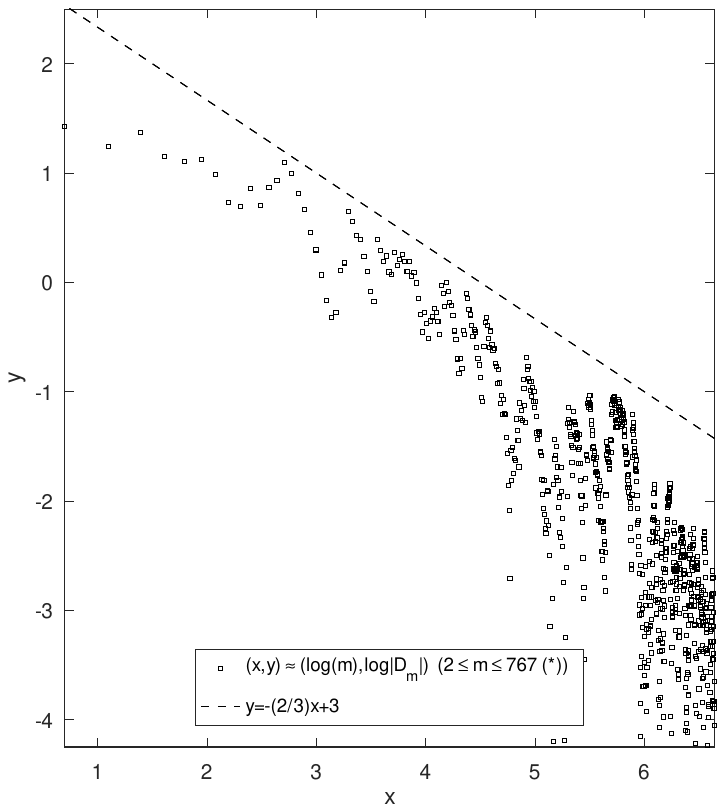} 	
\end{center}
\caption{\small A plot of $\log |D_m|$ versus $\log (m)\,$ [${}^{\star}$points with $y=\log |D_m| < -4.25$ omitted].} 
\label{fig:MC3}
\end{figure}

 
What we see in Figure~\ref{fig:MC2} appears to fit quite well with the conjecture \eqref{MC-2}. 
Seeking further evidence in support of that conjecture, we consider the differences 
\begin{equation}\label{MC-3}
D_m := C_m - C^{*}_{m-1} = m\cdot\left( C^{*}_m - C^{*}_{m-1}\right)\quad \text{($m\geq 2$)} .
\end{equation} 
Using the data ${\mathcal C}_2,\ldots ,{\mathcal C}_{2M-1}$ and ${\mathcal C}^{*}_1,\ldots ,{\mathcal C}^{*}_{2M-2}$,  
we compute certain short closed intervals ${\mathcal D}_2,\allowbreak\ldots , {\mathcal D}_{2M-1}$ satisfying 
${\mathcal D}_m\ni D_m\,$ ($2\leq m\leq 2M-1$). It turns out that none of these intervals contains $0$: 
in fact $560$ of these intervals are compact subsets of $(0,\infty)$, while the other $206$ are compact subsets of 
$(-\infty , 0)$. The data ${\mathcal D}_2,\ldots , {\mathcal D}_{2M-1}$ therefore enables us to 
compute intervals $[a_2,b_2],\ldots , [a_{2M-1},b_{2M-1}]\subset{\mathbb R}$ with 
$a_m \leq \log |D_m|\leq b_m$ for $2\leq m\leq 2M-1$. We find that, while the median length of these $2M-2$ 
intervals is very small (less than $3.1\times 10^{-4}$), some of these intervals are quite long (in particular, 
$b_{597}-a_{597}\approx 0.48$). There are, nevertheless, only $10$ of the intervals whose length exceeds $0.025$, and each 
of these is a subset of $(-\infty , -4.28)$. Happily, it turns out that cases with  $\log |D_m| < -4.28$ are   
not of critical importance in the discussion that is to follow: all that really matters there is 
our data in respect of cases with (say) $|D_m| > 3 m^{-2/3}$, which implies $\log |D_m| > \log(3) - \frac23 \log(2M) > -3.34$. 
For this reason Figure~\ref{fig:MC3} (below) shows only points $(x,y)\approx (m, \log |D_m|)$ that have $y\geq -4.25$: it should be noted that, 
in plotting this figure, we are able to use the number $\gamma_m := \frac 12 (a_m + b_m)$ as a satisfactory approximation to $\log |D_m|\,$ 
(since we know that $|\gamma_m - \log |D_m||\leq \frac12 (b_m - a_m)\leq 0.0125$ whenever $2\leq m\leq 2M-1$ and $\gamma_m \geq -4.25$). 

What we see in Figure~\ref{fig:MC3} persuades us it is reasonable to conjecture that 
\begin{equation}\label{MC-4}
\left| D_m\right| < e^3 m^{-\frac23}\quad\text{($m \geq 2$)} .
\end{equation}
It is assumed henceforth (until the end of the subsection) that this conjecture is correct. 
\par
By \eqref{MC-4} and \eqref{MC-3}, one has  
$| C^{*}_m - C^{*}_{m-1}| < e^3 m^{-5/3}$ for $m\geq 2$. 
It follows that
\begin{equation*} 
\left| C^{*}_{\ell} - C^{*}_{m-1}\right| \leq \sum_{k=m}^{\ell} \left| C^{*}_k - C^{*}_{k-1}\right| 
< e^3 \sum_{k=m}^{\infty} k^{-\frac53} < e^3 \int_{m-1}^{\infty} x^{-\frac53} dx = \frac{3 e^3}{2 (m-1)^{2/3}}\;, 
\end{equation*}
whenever $\ell\geq m > 1$. The sequence $C^{*}_1,C^{*}_2,\ldots\ $ is therefore a Cauchy sequence, and so 
converges to some limit, $c_0\,$ (say), in ${\mathbb R}$. The same inequalities (displayed above) imply 
also that one has 
\begin{equation}\label{MC-5} 
\left| c_0 - C^{*}_{m-1}\right| < {\textstyle\frac32} e^3 (m-1)^{-\frac23}\quad\text{($m\geq 2$)} .
\end{equation} 
Recalling that $D_m := C_m - C^{*}_{m-1}$, we find by \eqref{MC-4} and \eqref{MC-5} that 
\begin{equation}\label{MC-6} 
\left| C_m - c_0\right| < {\textstyle\frac52} e^3 (m-1)^{-2/3}\quad\text{($m\geq 2$)} . 
\end{equation} 
We can therefore conclude 
that if our conjecture \eqref{MC-4} is correct then one has 
both \eqref{MC-1} and \eqref{MC-2} (with $c_0 := \lim_{m\rightarrow\infty} C^{*}_m$).  
Thus the data that was mustered in support of our conjecture \eqref{MC-4} 
serves also to support the conjecture \eqref{MC-1}. 
\par
It turns out that ${\mathcal C}^{*}_{2M-1}$  
is contained in the interval with midpoint $9.916858$ and diameter $4.6\times 10^{-5}$. 
Since ${\mathcal C}^{*}_{2M-1}\ni C^{*}_{2M-1}$, it follows (using \eqref{MC-5}) that one 
must have $|c_0 - 9.916858|< \frac32 e^3 (2M-1)^{-2/3} + 2.3\times 10^{-5} < 0.35959$. 
Thus, subject to our working hypotheses (including \eqref{MC-4}), 
we find that there exists a  constant 
\begin{equation}\label{MC-7} 
c_0 \in (9.557 , 10.277) 
\end{equation}
such that \eqref{MC-1} holds.   
\par
In order to test the credibility of our conjecture \eqref{MC-4} we examine certain 
implications of one of its consequences: the set of bounds in \eqref{MC-6}. 
We consider, for selected values of $c\in (0,\infty)$, 
the implications for the sum 
\begin{equation*} 
F_{2M}' := \sum_{m\geq 2M} \lambda_m^{-2} 
\end{equation*} 
if one has \eqref{MC-6} with $c_0 = c$. 
In view of Theorem~\mbox{1.1}, we consider only cases with $c_0 = c\geq 1$. 
In such cases \eqref{MC-6} gives us 
\begin{equation}\label{MC-8} 
F_{2M}' \leq \sum_{m=2M}^L \frac{\log^3 m}{\left( c m - B m^{1/3}\right)^2} 
+ \frac{L^2}{\left( c L  - B L^{1/3}\right)^2} \sum_{m=L+1}^{\infty}  \frac{\log^3 m}{m^2} 
\end{equation}
when $B = \frac52 e^3 (1 - \frac{1}{2M})^{-2/3} \approx 50.26$ and $L$ is an arbitrary integer greater than $2M$. 
One obtains also an analogous lower bound for $F_{2M}'\,$ (similar to the above upper bound, but with $-B$ in place of $B$).
We work with the case $L=10^6$ of these bounds, since this yields results that are sharp enough for our purposes. 
A numerical bound for the first of the sums over $m$ occuring in \eqref{MC-8}  can easily be computed:  we use  
Octave's {\tt interval} package to accomplish this. 
The other (infinite) sum over $m$ in \eqref{MC-8} lies between $G$ and $(1-L^{-1})G$, 
where $G := \int_L^{\infty} (\log x)^3 x^{-2} dx = (\log^3 L + 3\log^2 L + 6\log L + 6)/L$. 
Using these bounds, together with \eqref{MC-8} and the analogous lower bound for $F_{2M}'$,  
we compute certain numbers $\tau_c' > \tau_c > 0$ such that  
$[\tau_c , \tau_c']\ni F_{2M}'$ if one has \eqref{MC-6} with $c_0 = c$. 
\par
The preceding only becomes really useful if we have some alternative (and independent) way of estimating $F_{2M}'$. 
The case $k=K$ of the result \eqref{TraceLemma-k} (from Lemma~\mbox{6.1}) gives us this, since it implies that we have 
\begin{equation}\label{MC-9}
F_{2M}' =  \| K\|^2 - \sum_{m=1}^{2M-1} \frac{1}{\lambda_m^2}\;. 
\end{equation}
In view of the relations \eqref{abs(lam)_intervals},  our numerical data ${\mathcal B}_1,\ldots ,{\mathcal B}_{2M-1}$ enables 
us to estimate the sum over $m$ in \eqref{MC-9} quite accurately: we find, in fact, that 
this sum lies in the closed interval with midpoint $0.07513436545$ and diameter $7.09\times 10^{-8}$. 
By combining this result with the partial decimal expansion of $\| K\|^2$ from Remarks~\mbox{6.11}~\mbox{(2)}, we 
deduce, using \eqref{MC-9}, that  
\begin{equation}\label{MC-10}
F_{2M}'\in \left\{ x\in{\mathbb R} : \left| x - 0.006386245\right| \leq 3.6\times 10^{-8}\right\} = {\mathcal S}_{2M}'\quad\text{(say)} . 
\end{equation} 
\par 
Given an arbitrary number $c\in [1,\infty)$,  the hypothesis that one has \eqref{MC-6} with $c_0 = c$ can be tested 
by checking if the intervals $[\tau_c , \tau_c']$ and ${\mathcal S}_{2M}'$ overlap: if these two sets are disjoint then
it follows from \eqref{MC-10} that $[\tau_c , \tau_c']\not\ni F_{2M}'$, so that 
one cannot have \eqref{MC-6} with $c_0 = c$. We apply this test for certain numbers $c\geq 1$ 
satisfying $256 c\in{\mathbb Z}$. What we find is that, while $[\tau_c , \tau_c']\cap  {\mathcal S}_{2M}'$ is non-empty 
for all $c\in\{ b/256 : b\in{\mathbb Z}\ {\rm and}\ 2436\leq b\leq 2587\}$, one has 
$\tau_c \geq 0.0063909 > \sup {\mathcal S}_{2M}'$ when $c=\frac{2435}{256}=9.51171875$, 
and $\tau_c' \leq 0.0063844 < \inf {\mathcal S}_{2M}'$ when $c=\frac{2588}{256}=10.109375$. 
Since the upper bound in \eqref{MC-8} and the analogous lower bound for $F_{2M}'$ 
are both monotonic decreasing functions of $c$, it follows that if \eqref{MC-6} holds then the 
relevant constant $c_0$ satisfies 
\begin{equation}\label{MC-11}
9.511 < c_0 < 10.110\;.
\end{equation} 
\par 
We compare \eqref{MC-11} with the range for $c_0$ given in \eqref{MC-7}. 
It is encouraging that the two ranges intersect quite substantially  
(both are less than $31\%$ longer than their intersection), even though the method used to get \eqref{MC-11} 
is different in principle from that which gave us \eqref{MC-7}. 
At the same time, since neither range contains the other, it is worth noting that the combination 
of \eqref{MC-7} and \eqref{MC-11} gives us $c_0\in (9.557, 10.11)$.
In remarks, below, we speculate about a further shortening of the range for $c_0$ . 

\begin{remarks}

Our conjecture \eqref{MC-4} may be too cautious: 
certain trends observed in our numerical data lead us to speculate that 
\begin{equation}\label{MC-12}
m \left| D_m\right| < e^3 \log(m)\quad\text{($m\geq 2$)} , 
\end{equation}
and that 
\begin{equation}\label{MC-13} 
\sum_{\ell = m}^{\infty} \frac{|D_{\ell}|}{\ell} < \frac{50}{m}\quad\text{($m\geq 2$)} .
\end{equation}
We would need more data in order to make a truly convincing case in support of 
(or against) these hypotheses, 
so here we shall only briefly discuss what  
leads us to suggest that \eqref{MC-12} and \eqref{MC-13} might be true. 
In formulating \eqref{MC-12} we were primarily influenced by 
a visual examination of the distribution of points $(\log(m), m|D_m|)\in{\mathbb R}^2\,$ 
(particularly in the cases with $2\leq m\leq 384$). 
The formulation of \eqref{MC-13} was informed by an examination 
of estimates obtained for the sums of the form  $\sum_{\ell = m+1}^{6m} |D_{\ell}|/\ell$ with $m\in\{ 1,2,\ldots  ,127\}$. 
We doubt that either of the constants ($e^3$ and $50$) occurring in \eqref{MC-12} 
and \eqref{MC-13} is optimal. 
\par 
Let us suppose now that the hypotheses \eqref{MC-12} and \eqref{MC-13} are valid. 
In view of the definition \eqref{MC-3}, it follows immediately from \eqref{MC-13} 
that there exists  a unique real constant $c_0$ satisfying 
\begin{equation}\label{MC-14} 
\left| c_0 - C^{*}_{m-1}\right| < \frac{50}{m}\quad\text{($m\geq 2$)} . 
\end{equation}
By this, \eqref{MC-12}, \eqref{MC-3} (again) and the triangle inequality, we get: 
\begin{equation}\label{MC-15}
\left| C_m - c_0\right| < \frac{50 + e^3 \log m}{m} 
\quad\text{($m\geq 2$)} . 
\end{equation}
Similarly to how \eqref{MC-7} was  deduced from \eqref{MC-5}, it follows readily from 
\eqref{MC-14} (and the data ${\mathcal C}^{*}_{2M-1}$) that one has $c_0\in (9.8517 , 9.9820)$. 
At the same time, similarly to how we got \eqref{MC-11} from \eqref{MC-6} and \eqref{MC-10},  
we find that it follows from \eqref{MC-15} and \eqref{MC-10} that one has 
$c_0\in [\frac{4970}{512} , \frac{5073}{512}] \subset (9.7070 , 9.9083)$. 
We conclude that if the hypotheses \eqref{MC-12} and \eqref{MC-13} are valid then 
one has \eqref{MC-1} for some constant $c_0\in (9.8517 , 9.9083)$. 
\par 
The above is not solely in pursuit of a shorter range for $c_0$: 
it also a test of \eqref{MC-12} and \eqref{MC-13},  passed  
by virtue of the fact that 
$(9.7070 , 9.9083)\cap(9.8517 , 9.9820)\neq\emptyset$. 
\end{remarks}

\subsection{The series $\sum_{m=1}^{\infty} \varkappa_m$ and some conjectures}

By Corollary~\mbox{4.2}, the series $\varkappa_1 + \varkappa_2 + \varkappa_3 + \ \ldots\ = \sum_{m=1}^{\infty}\varkappa_m$ is not absolutely convergent. 
If our conjecture \eqref{MC-1} is correct, then one has  $\sum_{\ell = 1}^m |\varkappa_{\ell}| \sim \frac25 c_0^{-1} \log^{5/2} m$ 
as $m\rightarrow\infty$. Nevertheless, subject to the validity of the following conjecture 
(which is supported by evidence mentioned in Section~\mbox{7.2}), it can be seen to follow that 
$\sum_{m=1}^{\infty}\varkappa_m$ is a convergent series. 

\begin{conjecture} 
There exists a function $f : {\mathbb N}\rightarrow {\mathbb N}$ satisfying both 
\begin{equation}\label{Races-Conjecture-1}
\lim_{m\rightarrow\infty} m^{-\varepsilon} f(m) = 0\quad\text{($\varepsilon >0$)} 
\end{equation}
and 
\begin{equation}\label{Races-Conjecture-2}
\max\left\{  \varkappa^{+}_{m+f(m)} ,  |\varkappa^{-}_{m+f(m)}|\right\} 
< \min\left\{  \varkappa^{+}_{m} ,  |\varkappa^{-}_{m}|\right\} 
\quad\text{($m\in{\mathbb N}$)} . 
\end{equation}
\end{conjecture} 

See \eqref{LeapFrog}, and  the second paragraph of Section~\mbox{7.2}, for the empirical evidence that leads us to make the above conjecture: 
what is noted there might even be considered supportive of the stronger conjecture that 
one has \eqref{Races-Conjecture-2} when $f(m)=2$ for all $m\in{\mathbb N}$. 
\par
Supposing now that Conjecture~\mbox{7.3} is correct, we take $f$ to be a function satisfying the conditions 
\eqref{Races-Conjecture-1} and \eqref{Races-Conjecture-2}. Those conditions 
will still be satisfied if we substitute $\max_{\ell\leq m} f(\ell)$ in place of $f(m)$:  
we may therefore assume henceforth that $f$ is monotonic increasing. 
We first apply \eqref{Races-Conjecture-1} and \eqref{Races-Conjecture-2} to the sums 
\begin{equation*}
T_m := \sum_{\ell = 1}^m \left( \varkappa^{+}_m + \varkappa^{-}_m\right)\quad \text{($m\in{\mathbb N}$)} . 
\end{equation*} 
Given an arbitrary choice of integers $\ell >m\geq1$, one has 
\begin{align*} 
T_{\ell} - T_m = \sum_{j=m+1}^{\ell} \varkappa^{+}_j +  \sum_{j=m+1}^{\ell} \varkappa^{-}_j 
 &< \sum_{j=m+1}^{\ell + f(\ell)} \varkappa^{+}_j +  \sum_{j=m+1}^{\ell} \varkappa^{-}_j \\
 &= \sum_{j=m+1}^{m + f(\ell)} \varkappa^{+}_j +  \sum_{j=m+1}^{\ell} \left( \varkappa^{+}_{j+f(\ell)} + \varkappa^{-}_j\right) \;.
\end{align*} 
We observe, firstly, that the penultimate sum over $j$ here is not greater than $f(\ell)\varkappa^{+}_m$ and, secondly,    
that for $j\leq \ell$ one has $\varkappa^{+}_{j+f(\ell)} \leq \varkappa^{+}_{j+f(j)} < -\varkappa^{-}_j\,$ 
(by virtue of the monotonicity of $f$ and the relations \eqref{Races-Conjecture-2}). 
We therefore have $T_{\ell} - T_m \leq f(\ell)\varkappa^{+}_m$ when $\ell >m\geq1$. 
One can show (similarly) that $T_{\ell} - T_m \geq f(\ell)\varkappa^{-}_m$ when $\ell >m\geq1$.
It follows that for $\ell >m\geq1$ one has 
$|T_{\ell} - T_m|\leq f(\ell) \max\{ \varkappa^{+}_m  , |\varkappa^{-}_m|  \}\leq  f(\ell) |\varkappa_m|$.
By this, Theorem~\mbox{1.1} and \eqref{Races-Conjecture-1}, we get 
\begin{equation}\label{RC-consequence-1}
\left| T_{\ell} - T_m\right| = O_{\varepsilon}\left( \ell^{\varepsilon} m^{-1}\log^{3/2} m\right) 
\quad\text{($\ell > m\geq 2$)} , 
\end{equation} 
with $\varepsilon$ here denoting an arbitrarily small positive constant. 
\par 
For our present purposes the case $\varepsilon = \frac12$ of \eqref{RC-consequence-1} is more than sufficient. 
We observe that it implies that  $|T_{\ell} - T_m| \ll \ell^{1/2} m^{-3/4} \ll m^{-1/4}$ whenever $m$ and $\ell$ 
are positive integers satisfying $2\leq m < \ell \leq 2m$. Thus, for $\ell > m \geq 2$, we have: 
\begin{equation*}
\left| T_{\ell} - T_m\right| \leq \sum_{k=1}^{\lceil \log_2 (\ell / m)\rceil} \left| T_{\min\{ \ell , 2^k m\}} - T_{2^{k-1} m}\right| 
\ll \sum_{k\in{\mathbb N}} \left( 2^{k-1} m\right)^{-1/4} \ll m^{-1/4}\;.
\end{equation*} 
It follows from this that $T_1,T_2,T_3,\ldots\ $ is a Cauchy sequence, and so  there must 
exist a unique real number $c_1\,$ (say) such that 
\begin{equation}\label{RC-consequence-2}
\lim_{m\rightarrow\infty} T_m = c_1\;. 
\end{equation}
Since one has both $\varkappa^{+}_m\rightarrow 0$ and $\varkappa^{-}_m\rightarrow 0$ as $m\rightarrow\infty$, 
it is therefore a corollary of \eqref{RC-consequence-2} that the series 
$\varkappa^{+}_1 + \varkappa^{-}_1 + \varkappa^{+}_2 + \varkappa^{-}_2 + \ \ldots\ $ is convergent, with sum $c_1$. 
\par 
Finally, regarding the series $\sum_{m=1}^{\infty} \varkappa_m$, we find it helpful to observe 
that \eqref{Races-Conjecture-2} implies that one has 
\begin{equation}\label{converter-LB}
\left| \varkappa_{2m + 2f(m) - 1}\right| < \min\left\{  \varkappa^{+}_{m} ,  |\varkappa^{-}_{m}|\right\} 
\quad\text{($m\in{\mathbb N}$)} . 
\end{equation}
Therefore when $m\in{\mathbb N}$ one must have 
\begin{equation*} 
\sum_{\ell =1}^{2m+2f(m)-1} \varkappa_{\ell} 
= \sum_{\ell = 1}^{m+a} \varkappa^{+}_{\ell} + \sum_{\ell = 1}^{m+b} \varkappa^{-}_{\ell} 
= T_m + \sum_{m<\ell\leq m+a} \varkappa^{+}_{\ell} + \sum_{m<\ell\leq m+b} \varkappa^{-}_{\ell} \;, 
\end{equation*} 
where  $a=a(m)$ and $b=b(m)$ are certain non-negative integers satisfying $a+b=2f(m)-1$. 
This gives us $|\sum_{\ell = 1}^{2m} \varkappa_m - T_m| < 2f(m) (|\varkappa_{2m}| + \max\{ \varkappa^{+}_m , |\varkappa^{-}_m|\})\leq 4f(m)|\varkappa_m|$, 
so that, by \eqref{RC-consequence-2}, \eqref{Races-Conjecture-1} and Theorem~\mbox{1.1}, one can deduce that 
$\lim_{m\rightarrow\infty}\sum_{\ell = 1}^{2m} \varkappa_m = c_1$. 
By this and Theorem~\mbox{1.1} (again), it follows that the  series $\sum_{m=1}^{\infty} \varkappa_m$ is convergent, 
and has the same sum, $c_1$, as the series $\varkappa^{+}_1 + \varkappa^{-}_1 + \varkappa^{+}_2 + \varkappa^{-}_2 + \ \ldots\ $. 

\begin{figure}[th]
\begin{center}
\includegraphics*[scale=0.98, bb=48 200 549 593]{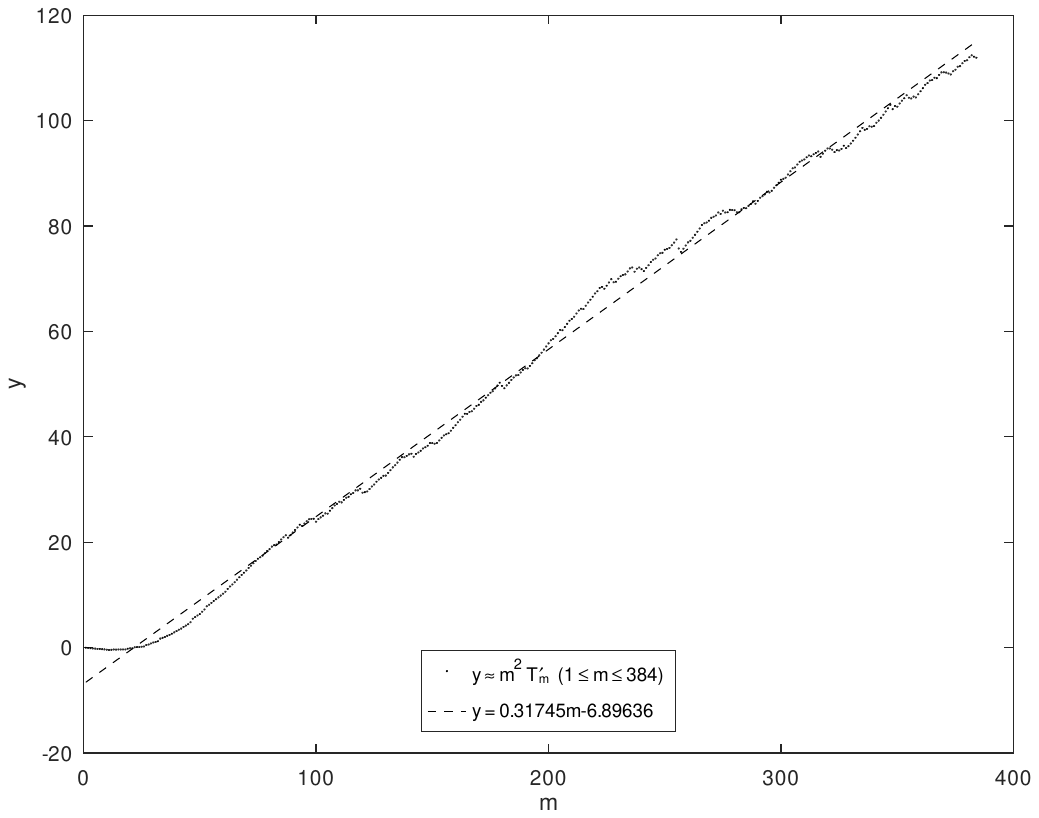} 	
\end{center}
\caption{\small A plot of $m^2 T_m' = m^2 (T_m - \int_0^1 K(x,x)dx)$ versus $m$.} 
\label{fig:PSS3}
\end{figure}

\begin{remarks} 

\item{\it 1)}\quad Supposing that Conjecture~\mbox{7.3} is correct, and assuming that  
$f : {\mathbb N}\rightarrow{\mathbb N}$ is a monotonic increasing function 
satisfying the conditions \eqref{Races-Conjecture-1} and \eqref{Races-Conjecture-2},  
we find, by \eqref{Races-Conjecture-2}, that one has \eqref{converter-LB}, and also 
\begin{equation*} 
\left|\varkappa^{\pm}_m\right| 
\leq \bigl|\varkappa^{\pm}_{m - f(m) +f(m-f(m))}\bigr| 
< \min\left\{ \varkappa^{+}_{m-f(m)} , \bigl| \varkappa^{-}_{m-f(m)}\bigr|\right\} 
\leq  \left| \varkappa_{2(m-f(m))}\right|\;, 
\end{equation*}
whenever $\pm\in\{ + , -\}$ and $m$ is a positive integer with $m>f(m)$. 
By this and \eqref{Races-Conjecture-1}, there exists an $m_0\in{\mathbb N}$ 
such that 
$\{\lambda^{+}_m,|\lambda^{-}_m|\}\subset (|\lambda_{2(m-f(m))}| , |\lambda_{2(m+f(m))}|)$ when  $m > m_0$.  
Using this, one can show that if Conjecture~\mbox{7.3} is correct, and if $\eta\in (0,1)$ is such that 
\begin{equation}\label{MC-powerbound} 
|\lambda_m| = c_0 m \log^{-3/2} (m) + O\left(m^{\eta}\right) \quad\text{as $m\rightarrow\infty$} ,
\end{equation}
then 
\begin{equation*}
|\lambda^{\pm}_m| = 2 c_0 m \log^{-3/2} (2m) + O\left(m^{\eta}\right) \quad\text{as $m\rightarrow\infty$} . 
\end{equation*}
Here it is worth recalling (from the previous subsection) the conditional bounds \eqref{MC-6} 
implied by our conjecture \eqref{MC-4}: if those bounds are valid (as empirical evidence 
presented in Section~\mbox{7.3} seems to indicate) then 
one has \eqref{MC-powerbound} when $\eta = \frac13$. 
The hypothetical bounds \eqref{MC-15} would imply 
that one has \eqref{MC-powerbound} for $\eta = 0$. 

\item{\it 2)}\quad By revisiting some of the steps in our reasoning (above) and taking more care there, 
we can show that if Conjecture~\mbox{7.3} is correct then, for some $c_1\in{\mathbb R}$ and all $\varepsilon > 0$, 
one has both $|T_m - c_1| = O(m^{\varepsilon -1})$ 
and $|\sum_{\ell\leq m} \varkappa_m - c_1| = O(m^{\varepsilon -1})$ as $m\rightarrow\infty$. 
\par 
We conjecture that \eqref{RC-consequence-2} holds with 
\begin{equation}\label{Mercer-Conjecture}
c_1 = \int_0^1 K(x , x) dx = {\textstyle\frac12}\int_1^{\infty} \left( {\textstyle\frac12} - \{ y\} \right) y^{-3/2} dy 
= {\textstyle\frac32} + \zeta({\textstyle\frac12})  \approx 0.03964549
\end{equation} 
(\eqref{DefK} and \eqref{LaplaceB1} giving us the second and third equalities here), and that one in fact has: 
\begin{equation}\label{Merc-2}
T_m = \int_0^1 K(x , x) dx + O\left( \frac{\log^2 (m)}{m}\right)\qquad \text{($m\geq 2$)} . 
\end{equation} 
In proposing the first of these conjectures we are guided by Mercer's theorem \cite[Section~3.12]{Tr1957}, 
which implies that 
one has $\sum_{m=1}^{\infty} \varkappa_m (k) = \int_0^1 k(x,x)dx$ for any continuous real symmetric integral kernel 
$k : [0,1]\times [0,1]\rightarrow{\mathbb R}$ that has only finitely many negative eigenvalues.  
Our kernel $K$, however, is not continuous: and we also know (see Section~1) that both $K$ and $-K$ have infinitely many negative 
eigenvalues. The second conjecture, \eqref{Merc-2}, is implied by the first if 
\eqref{Races-Conjecture-2} holds with  $f(m) := \lceil \log^{1/2} (3m)\rceil $. 
\par 
Through an analysis of our numerical data 
${\mathcal B}^{+}_1,\allowbreak {\mathcal B}^{-}_1,\allowbreak \ldots ,\allowbreak {\mathcal B}^{+}_M,\allowbreak {\mathcal B}^{-}_M$ 
we obtain empirical evidence that appears supportive of the strong conjecture \eqref{Merc-2}. 
This analysis begins with the  computation of short intervals 
${\mathcal T}_1',\ldots ,{\mathcal T}_M'$ satisfying 
\begin{equation}\label{T_m_prime}
{\mathcal T}_m'\ni T_m' := T_m - \int_0^1 K(x,x) dx\quad\ \text{($1\leq m\leq M$)} .
\end{equation} 
Recall that $T_m := \sum_{j=1}^m (1/\lambda^{+}_j + 1/\lambda^{-}_j)$, 
and that we are assuming here that \eqref{Work-Horse} holds in all relevant cases.    
Therefore the computation of ${\mathcal T}_1',\ldots ,{\mathcal T}_M'$ from the data 
${\mathcal B}^{+}_1,\allowbreak {\mathcal B}^{-}_1,\allowbreak \ldots ,\allowbreak {\mathcal B}^{+}_M,\allowbreak {\mathcal B}^{-}_M$ 
is a straightforward task, made easy by the use  of Octave's {\tt interval} package  
(using the estimate $\zeta(\frac12) \approx -1.4603545088095868$, from \cite[A059750]{Sl2022}, 
we get a better approximation to $\int_0^1 K(x,x)dx$ than is shown in \eqref{Mercer-Conjecture}). 
We find that ${\mathcal T}_m' \subset (-\infty , 0)$ for $m\leq 21$, 
that  ${\mathcal T}_m' \subset (0, \infty)$ for $22\leq m\leq M=384$, and that 
$|\log |\sup{\mathcal T}_m'| - \log |\inf{\mathcal T}_m'|| < ( 1  + 10^{-6} m^3)\times 2.8\times 10^{-4} < 1/60$ 
for $1\leq m\leq M$. 
\par 
Using the midpoints of ${\mathcal T}_1',\ldots ,{\mathcal T}_M'$ 
as approximations to the numbers $T_1',\ldots ,T_m'$, we obtain 
the approximation to the graph of the function $m\mapsto m^2 T_m'\,$ ($1\leq m\leq M$) shown in Figure~\ref{fig:PSS3} (below). 
Also indicated in Figure~\ref{fig:PSS3} is the line $y = 0.31745 m - 6.89636\approx 0.317 (m - 21.7)$, found using the method of least squares.   
In view of how uniformly close this line is to our approximations to the points $(m, m^2 T_m')$ ($1\leq m\leq M$), 
we  think it reasonable to conjecture that there exists an approximately linear relationship between 
$m^2 T_m'$ and $m$. 
Indeed, based on what one can see in Figure~\ref{fig:PSS3}, it appears plausible that one might even have $T_m' \sim c_2 / m$ 
as $m\rightarrow \infty\,$  (with $c_2$ some absolute constant lying between $\frac14$ and $\frac13$). 
Thus, recalling \eqref{T_m_prime}, 
we have (in Figure~\ref{fig:PSS3}) some evidence for a conjecture even stronger than \eqref{Merc-2}. 
It should however be noted that,  
although the straight line plotted in Figure~\ref{fig:PSS3} is quite a good fit for the points  $(m, m^2 T_m')$ with $65\leq m \leq 384$, 
the fit is less satisfactory for smaller values of $m$. The empirical evidence for our conjecture \eqref{Merc-2} is 
therefore not as solid as we would wish: our data set (restricted to cases where $m\leq 384$) is simply not big enough 
to yield a convincing weight of evidence. 
\par 
Although we have no proof that $\lim_{m\rightarrow\infty} T_m = \int_0^1 K(x,x)dx$, we 
do (at least) possess an unconditional proof that, for all integers $N\geq 3$, one has
\begin{align*}
{\rm Trace}(H(N)) - \int_0^1 K(x,x)dx 
 &=O\left( \frac{\log N}{N}\right) + \sum_{i=1}^N {\textstyle\frac{1}{(x_i - x_{i+1})}}  \int_{x_{i+1}}^{x_i}\int_{x_{i+1}}^{x_i} 
\left( K(x,y) - K(x,x)\right) dx dy \\ 
 &= O\biggl(\sqrt{\frac{\log N}{N}}\;\biggr) \;,
\end{align*}
with $x_i = x_i (N)\,$ ($1\leq i\leq N+1$) given by the definitions  \eqref{epsilonDef}--\eqref{x_iDef}. 
Our proof of this does not merit inclusion here, though it may be worth mentioning 
that it relies on a lemma somewhat similar (in both statement and proof) to our Lemma~\mbox{3.2}. 

\end{remarks} 

\section*{Acknowledgements} 
\addcontentsline{toc}{section}{\protect\numberline{}Acknowledgements}
This preprint was prepared using \LaTeX. 
The author acknowledges his use of both GNU~Octave \cite{GNU_OCTAVE} 
and the GNU~Octave interval package for `set-based' interval arithmetic \cite{octave-interval}.

\addcontentsline{toc}{section}{\protect\numberline{}References}
\printbibliography 	

\end{refsection} 

\begin{appendices}

\begin{refsection}[refsA.bib]

\section{Proving results stated in Section~4}

\subsection{Proof of the first equality in \eqref{Gamma_h-HSnormEval}} 

By \eqref{kernel-h} and \eqref{DefK}, we have 
$h^2 (t) = e^{-t} \widetilde B_1^2 (e^t) = e^{-t} ( \widetilde B_2 (e^t)  + \frac{1}{12})$ for $t>0$, 
where $\widetilde B_2 (x) := \{ x\}^2 -\{ x\} + \frac16$.
It follows by \eqref{GammaHSnorm} that 
\begin{equation}\label{BWc1} 
\| \Gamma_h \|_{\rm HS}^2 = 
\mathscr{L}\left( \left(  \widetilde B_2 (e^t)  + \textstyle{\frac{1}{12}}\right) t ; 1\right) . 
\end{equation} 
The periodic Bernoulli function $\widetilde B_2 (x)$ is bounded and continuous  on ${\mathbb R}$, and 
it satisfies $\widetilde B_2 (1) = \frac16$ and $\widetilde B_2' (x) = 2 \widetilde B_1 (x)$, 
for $x\in{\mathbb R}\backslash{\mathbb Z}$.  
We therefore find, by the application of well-known properties of the Laplace transform 
(for which see \cite[Section~1.14~(iii)]{Ol2010})  that, for ${\rm Re}(s) >0$, one has both 
\begin{equation*} 
\mathscr{L}\left( \left(  \widetilde B_2 (e^t)  + \textstyle{\frac{1}{12}}\right) t ; s\right)  
= - \frac{\partial}{\partial s} \mathscr{L}\left( \widetilde B_2 (e^t)  + \textstyle{\frac{1}{12}} ; s\right) 
\end{equation*} 
and
\begin{align*} 
\mathscr{L}\left( \widetilde B_2 (e^t)  + \textstyle{\frac{1}{12}} ; s\right)  &= 
\mathscr{L}\left( \widetilde B_2 (e^t)  - \widetilde B_2 (1)  ; s\right) + \mathscr{L}(\textstyle{\frac14} ; s)  \\ 
 &= \mathscr{L}\left( \int_0^t 2\widetilde B_1 (e^{\tau}) e^{\tau} d\tau ; s\right) + \frac{1}{4s} \\ 
 &= \frac1s  \mathscr{L}\left( 2\widetilde B_1 (e^t) e^t ; s\right) + \frac{1}{4s}  
= \frac2s  \mathscr{L}\left( \widetilde B_1 (e^t) ; s - 1\right) + \frac{1}{4s} \,.   
\end{align*} 
By this and \eqref{LaplaceB1}, we see that when ${\rm Re}(s) > 1$ and $s\neq 2$ one has 
\begin{equation}\label{BWc2} 
\mathscr{L}\left( \left(  \widetilde B_2 (e^t)  + \textstyle{\frac{1}{12}}\right) t ; s\right)   
= \frac{1}{4s^2} + 2\frac{\partial}{\partial s}\left( \frac{\zeta(s-1) - \frac{1}{s-2} - \frac12}{s(s-1)}\right) . 
\end{equation} 
\par 
We put now $s=1+w$. For $0<|w|<1$, one has both 
\begin{equation*} 
\frac{\zeta(w) - \frac{1}{w-1} - {\textstyle{\frac12}}}{w} 
= \frac{\zeta(0) +  {\textstyle{\frac12}}}{w} + \sum_{n=1}^{\infty} \left( \frac{\zeta^{(n)} (0)}{n!} + 1\right) w^{n-1}   
\end{equation*} 
and $(1+w)^{-1} = \sum_{m=0}^{\infty} (-w)^m$. 
Since $\zeta (0) = -\frac12$, we deduce that there is a power series of the form 
\begin{equation*} 
\frac{\zeta' (0)}{1!} + 1 + \left( \frac{\zeta'' (0)}{2!} - \frac{\zeta' (0)}{1!} \right) w + c_2 w^2 + c_3 w^3  +\ \ldots \ 
\end{equation*} 
summing to $(\zeta(w) - \frac{1}{w-1} - \frac12)/((1+w)w)$ when $0<|w|<1$. 
This implies the result   
\begin{equation*} 
\lim_{w\rightarrow 0} \frac{\partial}{\partial w} \left( \frac{\zeta(w) - \frac{1}{w-1} - \textstyle{\frac12}}{(1+w)w}\right) 
= \frac{\zeta'' (0)}{2} - \zeta' (0)\;,
\end{equation*} 
which, taken together with the fact that 
$\mathscr{L}((\widetilde B_2 (e^t) + \frac{1}{12})t ; s)$ is analytic for 
${\rm Re}(s) >0$ (and so is continuous at $s=1$), enables us to deduce from \eqref{BWc2} that 
$\mathscr{L}((\widetilde B_2 (e^t) + \frac{1}{12})t ; 1) 
= \frac14 + \zeta'' (0) - 2\zeta' (0)$. By this and \eqref{BWc1}, we  obtain 
the first equality in \eqref{Gamma_h-HSnormEval}. \hfill $\square$ 

\subsection{Proof of \eqref{Parseval-A}, \eqref{H2(D)HankelAction} and  \eqref{c_m-in-series}} 

We adopt a method  
outlined in \cite[Section~1.5]{Ya2019}: one alternative would have been to use instead 
Theorem~\mbox{4.6} of the book of J.R. Partington that is cited in Section~4. 
\par 
Recall that the family $\{ p_m\}_{m=0}^{\infty}$ is an orthonormal basis for $H^2 ({\mathbb D})$. 
Using the case $\alpha =0$ of \cite[Formula~(18.17.34)]{Ol2010}  one may verify that, 
when $m$ is a positive integer,  
the surjective isomorphism $L^{-1} V : H^2 ({\mathbb D}) \rightarrow L^2 (0,\infty)$ 
maps $p_m$ to the function 
\begin{equation}\label{Def_varpi_m} 
\varpi_m (t) := (-1)^m \sqrt{2}\cdot L_m (2t) e^{-t}\qquad\text{($0<t<\infty$)} , 
\end{equation} 
where $L_m (x) = L_m^{(0)} (x)$ denotes the Laguerre polynomial of degree $m$ 
(defined as in \cite[Chapter~18]{Ol2010}). It follows that  
$\{ \varpi_m\}_{m=0}^{\infty}$ is an orthonormal basis for $L^2 (0,\infty)$, 
and that the isomorphism $U := L^{-1} V$ is fully described by the relations: 
\begin{equation}\label{U-action}
U\sum_{m=0}^{\infty} a_m p_m 
= \sum_{m=0}^{\infty} a_m \varpi_m \qquad 
\text{($\{ a_m\}_{m=0}^{\infty} \in \ell^2$)} . 
\end{equation} 
\par 
The family $\{ \varpi_m\}_{m=0}^{\infty}$ is (of course) not the only basis for 
$L^2 (0,\infty)$: for example, when $\alpha > -1$ and $L^{(\alpha)}_0 (x),L^{(\alpha)}_1 (x),L^{(\alpha)}_2 (x),\ldots\ $ 
are Sonine's generalized Laguerre polynomials (defined as in \cite[Chapter~18]{Ol2010}), 
the family of functions 
\begin{equation}\label{Def_Sonine} 
\varpi_{\alpha , m} (t) := (-1)^m\sqrt{\textstyle\frac{m! 2^{1+\alpha}}{\Gamma(m+\alpha + 1)}} \cdot L^{(\alpha)}_m (2t) t^{\alpha /2} e^{-t} 
\qquad\text{($m=0,1,2,3,\ldots\ $)} 
\end{equation} 
is an orthonormal basis for $L^2 (0,\infty)$. By this and \eqref{GammaHSnorm}, 
it follows that one has, in $L^2 (0,\infty)$, the orthogonal decomposition:  
\begin{equation}\label{root(t)h(t)}
\eta(t) := h(t) t^{1/2}  = \sum_{m=0}^{\infty} (-1)^m c_m \cdot 2 L^{(1)}_m (2t) t^{1/2} e^{-t} 
\end{equation} 
with  
\begin{equation}\label{Def-c_m}
c_m := {\textstyle\frac{1}{\sqrt{m+1}}} \left\langle \eta , \varpi_{1 , m}\right\rangle_{L^2 (0,\infty)} 
=  \left\langle h, \varphi_m\right\rangle_{L^2 (0,\infty)} \;, 
\end{equation} 
where, for $0<t<\infty$, 
\begin{align}\label{REphi_m} 
\varphi_m (t) := {\textstyle\frac{1}{\sqrt{m+1}}}\cdot t^{1/2} \varpi_{1 , m}(t) &= {\textstyle\frac{(-1)^m 2}{m+1}} L^{(1)}_m (2t) t e^{-t} \nonumber\\ 
 &=(-1)^m \left( L_m (2t) - L_{m+1} (2t)\right) e^{-t} \nonumber\\ 
 &= {\textstyle\frac{1}{\sqrt{2}}}\left( \varpi_m (t) + \varpi_{m+1} (t)\right)  
\end{align} 
(the last two equalities following by \cite[Formula~(18.9.14)]{Ol2010} and \eqref{Def_varpi_m}). 
By virtue of Parseval's identity, one has also: 
\begin{equation*}
\sum_{m=0}^{\infty} (m+1) \left| c_m\right|^2 = \| \eta\|_{L^2 (0,\infty)}^2 = \int_0^{\infty} h^2 (t) t dt \;. 
\end{equation*}
By this and \eqref{GammaHSnorm}, we obtain the result \eqref{Parseval-A}: note that 
in proving (below) the result \eqref{H2(D)HankelAction} we shall, at the same time, be showing that \eqref{Def-c_m} defines exactly the 
same sequence $\{ c_n\}_{n=0}^{\infty}\in\ell^2$ as that which appears in Section~4. 
\par 
It is implicit in \eqref{root(t)h(t)} that the `remainder term' functions $R_0,R_1,R_2,\ldots\ $ 
defined on $(0,\infty)$ by 
\begin{equation*} 
R_N(t) = \left( h(t) - \sum_{n=0}^N (-1)^n c_n \cdot 2 L^{(1)}_n (2t) e^{-t}\right) \cdot t^{1/2} 
\end{equation*} 
satisfy $\| R_N\|_{L^2 (0,\infty)} \rightarrow 0$ as $N\rightarrow \infty$. 
Although it does not follow that this remains true when $R_N(t)$ is replaced by $t^{-1/2} R_N(t)$, 
one can at least deduce that, for each constant $t>0$, one has 
$\lim_{N\rightarrow\infty} \int_t^{\infty} w^{-1}R_N^2(w)dw = 0$. 
Suppose now that $g\in L^2 (0,\infty)$ and $0<t<\infty$. 
By combining the last observation with \eqref{DefGamma_h}, we find that 
\begin{equation*} 
\left( \Gamma_h g\right) (t) 
= \sum_{n=0}^{\infty}  (-1)^n c_n\int_0^{\infty} 2 L^{(1)}_n (2t + 2u) e^{-(t+u)} g(u) du\; . 
\end{equation*} 
We apply the special case $r=2$, $\alpha_1 = \alpha_2 = 0$ of the Addition Theorem \cite[(18.18.10)]{Ol2010} to 
expand the term $L^{(1)}_n (2t + 2u)$ in the above equation. This (if we recall the definitions \eqref{Def_varpi_m}) yields 
the result:  
\begin{equation*} 
\left( \Gamma_h g\right) (t) 
= \sum_{n=0}^{\infty} c_n\sum_{r=0}^n \left\langle g , \varpi_{n-r}\right\rangle_{L^2 (0,\infty)} \cdot \varpi_r (t) \;. 
\end{equation*} 
\par 
Considering in particular the case where  (for an arbitrary $m\in\{ 0\}\cup{\mathbb N}$) one has $g=\varpi_m$, 
we find (by the result just obtained, and the orthonormality of $\{ \varpi_n\}_{n=0}^{\infty}$) 
that one has $\Gamma_h \varpi_m (t) = \sum_{n=m}^{\infty} c_n \varpi_{n-m} (t) =\sum_{n=0}^{\infty} c_{m+n} \varpi_n (t)$ 
for $0<t<\infty$. 
This establishes that the series $\sum_{n=0}^{\infty} c_{m+n} \varpi_m$ 
converges pointwise on $(0,\infty)$, and has the pointwise sum $\Gamma_h \varpi_m$. By \eqref{Parseval-A} and 
the orthonormality of $\{ \varpi_n\}_{n=0}^{\infty}$, the same series $\sum_{n=0}^{\infty} c_{m+n} \varpi_m$ 
is convergent in the Hilbert space $L^2 (0,\infty)$,  so it has there a sum $\varsigma_m$ (say),   
and there exists some strictly increasing sequence of positive integers 
$\{ N_k\}_{k=0}^{\infty}$ such that the sequence $\{\sum_{n=0}^{N_k} c_{m+n} \varpi_m\}_{k=0}^{\infty}$ 
converges pointwise to  that sum $\varsigma_m$ almost everywhere in the interval $(0,\infty)$ 
(see the discussion of the Riesz-Fisher theorem in \cite[Section~10.25]{Ap1974} regarding these assertions). 
Therefore, recalling it was found (above) that  
$\lim_{N\rightarrow\infty} \sum_{n=0}^N c_{m+n} \varpi_m (t) = \Gamma_h \varpi_m(t)\,$  
for $0<t<\infty$, we may conclude now that $\varsigma_m = \Gamma_h \varpi_m$ almost everywhere in $(0,\infty)$. 
It follows that, in the Hilbert space $L^2 (0,\infty)$, one has:  
$\Gamma_h \varpi_m =\varsigma_m = \sum_{n=0}^{\infty} c_{m+n} \varpi_n$. 
\par 
By \eqref{Parseval-A}, \eqref{U-action} and the conclusion reached in the last paragraph, 
one has $U\sum_{n=0}^{\infty} c_{m+n} p_n = \sum_{n=0}^{\infty} c_{m+n} \varpi_n = \Gamma_h \varpi_n = \Gamma_h U p_m$, 
for all non-negative integers $m$. 
Since $U^{-1} \Gamma_h U = \widetilde \Gamma_h$, it follows that 
\begin{equation*}
\sum_{n=0}^{\infty} c_{m+n} p_n = \widetilde \Gamma_h p_m 
\qquad\text{($m=0,1,2,\ldots\ $)} . 
\end{equation*} 
Therefore, if $\{ a_m\}_{m=0}^{\infty}\in\ell^2$ then   
\begin{equation*} 
\widetilde \Gamma_h \sum_{m=0}^{\infty} a_m p_m  
=\lim_{M\rightarrow\infty} \sum_{m=0}^M a_m \sum_{n=0}^{\infty} c_{m+n} p_n  
=\lim_{M\rightarrow\infty} \sum_{n=0}^{\infty} \left( \sum_{m=0}^M a_m c_{m+n}\right) p_n\; ,
\end{equation*} 
and so we obtain \eqref{H2(D)HankelAction} (since it follows from \eqref{Parseval-A}, via the Cauchy-Schwarz inequality, that 
$\|\{\sum_{m>M} a_m c_{m+n}\}_{n=0}^{\infty}\|_{\ell^2} 
\leq  (\sum_{m>M} |a_m|^2)^{1/2} \cdot \| \Gamma_h\|_{\rm HS}  < \infty$ when $M\geq 0$). 
\par 
Given that $U : H^2 ({\mathbb D}) \rightarrow L^2 (0,\infty)$ is an isomorphism, it 
follows from \eqref{Def-c_m} that one has $c_m = \langle U^{-1} h , U^{-1}  \varphi_m \rangle_{H^2 ({\mathbb D})}$. 
Now $U^{-1}  \varphi_m =  (p_m + p_{m+1})/\sqrt{2}$, by \eqref{REphi_m} and \eqref{U-action}, 
and $U^{-1} h = V^{-1} L h  = \frac{1}{\sqrt{2\pi}} V^{-1} H$. 
Thus, with  
\begin{equation*} 
q(z) := \left( (2\sqrt{\pi} V)^{-1} H\right) (z) = (1 + z)^{-1} H(Mz)\qquad 
\text{($z\in{\mathbb D}$)} , 
\end{equation*} 
one has:  
\begin{align*} 
c_m &= \left\langle q(z) , p_m (z)\right\rangle_{H^2 ({\mathbb D})} + \left\langle q(z) , p_{m+1} (z)\right\rangle_{H^2 ({\mathbb D})}  \\ 
 &=  \left\langle q(z)z , p_{m+1} (z)\right\rangle_{H^2 ({\mathbb D})} + \left\langle q(z) , p_{m+1} (z)\right\rangle_{H^2 ({\mathbb D})}  \\ 
 &=  \left\langle (z + 1) q(z) , p_{m+1} (z)\right\rangle_{H^2 ({\mathbb D})} \\ 
 &=  \left\langle  H(Mz) , p_{m+1} (z)\right\rangle_{H^2 ({\mathbb D})} \;.  
\end{align*} 
It follows that $c_m$ is the coefficient of $z^m$ in the Laurent series expansion of the function 
$z^{-1} H(Mz)$ at the simple pole $z=0$. We  therefore have \eqref{c_m-in-series}. \hfill $\square$

\printbibliography 

\end{refsection}

\section{Computing the elements of $H'(N)$} 

In this appendix we describe how we compute   
our estimates $H_1',\ldots , H_{2N - 1}'$ for 
the numbers $H_1,\ldots , H_{2N - 1}$  defined in Section~\mbox{5.2}. 

\subsection{Lemmas}

Let $\delta = \delta(N)$ and $\varepsilon = \varepsilon(N)$ 
be defined as in Section~\mbox{3.1}. 
With \eqref{H_nDef} in mind, we start with a formula giving the exact 
values of the numbers $M_1,\ldots ,M_{2N - 1}$ that are defined in \eqref{m_ijHankelformed}. 

\begin{lemma} 
Let $n,N\in{\mathbb N}$ satisfy $N\geq 3$ and $1\leq n\leq 2N - 1$. Then 
\begin{align}\label{M_nXACT-1} 
M_n &= {\textstyle\frac12} - \frac{\delta}{(1 - \delta)^2} \cdot 
\Biggl( \delta \sum_{e^{n\varepsilon} < m \leq e^{(n+1)\varepsilon}}   
\left( \frac{\left( 1 - (n+1)\varepsilon\right) + \log m}{\delta^{n+1} m} - 1\right) \nonumber\\ 
 &\phantom{{ = {\textstyle\frac12} + \frac{\delta}{(1 - \delta)^2} \cdot \Biggl( }}  
- \delta^{-1}\sum_{e^{(n-1)\varepsilon} < m \leq e^{n\varepsilon}}  
\left( \frac{\left( 1 - (n-1)\varepsilon\right) + \log m}{\delta^{n-1} m} - 1\right)  \nonumber\\ 
 &\phantom{{ = {\textstyle\frac12} + \frac{\delta}{(1 - \delta)^2} \cdot \Biggl( }}  
- 4\sinh^2 (\varepsilon / 2)  \left\lfloor e^{n\varepsilon}\right\rfloor 
+ \varepsilon^2 e^{n\varepsilon} \Biggr) \;. 
\end{align} 
\end{lemma} 

\begin{proof} 
By \eqref{DefK}, the integrand $K\left( 1 , \delta^{n-1} uv\right)$ 
that occurs in \eqref{m_ijHankelformed} 
is equal to  $\frac12 - \{ \delta^{1-n} u^{-1} v^{-1}\}$. 
The term $\frac12$ here contributes the term $\frac12$ 
occurring on the right-hand side of  \eqref{M_nXACT-1}, so it only remains for us to show that 
$\delta^{-1} \int_{\delta}^1 \int_{\delta}^1   \{ \delta^{1-n} u^{-1} v^{-1}\} du dv$ is equal 
to the expression contained between the largest pair of brackets in \eqref{M_nXACT-1}. 
\par 
By \eqref{deltaDef}, the substitutions 
$u = \exp(-\frac12 \varepsilon \tau - \frac12 \varepsilon (\phi + 1))$ and   
$v = \exp(\frac12 \varepsilon \tau - \frac12 \varepsilon (\phi + 1))$, 
and the observations that 
$\int_{|\phi| - 1}^{1-|\phi|} d\tau = 2-2|\phi|\,$ ($-1 < \phi < 1$) and $\int_{-1}^1 ( 1 - |\phi|) d\phi = 1$, 
we find that 
\begin{multline}\label{BWg1} 
\delta^{-1} \int_{\delta}^1 \int_{\delta}^1   \{ \delta^{1-n} u^{-1} v^{-1}\} du dv 
= \varepsilon^2 \int_{-1}^1 \left( 1 - |\phi|\right) \left\{ e^{(n + \phi)\varepsilon}\right\} e^{-\varepsilon \phi} d\phi \\ 
=  \varepsilon^2 e^{n\varepsilon}  
- \varepsilon^2 \int_{-1}^1 ( 1 - |\phi|) \lfloor e^{(n + \phi)\varepsilon}\rfloor e^{-\varepsilon \phi} d\phi  
= \varepsilon^2 e^{n\varepsilon} -  \varepsilon^2  I_n \quad\text{(say)} . 
\end{multline} 
Now $I_n = \int_{-1}^0 \lfloor e^{(n + \phi)\varepsilon}\rfloor \phi  e^{-\varepsilon \phi} d\phi 
- \int_0^1 \lfloor e^{(n + \phi)\varepsilon}\rfloor \phi  e^{-\varepsilon \phi} d\phi 
+ \int_{-1}^1 \lfloor e^{(n + \phi)\varepsilon}\rfloor e^{-\varepsilon \phi} d\phi = A_n -B_n + C_n\,$ 
(say), and, regarding the first of these three integrals, we have 
\begin{multline*} 
A_n =  \int_{-1}^0 \Biggl( \,\sum_{0 < m\leq  e^{(n + \phi)\varepsilon}} 1 \Biggr) \phi  e^{-\varepsilon \phi} d\phi  
= \sum_{0 < m\leq  e^{n \varepsilon}} 
  \int_{\max\left\{ -1\, ,\, \varepsilon^{-1} \log m - n\right\}}^0 \phi  e^{-\varepsilon \phi} d\phi \\ 
= \left\lfloor e^{(n-1)\varepsilon}\right\rfloor \alpha (-1) 
+ \sum_{ e^{(n-1)\varepsilon} < m\leq  e^{n \varepsilon}} \alpha \left(  \varepsilon^{-1} \log m - n\right) \;, 
\end{multline*} 
where $\alpha (\phi) := \varepsilon^{-2} (( 1 + \varepsilon \phi )  e^{-\varepsilon \phi} - 1)$:  
similar calculations yield broadly similar formulae for the integrals $B_n$ and $C_n$ 
(we omit the details). By these formulae for $A_n$, $B_n$ and $C_n$, we find (after a short calculation) 
that $e^{n\varepsilon} - I_n$ is equal to $\varepsilon^{-2}$ times the 
expression contained between the largest pair of brackets in \eqref{M_nXACT-1}: 
upon combining this with the relations \eqref{BWg1} and our observations in the previous paragraph, 
we have a  proof of the lemma.
\end{proof} 

\begin{remarks} 
The formula \eqref{M_nXACT-1} is, in practice, unsuited to the task of 
obtaining good numerical approximations to the number $M_n$.  
In particular, since Lemma~\mbox{3.4} implies that the factor $\delta / (1 - \delta)^2$ that occurs in \eqref{M_nXACT-1} 
must grow like $(N/\log N)^2$ as $N\rightarrow\infty$, it follows that one must take care to minimise any rounding 
errors when estimating terms occurring in the sums over $m$ in \eqref{M_nXACT-1}. For this reason we 
need the following alternative formulation of the result  \eqref{M_nXACT-1}. 
\end{remarks} 

\begin{lemma} 
Let $n,N\in{\mathbb N}$ satisfy $N\geq 3$ and $1\leq n\leq 2N - 1$. 
Let $R\geq 2$ be even. 
Then 
\begin{equation}\label{M_nXACT-1alt} 
M_n = D_R (n) + E_R (n)\;,
\end{equation} 
where $E_R (n)$ satisfies 
\begin{equation}\label{TruncErrBound} 
\left| E_R (n)\right| < 
e^{2\varepsilon} \left( 1 + 2  \varepsilon e^{\varepsilon n} \left( 1 + \frac{1}{R + 2}\right)\right) 
\cdot\frac{\varepsilon^{R - 1}}{(R + 1)!} \;, 
\end{equation} 
while 
\begin{equation*}
D_R (n) = K\left( 1 , \delta^n\right) 
+ \frac{e^{\varepsilon n}}{4\sinh^2 (\varepsilon / 2) } \sum_{r=2}^R 
\frac{\left( a_r(n) - b_r(n) + c_r\right) \varepsilon^r}{r!}\;, 
\end{equation*} 
with:  
\begin{equation*} 
a_r (n) := \sum_{e^{n\varepsilon} < m \leq e^{(n + 1)\varepsilon}}  
\frac{\left( \frac{\log m}{\varepsilon} - (n + 1)\right)^r}{m}\;,
\end{equation*} 
\begin{equation*} 
b_r (n) := \sum_{e^{(n - 1)\varepsilon} < m \leq e^{n\varepsilon}}   
\frac{\left( \frac{\log m}{\varepsilon} - (n - 1)\right)^r}{m}\quad\text{and}\quad    
c_r := 
\begin{cases} 
2 & \text{if $2 \mid r$ and $r>2$} , \\
0 & \text{otherwise} . 
\end{cases} 
\end{equation*} 
\end{lemma} 

\begin{proof} 
We shall show firstly that, provided we define $E_{\infty} (n)$ to equal zero, 
the equation \eqref{M_nXACT-1alt} holds for $R=\infty$. 
\par 
Let $S_1$ and $S_2$ denote (respectively) the first and second of the sums over $m$ occurring in
\eqref{M_nXACT-1}. Recalling \eqref{deltaDef}, we find that the summand in $S_1$ can be 
expressed as $g(\rho_1(m)\varepsilon)$, where $\rho_1(m) := \varepsilon^{-1} \log(m) - (n + 1)$ and 
$g(z) := (z + 1) e^{-z}$. Similarly, the  summand in $S_2$ equals 
$g(\rho_2(m)\varepsilon)$, where $\rho_2(m) := \varepsilon^{-1} \log(m) - (n - 1)$. 
Therefore, by observing that 
$g(z) - 1 = -e^{-z} \sum_{r=2}^{\infty} z^r /r!\,$ ($z\in{\mathbb C}$), 
and that $\exp (-\varepsilon - \rho_1(m)\varepsilon) =  e^{n\varepsilon} /m = \exp (\varepsilon - \rho_2(m)\varepsilon)$, 
we find that one has $\delta S_1 - \delta^{-1} S_2 
= e^{n\varepsilon} \sum_{r=2}^{\infty} (b_r (n) - a_r (n)) \varepsilon^r /r!$.   
It follows, since $4\sinh^2 (\varepsilon /2)$ equals  $(1 - \delta)^2 / \delta$,  
that completion of the proof of the case  $R=\infty$ of \eqref{M_nXACT-1alt}  (with $E_{\infty}(n):=0$) 
depends only on our being able to show that  one has: 
\begin{equation*} 
\frac12 + \left\lfloor e^{n\varepsilon}\right\rfloor 
-\frac{\varepsilon^2 e^{n\varepsilon}}{4\sinh^2 (\varepsilon /2)} 
= K\left( 1 , \delta^n\right) + \frac{e^{n\varepsilon}}{4\sinh^2 (\varepsilon /2)}
\sum_{r=1}^{\infty} \frac{c_r \varepsilon^r}{r!}\;.
\end{equation*} 
This we can do by noting that $\frac12 + \lfloor e^{n\varepsilon}\rfloor = K( 1 , \delta^n ) + e^{n\varepsilon}$, 
and that one has 
$4\sinh^2 (\varepsilon /2) - \varepsilon^2 = 2(\cosh(\varepsilon) - 1 - \frac12 \varepsilon^2) 
= 2(\frac{1}{4!}\varepsilon^4 + \frac{1}{6!}\varepsilon^6 + \ldots\ )$. 
\par 
It follows that for $R\geq 2$ one has \eqref{M_nXACT-1alt} with 
$E_R (n) := (2\sinh(\varepsilon / 2))^{-2} e^{\varepsilon n} \sum_{r>R}  
(a_r(n) - b_r(n) + c_r)\varepsilon^r /r!$, so that if  $R\geq 2$ is even then 
$\varepsilon^2 e^{-\varepsilon n} |E_R (n)| 
< \sum_{r = R+1}^{\infty} | a_r (n) - b_r (n)  | \varepsilon^r /r!  + \sum_{r=R+2}^{\infty} 2  \varepsilon^r /r!$. 
The estimate \eqref{TruncErrBound} follows from this by noting that one has  the uniform bound 
$|a_r (n) - b_r (n)|\leq\sum_{|\log(m) - \varepsilon n|<\varepsilon} m^{-1}\leq e^{-(n-1)\varepsilon} + 2\varepsilon\,$ 
($r\in{\mathbb N}$), 
and then applying the inequality  
$(s + t)! \geq (s!)(t!)\,$ (for integers $s,t\geq 0$) and the equality 
$\sum_{s=0}^{\infty} \varepsilon^s /s! = e^{\varepsilon}$. 
\end{proof} 

\begin{remarks} 
Lemma~\mbox{B.3} gives us much better control of rounding errors than Lemma~\mbox{B.1}, but  
this control does become unsatisfactory when $N$ is large and 
$n$ is near the top of its range, $\{1,2,\ldots, 2N -1\}$. 
We improve control of the rounding error in such cases by using 
Lemmas~\mbox{B.5} and~\mbox{B.6} (below), in place of Lemma~\mbox{B.3}: see Section~\mbox{B.2} for further details. 
\end{remarks} 

\begin{lemma}
Let $n,N\in{\mathbb N}$ satisfy $N\geq 3$ and $1\leq n\leq 2N - 1$. 
Put 
\begin{equation*} 
\alpha_r := \sum_{2\leq m\leq r}\frac1m \qquad\text{($r \in {\mathbb N}$)} ,  
\end{equation*} 
so that, in particular, $\alpha_1 = 0$.
Then, for $R\in{\mathbb N}$, one has: 
\begin{equation}\label{M_n-IBP} 
- \left( 2 \sinh (\varepsilon / 2)\right)^2  M_n = 
\sum_{3\leq r\leq R} \frac{\alpha_{r-1}}{r}\cdot t_r (n) e^{-(r-1)\varepsilon n} + \varepsilon\cdot U_R (n) e^{-(R-1)\varepsilon n} \;, 
\end{equation} 
where, for $r\in{\mathbb N}$,  
\begin{equation*} 
t_r (n) := e^{r\varepsilon}  \widetilde B_r \left( e^{(n - 1)\varepsilon}\right) 
+ e^{-r\varepsilon}  \widetilde B_r \left( e^{(n + 1)\varepsilon}\right) 
- 2 \widetilde B_r \left( e^{n\varepsilon}\right)  
\end{equation*} 
and 
\begin{equation*} 
U_r (n)  := \int_{-1}^1 \left( \alpha_r \cdot\frac{\phi}{|\phi|} + \varepsilon (1 - |\phi|)\right) 
e^{-r\varepsilon\phi} \widetilde B_r \left( e^{(n + \phi)\varepsilon}\right) d\phi \;, 
\end{equation*} 
with $\widetilde B_r (x)$ denoting the $r$-th periodic Bernoulli function. 
\end{lemma} 

\begin{proof} 
Similarly to how \eqref{BWg1} was obtained, one can show that 
\begin{equation*} 
\delta^{-1} \int_{\delta}^1 \int_{\delta}^1  K\left( 1 , \delta^{n - 1} u v\right)  du dv 
= \varepsilon^2 \int_{-1}^1 \left( 1 - |\phi|\right) \left( {\textstyle\frac12} - \left\{ e^{(n + \phi)\varepsilon}\right\}\right) e^{-\varepsilon \phi} d\phi \;. 
\end{equation*}
Thus, recalling \eqref{m_ijHankelformed}, we find that 
\begin{equation*} 
-\frac{(1 - \delta)^2}{\delta} \cdot M_n 
= \varepsilon^2 \int_{-1}^1 \left( 1 - |\phi|\right) e^{-\varepsilon \phi} \widetilde B_1\left( e^{(n + \phi)\varepsilon}\right)  d\phi 
= \varepsilon\cdot U_1 (n)\;. 
\end{equation*} 
This shows (since $\delta = e^{-\varepsilon}$) that we have \eqref{M_n-IBP} for $R=1$. 
\par 
Suppose now that $R\in{\mathbb N}$. Through a straightforward application of  integration by parts 
(splitting up the integral at the point $\phi = 0$ 
and then exploiting the fact that $\frac{d}{dx} \widetilde B_{R+1} (x) = (R + 1)\widetilde B_R (x)$ for $x\in{\mathbb R}$) 
one may establish that 
\begin{equation*} 
U_R (n) =  \left( \frac{\alpha_R}{(R + 1)\varepsilon} \cdot t_{R+1} (n) + U_{R + 1} (n) \right) e^{-\varepsilon n} \;. 
\end{equation*} 
The cases $R=2,3,4,\ldots\ $ of \eqref{M_n-IBP} therefore follow by induction (and the observation 
that $\alpha_1 = 0$) from the case 
$R=1$ that was verified above.
\end{proof} 

\begin{lemma} 
Let $n,N\in{\mathbb N}$ satisfy $N\geq 3$ and $1\leq n\leq 2N - 1$. 
Let the sequences $\alpha_1,\alpha_2,\ldots\ $ and $U_1 (n), U_2 (n),\ldots\ $ be as defined in Lemma~\mbox{B.5}. 
Then, for each integer $R \geq 2$, one has  
\begin{equation*} 
\left| U_R (n)\right| 
\leq  e^{R\varepsilon} \cdot \min\left\{ \left( 2\alpha_R + \varepsilon\right) \beta_R \,,\,  
\left( e^{(n-1)\varepsilon} \alpha_R \beta_{R-1} + \alpha_{R-1} \beta_R\right) R \varepsilon\right\} \;, 
\end{equation*} 
where $\beta_r := \max_{0\leq x <1} |\widetilde B_r (x)|\,$ ($r\in{\mathbb N}$). 
\end{lemma} 

\begin{proof} 
Let $R$ be a positive integer. 
We have 
\begin{equation*} 
U_R (n) = \alpha_R 
\int_{-1}^1  (\phi /|\phi|) f_R (\phi) d\phi +  \varepsilon \int_{-1}^1  (1 - |\phi|)  f_R (\phi) d\phi  
= \alpha_R I_1 + \varepsilon I_2\quad\text{(say)} ,  
\end{equation*} 
where $f_R (\phi) := e^{-R\varepsilon\phi} \widetilde B_R (e^{(n + \phi)\varepsilon})$. 
Since $|f_R(\phi)| \leq e^{R\varepsilon} \beta_R$ for $\phi\geq -1$, one trivially has 
both $|I_1|\leq 2 e^{R\varepsilon} \beta_R$ and $|I_2|\leq C e^{R\varepsilon} \beta_R$, 
with $C = \int_0^1 |1 -|\phi|| d\phi =  \int_{-1}^1  (1 - |\phi|) d\phi = 1$. 
The bound $|U_R (n)| \leq  (2\alpha_R + \varepsilon) e^{R\varepsilon} \beta_R$ 
follows (given that we have $\alpha_R, \varepsilon\geq 0$). 
\par 
Suppose now that $R \neq 1$. 
We have $I_1 = \int_0^1 f_R (\phi) d\phi - \int_{-1}^0 f_R (\phi) d\phi = \int_0^1 (f_R (\phi) - f_R (-\phi)) d\phi 
= \int_0^1 \int_{-\phi}^{\phi} f_R' (\theta) d\theta d\phi = \int_{-1}^1 (1 - |\theta|) f_R' (\theta) d\theta$. 
Since $R\geq 2$, we have here $f_R' (\theta) =  (e^{n\varepsilon} f_{R-1}(\theta) - f_R (\theta)) R\varepsilon\,$ 
almost everywhere in the interval $[-1, 1]$. We therefore find that 
$I_1 = (e^{n\varepsilon} I_3 - I_2) R\varepsilon$, where $I_2$ is as above, while 
$I_3$ is the similar integral $\int_{-1}^1 (1 - |\phi|) f_{R-1} (\phi) d\phi$. 
It follows that we have 
\begin{equation*} 
U_R (n) = \alpha_R  (e^{n\varepsilon} I_3 - I_2) R\varepsilon 
+ \varepsilon I_2 = (e^{n\varepsilon}  \alpha_R I_3 - (\alpha_R - \frac1R) I_2 ) R\varepsilon\;. 
\end{equation*}  
Here $\alpha_R - \frac1R = \alpha_{R-1}$, since $R\geq 2$. 
We showed earlier that $|I_2|\leq e^{R\varepsilon} \beta_R$: 
one has, similarly,  $|I_3|\leq e^{(R - 1)\varepsilon} \beta_{R - 1}$. 
Thus are able to deduce that 
$|U_R (n)| \leq 
( e^{n\varepsilon}  \alpha_R  e^{(R - 1)\varepsilon} \beta_{R - 1} + \alpha_{R-1} e^{R\varepsilon} \beta_R) 
R\varepsilon$. This, combined with what was found in the last paragraph, 
gives the desired result. 
\end{proof} 

\begin{remarks} 
For our numerical work  we need only the case $R=4$ of the results in the last two lemmas. 
In applying \eqref{M_n-IBP} we note that \eqref{epsilonDef} gives us 
\begin{align*} 
\left( 2 \sinh (\varepsilon / 2)\right)^2 
 &= {\textstyle\frac12} e^{-2N\varepsilon} \left( 1 + \sqrt{ 1 + {\textstyle\frac14} e^{-2N\varepsilon}}\right)^{-1} \\ 
 &= {\textstyle\frac12} e^{-2N\varepsilon} \left( 2 + {\textstyle\frac18} e^{-2N\varepsilon}\right)^{-1} \cdot 
 \left( 1 + {\textstyle\frac{\theta}{256}} e^{-4N\varepsilon}\right) \;, 
\end{align*} 
for some $\theta = \theta (N) \in (0, 1)$. Noting also that 
$|M_n|\leq \frac12\,$  (by \eqref{m_ijHankelformed} and \eqref{DefK}), 
we find that \eqref{M_n-IBP}  (for $R=4$) yields: 
\begin{equation}\label{M_n-IBPshort} 
M_n = - \left( 1 + {\textstyle\frac{1}{16}} e^{-2N\varepsilon}\right) \cdot {\textstyle\frac23} t_3(n)  e^{-(2n - 2N)\varepsilon}  
-{\textstyle\frac56} t_4(n) e^{-(3n - 2N)\varepsilon} + E_4^{*} (n)\;, 
\end{equation} 
where the term $E_4^{*}(n)$ satisfies  
\begin{equation}\label{M_n-IBPerr}   
\left| E_4^{*}(n)\right| 
< {\textstyle\frac{5}{96}} \left| t_4(n)\right| e^{-3n\varepsilon} 
+ \left( 1 + {\textstyle\frac{1}{16}} e^{-2N\varepsilon}\right)  \cdot 4\varepsilon\left| U_4 (n)\right| e^{-(3n - 2N)\varepsilon}  
+ {\textstyle\frac{1}{512}} e^{-4N\varepsilon} 
\end{equation} 
for $n = 1,\ldots , 2N - 1$. 
\end{remarks} 

\subsection{The machine computations} 

We have applied Lemmas~\mbox{B.3}, \mbox{B.5} and~\mbox{B.6} in computing, 
for each $N$ in the geometric sequence $16,32,\ldots, 2^{21}$, and 
each positive integer $n\leq 2N-1$, a number $H_n'$ that us a useful approximation 
to $H_n$.
We got a desktop computer (running GNU~Octave) to perform this computation 
for us. The main steps of the computation can be summarised as follows. 
\begin{description}
\item[Step~1.] A simple iterative algorithm and Octave's {\tt interval} package  are used 
to compute a short real interval  $[\varepsilon', \varepsilon'']$  containing the number $\varepsilon = \varepsilon(N)$. 
\item[Step~2.] Bounds for the numbers $M_n\,$  ($1\leq n\leq 2N-1$) are computed as follows. 
Firstly, after making a suitable choice of $R$, we use \eqref{TruncErrBound}, 
the {\tt interval} package and the result of Step~$1$, in order to compute, for $1\leq n < 2N$, 
real intervals ${\mathcal D}_R (n)$ and ${\mathcal E}_R (n)$ containing (respectively) 
the numbers  $D_R (n)$ and $E_R (n)$ occurring on the right-hand side of equation \eqref{M_nXACT-1alt}. 
Using this data and the {\tt interval} package, we compute $A_n,B_n$ with 
$A_n \leq M_n \leq B_n\,$ ($1\leq n\leq 2N-1$). 
For $N=2^{19}$, and for $N=2^{21}$, there are a few values of $n$ that require special 
treatment, due to there being at least one of 
the three integers $\lfloor e^{(n-1)\varepsilon}\rfloor,\lfloor e^{n\varepsilon}\rfloor,\lfloor e^{(n+1)\varepsilon}\rfloor$ 
whose value is not uniquely determined by our computation. 
This makes it more troublesome than in other cases to estimate the terms $a_r (n)$, $b_r(n)$ and $K(1,\delta^n)$ occurring in Lemma~\mbox{B.3}.  
In such `exceptional' cases we simply put $A_n = -\frac12$ and $B_n = \frac12\,$ 
(the bounds $-\frac12\leq M_n \leq \frac12$ being trivially valid). 
\item[Step~2$'$.] Using \eqref{M_n-IBPshort}, \eqref{M_n-IBPerr}, Lemma~\mbox{B.6}  
(combined with the relations $\beta_3 =  \sqrt{3} / 36$ and 
$\beta_4 = \frac{1}{30}$), the {\tt interval} package and the result of Step~$1$,
we compute, for $1 \leq n < 2N$, supplementary bounds 
$A_n', B_n'$ satisfying $A_n'\leq M_n\leq B_n'$  . 
This does not require any special treatment of the above mentioned exceptional cases, 
as the functions $\widetilde B_3 (x)$ and $\widetilde B_4 (x)$ are continuously differentiable 
and periodic on ${\mathbb R}$. 
\item[Step~3.] The results of Steps~$2$ and~$2'$ are combined, by computing both 
$A_n'' := \max\{ A_n , A_n' \}$ and $B_n'' := \min\{ B_n , B_n' \}$. We then have   
$[A_n'' ,B_n'']\ni M_n\,$ ($1\leq n\leq 2N-1$). 
\item[Step~4.] Using the relations \eqref{H_nDef} and \eqref{deltaDef}, 
the {\tt interval} package 
and the results of Steps~$1$ and~$3$, we compute a short real interval ${\mathcal H}_n$ containing $H_n$. We then put 
\begin{equation}\label{H_n'-numeric}
H_n' := {\rm mid}\left( {\mathcal H}_n\right)\qquad\text{($1\leq n\leq 2N-1$)} ,
\end{equation} 
where `mid()' is the function that the {\tt interval} package provides for computing a
double precision approximation to the midpoint of an interval. 
\end{description} 

\begin{remarks}[assuming $N\in\{2^4,2^5,\ldots,2^{21}\}$ and $2N > n\in{\mathbb N}$] 

\item{\it 1)}\quad From the result of Step~$1$ one can obtain a double precision 
(binary64) approximation to $\varepsilon(N)$ involving a relative error 
not exceeding $5{\tt u} = 5\times 2^{-53}$: for $N\geq 128$ the relative error does not exceed 
$2{\tt u}$. Note however that it is the interval $[\varepsilon', \varepsilon'']$  itself 
(and not any numerical approximation derived from it) that serves as the basis for Steps~\mbox{$2$--$4$}.

\item{\it 2)}\quad It turns out that, after Steps~\mbox{$1$--$3$}, we have $[A_n',B_n']\subseteq [A_n,B_n]\,$  
(and so $A_n'' = A_n'$ and $B_n'' = B_n'$) 
whenever $N\geq 2^{13}$ and $n>\frac32 N$. 
Therefore, given that the exceptional cases (in Step~$2$) occur only when one has both  
$N\geq 2^{19}$ and $n > \frac74 N$, 
we are confident that our suboptimal treatment 
of those cases does not adversely affect 
the bounds $A_n'',B_n''$ ultimately obtained. 

\item{\it 3)}\quad Steps~$2$ and~$2'$ are independent of one another, so that these 
two steps may be completed simultaneously.

\begin{table}[hb]
\centering
\begin{minipage}{14cm} 
\resizebox{14cm}{!}{ 
\begin{tabular}[c]{|r|c|c|c|} 
\hline 
$N$ & $\| k_N'\|$ & $E(N)$ & $F(N)$ \\ \hline 
$2^4$ & $0.119604480526698$ & $1.126430816675218\times 10^{-14}$ & $2.592592886256919\times 10^{-1}$ \\ \hline
$2^5$ & $0.156662714376511$ & $1.637762169439114\times 10^{-14}$ & $2.386994018110340\times 10^{-1}$ \\ \hline
$2^6$ & $0.187352381209700$ & $3.005587409809924\times 10^{-14}$ & $2.154523050603484\times 10^{-1}$ \\ \hline
$2^7$ & $0.210787762259596$ & $4.673947250048874\times 10^{-14}$ & $1.925853830963158\times 10^{-1}$ \\ \hline
$2^8$ & $0.231407573306031$ & $8.333956270747470\times 10^{-14}$ & $1.672457638249101\times 10^{-1}$ \\ \hline
$2^9$ & $0.246906311948199$ & $1.994182752555530\times 10^{-13}$ & $1.433802065171640\times 10^{-1}$ \\ \hline
$2^{10}$ & $0.258849641090915$ & $3.068598934835102\times 10^{-13}$ & $1.204884799799176\times 10^{-1}$ \\ \hline
$2^{11}$ & $0.267380429163074$ & $5.635821559916214\times 10^{-13}$ & $1.001414829196205\times 10^{-1}$ \\ \hline 
$2^{12}$ & $0.273492166462582$ & $5.808121879145198\times 10^{-13}$ & $8.199173973257347\times 10^{-2}$ \\ \hline
$2^{13}$ & $0.277662301306867$ & $7.381750011358651\times 10^{-13}$ & $6.651508801579206\times 10^{-2}$ \\ \hline
$2^{14}$ & $0.280486652055700$ & $8.943962706048118\times 10^{-13}$ & $5.336523699339534\times 10^{-2}$ \\ \hline
$2^{15}$ & $0.282343337060482$ & $1.030037000364939\times 10^{-12}$ & $4.245998726249360\times 10^{-2}$ \\ \hline
$2^{16}$ & $0.283543576117694$ & $1.938335491393713\times 10^{-12}$ & $3.352090307840420\times 10^{-2}$ \\ \hline
$2^{17}$ & $0.284302064005348$ & $1.765470239629683\times 10^{-12}$ & $2.632388465054856\times 10^{-2}$ \\ \hline
$2^{18}$ & $0.284779315906853$ & $2.237110376240634\times 10^{-12}$ & $2.052685393399503\times 10^{-2}$ \\ \hline
$2^{19}$ & $0.285073640586082$ & $4.031995570417358\times 10^{-12}$ & $1.592576352417425\times 10^{-2}$ \\ \hline
$2^{20}$ & $0.285253128240220$ & $4.370383787783687\times 10^{-12}$ & $1.229891580548921\times 10^{-2}$ \\ \hline
$2^{21}$ & $0.285361694024823$ & $6.222827840613366\times 10^{-12}$ & $9.450612896218688\times 10^{-3}$ \\ \hline
$+\infty$ & $0.285518143908160$ & {\bf --} & {\bf --} \\ 
\hline 
\end{tabular} 
}
\vskip 3mm
\ Table~B-1 
\end{minipage} 
\end{table}

\item{\it 4)}\quad  Since the interval ${\mathcal H}_n$ computed in Step~$4$  always contains $H_n$, it follows that 
we always have $|H_n' - H_n|\leq \sup\{ | x - H_n' | : x\in {\mathcal H}_n\}$. Thus, with the help of the 
{\tt interval} package, a useful upper bound for $|H_n' - H_n|$ can be computed. 
This, together with what is noted in Remarks~\mbox{6.11}~\mbox{(1)--(3)} (including, in particular, the 
inequality \eqref{HSnormBound}), enables us to compute 
both a fairly sharp upper bound $F(N)$ for $\| k_N' - K\|$ and  
a useful upper bound $E(N)$ for $\| H'(N) - H(N)\| = \| k_N' - k_N\|$. 
Using just the data $H_1',\ldots ,H_{2N-1}'$, and the {\tt interval} package, we get also 
a good approximation to the number $\| H'(N) \| = \| k_N'\|\,$ 
(one accurate enough to determine the result of rounding $\| k_N'\|$ to 15 significant digits).  
The results obtained 
are shown in Table~\mbox{B-1} (above). 

\item{\it 5)}\quad The number $0.2855\ldots\ $ at the bottom of the second column  of Table~\mbox{B-1} 
is $\| K\|$ rounded to  15 significant digits (see \eqref{Gamma_h-HSnormEval} and Remarks~\mbox{6.11}~\mbox{(2)}). 
Other than that, what is shown 
in the second column of Table~\mbox{B-1} is  $\| k_N'\|$ rounded to 15 significant digits. 

\end{remarks}

\begin{refsection}[refsC.bib]

\section{Fast estimation of $H'(N) X$} 

In this appendix we describe how we compute  
matrix-vector or matrix-matrix products of the form $H'(N) X$,  
where $H'(N)$ is the $N\times N$ Hankel matrix discussed in Section~\mbox{5.2}   
(note that we shall assume that the numerical values of the elements of 
$H'(N)$ and $X$ are already known: 
see Appendix~B regarding the computation of  the elements of $H'(N)$). 
In order to rapidly  compute these products we have used 
an Octave function {\tt fast\_hmm()} that 
implements an algorithm of Luk and Qiao \cite[Section~4, Algorithm~2]{LQ2000}.  
We obtained this function by converting a Matlab function that was authored by S. Qiao. 
The Luk-Qiao algorithm utilises fast discrete Fourier transforms (FFTs)  
and makes it possible to compute 
the product of a vector and an $N\times N$ real or complex Hankel matrix in time $O(N\log N)$. 
Without this sort of speed much of our computational work 
(particularly on the cases with $N \geq 2^{14}$) would not have been feasible. 
\par 
The function {\tt fast\_hmm()} takes the following as input data:  
the matrix $X\,$  (which is passed to it as an argument), and the column vector $(H_1',\ldots , H_N')^{\rm T}$ and row  
vector $(H_N',\ldots ,H_{2N-1}')\,$  (passed to it as global variables {\tt C} and {\tt R}, respectively). 
To denote the result returned by this function 
we use the notation ${\tt fast\_hmm}(X)\,$  (i.e. we treat both 
{\tt C} and {\tt R} as constants, since they depend only on $N$). 
If the columns of $X$ (ordered from left to right) are ${\bf x}_1,\ldots , {\bf x}_w\in{\mathbb C}^N$, then 
the columns of the computed result ${\tt fast\_hmm}(X)$ are: ${\tt fast\_hmm}({\bf x}_1),\ldots , {\tt fast\_hmm}({\bf x}_w)$. 
We therefore confine further discussion of the function {\tt fast\_hmm()} to cases in which the input data $X$ 
is some vector `${\bf x}$' (i.e. an $N\times 1$ matrix). 
Since the computations carried out (internally) by the function {\tt fast\_hmm()} involve the use 
of double precision floating-point arithmetic, this function cannot be relied upon to 
compute the product ${\bf y}=H'(N){\bf x}$ exactly: the result $\tilde{\bf y} = {\tt fast\_hmm}({\bf x})$ 
that it returns is affected by rounding errors, and so is (in general) only an approximation to ${\bf y}$. 

\subsection{Fast Hankel matrix-vector multiplication} 

Like S. Qiao's Matlab function, from which it derives,  
the function {\tt fast\_hmm()} is an application of the following result 
concerning Hankel matrix-vector products $H {\bf x}$ and the discrete Fourier transform 
${\mathcal F}$ defined on ${\mathbb C}^{2N}$ by putting 
${\mathcal F}((z_1,\ldots ,z_{2N})^{\rm T}) = (\hat z_1,\ldots ,\hat z_{2N})^{\rm T}$, where    
$\hat z_j := \sum_{k=1}^{2N} z_k \exp(-2\pi i (j-1)(k-1)/(2N))$ for $1\leq j\leq 2N$.
\begin{lemma}[The Luk-Qiao Algorithm] 
Let $N\in{\mathbb N}$. Suppose that ${\bf x} = (x_1,\ldots , x_N)^{\rm T}\in {\mathbb C}^N$, that 
$H$ is a (real or complex) Hankel matrix with 
last row $(R_1,\ldots ,R_N)$
and first column 
$(C_1,\ldots ,C_{N-1}, R_1)^{\rm T}$, and that 
\begin{equation*}
{\mathcal F}^{-1}( {\mathcal F}({\bf p}) \mathop{.*} {\mathcal F}({\bf q})) = (Y_1,\ldots , Y_{2N})^{\rm T} \in{\mathbb C}^{2N}\;, 
\end{equation*} 
where ${\bf p} = (R_1,\ldots , R_N, 0, C_1, \ldots , C_{N-1})^{\rm T}\in {\mathbb C}^{2N}$,  
${\bf q} = (x_N,x_{N-1}, \ldots , x_1, 0, 0, \ldots , 0)^{\rm T}\in {\mathbb C}^{2N}$,   
and ${\bf u}\mathop{.*}{\bf v} := (u_1 v_1, \ldots , u_{2N} v_{2N})^{\rm T}\in{\mathbb C}^{2N}\,$ (${\bf u},{\bf v}\in{\mathbb C}^{2N}$).  
Then the product $H{\bf x}$ and the vector ${\bf y} := (Y_1,\ldots , Y_N)^{\rm T}\in{\mathbb C}^N$ are equal. 
\end{lemma}
\begin{proof} 
In \cite[Section~4]{LQ2000} Luk and Qiao sketch the proof of a result that is similar to this, but 
involves discrete Fourier transforms on ${\mathbb C}^{2N-1}$. 
We need only a minor modification of that proof, in which ${\bf p}$ is substituted for 
the vector `$\hat{\bf c}$' specified in \cite[Equation~(10)]{LQ2000} 
(so that the circulant `$C(\hat{\bf c})$' defined in \cite{LQ2000}  becomes a $2N\times 2N$ matrix). 
\end{proof}
In order to achieve an efficient algorithmic application of Lemma~\mbox{C.1} one wants a means 
of computing the discrete Fourier transform $\mathcal F$ and its inverse both quickly and accurately: 
our function {\tt fast\_hmm()} utilises, for this purpose,   the `fast Fourier transform' function {\tt fft()}    
and its inverse {\tt ifft()}, which are built-in functions of Octave 
(in this we have again followed Qiao, who used the corresponding Matlab functions). 
\par 
The Hankel matrix-vector products $H{\bf x}$ that concern us are real 
(we always have ${\bf x}\in{\mathbb R}^N$ and $H = H'(N)$, which is a real Hankel matrix), 
so that in each case the corresponding vector 
${\bf y} = (Y_1,\ldots , Y_N)^{\rm T}\,$ (defined as in Lemma~\mbox{C.1}) is real also. 
The same is not quite true of the approximations to $Y_1,\ldots , Y_N$ 
computed by the function {\tt fast\_hmm()}, 
for amongst these approximations, $\tilde Y_1,\ldots , \tilde Y_N\,$   (say), 
there may be some with a non-zero imaginary part. 
Therefore, 
instead of simply returning the complex vector $(\tilde Y_1, \ldots , \tilde Y_N)^{\rm T}$ as its result, 
our Octave function  {\tt fast\_hmm()} returns the real vector 
$\tilde{\bf y} = ( {\rm Re}(\tilde Y_1), \ldots , {\rm Re}(\tilde Y_N))^{\rm T}$. 
Since ${\bf y}$ is real, this last adjustment does not make the 
function {\tt fast\_hmm()} any less accurate than it would otherwise be:
it does, however, mean that this function requires that all of its input data 
(the matrix $X$, or vector ${\bf x}$, and global variables {\tt C} and {\tt R}) be real, 
whereas the original Matlab function of S.~Qiao is not so restricted in scope.

\subsection{A rounding error analysis} 

In this section we discuss how accurate the function {\tt fast\_hmm()} is,   
when used  to compute a Hankel matrix-vector product of the special form $H'(N){\bf x}$, with  
${\bf x}\in{\mathbb R}^N$. The relevant global variables {\tt C} and {\tt R} are assumed to be equal to 
$(H_1',\ldots , H_N')^{\rm T}$ and $(H_N',\ldots ,H_{2N-1}')$, respectively.  
We suppose also that the vector ${\bf p}$ is defined as in the relevant case of Lemma~\mbox{C.1}, 
so that one has ${\bf p} = (H_N',\ldots ,H_{2N-1}', 0,H_1',\ldots , H_{N-1}')^{\rm T}\in{\mathbb R}^{2N}$.
\par 
We shall assume, initially, that $N = 2^n$ for some integer $n$ with $4\leq n\leq 41$.  
Note that we need $N$ to be a (positive integer) power of $2$ 
when applying \eqref{Higham-1} (below), although there are of course some weaker alternatives to 
\eqref{Higham-1} that could be applied, were $N$ not of this form. 
Ultimately we consider just the cases with $n := \log_2(N)\in \{ 12,11,\ldots , 21\}$, which 
are all that we need for the computational work described in Section~\mbox{D.5} (see, in particular, Remarks~\mbox{D.6} there). 
We need the error analysis of {\tt fast\_hmm()} just so that we can compute an appropriate value for 
the factor $\Delta(A,g) > 0$ occurring in the bound \eqref{DieCast-2}.  
\par 
Before considering the function {\tt fast\_hmm()}, we first need to discuss the 
accuracy of Octave's built-in functions {\tt fft()} and {\tt ifft()}.
\par 
Let ${\mathcal B}(D)$ denote (when $D\in{\mathbb N}$)  the set of all vectors 
${\bf r} = (r_1,\ldots ,r_D)^{\rm T}\in{\mathbb C}^D \backslash\{ {\bf 0}\}$ 
such that the real and imaginary parts of $r_1,\ldots ,r_D\in{\mathbb C}$ 
are representable in binary64 format. 
We shall assume that the error bound \cite[Theorem~24.2]{Hi2002} 
for the Cooley-Tukey FFT algorithm is applicable to the results that we get from Octave's {\tt fft()} function.   
That is, we assume that, for the relevant vectors ${\bf r}\in{\mathcal B}(2N) = {\mathcal B}(2^{n+1})$, one has
\begin{equation}\label{Higham-1}
\| {\tt fft}({\bf r}) - {\mathcal F}({\bf r})\| / \| {\mathcal F}({\bf r})\| 
\leq d \eta_d / (1 - d\eta_d)\;, 
\end{equation}
with 
\begin{equation*} 
d := n+1\quad \text{and}\quad \eta_d := \mu_d + \gamma_4\cdot  (\sqrt{2} + \mu_d)\;, 
\end{equation*}
where 
$\mu_d$ is the maximum modulus of 
the errors in the estimates for roots of unity that are used (as `weights' or `twiddle factors') in computing the term 
${\tt fft}({\bf r})\in{\mathcal B}(2^d)$, while  
$\gamma_4 = 4{\tt u} / (1-4{\tt u})\,$ (with ${\tt u} := 2^{-53}$). 
Although we have not been able to verify this assumption, we do have reasons for   
believing that it is not far from the truth. 
In particular, Octave's {\tt ver()} function tells us that the results returned by 
the function {\tt fft()} are computed using version 3.3.8-sse2 of the FFTW library, which  
(according to \cite[Section~II]{FJ2005}) utilises primarily the Cooley-Tukey FFT algorithm. 
\par
Methods of computing approximations to roots of unity  (for use in the Cooley-Tukey algorithm) are discussed 
in \cite[Section~3]{TZ2001}: see there the `Algorithms 3.1--3.3', in particular. 
It appears, from tests we have carried out on the {\tt fft()} function, 
that the function $d\mapsto\mu_d$ is approximately linear. 
Of the algorithms considered in \cite{TZ2001} only one (`Algorithm~3.3': also known as `subvector scaling') 
is consistent with this linearity. 
Assuming that subvector scaling is indeed the method used to precompute approximate roots 
of unity  for the {\tt fft()}  function, 
it follows (see \cite[(4.7)--(4.9), or (3.3)]{TZ2001}) that one has 
\begin{equation}\label{UnityErrBound}
\mu_d \leq c_{\scriptscriptstyle RSS}  {\tt u} d\quad\ \text{($0\leq d\leq 64$)} 
\end{equation}
with $c_{\scriptscriptstyle RSS} :={\textstyle\frac{4}{\sqrt{3}} + \frac{1}{\sqrt{2}}}$.  
\par
We combine \eqref{UnityErrBound} with Higham's theorem \eqref{Higham-1}, getting: 
\begin{equation}\label{Higham-2}
\| {\tt fft}({\bf r}) - {\mathcal F}({\bf r})\| / \| {\mathcal F}({\bf r})\| 
\leq \left( c_{\scriptscriptstyle RSS} d^2 + 6 d\right) {\tt u}
\end{equation}
for the relevant vectors ${\bf r}\in{\mathcal B}(2^d)$.  
This implies a similar bound for  $\| {\tt ifft}({\bf s}) - {\mathcal F}^{-1}({\bf s})\| / \| {\mathcal F}^{-1}({\bf s})\|$.  
The case $d=n+1$ of these `mean square' error bounds (for {\tt fft()} and {\tt ifft()}) can be shown to imply that  
\begin{equation}\label{fast_hmmErrBound-1}
\left\| \tilde{\bf y}  - H'(N){\bf x}\right\| 
\leq c_{\scriptscriptstyle RSS}{\tt u}\cdot\left( {\textstyle\frac32} + c_{\scriptscriptstyle RSS}{\tt u}\cdot (n+2)^2\right) (n+2)^2 
\left( \| {\bf p}\|_1 \| {\bf x}\| + \| {\bf p}\| \| {\bf x}\|_1\right) 
\end{equation}
when one has $N = 2^n$, ${\bf x}\in{\mathcal B}(N)\cap{\mathbb R}^N$ and $\tilde{\bf y} := {\tt fast\_hmm}({\bf x})$. 
\par 
It is convenient to simplify \eqref{fast_hmmErrBound-1}, using the fact that  
$\| {\bf z}\|_1 \leq m^{1/2} \| {\bf z}\|\,$ ($m\in{\mathbb N}$, ${\bf z}\in{\mathbb C}^m$). 
Using also the (empirical) observation that one has 
\begin{equation*} 
\| {\bf p}\| / \|  H'(N)\| < {\textstyle \frac{21}{125}} N^{-\frac12} \left(  {\textstyle \frac{28}{3}} + \log N\right) 
\quad\ \text{($4\leq n\leq 21$)} , 
\end{equation*}
we deduce from \eqref{fast_hmmErrBound-1} that 
if $n\in\{12, 13,\ldots ,21\}$ and $N=2^n$ then    
\begin{align} 
\frac{\left\| \tilde{\bf y} - {\bf y}\right\|}{\| {\bf x}\|} &\leq 
(1 + 1064{\tt u})(1 + \sqrt{2}) \cdot {\textstyle\frac{63}{250}} c_{\scriptscriptstyle RSS} {\tt u} \cdot (2+n)^2 
\left( {\textstyle\frac{28}{3}} + \log N\right) \| H'(N)\| \nonumber \\ 
 &< {\textstyle \frac53}  {\tt u} N  \| H'(N)\| \label{fast_hmmErrBound-2}
\end{align}
when  ${\bf x}\in{\mathcal B}(N)\cap{\mathbb R}^N$, 
${\bf y} =  H'(N) {\bf x}\in{\mathbb R}^N$ and 
$\tilde{\bf y} = {\tt fast\_hmm}( {\bf x})\in{\mathcal B}(N)$. 
\begin{remarks}

\item{\it 1)}\quad In view of the explicitly stated assumptions, \eqref{Higham-1} and \eqref{UnityErrBound}, it is  
clear that we do not have a complete proof of \eqref{fast_hmmErrBound-2}. 
We fall short of this because we cannot be certain that the error-bound \eqref{Higham-2} for the {\tt fft()} function is correct.
However, the results of a benchmarking of the FFTW~3.1 library by the authors of \cite{FJ2005} 
include nothing that would suggest \eqref{Higham-2} is false: see, in particular \cite{benchFFT} for the accuracy benchmark results obtained  
with a 3.1 GHz Intel Xeon E3-1220v3 four core processor (similar to the one in our computer).  Therefore we are reasonably confident 
of the correctness of \eqref{Higham-2} and \eqref{fast_hmmErrBound-2}, in the context of our work. 
\par 
We rely on \eqref{fast_hmmErrBound-2} only when computing our `probable bounds':  
it is used in Step~3 of the statistics-based algorithm described in Section~\mbox{D.5}, and nowhere else. 
In particular, with the exception of the conditional bounds $\tilde{\mathcal U}_m$, $\tilde{\mathcal U}^{\pm}_m$ and $\tilde\mu_m\,$ 
(briefly discussed in Remarks~\mbox{6.14}~\mbox{(5)}), the results we describe in Section~\mbox{6.5} 
are obtained independently of \eqref{fast_hmmErrBound-2} and \eqref{Higham-2}. 

\item{\it 2)}\quad The above error analysis of the function {\tt fast\_hmm()} depends also on 
two other assumptions that have not yet been mentioned:  
we have assumed that `overflows' (where a binary64 `exponent' exceeds $1023$) do not 
occur in the course of the relevant computations, 
and that `underflows' (where an exponent is less than $-1022$) 
are a negligible source of rounding errors.  
\par 
When overflows occur at some point in a computation, it is (usually) easy to spot that this 
has happened: one need only check the result of the computation for the presence 
of `infinite' or `undefined' quantities (in Octave these appear as an `{\tt Inf}', a `{\tt -Inf}', 
or a `{\tt NaN}'). 
\par 
In considering the effect of underflows on accumulated rounding error, we first 
 examine whether they might make a significant revision of \eqref{Higham-2} necessary. 
Assuming that ${\tt fft}({\bf r})$ is computed by the Cooley-Tukey FFT algorithm 
(or a method not too different from this), one can show that if $ {\bf r}\in {\mathcal B}(2^d)$ and 
$\| {\bf r}\| \geq {\tt u}^{-1} \cdot 2^{d - 1022}\,$ (say) 
then the contribution to $\| {\tt fft}({\bf r}) - {\mathcal F}({\bf r})\|$ coming from underflow errors 
will be insignificant, when compared with the upper bound for  
$\| {\tt fft}({\bf r}) - {\mathcal F}({\bf r})\|$ implied by \eqref{Higham-2}. 
As a result of this finding, together with the fact that $\| {\bf p}\| \geq N^{-1/2} \| H'(N)\|$ and the observation that 
$\| H'(N)\|\in (\frac14 , \frac27)\,$ (when $12\leq n\leq 21$), 
we can conclude that if  ${\bf x}\in {\mathcal B}(N)\cap {\mathbb R}^N$ satisfies  
$\| {\bf x} \| \geq 2^{-968} N\,$  (say)  then any effect that underflow errors might have on the result  
$\tilde{\bf y} = {\tt fast\_hmm}( {\bf x})$ 
will be too insignificant to cause \eqref{fast_hmmErrBound-2} to become invalid. 

\item{\it 3)}\quad We chose to make \eqref{fast_hmmErrBound-2} similar in form to  
the case $p=F\,$ (for `Frobenius norm') of the result stated 2 lines below \cite[(3.13)]{Hi2002}, which is 
a commonly used bound for the rounding error when matrices are multiplied together 
using finite precision floating-point arithmetic. 
Thus \eqref{fast_hmmErrBound-2} is a compromise that favours simplicity over strength:  
in the case $N=2^{21}$, for example, it is a much weaker bound than 
that which precedes it. This compromise is ultimately quite harmless:  
for it turns out that we always have $(\sum_{j=1}^{2M} \alpha_j^2)^{1/2} > \frac78 \| H'(N)\|$ when applying \eqref{tilde_y_errbound}--\eqref{Delta(A,g)Formula}. 

\end{remarks} 

\subsection{Estimating $\| (H'(2^{21}))^2 \|$}

We require, for an application discussed in Section~\mbox{D.3}, 
a sharp numerical upper bound for the Frobenius norm of the matrix $(\eta_{i,j}) := \left( H'(2^{21})\right)^2$.  
We describe here the iterative method (suited to our limited computing resources) 
by which we obtain the required result.  Note that what is actually obtained is an estimate for 
$\| (H'(N))^2 \|\,$ ($N=2^{21}$), along with a bound for the size of the error: the upper bound we that we seek 
is a corollary. Although our method does apply more generally, what mainly concerns us here is its application 
to the single case where $N=2^{21}$. Thus, 
until we come to the (more general) Remarks~\mbox{C.3}, below, 
the letters $n$ and $N$ will denote the numbers $21$ and $2^{21}$, respectively.  
\par 
Our method involves a total of $T=2^{14}$ steps, with the $t$-th step 
yielding an $N\times N/T$ matrix $P(t)$ of double precision numbers $p_{i,j}(t)$
such that, for $r=1,\ldots,N/T$, 
the $r$-th column of $P(t)$ is approximately equal to the $(t+(r-1)T)$-th column of 
the matrix $(\eta_{i,j})$. 
The computed data $P(t)$ gets used (immediately) in 
computing an approximate value $\sigma(t)$ for the sum 
$\sum_{i=1}^N\sum_{j=1}^{N/T} (p_{i,j}(t))^2$: aside from that it is used only in computing $P(t+1)$, 
so that the same working memory ($2$ gigabytes) can be used 
to store each of $P(2), P(3), \ldots , P(T)$  in turn.  
The matrix $P(1)$ is an exception here, since the data 
contained in its first column is utilised in every one of the $T$ steps of the method. 
\par 
We compute $P(1)$ using the FFT-based Octave function 
{\tt fast\_hmm()} discussed in the previous sections of this appendix. 
Then, using a modified `interval arithmetic' version of 
{\tt fast\_hmm()} that makes extensive use of the {\tt interval} package,   
upper bounds $W_{i,j}$  for the absolute values of the associated approximation errors 
$p_{i,j}(1) - \eta_{i,1+(j-1)T}\,$  ($1\leq i\leq N$, $1\leq j\leq N/T$) 
are obtained. 
\par
In computing $P(2), P(3),\ldots, P(T)$ we exploit the structure 
of the matrix $(\eta_{i,j})$. We use, in particular, 
the facts that one has 
\begin{equation*} 
\eta_{i, j} = \eta_{i-1, j-1} + H_{i+N-1}' H_{j+N-1}' - H_{i-1}' H_{j-1}'\quad\ \text{($2\leq i,j\leq N$)} 
\end{equation*} 
and
\begin{equation*} 
\eta_{1, j} = \eta_{j ,1}\quad\ \text{($2\leq j\leq N$)} . 
\end{equation*}
This enables rapid computation of $P(t)\,$ (for $t>1$), 
using the previously computed data $P(t-1)$, $P(1)$ 
and $H_1',\allowbreak\ldots,\allowbreak H_{2N-1}'$. 
\par 
After the final ($T$-th) step, a number $\Sigma =\Sigma(N,T)$ 
approximating the sum $\sum_{t=1}^T \sigma(t)$ is computed: 
this $\Sigma$ serves as our best approximation to the true value 
of $\| \left( H'(N)\right)^2\|^2$. 
In computing $P(1), \sigma(1),\ldots , P(T), \sigma(T)$ and $\Sigma$ 
we use only double precision arithmetic: use of the {\tt interval} package 
there would  greatly increase the running time needed. 
We  therefore take care that the computation of the sums $\sigma(1),\ldots , \sigma(T)$ 
and $\Sigma$ (in particular) is done in a way that gives us good control of rounding errors. 
\par 
We have carried out an analysis of the potential effect of rounding errors 
in the above computations. We omit the details of this analysis, noting only that 
it shows that, for some $\theta\in [-1,1]$, one has  
\begin{equation}\label{RoundErrControl}
\| \left( H'(N)\right)^2\| 
= \sqrt{\Sigma}\cdot\exp\left( (1 + n){\tt u}\theta\right)  
+ \left(\| X\| + \| Y\|\right)\theta\;, 
\end{equation} 
where $X$ and $Y\,$ (certain $N\times N$ matrices of error terms) are such that  
\begin{equation} \label{FrobXbound}
\| X\|^2 < T \sum_{i=1}^N {\sum_{1\leq j\leq N/T}}^{\!\!\!\!\!\!*}\ \left( W_{i,j} + 
\left( \frac{2T{\tt u}}{1 - 2T{\tt u}} \right) \left| p_{i,j}(1)\right| \right)^2 
\end{equation}
(with the asterisk signifying that the summand is to be doubled when $j=1$) and 
\begin{align}\label{FrobYbound}
\| Y\|^2 &< \frac{e^{\varepsilon / 2}}{4} \left( \frac{2T{\tt u}}{1 - 2T{\tt u}} \right)^2 
\left( 1 + \frac{1}{2(e^{T\varepsilon} - 1)}\right) \nonumber \\ 
 &\phantom{{<}} \quad \times \left( \sinh\left( \frac{(2T - 1)\varepsilon}{2}\right) 
- (2T - 1) \sinh\left( \frac{\varepsilon}{2}\right) \right)  
\end{align} 
(with $\varepsilon = \varepsilon(N)$ as defined in Section~\mbox{3.1}).  
The bound for $\| Y\|^2$ depends on the fact that  
$|H_j'| \leq \frac12 (1 - e^{-\varepsilon}) e^{-(j-1)\varepsilon /2}\,$ 
($1\leq j < 2N$). 
\par 
Using \eqref{FrobXbound}, \eqref{FrobYbound} and 
the relevant numerical data obtained while computing $\Sigma$, 
we find that $\| X\| < 1.4336964\times 10^{-12}$ and 
$\| Y\| <  4.994\times 10^{-14}$.   
The computed value of $\Sigma$, rounded to 16 significant digits, is:  
$1.443938638781244\times 10^{-4}$. 
Using \eqref{RoundErrControl}, and our data on $\Sigma$, $X$ and $Y$, we deduce that  
$\| \left( H'(N)\right)^2\|$ lies in the real interval with 
midpoint $1.20163997885\times 10^{-2}$ and 
length $3.1\times 10^{-12}$. We obtain, in particular, the bound  
\begin{equation}\label{FroHsquared21}
\| \left( H'(2^{21})\right)^2\|^2 < 1.4439386391378115481049 \times 10^{-4} \;, 
\end{equation} 
for use in Section~\mbox{D.3}. 
\begin{remarks}

\item{\it 1)}\quad Note that $\| (H'(2^{21}))^2\|$ is quite a good approximation for $\| (H(2^{21}))^2\|$. Indeed, 
by virtue of the inequality $\| PQ\| \leq \min\{ \| P\|_2 \| Q\| , \| P\| \| Q\|_2\}\,$ (valid when $PQ$ exists), 
one has 
$\| (H'(N))^2 - (H(N))^2\| \leq ( \| H'(N)\|_2 +  \| H(N)\|_2  ) \| H'(N) - H(N)\|$ 
for $N\geq 3\,$ 
(regardless of whether or not $H'(N)$ and $H(N)$ commute);   
since we have here $\| H'(N)\|_2 = |\varkappa_1 (k_N')|$,  $\| H(N)\|_2 = |\varkappa_1 (k_N)|$ 
and $\| H'(N) - H(N)\|_2 \leq \| H'(N) - H(N)\| = \| k_N' - k_N\|$, 
we therefore can deduce (by \eqref{RELBlemma-v2}, \eqref{Weyl-4} and the triangle inequality) that 
\begin{equation}\label{H'-to-H_transfer}
\left| \| (H'(N))^2\|  - \| (H(N))^2\| \right| \leq \left( 2\,{\mathcal U}_1 + E(N)\right) E(N)\quad\ 
\text{($N\in\{ 2^4, 2^5, \ldots , 2^{21}\}$)} ,
\end{equation}
with ${\mathcal U}_1$ and $E(N)$ here denoting, respectively,  
the upper bounds for  $|\varkappa_1 |$ and $\| k_N' - k_N\|$ 
whose computation is discussed in Section~\mbox{6.5} and Remarks~\mbox{B.8}~\mbox{(4)}. 
Noting (see Tables~1 and~\mbox{B-1}) that ${\mathcal U}_1 < 1/12.4611 < 0.0802498$, 
and that $E(N) < 6.23\times 10^{-12}$ in the relevant cases, 
we deduce that the factor $2\,{\mathcal U}_1 + E(N)$ in \eqref{H'-to-H_transfer} 
is never greater than  $0.1605$. In particular, \eqref{H'-to-H_transfer} 
gives us: 
$\bigl| \| (H'(N))^2\|  - \| (H(N))^2\| \bigr| \leq 0.1605 \times 6.23\times 10^{-12} < 10^{-12}$. 

\item{\it 2)}\quad In cases where $N\leq 2^{13}$ the computation  
of an accurate estimate for the Frobenius norm of $(H'(N))^2$ is a relatively easy task, 
compared to what we found it necessary to do for $N=2^{21}$. In such cases we are able to store 
the entire $N\times N$ matrix $H'(N)$ in our machine's 16 gigabytes of RAM,  
as a matrix of intervals (all of length $0$), so that the computation of upper and lower bounds 
for $\| (H'(N))^2\|$ requires nothing more than a 
straightforward utilisation of the {\tt interval} package. The same is true 
even in the case $N=2^{14}$, but in this particular case we find it 
convenient to shorten the time required for the computation by using the  
interval arithmetic version of the function {\tt fast\_hmm()} that was mentioned earlier in this section. 
For $N\in\{ 2^{17}, 2^{18}, 2^{19}, 2^{20}\}$ we resort to the  
method whose use in the case $N=2^{21}$ is described above 
(only changing the parameters $N$ and $T$).
For $N=2^{15}$ and $N=2^{16}$ we get the best results with a version of the same method (that used for $N=2^{21}$)  
in which there is additional use made of the {\tt interval} package 
and of our interval arithmetic version of the function {\tt fast\_hmm()}:  
the bounds \eqref{FrobXbound} and \eqref{FrobYbound} are not relevant here, 
since each step of the computation takes automatic account of rounding errors. 
\par 
For $4\leq n\leq 21$ and $N=2^n$ we obtain, as a result of the work just mentioned (and 
the work, described earlier, on the case $n=21$), a specific 
short interval $[\phi_N, \phi_N']$ in which $\| (H'(N))^2 \|$ lies. Combining these 
results with \eqref{H'-to-H_transfer}, we get: 
\begin{equation}\label{HH_FroNorm_interval}
\big| \| (H(N))^2 \| - \Phi_N\bigr| \leq \Delta_N\quad\  \text{($4\leq n\leq 21$, $N=2^n$)} , 
\end{equation} 
where $\Phi_N$ and $\Delta_N$ are (respectively) 
a double precision approximation to $\frac12 (\phi_N' + \phi_N)$ and an upper bound for 
$\max\{ \Phi_N - \phi_N \,,\, \phi_N' - \Phi_N\} +  \left( 2\,{\mathcal U}_1 + E(N)\right) E(N)\,$  
(both computed using Octave's {\tt interval} package). 
In Tables~\mbox{C-1a} and~\mbox{C-1b} (below) the number $f_n$ is a decimal approximation 
to $\Phi_N$, while $d_n$ satisfies $d_n \geq \Delta_N + | f_n - \Phi_N|$. 
Thus, by \eqref{HH_FroNorm_interval}, we have $f_n + d_n \geq \| (H(2^n))^2\| \geq f_n - d_n$ for each 
$n$ in these two tables. For $n\geq 15$, where the method employed in computing $[\phi_N , \phi_N']$ 
is similar to that used for $N=2^{21}$, we include in Table~\mbox{C-1b} the relevant 
specification of the parameter $T=T(N)=T(2^n)$. 

\begin{table}[b!]
\centering 
\begin{minipage}{79mm}
\resizebox{!}{28mm}{
\begin{tabular}[c]{|r|c|c|} 
\hline 
$n$ & $f_n$ & $d_n$ \\ \hline 
$4$ & $0.007134544407958$ & $2.3\times 10^{-15}$ \\ \hline 
$5$ & $0.009031755502642$ & $2.8\times 10^{-15}$ \\ \hline 
$6$ & $0.010339967816155$ & $5.4\times 10^{-15}$ \\ \hline 
$7$ & $0.011144642216076$ & $7.7\times 10^{-15}$ \\ \hline 
$8$ & $0.011604013577423$ & $1.4\times 10^{-14}$ \\ \hline 
$9$ & $0.011837726495411$ & $3.3\times 10^{-14}$ \\ \hline 
$10$ & $0.011945067130672$ & $5.0\times 10^{-14}$ \\ \hline 
$11$ & $0.011989683583481$ & $9.1\times 10^{-14}$ \\ \hline 
$12$ & $0.012006969963481$ & $9.4\times 10^{-14}$ \\ \hline 
\end{tabular}
} 
\vskip 3mm
\ Table~C-1a 
\end{minipage}
\begin{minipage}{84mm}
\resizebox{!}{28mm}{ 
\begin{tabular}[c]{|r|c|c|c|}
\hline  
$n$ & $f_n$ & $d_n$ & $T(N)$ \\ \hline 
$13$ & $0.012013220211127$ & $1.2\times 10^{-13}$ & {\bf --} \\ \hline 
$14$ & $0.012015371155722$ & $1.7\times 10^{-13}$ & {\bf --} \\ \hline 
$15$ & $0.012016077594365$ & $2.0\times 10^{-13}$ & $2^5$ \\ \hline 
$16$ & $0.012016301628868$ & $3.6\times 10^{-13}$ & $2^7$ \\ \hline 
$17$ & $0.012016370622190$ & $5.8\times 10^{-13}$ & $2^{10}$ \\ \hline 
$18$ & $0.012016391379146$ & $8.4\times 10^{-13}$ & $2^{12}$ \\ \hline 
$19$ & $0.012016397501435$ & $1.4\times 10^{-12}$ & $2^{13}$ \\ \hline 
$20$ & $0.012016399279054$ & $1.8\times 10^{-12}$ & $2^{14}$ \\ \hline 
$21$ & $0.012016399788544$ & $2.5\times 10^{-12}$ & $2^{14}$ \\ \hline 
\end{tabular} 
} 
\vskip 3mm
\ Table~C-1b 
\end{minipage} 
\end{table} 

An analysis if the data $\Phi_N,\Delta_N\,$ ($N\in\{ 2^{4}, 2^{5}, \ldots , 2^{21}\}$) 
leads us to conjecture that  
\begin{equation}\label{T4hypothesis_1}
D(N) := \| (H(2N))^2\|^2 - \| (H(N))^2\|^2 \sim \eta_0 N^{-2} \log^{3/2} N 
\quad\text{as $N\rightarrow \infty$} ,
\end{equation} 
where $\eta_0$ is a positive constant (approximately $\frac{21}{34}$). 
Considering, in particular, the sequence $V_4, V_5, \ldots \ $ given by  
\begin{equation*} 
V_n := \left( N^2 D(N)\right)^2 \quad\ \text{($N=2^n$, $n\geq 4$)} , 
\end{equation*}
we make use of the data $\Phi_N,\Delta_N$  in computing real intervals 
${\mathcal V_n}\ni V_n\,$ ($4\leq n\leq 20$). 
An initial examination of the results
(including some work on finite differences of up to the fourth order) reveals that,  
for a suitably chosen constant $m \approx\frac18$, the polynomial $m(n-9)^3$ serves as a 
fairly good approximation to $V_n$ when $10\leq n\leq 20$. Further work (studying a consistent bias 
in the error associated with this approximation) leads us to conjecture that 
there exist constants $m_1\approx 0.127$ and $c_1\approx 0.37$ such that one has 
\begin{equation}\label{Conjecture_RE_V} 
V_n = m_1 (n - 9)^3 + c_1 (n-9)\log(n) + O(n) 
\end{equation}
for all integers $n\geq 10$. This conjecture, if correct, would imply immediately that one has  
\eqref{T4hypothesis_1}, with $\eta_0 = m_1^{1/2} / (\log 2)^{3/2} \approx\frac{21}{34}$. 
\par 
The conjecture \eqref{Conjecture_RE_V}  is in one respect guesswork, since 
we have no theoretical or heuristic explanation for the factors of form $n-9$ that appear in it: we 
should perhaps conjecture only that $V_n = m_1 (n^3 - 3 d_1 n^2) + c_1 n\log(n) + O(n)$, 
where $d_1$ is some constant close to $9$.  Our confidence 
in the conjecture \eqref{Conjecture_RE_V}  derives in part from 
similarities between it and the conditional asymptotic formula 
for $\| K - k_N\|^2 = \| K\|^2 - \| k_N\|^2 = \| K\|^2 - \| H(N)\|^2$ 
that is mentioned in Remarks~\mbox{3.6}, in the lines below \eqref {Uncond_Frob_Asymp}. 

\renewcommand{\thefigure}{C-1} 

\begin{figure}[b!]
\begin{center}
\includegraphics*[scale=0.9, bb=48 200 549 593]{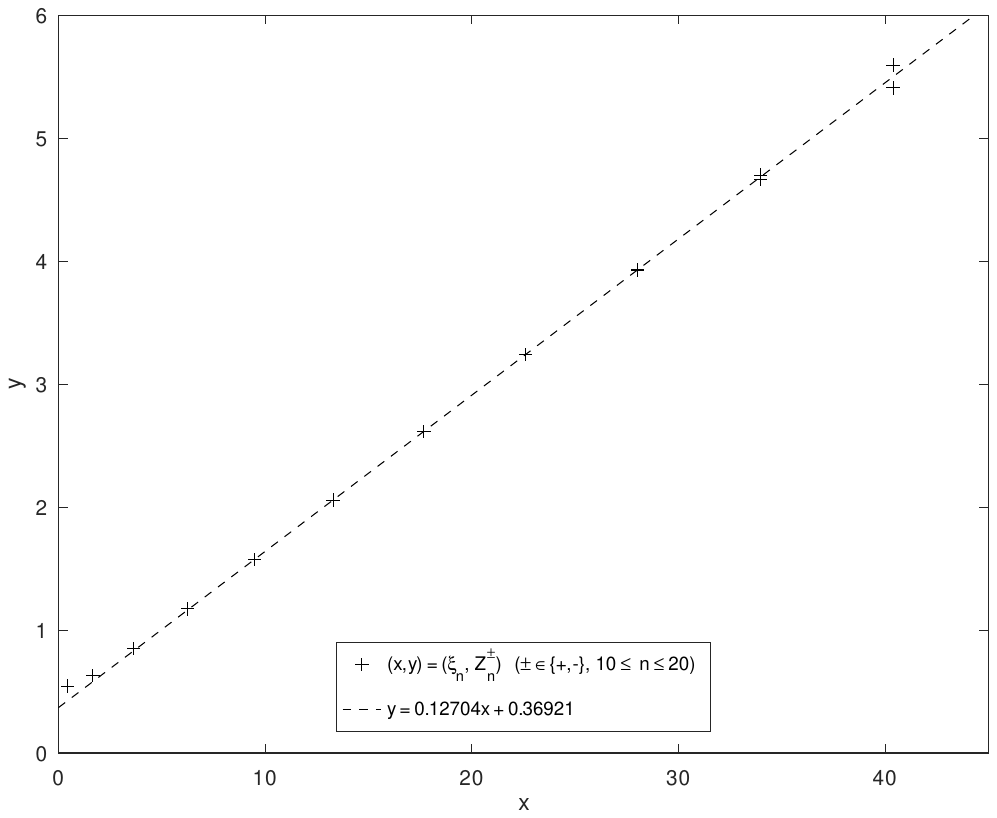} 	
\end{center}
\caption{\small Evidence supporting Conjecture \eqref{Conjecture_RE_V}} 
\label{fig:HH-FroNorm}
\end{figure}

Our strongest evidence in favour of the conjecture \eqref{Conjecture_RE_V} 
is depicted in the graph in Figure~\ref{fig:HH-FroNorm} (below), where we have plotted points 
$\left( \xi_n , Z^{-}_n\right)$ and 
$\left( \xi_n , Z^{+}_n\right)\,$ ($10\leq n\leq 20$), with: 
\begin{equation}\label{Z-braces} 
Z^{-}_n := \frac{\inf{\mathcal V}_n}{(n-9)\log(n)}\leq Z_n := \frac{V_n}{(n-9)\log(n)} 
\leq \frac{\sup{\mathcal V}_n}{(n-9)\log(n)} =: Z^{+}_n 
\end{equation} 
and  $\xi_n := (n-9)^2 / \log(n)$. The line $y=0.12704 x + 0.36921$ shown in Figure~\ref{fig:HH-FroNorm} 
is an approximation to the `line of best fit' (as determined by the ordinary least squares method) 
for the 7 points $\left( \xi_{14} , \frac12(Z^{-}_{14}+Z^{+}_{14})\right) , \ldots , \left( \xi_{20} , \frac12(Z^{-}_{20}+Z^{+}_{20})\right)$,   
which were selected (for this purpose) 
because the sequences $\frac12(Z^{-}_{n}+Z^{+}_{n})\,$ ($14\leq n\leq 20$) 
and $\xi_n\,$ ($14\leq n\leq 20$) are highly correlated, with a correlation coefficient greater than $0.9999997$: we 
remark that, by using Octave's {\tt interval} package,  one can establish that the correlation coefficient of the 
pair of sequences $\xi_{n}\,$ ($14\leq n\leq 18$) and $Z_n\,$ ($14\leq n\leq 18$) is greater than $0.9964$. 
After subjecting the data set $\Phi_N,\Delta_N\,$ ($N\in\{ 2^{10}, 2^{11}, \ldots , 2^{21}\}$) to relevant further 
processing and analysis 
(the details of which we omit, for the sake of brevity), we are satisfied that 
the distribution of the points $(\xi_n, Z^{\pm}_n)$ shown in Figure~\ref{fig:HH-FroNorm} is consistent with  
the hypothesis that there are constants $m_1\approx 0.127$ and $c_1\approx  0.37$ such that 
\begin{equation}\label{Conjecture_RE_V_Equiv} 
Z_n = m_1\xi_n + c_1 + O\left(\frac{1}{\log(n)}\right) \quad\ \text{($n\geq 10$)}  , 
\end{equation}
when $Z_n$ and $\xi_n$ are defined as in, and below, \eqref{Z-braces}. 
This hypothesis is, of course, equivalent to the conjecture \eqref{Conjecture_RE_V}. 
Regrettably, our data set $\Phi_N,\Delta_N\,$ ($N\in\{ 2^{10}, 2^{11}, \ldots , 2^{21}\}$ 
is not of sufficient scope or quality to provide us with truly convincing evidence in support 
of the hypothesis \eqref{Conjecture_RE_V_Equiv}, so all we can say with any confidence is that 
the conjecture \eqref{Conjecture_RE_V} does fit our limited empirical data remarkably well 
(considering the simplicity of the conjecture). 
 
\item{\it 3)}\quad With regard to our conjecture \eqref{T4hypothesis_1}, note that one certainly does have:   
\begin{equation}\label{K2monotonic_approx} 
0 < \| K_2\|^2 - \| (H(N))^2\|^2 \ll N^{-1}\log^3 N  
\quad\text{($N\geq 3$)} ,
\end{equation}
where $K_2$ is the `iterated kernel' given by 
\begin{equation*} 
K_2(x,y) = \int_0^1 K(x,z)K(z,y) dz\qquad \text{($0\leq x,y\leq 1$)} . 
\end{equation*} 
Indeed, since the reciprocal eigenvalues of $K_2$ are precisely 
$\varkappa_1^2, \varkappa_2^2, \varkappa_3^2, \ldots\ $, 
while the eigenvalues of $(H(N))^2$ are 
$(\varkappa (k_N))^2, \ldots , (\varkappa_N (k_N))^2\,$ 
(we use here the notation defined in Section~\mbox{5.1}), 
it follows by \eqref{TraceLemma-k} and \eqref{RELBlemma-v2}  
that one has 
\begin{align}\label{TrDiffIdentity}
0 < \| K_2\|^2 - \| (H(N))^2\|^2 &= \sum_{m=1}^{\infty} \left( \varkappa_m^4 - (\varkappa_m (k_N))^4 \right) \\ 
&\leq \sum_{m=1}^{\infty} \left(  \varkappa_1^2 + (\varkappa_1 (k_N))^2\right) 
\left( \varkappa_m^2 - (\varkappa_m (k_N))^2 \right) \nonumber\\ 
&= \left(  \varkappa_1^2 + (\varkappa_1 (k_N))^2\right)  
\left( \| K \|^2 - \| k_N\|^2\right) \;, \nonumber
\end{align}
which implies \eqref{K2monotonic_approx}, as we have here 
$|\varkappa_1 (k_N)|\leq |\varkappa_1|$ and (by \eqref{ParallelAxis-2} and 
what is established in the proof of Proposition~\mbox{1.2}, in Section~3)  
$\,\| K \|^2 - \| k_N\|^2 = \| K - k_N\|^2<\frac14 N^{-1}\log^3 N$ for all sufficiently large positive integers $N$. 
\par 
In view of \eqref{K2monotonic_approx},   
it follows that if our conjecture \eqref{T4hypothesis_1} is correct then 
\begin{equation*}
\| K_2\|^2 - \| (H(N))^2\|^2 \sim {\textstyle\frac43}\eta_0 N^{-2} \log^{\frac32} N \quad\ \text{as $N\rightarrow\infty$} , 
\end{equation*} 
so that (by  \eqref{RELBlemma-v2} and \eqref{TrDiffIdentity}) one has:  
\begin{align}
0 \leq |\varkappa_m| - |\varkappa_m (k_N)|
 &\leq \frac{\| K_2\|^2 - \| (H(N))^2\|^2}{\left(  \varkappa_m^2 + (\varkappa_m (k_N))^2 \right)\left(  |\varkappa_m| + |\varkappa_m (k_N)| \right) } \label{pre-HyperConverge} \\ 
 &\ll \left( |\varkappa_m| + |\varkappa_m (k_N)| \right)^{-3} N^{-2} \log^{\frac32} N \nonumber 
\end{align} 
for $m\geq 1$, $N\geq 3$. Finally, it is worth noting that the inequalites in the single line \eqref{pre-HyperConverge} are independent of 
our conjectures (i.e. they hold unconditionally) and might be used to 
get new (or improved) numerical estimates for singular values of $K$, provided that one had sufficiently accurate 
numerical estimates for $\| K_2\|$, $\| (H(N))^2\|$ and relevant eigenvalues of $H(N)$.

\end{remarks}

\printbibliography 

\end{refsection}

\begin{refsection}[refsD.bib]

\section{Intervals containing eigenvalues of $H'(N)$}

In this appendix we discuss certain numerical computations enabling us to check the accuracy 
of the approximate eigenvalues  
$\alpha_j\in{\mathbb R}\,$ ($1\leq j\leq N'$) first mentioned in Section~\mbox{5.3};   
we employ much of the notation introduced in Section~5,  
but shall also use some alternative notation (where that is more convenient): 
recall, in particular, that each of  
$\varkappa^{+}_1(k_N'), \ldots , \varkappa^{+}_{P'}(k_N')$ 
and $\varkappa^{-}_1(k_N'), \ldots , \varkappa^{-}_{Q'}(k_N')$ is both  
a reciprocal eigenvalue of $k_N'$ and an eigenvalue of $H'(N)$. 
\par 
We deal exclusively with cases in which $N$ and $n$ are integers with $N = 2^n$ 
and $4\leq n\leq 21$. 
The data utilised includes both the approximate eigenvalues themselves and 
the corresponding approximate eigenvectors ${\bf v}_j\in{\mathbb R}^N$ ($1\leq j\leq N'$).
\par 
As a consequence of Weyl's inequalities for eigenvalues of Hermitian matrices, we have: 
\begin{equation}\label{CrudeWeyl}
\left| \Lambda_j' - \Lambda_j'' \right| \leq \| A\|_2
\quad \text{($1\leq j\leq N$)} , 
\end{equation}
where $\Lambda_1'\geq\Lambda_2'\geq  \ldots \geq \Lambda_N'\,$ 
(resp.  $\Lambda_1''\geq\Lambda_2''\geq  \ldots \geq \Lambda_N''$)  
are the eigenvalues of $H'(N)\,$ (resp. $H''(N,N')$), 
while $H''(N,N')$ and $A$ are the real symmetric matrices defined in \eqref{Def-H''} and \eqref{Def-A(N)}. 
The sequence $\Lambda_1',\ldots ,\Lambda_N'$ here is a  permutation of $\varkappa_1(k_N'), \ldots , \varkappa_N(k_N')$. 
Since $A$ is real and symmetric, we have also: 
\begin{equation}\label{A_to_Asquared}
\sqrt{\| A^2\|} \geq\| A\|_2 \;. 
\end{equation}

\subsection{The cases with $n\leq 9$} 

Recall from Section~\mbox{5.3} that when $4\leq n\leq 9$ we have $N' := p' + q' = N$, so that there are  
precisely as many terms in the sequence $\alpha_1,\ldots ,\alpha_{N'}$  
as there are eigenvalues of $H'(N)$. 
Using \eqref{A_to_Asquared} and a numerical bound for $\| A^2\|$, 
we can get an upper bound for $\| A\|_2$ that is adequate for our purposes 
(specifically, its application in \eqref{CrudeWeyl}) when $n\leq 9$. 
We first compute an $N\times N$ matrix of intervals, ${\mathcal B}\,$ (say), 
such that each element of ${\mathcal B}$ contains the corresponding element of the matrix $A^2$;  
we use the data ${\mathcal B}$ to compute a sharp upper bound for $\| A^2\|$, and so,  
through \eqref{A_to_Asquared}, we obtain our bound for $\| A\|_2$. 
In carrying out these computations we utilise Octave's {\tt interval} package:  
owing to the moderate sizes of the relevant matrices, 
run times and usage of working memory fall within acceptable limits. 
We shall not give the results in full, but will note  
that they imply $\| A\|_2 < {\tt u}\cdot (n-3)  < 10^{-15}$ for $4\leq n\leq 9$. 
\par 
Since $\| A\|_2$ is so very small (when $n\leq 9$), it follows by \eqref{CrudeWeyl} that 
we may get a good idea of the locations of the eigenvalues of $H'(N)$ by 
determining (as best we can) where the eigenvalues of $H''(N, N')$ lie. 
\par 
By its definition, in \eqref{Def-H''}, the matrix $H''(N, N')$ is real and symmetric: 
one has, in particular, $H''(N, N') = V D V^{\rm T}$, where (since we have $N' = N$ when $4\leq n\leq 9$)
$\, V$ is the $N\times N$ matrix with columns ${\bf v}_1,\ldots , {\bf v}_N$, 
while $D$ is the $N\times N$ diagonal matrix $(d_{i,j})$ with $d_{i,i}=\alpha_i\,$ 
($1\leq i\leq N$). It turns out that when $4\leq n\leq 9$ the matrix $V$ is nearly orthogonal, and so,   
by the computation of an upper bound for $\| V^{\rm T} V - I_N\|_1$, it 
can readily be established that $V$ is invertible. 
The matrices $H''(N, N')=V D V^{\rm T}$ and $DV^{\rm T}V$ are therefore similar,   
and so have the same $N$ eigenvalues: $\Lambda_1''\geq\Lambda_2''\geq  \ldots \geq \Lambda_N''$. 
\par 
By virtue of the near-orthogonality of $V$, the matrix $DV^{\rm T}V$ 
is nearly diagonal when $4\leq n\leq 9\,$ (i.e. those of its elements that lie on the main diagonal are large in comparison with 
those that do not). Because of this, we can determine explicitly (via a computation 
implementing the Gershgorin disc theorem \cite{Ge1931}) 
a union of $N$ very short intervals, $[x_i, y_i]\ni \alpha_i\,$ ($1\leq i\leq N$), 
that contains every eigenvalue of the matrix $DV^{\rm T}V$. 
It can be checked that $x_i > y_{i+1}$ for $1\leq i < N$, so that (by Gershgorin's theorem, again) 
we may conclude that $[x_i,y_i]\ni \Lambda_i''$ for $1\leq i \leq N$. 
We therefore find, using also \eqref{CrudeWeyl}, that for $1\leq i\leq N$ the eigenvalue
$\Lambda_i'\,$ (belonging to $H'(N)$) lies in the interval 
$[a_i, b_i] = [x_i - R, y_i + R]$, where $R$ is our numerical upper bound for $\| A\|_2$.
\par
One can check that the widening of $[x_i, y_i]$ to $[x_i - R, y_i + R]\,$ ($1\leq i\leq N$) 
does not lead to overlaps: 
one always ends up having $a_i > b_{i+1}$ for $1\leq i < N$, so that the 
`strict ordering' of the intervals is preserved. We therefore conclude that, 
for $4\leq n\leq 9$, all eigenvalues of $H'(N)$ are simple. 
One can check also that $a_{p'} > 0 > b_{p' + 1}$, so that $H'(N)$ is non-singular, and  
has $p'$ positive eigenvalues, $\Lambda_1' > \Lambda_2' > \ldots > \Lambda_{p'}'$ and 
$N - p' = N' - p' = q'$ negative  eigenvalues, $\Lambda_{N}' < \Lambda_{N-1}' < \ldots < \Lambda_{p' + 1}'$. 
Thus, for $4\leq n\leq 9$, 
we have: $P'(N) = p'$, and 
$\varkappa^{+}_i (k_N') \in [a_i, b_i]\,$ ($1\leq i\leq p'$);  
$Q'(N) = q' = N - p'$, and 
$\varkappa^{-}_j (k_N') \in [a_{N+1-j}, b_{N+1-j}]\,$ ($1\leq j\leq q'$).

\par 
Table~\mbox{D-1} (below) presents a sample of our numerical results, focussing on those connected 
with the `middle' eigenvalue $\Lambda_{N/2}'$, which (it turns out) is always one of the smaller eigenvalues of $H'(N)$. 
The values given for $a_{N/2}$, $b_{N/2}$ and $\alpha_{N/2}$ are rounded to $12$ significant digits; 
those given for $R$ and $y_{N/2} - x_{N/2}$ are rounded to just $2$ significant digits. 
Although each `Gershgorin interval' $[x_{N/2}, y_{N/2}]$ in Table~\mbox{D-1}  
has length less than $\frac{1}{20}R$, the same is not true of all the intervals $[x_i, y_i]\,$ ($1\leq i\leq N$): 
one has, in fact, $\max_{1\leq i\leq N} (y_i - x_i) \in (\frac43 R, 2R)$ for $4\leq n\leq 9$. 

\begin{table}[hb] 
\centering
\begin{minipage}{165mm} 
\resizebox{165mm}{!}{  
\begin{tabular}[c]{|r|r|r|c|c|c|c|c|} 
\hline 
$n$ & $ p' $ & $q'$ & $a_{N/2}$ & $\alpha_{N/2}$ & $b_{N/2}$ & $ R$ & $y_{N/2} - x_{N/2}$ \\ \hline 
$4$ & $8$ & $8$ & $2.19270977044\times 10^{-5}$ & $2.19270977045\times 10^{-5}$ & $2.19270977046\times 10^{-5}$ & $1.0\times 10^{-16}$ & $2.3\times 10^{-18}$ \\ \hline
$5$ & $16$ & $16$ & $1.31363689321\times 10^{-5} $ & $1.31363689323\times 10^{-5}$ & $1.31363689325\times 10^{-5}$ & $2.1\times 10^{-16}$ & $1.0\times 10^{-17}$ \\ \hline
$6$ & $32$ & $32$ & $1.07804386572\times 10^{-6}$ & $1.07804386597\times 10^{-6}$ & $1.07804386622\times 10^{-6}$ & $2.5\times 10^{-16}$ & $4.9\times 10^{-18}$ \\ \hline
$7$ & $63$ & $65$ & $-7.30264237532\times 10^{-8}$ & $-7.30264234109\times 10^{-8}$ & $-7.30264230686\times 10^{-8}$ & $3.4\times 10^{-16}$ & $1.4\times 10^{-18}$ \\ \hline
$8$ & $128$ & $128$ & $9.26867206315\times 10^{-9}$ & $9.26867249924\times 10^{-9}$ & $9.26867293533\times 10^{-9}$ & $4.3\times 10^{-16}$ & $5.3\times 10^{-18}$ \\ \hline
$9$ & $255$ & $257$ & $-1.23880253802\times 10^{-9}$ & $-1.23880199194\times 10^{-9}$ & $-1.23880144586\times 10^{-9}$ & $5.5\times 10^{-16}$ & $8.3\times 10^{-19}$ \\ 
\hline
\end{tabular} 
} 
\vskip 3mm
\ Table~D-1  
\end{minipage}
\end{table} 

\par 
Using the data $\alpha_i, a_i, b_i\,$ ($1\leq i\leq N$), we find that  
\begin{equation*}
\max_{1\leq i\leq N} \left| \alpha_i - \Lambda_i'\right| < {\textstyle \frac85} {\tt u}\cdot (n-2) 
\qquad\text{($4\leq n\leq 9$)}  ,
\end{equation*}
and that 
$| \alpha_i - \Lambda_i'|/|\alpha_i| < 4.5\times 10^{-7}\,$ ($1\leq i \leq N$, $4\leq n\leq 9$). 
The ARPACK software library (on which Octave's {\tt eigs()} function relies) must be given the credit for the accuracy of the approximations $\alpha_1,\ldots ,\alpha_N$ when $4\leq n\leq 9$. 

\subsection{Lower bounds for the sizes of eigenvalues, when $n\geq 10$} 

For the remainder of this appendix we focus solely on the cases with $10\leq n\leq 21$. 
Thus, recalling the relevant details from Section~\mbox{5.3}, we may assume  
henceforth that we have $p' = q' = M\,$ (so that 
$\alpha^{+}_j > 0 > \alpha^{-}_j$ for $1\leq j\leq M$).  
\par 
In this section we explain how we have gone about computing satisfactory numerical values 
for the lower bounds $L^{\pm}_1,\ldots , L^{\pm}_M$ in \eqref{LBs_for_eigenmoduli}. 
\par 
Put 
\begin{equation}\label{DefDeltaPlus}
{\bf \Delta}^{+}_j := (\alpha^{+}_j I_N - H'(N))  {\bf v}^{+}_j\quad\text{($1\leq j\leq M$)}  ,
\end{equation}
and let $E^{+} = (e^{+}_{i,j})$ and $F^{+} = (f^{+}_{i,j})$ be the $M\times M$ matrices given by:
\begin{equation*}
E^{+} =  I_M - (V^{+})^{\rm T} V^{+}\qquad\text{and}\qquad F^{+}  
=  D^{+} - (V^{+})^{\rm T} H'(N) V^{+}\;,
\end{equation*}
where $D^{+}$ is the diagonal matrix that has 
$\alpha^{+}_1 , \alpha^{+}_2 , \ldots , \alpha^{+}_M$ 
as the elements along its principal diagonal, while  
$V^{+}$ is the $N\times M$ matrix with columns  ${\bf v}^{+}_1,{\bf v}^{+}_2,\ldots ,{\bf v}^{+}_M$.
For $1\leq m\leq M$ let $E^{+}_m\,$ (resp. $F^{+}_m$) denote the $m$-th leading principal submatrix 
of $E^{+}\,$  (resp. $F^{+}$). 

\begin{lemma}
Let $m\leq M$ be a positive integer. 
Put 
\begin{equation}\label{DefPhi(m)}
\Phi^{+} = \Phi^{+}(m) := \frac{\| F^{+}_m\|_1}{\alpha^{+}_m} + \sum_{1\leq j < m} \frac{\| {\bf \Delta}^{+}_j\|^2}{(\alpha^{+}_j - \alpha^{+}_m)^2} \;.  
\end{equation}
Then one has both 
\begin{equation}\label{New_Lower} 
\varkappa^{+}_m (k_N')\geq \alpha^{+}_m - \frac{\| {\bf  \Delta}^{+}_m \|}{\sqrt{1 - \Phi^{+}}} 
\quad\text{if $\Phi^{+} < 1$} 
\end{equation}
and 
\begin{equation}\label{CrudeLower} 
\varkappa^{+}_m (k_N')\geq \frac{\max\{ 0\,,\,\alpha^{+}_m - \| F^{+}_m\|_1\}}{1 + \| E^{+}_m\|_1} \;. 
\end{equation}
\end{lemma} 
\begin{proof}
Since $\varkappa^{+}_m (k_N')$ is (by definition) positive or zero,  
the results \eqref{New_Lower} and \eqref{CrudeLower} are trivially valid 
in cases where $\varkappa^{+}_m (k_N')\geq \alpha^{+}_m$. 
We therefore have only to show that \eqref{New_Lower} and \eqref{CrudeLower} are 
valid when $0\leq\varkappa^{+}_m (k_N') < \alpha^{+}_m\,$ 
(we assume henceforth that this is the case). 
\par
Let ${\bf u}^{+}_1,\ldots ,{\bf u}^{+}_{P'}$ be some 
orthonormal system of eigenvectors for the positive eigenvalues of $H'(N)$. 
By a suitable choice of ${\bf b}\in{\mathbb R}^m$ with 
$\| {\bf b}\|^2 = b_1^2 + \cdots + b_m^2 = 1$, one can ensure that
the vector ${\bf w} := \sum_{j=1}^m b_j {\bf v}^{+}_j\in{\mathbb R}^N$ 
satisfies  
${\bf w}\cdot {\bf u}^{+}_j = 0\,$ ($1\leq j \leq \min\{ m-1, P'\}$). 
Following this, we put $J_{\bf w}:={\bf w}^{\rm T} H'(N) {\bf w}$. 
Using the spectral decomposition of $H'(N)$, we get: 
\begin{equation}\label{AppSpecDecomp}
J_{\bf w} \leq \sum_{m\leq i\leq P'} \varkappa^{+}_i (k_N') \left( {\bf w} \cdot  {\bf u}^{+}_i\right)^2 
\leq \varkappa^{+}_m (k_N') \sum_{m\leq i\leq P'} \left( {\bf w} \cdot  {\bf u}^{+}_i\right)^2 \;. 
\end{equation} 
We have, at the same time, 
$J_{\bf w} = \sum_{j=1}^m \alpha^{+}_j b_j^2 - {\bf b}^{\rm T} F^{+}_m {\bf b} 
\geq \alpha^{+}_m - \| F^{+}_m\|_2\,$ 
(the last inequality following by virtue of \eqref{alphasOK}, and our choice of ${\bf b}$).  
Since the matrix $F^{+}$ is real and symmetric, it follows that 
one has $J_{\bf w} \geq \alpha^{+}_m - \| F^{+}_m\|_1$. 

\par 
By Bessel's inequality, the last of the sums occurring in \eqref{AppSpecDecomp} is not greater than $\| {\bf w}\|^2$. 
Since $\| {\bf w}\|^2 = {\bf b}^{\rm T} (I_m - E^{+}_m) {\bf b} = 1 -  {\bf b}^{\rm T} E^{+}_m {\bf b} \leq 1 + \| E^{+}_m\|_2$, 
and since $E^{+}$ is real and symmetric, it therefore follows from \eqref{AppSpecDecomp} that 
\begin{equation}\label{EasyBesselApp} 
J_{\bf w} \leq (1 + \| E^{+}_m\|_1)\varkappa^{+}_m (k_N')\;. 
\end{equation} 
This, when combined with our lower bound for $J_{\bf w}$, 
yields the result \eqref{CrudeLower}. 
\par
Our proof of \eqref{New_Lower} rests on the observation that, as a consequence of 
$H'(N)$ being symmetric, one has  
$\left(\alpha^{+}_j - \varkappa^{+}_i (k_N')\right) {\bf v}^{+}_j \cdot {\bf u}^{+}_i = {\bf \Delta}^{+}_j \cdot {\bf u}^{+}_i\,$   
($1\leq i\leq P'$, $1\leq j\leq m$). 
Since $\varkappa^{+}_m (k_N') < \alpha^{+}_m$, it follows that one has 
$|{\bf v}^{+}_j \cdot {\bf u}^{+}_i| \leq |{\bf \Delta}^{+}_j \cdot {\bf u}^{+}_i| / \left(\alpha^{+}_j - \varkappa^{+}_m (k_N')\right)$ 
when $1\leq j\leq m$ and $m\leq i\leq P'$. Using this, one finds that, for $m\leq i\leq P'$, one has 
$| {\bf w} \cdot {\bf u}^{+}_i | \leq \sum_{j=1}^m |b_j {\bf \Delta}^{+}_j \cdot {\bf u}^{+}_i | / (\alpha^{+}_j - \varkappa^{+}_m (k_N'))$, 
and so 
\begin{equation*} 
({\bf w} \cdot {\bf u}^{+}_i )^2 \leq \sum_{j=1}^m \frac{({\bf \Delta}^{+}_j \cdot {\bf u}^{+}_i)^2}{\left(\alpha^{+}_j - \varkappa^{+}_m (k_N')\right)^2}
\end{equation*}
(this last inequality following by the Cauchy-Schwarz inequality, given that $\| {\bf b}\| = 1$). 
By this, \eqref{AppSpecDecomp} and Bessel's inequality, we get: 
\begin {align*}
J_{\bf w} &\leq \varkappa^{+}_m (k_N') \sum_{j=1}^m \left(\alpha^{+}_j 
  - \varkappa^{+}_m (k_N')\right)^{-2} \sum_{m\leq i\leq P'}   ({\bf \Delta}^{+}_j \cdot {\bf u}^{+}_i)^2 \\ 
 &\leq \varkappa^{+}_m (k_N') \sum_{j=1}^m \frac{\| {\bf \Delta}^{+}_j \|^2}{\left(\alpha^{+}_j - \varkappa^{+}_m (k_N')\right)^2} \\ 
 &\leq \alpha^{+}_m \cdot \left( \frac{\| {\bf \Delta}^{+}_m \|^2}{\left(\alpha^{+}_m - \varkappa^{+}_m (k_N')\right)^2} 
 + \sum_{1\leq j < m} \frac{\| {\bf \Delta}^{+}_j\|^2}{(\alpha^{+}_j - \alpha^{+}_m)^2}\right) .
\end{align*} 
Since $\alpha^{+}_m > 0$, 
the result \eqref{New_Lower} follows directly from the combination of last of the above upper bounds for $J_{\bf w}$ 
with the lower bound $J_{\bf w} \geq \alpha^{+}_m - \| F^{+}_m\|_1$. 
\end{proof}

\begin{remarks} 
\item{\it 1)}\quad Similarly, for $1\leq m\leq M$, 
one can define in terms of $H'(N)$, $\alpha^{-}_1 ,  \ldots , \alpha^{-}_M$ 
and ${\bf v}^{-}_1, \ldots ,{\bf v}^{-}_M$ 
upper bounds (conditional and unconditional) for  
$\varkappa^{-}_m (k_N')$ that are analogous 
to the (conditional and unconditional) lower bounds for  
$\varkappa^{+}_m (k_N')$ in \eqref{New_Lower} and \eqref{CrudeLower}. 
\item{\it 2)}\quad It is clear (from the above proof) that \eqref{CrudeLower}  
would remain valid if the $1$-norms, there, were replaced by $2$-norms. 
This would, in fact, yield a stronger lower bound for $\varkappa^{+}_m (k_N')\,$       
(recall that when $S$ is a real symmetric matrix one has  
$\| S\|_1\geq \| S\|_2 = \rho(S)$, the spectral radius of $S$). 
However, we have found the bound \eqref{CrudeLower}  to be useful enough for our purposes,   
and very easy to use (as the relevant $1$-norms are simple to estimate).
\end{remarks} 

Using the {\tt interval} package and 
the interval arithmetic version  of the function {\tt fast\_hmm()} that is mentioned in Section~\mbox{C.3},    
we have computed sharp upper bounds for the Euclidean norms of the 
vectors ${\bf \Delta}^{+}_1,\ldots ,{\bf \Delta}^{+}_M$ defined in \eqref{DefDeltaPlus}. 
We have also computed (using the same tools) two $M\times M$ matrices of intervals with elements 
$[e_{i,j}', e_{i,j}'']\ni e^{+}_{i,j}$ and $[f_{i,j}', f_{i,j}'']\ni f^{+}_{i,j}\,$ ($1\leq i,j\leq M$). 
The data contained in these two matrices was used  in computing, for each $m\leq M$, 
certain sharp upper bounds for the $1$-norms occurring in \eqref{DefPhi(m)} and \eqref{CrudeLower}  
(the bound obtained for $\| E^{+}_m\|_1$, for example, being an  
upper bound for $\max_{1\leq j\leq m} \sum_{i=1}^m \max\{ |e_{i,j}'| , |e_{i,j}''|\}$). 
\par
The upper bounds found for $\| {\bf \Delta}^{+}_j\|$ and $\| F^{+}_j\|_1\,$ ($1\leq j\leq M$) were 
used, together with the data $\alpha^{+}_1 , \alpha^{+}_2 , \ldots , \alpha^{+}_M$, 
in computing sharp upper bounds for the non-negative real 
numbers $\Phi^{+}(1),\ldots , \Phi^{+}(M)$ 
defined in \eqref{DefPhi(m)}: it was found that 
we had $\max_{m\leq M} \Phi^{+}(m) < 1.9 \times 10^{-7}$ for $10\leq n\leq 21$,  
so that the lower bound for $\varkappa^{+}_m (k_N')$ in \eqref{New_Lower} 
was applicable in every case. 
\par 
Finally, for $m\leq M$, we used  the data $\alpha^{+}_m$ and the upper bounds found for 
$\| {\bf \Delta}^{+}_m\|$,  $\Phi^{+}(m)$, $\| E^{+}_m\|_1$ and $\| F^{+}_m\|_1$  
in computing, via \eqref{New_Lower} and \eqref{CrudeLower} (respectively), 
a pair of lower bounds $\ell^{+}_{j,1},\ell^{+}_{j,2}\in (0,\infty)$ for the eigenvalue $\varkappa^{+}_m (k_N')$ of $H'(N)$:  
by making similar use of 
the data $H'(N)$, $\alpha^{-}_1,\ldots ,\alpha^{-}_M$ and ${\bf v}^{+}_1,\ldots ,{\bf v}^{+}_M$,
and the analogues of \eqref{New_Lower} and \eqref{CrudeLower} mentioned in Remarks~\mbox{D.2}~\mbox{(1)}, 
we also obtained lower bounds $\ell^{-}_{j,1},\ell^{-}_{j,2}\in (0,\infty)$ for the modulus of $\varkappa^{-}_m (k_N')$. 
Thus we got our numerical bounds \eqref{LBs_for_eigenmoduli},  
with $L^{\pm}_m := \max\{  \ell^{\pm}_{j,1} , \ell^{\pm}_{j,2} \}$ for $m\leq M$ 
(and either consistent  choice of sign, $\pm$). 
Certain of the bounds $L^{+}_1,\ldots ,L^{+}_M\in (0,\infty)$ obtained for $n=21$ 
are shown in Table~\mbox{D-3}, which can be found at the end of Section~\mbox{D.7}. 
\par 
After completing all of the above computations, it was found that in each case (i.e. for $10\leq n\leq 21$) 
we had $L^{\pm}_m > |\alpha^{\pm}_m| - 2^{9/2} {\tt u} N^{1/2}$ for $1\leq m\leq M = 384$.
This finding, together with the bounds \eqref{LBs_for_eigenmoduli}, 
enables us to conclude that we in each case have  
$|\alpha^{\pm}_m| < |\varkappa^{\pm}_m(k_N')| + 2^{-38}$ 
for $1\leq m\leq M$ and $\pm\in\{ + , -\}$. 

\subsection{An upper bound for $\| A\|_2$, when $n=21$} 

In this section, and the next, we focus on the case $n=21$,  
and so have $N=2^{21}$ and $M=384\,$ (except where there is an indication to the contrary). 
This section deals with our computation of a useful numerical upper bound for $\| A\|_2\,$ 
($A$ being the $N\times N$ symmetric matrix given by \eqref{Def-A(N)}). 
Only the key steps and results are mentioned. 
\par 
In outline, our approach to bounding $\| A\|_2$ for $n=21$ resembles the 
approach taken for $4\leq n\leq 9$, which is described in Section~\mbox{D.1}. 
That is, we first compute a fairly sharp numerical upper bound for 
$\| A^2\|\,$ (the Frobenius norm of $A^2$), and 
then apply \eqref{A_to_Asquared} to get our upper bound for $\| A\|_2$.  
We have found, however, that the methods we 
used in computing  $\| A^2 \|$ when $n\leq 9$ are not well suited 
to the case $n=21\,$ (in which the relevant matrices are far larger). In particular, when $n=21$ we 
find it impractical to employ Octave's {\tt interval} package 
in the same way as is described in Section~\mbox{D.1}. The implementation of matrix multiplication 
that the {\tt interval} package provides is so slow (on our desktop computer) 
that the run-time required becomes prohibitive 
if one is dealing with matrices as large as $A$ is in the case $n=21$: 
indeed we estimate that, even for the simpler task of computing just the double precision product of two real matrices of this size,  
the run-time could be more than $3$ years (if using our computer and Octave). 
Therefore, for $n=21$, we have taken a circuitous approach to the estimation of $\| A^2\|$ 
that utilises the upper bound for $\| (H'(N))^2\|^2\,$ ($N=2^{21}$) obtained in Section~\mbox{C.3}. 
\par 
Recalling the relevant notation introduced in Section~\mbox{5.3}, we put 
${\bf \Delta}_j := (\alpha_j I_N - H'(N))  {\bf v}_j$ 
for $1\leq j\leq 2M$, 
and define $E = (e_{i,j})$ to be the $2M\times 2M$ matrix 
given by $E = I_{2M} - V^{\rm T} V$, 
where $V$ is the $N\times 2M$ matrix with columns 
${\bf v}_1,\ldots , {\bf v}_{2M}$ 
(note that we are assuming $n=21$ and $N=2^{21}$, so that, as discussed in Section~\mbox{5.3}, we have: $N' = 2M = 768$). 
By \eqref{Def-A(N)} and \eqref{Def-H''}, one has 
\begin{equation}\label{Asquared_expansion}
A^2 = \left( H'(N) - \sum_{j=1}^{2M} \alpha_j {\bf v}_j {\bf v}_j^{\rm T}\right)^{\!\!2}
= B + C_1 + C_1^{\rm T} - C_2\;,
\end{equation}
where 
\begin{equation*}
B = \left( H'(N)\right)^2 - \sum_{j=1}^{2M} \alpha_j^2 {\bf v}_j {\bf v}_j^{\rm T}\;, 
\end{equation*}
\begin{equation*}
C_1 = \sum_{j=1}^{2M} \alpha_j {\bf \Delta}_j {\bf v}_j^{\rm T}\quad \text{and} \quad 
C_2 = \sum_{i=1}^{2M} \sum_{j=1}^{2M} e_{i,j}\alpha_i \alpha_j  {\bf v}_i {\bf v}_j^{\rm T}\;. 
\end{equation*} 
A short calculation gives: 
\begin{align}\label{Fro_B_expansion}
\| B\|^2 &= \| \left( H'(N)\right)^2\|^2 - \sum_{j=1}^{2M} \alpha_j^4 \nonumber\\
&\phantom{{=}} +4 \sum_{j=1}^{2M} \alpha_j^3 {\bf v}_j^{\rm T} {\bf \Delta}_j 
-2 \sum_{j=1}^{2M}  \alpha_j^2 \| {\bf \Delta}_j \|^2 
+ \sum_{i=1}^{2M} \sum_{j=1}^{2M} e_{i,j}^2 \alpha_i^2 \alpha_j^2 \nonumber\\ 
&= \| \left( H'(N)\right)^2\|^2 - T_1 + 4T_2 -2T_3 + T_4 
\quad\text{(say)} . 
\end{align} 
At the same time, using the Cauchy-Schwarz inequality, one can show that  
\begin{equation}\label{C_1andC_2bounds}
\| C_1\|^2 \leq \left( 1 + \| E\|\right) T_3 \quad\text{and}\quad 
\| C_2\|^2 \leq   \left( 1 + \| E\|\right)^2 T_4\;. 
\end{equation} 
\par
Using the numerical data $\alpha_m, {\mathbf v}_m\,$ ($1\leq m\leq M$) 
it was found, with the help of Octave's {\tt interval} package,  that  we had here 
$\| E\| < 2.4 \times 10^{-13}$,  $|T_2| < 1.8353876 \times 10^{-15}$,  
$0\leq T_3 < 4.01\times 10^{-25}$ and  $0\leq T_4 < 2.7\times 10^{-33}$; 
using also the sharp upper bound \eqref{FroHsquared21} for  $\| (H'(N))^2 \|^2$, 
we got: 
\begin{equation}\label{Trace-Off-Bound}
\left\| (H'(N))^2 \right\|^2 - T_1 < 9.2694157127542684\times 10^{-9}\;. 
\end{equation} 
\par 
By the bounds obtained for $T_3$ and $T_4$, and the inequalities in \eqref{C_1andC_2bounds},  it follows that 
the numbers $\| C_1\|,\| C_2\|\geq 0$ are small (less than $10^{-12}$). 
Therefore not much is lost by using the bound 
$\|A^2\| \leq \| B\| + 2( 1 + \| E\|)^{1/2} T_3^{1/2} +  ( 1 + \| E\|) T_4^{1/2}$, 
which follows (via the triangle inequality) from \eqref{Asquared_expansion} and \eqref{C_1andC_2bounds}. 
This bound, when combined with the expansion \eqref{Fro_B_expansion} of $\| B\|^2$, the 
bound \eqref{Trace-Off-Bound} and the bounds obtained for $T_2$, $T_3$, $T_4$ and $\| E\|$, 
enables us to establish that 
$\| A^2\| < 9.6277844274708 \times 10^{-5}$. This, by virtue of \eqref{A_to_Asquared},  gives us: 
\begin{equation}\label{SpecNormA21}
\| A\|_2 < 9.81212740819788 \times 10^{-3} \quad\text{($N = 2^{21}$, $M=384$)} .  
\end{equation}

\begin{remarks} 

\item{\it 1)}\quad The upper bound for $\| A\|_2$ in \eqref{SpecNormA21} is  
approximately $4.37$ times the number $\max\{\alpha^{+}_M, |\alpha^{-}_M|\}$.  
In contrast, the results of Section~\mbox{D.5} (below) make it appear very plausible that $\| A\|_2$ is less than $\min\{ \alpha^{+}_M , |\alpha^{-}_M|\}$ when $n=21$: 
see in particular Table~\mbox{D-2} (there) for the `probable upper bound' ($\tilde R$) that we get for $\| A\|_2$ 
in this case. Despite these observations, we still think that our application of   
the bound \eqref{SpecNormA21} in Section~\mbox{D.4} (below) is worthwhile, since it produces results 
that are certain to be correct (whereas in using any `probable upper bound' for $\| A\|_2$ one risks making an error). 
\par 
The upper bound for $\| A^2\|$ from which we obtain \eqref{SpecNormA21} is sharp 
(a calculation shows that the factor by which it exceeds $\| A^2\|$ is less than $\exp(5\times 10^{-6})$). 
It follows that if  $\| A\|_2$ is actually more than $4$ times smaller 
than the bound in \eqref{SpecNormA21} then the principal reason for that bound being so weak is that the bound \eqref{A_to_Asquared}, from which it derives, 
is equally weak (in this instance). 

\item{\it 2)}\quad For the cases with $10\leq n < 21$, we have $\| (H'(N))^2\|\in [\phi_N,\phi_N']$, 
where $[\phi_N,\phi_N']$ is the explicit real interval whose computation is discussed in Remarks~\mbox{C.3}~\mbox{(2)}.  
Using just  the relevant data $\phi_N,\phi_N'$ and $\alpha_1, \ldots , \alpha_{2M}$, 
one can quickly calculate a set of estimates (i.e. a non-rigorous predictions) of  
the upper bounds for $\| A\|_2$ that would result from the application of \eqref{A_to_Asquared} 
in these cases. Our estimates for these hypothetical bounds are: 
$7.2591\times 10^{-4}$, $4.7993\times 10^{-3}$, $6.8020\times 10^{-3}$, $8.2061\times 10^{-3}$, $9.0758\times 10^{-3}$, 
$9.5135\times 10^{-3}$, $9.7036\times 10^{-3}$, $9.7756\times 10^{-3}$, $9.8007\times 10^{-3}$, $9.8089\times 10^{-3}$ 
and $9.8114\times 10^{-3}$, 
for the cases $n=10$, ..., $n=20\,$ (respectively). 
These `estimated bounds', and the result \eqref{SpecNormA21}, 
seem weak in comparison with the `probable upper bounds' $\tilde R = \tilde R(N)$ 
discussed in the Section~\mbox{D.5}, below: see Table~\mbox{D-2} there. 
We have therefore not thought it worthwhile pursuing any application to cases with $10\leq n\leq 20$ 
of the method used in this section (for the case $n=21$ only). 
We are (besides) doubtful that 
analogues of the bound \eqref{SpecNormA21}, for those cases, 
would facilitate any progress on where the eigenvalues of the kernel $K$ are located
(i.e. progress beyond what \eqref{SpecNormA21} itself enables us to achieve).  
Our thinking here is based on the empirical observation that one has 
$\| k_{2N}' - K\| < \|k_N' - K\|$ for $N\in\{ 2^4, 2^5, \ldots , 2^{20}\}$. 

\item{\it 3)}\quad We find that when $10\leq n\leq 21$ and $N=2^n$ 
the $N\times N$ matrix $A = A(N)$ is real and symmetric, and has rank greater than $0$.  
It follows that in these cases $1 \geq \| A\|_2^{2\ell} / \| A^{\ell}\|^2 \geq 1/N$ for all ${\ell}\in{\mathbb N}$. 
We (ideally) would have liked to make use of 
the implied bound $\| A\|_2 \leq \| A^{\ell}\|^{1/{\ell}}$ for some ${\ell}$ greater than $2$. 
This would have been feasible in certain cases  (i.e. those where $n$ is not too large), 
but our experience in computing the approximations to $\| H'(N)^2\|$ 
shown in Tables~\mbox{C-1a} and~\mbox{C-1b} convinces us that, 
in any case with both $n\geq 16$ and ${\ell}$ large enough 
to yield a worthwhile result (i.e. a bound for $\| A\|_2$ not much greater than the 
`probable upper bound' $\tilde R$ discussed in Section~\mbox{D.5}), 
the time required to estimate $\| A^{\ell}\|$ accurately would be excessive. 
\end{remarks} 

\subsection{Upper bounds for the sizes of eigenvalues, when $n=21$} 

As in the previous section, we focus here on the case $n=21$. 
\par 
In computing upper bounds for the moduli of the eigenvalues $\varkappa^{\pm}_1 (k_N') , \ldots , \varkappa^{\pm}_M (k_N')$
we use the bound for $\| A\|_2\,$ ($N=2^{21}$) stated in \eqref{SpecNormA21}, 
together with the following more general result (applicable 
when one has $10\leq n\leq 21$).

\begin{lemma}
Let $m\leq M$ be a positive integer. 
Suppose  that \eqref{alphasOK} holds, that $\alpha^{+}_m > 0$, and that $\tau > 0$ satisfies  
\begin{equation}\label{tau-critereon} 
(1 + \tau)\alpha^{+}_m \geq \| A\|_2 
+\sum_{j=m}^M 
\frac{\alpha^{+}_j \| {\bf \Delta}^{+}_j\|^2}{\left( (1 + \tau)\alpha^{+}_m - \alpha^{+}_j\right)^2}\;, 
\end{equation}
with  
${\bf \Delta}^{+}_j$ defined as in \eqref{DefDeltaPlus}. 
Then, provided that $\alpha^{+}_j\geq 0\geq \alpha^{-}_j\,$  ($1\leq j\leq M$), 
one will have $(1 + \tau)\alpha^{+}_m \geq \varkappa^{+}_m (k_N')$. 
\end{lemma} 
\begin{proof}
Since we certainly have $(1 + \tau)\alpha^{+}_m > \alpha^{+}_m > 0$, 
it will be enough to consider just the cases where one has $\varkappa^{+}_m (k_N') > \alpha^{+}_m> 0$.  
In particular, we may assume henceforth that $H'(N)$ has at least $m$ positive 
eigenvalues. 
\par
Let ${\bf u}^{+}_1,\ldots ,{\bf u}^{+}_m$ be some 
orthonormal system of eigenvectors for the $m$ greatest eigenvalues of $H'(N)$. 
By a suitable choice of ${\bf b}\in{\mathbb R}^m$ with 
$\| {\bf b}\|^2 = b_1^2 + \cdots + b_m^2 = 1$, one can ensure that
the unit length vector ${\bf w} := \sum_{j=1}^m b_j {\bf u}^{+}_j$ 
satisfies  
${\bf w}\cdot {\bf v}^{+}_j = 0\,$ ($1\leq j < m$). 
One then has ${\bf w}^{\rm T} H'(N) {\bf w} \geq \varkappa^{+}_m (k_N')\,$ 
(the $m$-th greatest eigenvalue of $H'(N)$). Therefore, given that 
$\alpha^{-}_j\leq 0\,$  ($1\leq j\leq M$), it follows from the 
definition of the matrix $A$ and our choice of ${\bf w}\in{\mathbb R}^N$ that  
${\bf w}^{\rm T} A {\bf w} \geq \varkappa^{+}_m (k_N') 
- \sum_{j=m}^M \alpha^{+}_j ( {\bf w}\cdot {\bf v}^{+}_j)^2$. 
At the same time, it follows by the Cauchy-Schwarz inequality 
that ${\bf w}^{\rm T} A {\bf w} = {\bf w}\cdot (A {\bf w}) \leq \| A {\bf w}\| \leq \| A\|_2$, so 
that $\varkappa^{+}_m (k_N') 
\leq \| A\|_2 +  \sum_{j=m}^M \alpha^{+}_j ( {\bf w}\cdot {\bf v}^{+}_j)^2$. 
\par 
By using the fact that $H'(N)^{\rm T} = H'(N)$, and then \eqref{alphasOK},  
the Cauchy-Schwarz inequality and Bessel's inequality, 
one finds that 
\begin{align*} 
\left| {\bf w}\cdot {\bf v}^{+}_j\right| 
 &\leq \sum_{i=1}^m 
\frac{| b_i {\bf u}^{+}_i \cdot {\bf \Delta}^{+}_j |}{\varkappa^{+}_i(k_N')  - \alpha^{+}_j} \\ 
 &\leq \left(\sum_{i=1}^m 
\frac{b_i^2}{\left(\varkappa^{+}_i(k_N')  - \alpha^{+}_j\right)^2}\right)^{\!\!1/2} 
\left\| {\bf \Delta}^{+}_j\right\| 
\leq \frac{\left\| {\bf \Delta}^{+}_j\right\|}{\varkappa^{+}_m(k_N')  - \alpha^{+}_j}
\end{align*}
for $m\leq j\leq M$. This, together with the conclusion from the preceding paragraph, 
shows that the inequality in \eqref{tau-critereon} is reversed if 
one substitutes for each $\tau$ there  the number $T>0$ such that 
$(1+T)\alpha^{+}_m  = \varkappa^{+}_m(k_N')$:  
in light of the monotonicity of the relevant expressions 
(viewed as functions of $\tau >0$), this shows that $T{\not >}\tau$, 
so that $\varkappa^{+}_m (k_N'){\not >}(1 + \tau)\alpha^{+}_m$. 
\end{proof}
\par 
In the case ($n=21$) that is our main concern here, all of the hypotheses and conditions 
of Lemma~\mbox{D.4} are satisfied: 
recall, in particular, that in all cases our computed data $\alpha^{\pm}_1,\ldots ,\alpha^{\pm}_M$ satisfies both \eqref{alphasOK} and 
the analogous inequalities $\alpha^{-}_1 < \cdots < \alpha^{-}_M$, 
and that for $10\leq n\leq 21$ one has (in the notation of Section~\mbox{5.3})  $\,p' = q' = M$, 
so that $\alpha^{+}_j > 0 > \alpha^{-}_j$  for $1\leq j\leq M$.
\par 
Let $R$ be the numerical upper bound for $\| A\|_2$ in \eqref{SpecNormA21},  
and let $\delta^{+}_1,\ldots ,\delta^{+}_M$ be (respectively) those sharp upper bounds for  
$\|{\bf \Delta}^{+}_1\|,\ldots , \|{\bf \Delta}^{+}_M\|\,$ (in the case $n=21$) 
whose computation is described in Remarks~\mbox{D.2}~\mbox{(2)}. 
In applying Lemma~\mbox{D.4}, for a given value of $m\leq M$, 
we have used the {\tt interval} package and an iterative 
method (involving bisection of intervals) to compute a sharp upper 
bound $\tau^{+}_m$ for the unique $\tau > 0$ such that 
the two sides of the inequality \eqref{tau-critereon} becomes equal when 
$\| A\|_2$ and $\| {\bf \Delta}^{+}_m\|, \ldots , \| {\bf \Delta}^{+}_M\|$ 
are replaced with the corresponding 
upper bounds, $R$ and $\delta^{+}_m,\ldots ,\delta^{+}_M$.    
The application of Lemma~\mbox{D.4} gives us, then, the bound:  
\begin{equation}\label{Def-U^{+}_m}
\varkappa^{+}_m (k_N') \leq (1 + \tau^{+}_m)\alpha^{+}_m = U^{+}_m\quad 
\text{(say)} .
\end{equation}   
Using instead the 
appropriate analogue of Lemma~\mbox{D.4}, together with the bound \eqref{SpecNormA21}, 
the data $\alpha_1^{-},\ldots ,\alpha_M^{-}$ and  some sharp 
upper bounds for $\|{\bf \Delta}^{-}_1\|,\ldots , \|{\bf \Delta}^{-}_M\|$, 
we compute (similarly)  a number $\tau_m^{-}>0$ such that 
$| \varkappa^{-}_m (k_N')| \leq (1 + \tau_m^{-})\alpha^{-}_m = U^{-}_m\,$ 
(say). By combining this and \eqref{Def-U^{+}_m} (for $1\leq m\leq M$)   
with the complementary lower bounds \eqref{LBs_for_eigenmoduli}, we arrived at the 
numerical results stated formally in \eqref{n=21intervals}. 
\par 
Certain of the bounds $U^{+}_1,\ldots ,U^{+}_M\in (0,\infty)$ obtained for $n=21$ 
are shown in Table~\mbox{D-3}, which at the end of Section~\mbox{D.7}: see 
also Remarks~\mbox{D.7}~\mbox{(2)}  
concerning how accurate an approximation to $\varkappa^{\pm}_m (k_N')$ the number $\alpha^{\pm}_m$ is when $N=2^{21}$. 

\subsection{Probable upper bounds for $\| A\|_2$, when $n\geq 10$} 

Recall that $A$, which we defined in \eqref{Def-A(N)},  is a real and symmetric $N\times N$ matrix:  
it therefore has $N$ real eigenvalues,  
$\mu_1, \ldots , \mu_N\,$ (say), and a corresponding orthonormal system of 
eigenvectors ${\bf w}_1,\ldots , {\bf w}_N\in{\mathbb R}^N$ such that 
\begin{equation}\label{Apx}
\| A^j {\bf x}\|^2 = \sum_{i=1}^N ({\bf w}_i \cdot {\bf x})^2  \mu_i^{2j} 
\quad\text{(${\bf x}\in{\mathbb R}^N$, $j=0,1,2,\ldots\ $)} .
\end{equation} 
We may assume here that $|\mu_1|\geq |\mu_2|\geq\ \cdots\ \geq |\mu_N|$, so that 
$\| A\|_2 = |\mu_1|$. By \eqref{Apx}, we have 
$\| A^j {\bf x}\|^2 \geq ({\bf w}_1 \cdot {\bf x})^2  \mu_1^{2j} 
=   ({\bf w}_1 \cdot {\bf x})^2  \| A\|_2^{2j}\,$ ($j\in{\mathbb N}$), and so 
\begin{equation}\label{DieCast-1}
\| A\|_2 \leq \left(\frac{\| A^j {\bf x}\|}{|{\bf w}_1 \cdot {\bf x}|}\right)^{1/j} 
\quad\text{(${\bf w}_1 \cdot {\bf x}\neq 0$, $j\in{\mathbb N}$)} .
\end{equation} 
\par
In applying \eqref{DieCast-1} we treat the eigenvector ${\bf w}_1$ as if it were 
just some unspecified unit vector in ${\mathbb R}^N$. 
Therefore, since any non-trivial application of \eqref{DieCast-1} requires some sort of information  
about the size of ${\bf w}_1 \cdot {\bf x}$, it is obvious that just one application of 
\eqref{DieCast-1} will not be sufficient.
With this in mind, we observe that when $X$ is an orthonormal basis for ${\mathbb R}^N$ one has  
$\max_{{\bf x}\in X} |{\bf w}_1 \cdot {\bf x}|\geq (N^{-1} \| {\bf w}_1\|^2)^{1/2} 
= 1/\sqrt{N}$, and so $\| A\|_2 \leq 
(N^{1/2} \max_{{\bf x}\in X} \| A^j {\bf x}\|)^{1/j}$. 
This, however, yields a fairly weak result when $j$ is small (i.e. when $j\leq 10$, say);  
and when $j$ is not small 
the time needed for us to compute $\| A^j {\bf x}\|$ for each member ${\bf x}$ of such a basis 
$X$ becomes excessive for large values of $N$. 
We opt instead to compute $\| A^j {\bf x}\|$ for each vector ${\bf x}$ occurring 
in some fairly short sequence, ${\bf x}_1,\ldots , {\bf x}_S$, of unit vectors, each 
randomly (and independently) selected from the set 
\begin{equation}\label{DefSetXn}
X_N := \left\{ {\bf y}\in{\mathbb R}^N : y_i^2 
= \textstyle{\frac1N}\ {\rm for}\ i=1,\ldots ,N\right\} 
\end{equation}
(with all $2^N$ elements of $X_N$ given an equal chance of being selected). 
The following lemma enables one to argue that, as $S$ increases, the probability 
of having $\max_{1\leq s\leq S} |{\bf w}_1 \cdot {\bf x}_s| \gg 1/\sqrt{N}$ will tend 
rapidly towards $1$, and so can be close to $1$ even when $S$ is  
smaller than $N$ by several orders of magnitude.  

\begin{lemma}
Let $C = \sqrt{\pi^3 /2}$ and let ${\bf w}=(w_1,\ldots , w_N)^{\rm T}$ be a unit vector in ${\mathbb R}^N$. 
Then 
\begin{equation}\label{RandomWalk-1}
\left| \left\{ {\bf x}\in X_N : |{\bf w}\cdot {\bf x}| \leq \textstyle{\frac{1}{C\sqrt{N}}}\right\}\right| 
\leq \textstyle{\frac12} | X_N|\;. 
\end{equation} 
\end{lemma} 

\begin{proof} 
Let $\Pi_N({\bf w})$ 
denote the expression on the left-hand side of the inequality \eqref{RandomWalk-1}. 
We just need to prove that $\Pi_N ({\bf w}) \leq 2^{N-1}$. 
Since $\Pi_N({\bf w}) = 0$ if $N=1$, we can assume that  $N\geq 2$.
\par 
We consider firstly the cases where $|w_N| > 1/C$. Note that $\Pi_N ({\bf w})$ 
lends itself to being reformulated as $N$ nested summations, with a variable summand 
whose value is always either $0$ or~$1$. Using a uniform 
bound for the innermost summation, we find that 
$\Pi_N ({\bf w}) \leq 2^{N-1} \max_{\beta\in{\mathbb R}} 
|\{ j\in\{ -1 , 1\} : | j w_N + \beta|\leq \textstyle{\frac1C}\}|$. 
The desired result follows by observing that, since $|w_N| > \frac1C$, 
there is no $\beta$ with both 
$|\beta - w_N|\leq \frac1C$ and $|\beta + w_N|\leq \frac1C$. 
\par 
Since $\Pi_N ({\bf w}) = \Pi_N (w_1,\ldots ,w_N)$ is invariant 
under any permutation of $w_1,\ldots ,w_N$, we also obtain the desired 
result if $\frac1C < \| {\bf w}\|_{\infty} := \max_{1\leq j\leq N} |w_j|$. 
\par 
In the remaining cases (those where $\| {\bf w}\|_{\infty} \leq \frac1C$) we 
begin by observing that 
\begin{equation*} 
\Pi_N ({\bf w}) \leq \sum_{{\bf x}\in X_N} 
\textstyle{\frac18}\pi^2 {\rm sinc}^2\left(\textstyle{\frac14} C \sqrt{N} {\bf x}\cdot {\bf w}\right)\;,
\end{equation*} 
where ${\rm sinc}(t):=(\pi t)^{-1} \sin(\pi t)$ for $t\neq 0$, and ${\rm sinc}(0) := 1$. 
The fact that ${\rm sinc}^2 (t)$ is the Fourier transform 
of the real function $\Lambda (s) := \max\{ 0 , 1 - |s|\}$ enables one to deduce that 
\begin{equation*} 
\Pi_N ({\bf w}) \leq \frac{2^{N-1}\pi^2}{C} 
\int_{-C/4}^{C/4} \Lambda\left( \frac{y}{C/4}\right) \left( \prod_{j=1}^N \cos(2\pi y w_j)\right) dy\;. 
\end{equation*}
We have here $|y w_j|\leq |y| / C \leq \frac14$ for $|y|\leq C/4$ and $1\leq j\leq N$, 
and so, by exploiting the fact that 
\begin{equation*} 
0\leq \cos(x)\leq \exp(-\textstyle{\frac12} x^2)\quad\text{for $-\frac{\pi}{2}\leq x\leq \frac{\pi}{2}$} , 
\end{equation*}
we find that 
\begin{equation*}
\Pi_N ({\bf w}) \leq \frac{2^{N-1}\pi^2}{C}  \int_{-\infty}^{\infty} 
\exp\left( -2\pi^2 y^2 \| {\bf w}\|^2\right) dy 
= \frac{2^N\pi^2}{C}  \int_{0}^{\infty} \exp\left( -2\pi^2 y^2\right) dy\;,  
\end{equation*}
given that $\| {\bf w}\| =1$. 
The last integral is $\Gamma(1/2) / (2\sqrt{2}\pi) = 1/\sqrt{8\pi}$, so that, 
since $\pi^2 / C = \sqrt{2\pi}$, we obtain the bound $\Pi_N ({\bf w}) \leq 2^{N-1} = \frac12 | X_N |$, 
as required.
\end{proof}

We now detail the four step algorithm that constitutes our application of the above lemma.  
It should be kept in mind that we assume $n\in{\mathbb N}$,  $10\leq n\leq 21$ and $N = 2^n$.
\begin{description}
\item[Step~1.] We choose a positive integer $S$ (with $30\leq S\leq 32$, in practice). 
Then we choose vectors ${\bf x}_1,\ldots ,{\bf x}_S$ `at random' from the 
set $X_N\,$ (defined in \eqref{DefSetXn}) and compute double precision approximations,
${\bf x}_1^{(0)},\ldots ,{\bf x}_S^{(0)}$, 
to these vectors.
We then have 
${\bf x}_s^{(0)} = (1+\delta){\bf x}_s\,$  ($1\leq s\leq S$), 
where $\delta$ depends only on $N$, satisfies $|\delta|\leq {\tt u}$, 
and equals $0$ if and only if $n$ is even. 
\item[Step~2.]  We choose positive integers $J$ and $G$ (with $2\leq J/2^{10}\leq 5$ 
and $9\leq G\leq 11$, in practice) 
and put $g = 2^G$. 
As $j$ runs through the sequence $1,2,\ldots , J$ we compute, for each value that $j$ takes, 
vectors ${\bf x}_1^{(j)},\ldots , {\bf x}_S^{(j)}$ that are approximations to 
the vectors $gA{\bf x}_1^{(j-1)},\ldots , gA{\bf x}_S^{(j-1)}$, and then 
use the {\tt interval} package to compute a 
number $\chi_j$ that is a sharp upper bound for 
the set $\{ \| {\bf x}_s^{(j)}\| : 1\leq s\leq S\}$. The purpose of the factor $g$ 
is to ensure $\max_{1\leq s\leq S} |\log_2 (\| {\bf x}_s^{(j)}\|)|$ is never 
greater than $1022$ (a prudent limit, since we work with binary64 floating-point 
arithmetic).  In practice we chose $G$ (i.e. $\log_2 (g)$)  
to be the integer nearest to $\log_2 (2/(\alpha^{+}_M + |\alpha^{-}_M|))$. 
This worked well, so long as the choice of $J$ was not overly ambitious. 
\item[Step~3.] We compute a positive number $\Delta(A, g)$ such that one has: 
\begin{equation}\label{DieCast-2} 
\| {\bf x}_s^{(j)} - gA{\bf x}_s^{(j-1)}\| 
\leq \Delta(A, g) \| {\bf x}_s^{(j-1)}\| 
\quad\text{($1\leq j\leq J$, $1\leq s\leq S$)} . 
\end{equation}
This computation depends on a rounding error analysis of 
the relevant part of Step~2. 
\item[Step~4.] We compute a number $\tilde R$ such that 
the data computed in Steps~2 and~3 
may, by virtue of \eqref{DieCast-1}, \eqref{RandomWalk-1},  and \eqref{DieCast-2}, 
reasonably be viewed as experimental evidence 
supporting the hypothesis that $\tilde R\geq\| A\|_2$.
This $\tilde R$ becomes our `probable upper bound' for $\| A\|_2$. 
\end{description} 
The above descriptions of Steps~\mbox{1--4} are intentionally cursory:   
for a discussion of  the `random' choice of vectors in Step~1,  see Appendix~E;  
and for further details of Steps~2 and~3, see Remarks~\mbox{D.6} (below). 
We now discuss Step~4 in more detail. 
\par
Let $Z = (\frac12 \pi^3 N)^{1/2}$. By \eqref{DieCast-1}, either 
\begin{equation}\label{DieCast-3}
Z \| (gA)^J {\bf x}_s\| \geq \| gA \|_2^J = g^J \| A\|_2^J  \quad\text{for some $s\in\{ 1,\ldots ,S\}$} ,  
\end{equation}
or else one has: 
\begin{equation}\label{DieCast-4}
| {\bf w}_1 \cdot {\bf x}_s | < Z^{-1}
\quad\text{($1\leq s\leq S$)} .
\end{equation}
At the same time, it follows from \eqref{DieCast-2} that one certainly has  
\begin{equation*}
\| {\bf x}_s^{(J)} \| \geq \| (gA)^J {\bf x}_s^{(0)}\| 
- \Delta(A, g) \sum_{k=1}^J \| gA\|_2^{J-k} \| {\bf x}_s^{(k-1)} \|  \quad\text{($1\leq s\leq S$)} ,
\end{equation*} 
and here (recalling Step~1) one may observe that  
\begin{equation*}
\| (gA)^J {\bf x}_s^{(0)}\| = (1+\delta)  \| (gA)^J {\bf x}_s\| \geq (1-{\tt u})  \| (gA)^J {\bf x}_s\|\;.
\end{equation*}   
\par 
Defining 
\begin{equation*}
\varrho = (1 - {\tt u}) \| gA\|_2\;,\qquad 
\xi_j = Z \chi_j \geq Z\cdot (\max_{1\leq s\leq S} \| {\bf x}_s^{(j)} \| )\quad\text{($0\leq j\leq J$)} , 
\end{equation*}
and 
\begin{equation}\label{Def_Xi_rho1}
\Xi(\varrho_1) =  \Delta(A,g)\sum_{k=1}^J \xi_{k-1} \varrho_1^{-k}\quad\text{($\varrho_1 > 0$)} ,
\end{equation}
one finds, by  virtue of the facts just noted above, that  \eqref{DieCast-3} implies  
$\xi_J \geq (1 - \Xi(\varrho))\varrho^J$. 
It follows that, with $\varrho_0 := \xi_J^{1/J}$, one will have 
\begin{equation*} 
\varrho\leq (1 - \Xi(\varrho_0))^{-1/J} \varrho_0\quad 
\text{if \eqref{DieCast-3} holds and $\Xi(\varrho_0) < 1$} . 
\end{equation*} 
We conclude that if \eqref{DieCast-3} holds, and if 
$\Xi(\varrho_0) < 1$, then $\tilde R\geq \| A\|_2$ if one has: 
\begin{equation*} 
\tilde R\geq R_0 := 2^{-G} (1-{\tt u})^{-1} (1 - \Xi(\varrho_0))^{-1/J} \varrho_0\;.
\end{equation*}  
Our Step~4 consists essentially of 
a check on the size of $\Xi(\varrho_0)$ and (provided it is confirmed that $\Xi(\varrho_0)<1$) 
the computation 
of a suitable number $\tilde R$, satisfying $\tilde R\geq R_0$. 
In order to get a near optimal value for $\tilde R$ we use the {\tt interval} package 
to compute a very short interval containing $R_0$, and set $\tilde R$ equal 
to the upper endpoint of that interval. 
We found that, in practice, the value of $\Xi(\varrho_0)$ 
never exceeded $2.7\times 10^{-3}$. 
\par 
Recall now that in every case where \eqref{DieCast-3} fails to hold the 
condition \eqref{DieCast-4} will be satisfied. 
Let us assume that in Step~1 we choose ${\bf x}_1, \ldots , {\bf x}_S$ by a method that 
can reasonably be considered random sampling (with replacement) 
from $X_N$. Then the probability of an outcome in which 
the condition  \eqref{DieCast-4} is satisfied  
equals $(|Y_N| / |X_N|)^S$, where 
$Y_N = \{ {\bf x}\in X_N : |{\bf w}_1\cdot {\bf x}| \leq Z^{-1}\}$. 
By  Lemma~\mbox{D.5}, this probability does not exceed $2^{-S}$. 
Therefore, for $S\geq 30$ (say), outcomes  
such that \eqref{DieCast-3} fails to hold should rarely occur: indeed, since $2^{-30} < 10^{-9}$, 
the odds of such an outcome should be less than one in a billion. 
The computations outlined in the last paragraph above yield either 
no result, or else a number $\tilde R$ that satisfies $\tilde R\geq \| A\|_2$ if \eqref{DieCast-3} holds.  
Therefore, when $S\geq 30$, the odds of obtaining a result $\tilde R$ that is not a valid  
upper bound for $\| A\|_2$ should also be less than one in a billion. 
This is why we call our $\tilde R$ a `probable upper bound for $\| A\|_2$'. 
\par
Use of the above four step algorithm (with $S$ set equal to $32$) yielded    
the probable upper bounds $\tilde R = \tilde R(N)\,$ ($10\leq n\leq 21$) shown in Table~\mbox{D-2} (below). 
The relevant `random' input data (required in Step~1 of the algorithm) was obtained 
via the RDRAND-based method described in Appendix~E. Although  we do not know, for certain,  
that samples from $X_N$ produced by this method are truly random, we can at least    
note that the RDRAND-based method leads to results 
(i.e. probable upper bounds for $\| A\|_2$) that are similar in strength to results that we have 
obtained (separately) by using instead Octave's built-in pseudorandom number generator, the {\tt rand()} function, for the sampling from $X_N$: 
Appendix~E gives further details of this. 
\par 
Note that Table~\mbox{D-2} shows $r$, $r'$ and $\tilde R$ rounded to 16 significant digits, while $J^{-1}\log Z$ 
has been rounded to just 2 significant digits. 
The final 3 columns of this table allow a worthwhile comparison to be made: 
for if $\alpha^{+}_1, \ldots , \alpha^{+}_M$ and 
$\alpha^{-}_1, \ldots , \alpha^{-}_M$ are (as one might anticipate) accurate approximations to 
the $M$ greatest and $M$ least of the eigenvalues of $H'(N)$, then one should expect to have 
\begin{equation*} 
\log\left( \frac{\tilde R}{r}\right) < \frac{\log Z}{J} = \frac{(n-1)\log(2) + 3\log(\pi)}{2J}\;.
\end{equation*}
The data in Table~\mbox{D-2} is in accordance with this.  Indeed, 
we even have $J^{-1}\log Z > \log(\tilde R/r')$ 
for every $n$ occuring in the table 
(and $\tilde R < r'$ for $n\in\{ 11, 12, 15, 21\}$).  
Only in one case ($n=13$) do we have $\tilde R > r$. 

\begin{table}[ht] 
\centering 
\begin{minipage}{165mm}
\resizebox{165mm}{!}{ 
\begin{tabular}[c]{|r|r|c|c|c|c|} 
\hline 
$n$ & $J$ & $\tilde R$ & $r := \max\{ \alpha^{+}_M , |\alpha^{-}_M|\}$ & $r' := \min\{ \alpha^{+}_M , |\alpha^{-}_M|\}$ & $J^{-1}\log Z$ \\ \hline 
$10$ & $2048$ & $3.784180848066032\times 10^{-4}$ & $3.972680567252615\times 10^{-4}$ & $3.778696884388128\times 10^{-4}$ & $2.4\times 10^{-3}$ \\ \hline
$11$ & $2048$ & $1.352835069802926\times 10^{-3}$ & $1.355139969660699\times 10^{-3}$ & $1.353634248413826\times 10^{-3}$ & $2.5\times 10^{-3}$ \\ \hline 
$12$ & $2048$ & $1.764081026332446\times 10^{-3}$ & $1.766508594121058\times 10^{-3}$ & $1.764691894167444\times 10^{-3}$ & $2.7\times 10^{-3}$ \\ \hline
$13$ & $2048$ & $2.023189046437590\times 10^{-3}$ & $2.021915127457577\times 10^{-3}$ & $2.021239822451941\times 10^{-3}$ & $2.9\times 10^{-3}$ \\ \hline
$14$ & $2048$ & $2.153344461025095\times 10^{-3}$ & $2.157309211841370\times 10^{-3}$ & $2.151723663344031\times 10^{-3}$ & $3.0\times 10^{-3}$ \\ \hline
$15$ & $2048$ & $2.210814012851814\times 10^{-3}$ & $2.216505476779971\times 10^{-3}$ & $2.215879099566874\times 10^{-3}$ & $3.2\times 10^{-3}$ \\ \hline
$16$ & $2048$ & $2.234411276710962\times 10^{-3}$ & $2.234618952358004\times 10^{-3}$ & $2.234264903531416\times 10^{-3}$ & $3.4\times 10^{-3}$ \\ \hline
$17$ & $2048$ & $2.240271002535091\times 10^{-3}$ & $2.243604211407860\times 10^{-3}$ & $2.238846592457597\times 10^{-3}$ & $3.5\times 10^{-3}$ \\ \hline
$18$ & $2048$ & $2.241893563765699\times 10^{-3}$ & $2.245291132772025\times 10^{-3}$ & $2.240144104968621\times 10^{-3}$ & $3.7\times 10^{-3}$ \\ \hline
$19$ & $2048$ & $2.241938318318294\times 10^{-3}$ & $2.245772758418342\times 10^{-3}$ & $2.240518884627916\times 10^{-3}$ & $3.9\times 10^{-3}$ \\ \hline
$20$ & $2048$ & $2.241950891002065\times 10^{-3}$ & $2.245909791156905\times 10^{-3}$ & $2.240625617294758\times 10^{-3}$ & $4.0\times 10^{-3}$ \\ \hline
$21$ & $5120$ & $2.240542901682991\times 10^{-3}$ & $2.245948557551722\times 10^{-3} $ & $2.240655836509706\times 10^{-3}$ & $1.7\times 10^{-3}$ \\ 
\hline 
\end{tabular} 
} 
\vskip 3mm
\ Table~D-2 (with $M=384$ and $N=2^n$) 
\end{minipage} 
\end{table}

\begin{remarks} 
\item\quad\ We have not yet discussed how, in Step~3, a suitable value for $\Delta(A,g)$ is determined.    
We shall remedy this omission, after first providing relevant information about Step~2 and the 
rounding errors that may occur there. 
\par 
When $1\leq j\leq J$ and $1\leq s\leq S$, the vector ${\bf x}^{(j)}_s$ computed in Step~2 is  
$2^G$ times an approximation to a difference, $\tilde{\bf y} - \tilde{\bf z}$, in which 
$\tilde{\bf y}$ and $\tilde{\bf z}$ are our computed approximations 
to the vectors ${\bf y} := H'(N) {\bf x}^{(j-1)}_s$ and 
${\bf z} := \sum_{i=1}^{2M} ( {\bf v}_i \cdot {\bf x}^{(j-1)}_s) \alpha_i  {\bf v}_i$.  
Note, in particular, that $\tilde{\bf y} = {\tt fast\_hmm}( {\bf x}^{(j-1)}_s)$, 
where {\tt fast\_hmm()} is the Octave function discussed extensively in Appendix~C. 
Thus, by \eqref{fast_hmmErrBound-2}, we have: 
\begin{equation}\label{tilde_y_errbound}
\left\| \tilde{\bf y} - {\bf y}\right\| / \| {\bf x}^{(j-1)}_s\|  < {\textstyle \frac53}  {\tt u} N  \| H'(N)\| 
\quad \text{if $\,12\leq n\leq 21$} . 
\end{equation}
We overlook, for now, the reservations concerning \eqref{fast_hmmErrBound-2} that are expressed in Remarks~\mbox{C.2}. 
We also postpone discussion of what replaces \eqref{tilde_y_errbound} when $n\in\{ 10, 11\}$. 
\par
The difference between $\tilde{\bf z}$ and ${\bf z}$ 
arises from our use of binary64 floating-point arithmetic (and the consequent rounding errors). 
Through an elementary analysis of rounding error (omitted here) 
and an application (also omitted) of \cite[Section~27, Theorem~1]{Da1980}, 
it can be shown that when $N+2M+1\leq (8/{\tt u}^2)^{1/5} = 2^{21.8}$ one has  
\begin{equation*}
\left\| \tilde{\bf z} - {\bf z}\right\| / \bigl\| {\bf x}^{(j-1)}_s \bigr\| 
\leq \gamma_{N+2+2M} \cdot \left( \sum_{i=1}^{2M} |\alpha_i|^2 \| {\bf v}_i\|^2 \right)^{\!\!\frac12} 
\left( \max_{1\leq k\leq 2M} \sum_{i = 1}^{2M} |{\bf v}_k \cdot {\bf v}_i|\right)^{\!\!\frac12} \;, 
\end{equation*}
where $\gamma_m := m{\tt u} / (1-m{\tt u})$, while 
$\alpha_i, {\bf v}_i\,$  ($1\leq i\leq 2M$) 
are as indicated in Section~\mbox{5.3}. 
From this we get the slightly simpler bound: 
\begin{equation}\label{tilde_z_errbound}
\left\| \tilde{\bf z} - {\bf z}\right\| / \bigl\| {\bf x}^{(j-1)}_s \bigr\| 
\leq \gamma_{N+2+2M} \cdot \left( 1 + \| E\|_1\right)  \left( \sum_{i=1}^{2M} |\alpha_i|^2 \right)^{\!\!\frac12} \;,
\end{equation} 
with $E$  defined as in Section~\mbox{D.3}. 
\par 
We find, in practice,  that $\| E\|_1 < 1.7\times 10^{-13}$ for $10\leq n\leq 21$:  
given that this is so, it follows by \eqref{tilde_y_errbound} and \eqref{tilde_z_errbound} that when $12\leq n\leq 21$ 
the bounds \eqref{DieCast-2} will hold (comfortably) if we have  
\begin{equation}\label{Delta(A,g)Formula}
\Delta (A, g) \geq 2 g N {\tt u}\cdot \left( \| H'(N)\| + \biggl(\sum_{j=1}^{2M} \alpha_j^2\biggr)^{\!\!1/2} \right) .
\end{equation}
Thus in Step~3 we only need $\Delta(A, g)>0$ to be large enough 
that \eqref{Delta(A,g)Formula} holds: 
this we ensure by first using the {\tt interval} package to compute a short interval $[a, b]$ containing the 
numerical value of the expression on the right-hand side of \eqref{Delta(A,g)Formula}, and then putting   
$\Delta(A,g) = b$. 
\par 
As already noted, some points of the above error analysis need further checks. 
In particular, in view of Remarks~\mbox{C.2}~\mbox{(2)}, we need to eliminate the possibility that 
overflows or underflows occurring in the computations carried out by the function {\tt fast\_hmm()} 
might make the error bound \eqref{tilde_y_errbound} invalid. 
\par 
There is no great difficulty in guarding against overflows. 
Indeed, by means of a suitable choice of parameters $G$ and $J$, in Step~2, 
we have (in practice) avoided having any overflows occur: this was checked 
by applying Octave's {\tt isfinite()} function to the vectors ${\bf x}^{(j)}_s\,$ ($0\leq j\leq J$, $1\leq s\leq S$). 
Note that such checks are not superfluous, since it is possible to have overflows occur without this 
being apparent from the numbers $\xi_0,\ldots ,\xi_J$ computed in Step~4   
(we used there Octave's {\tt max()} function, which will overlook any {\tt NaN} in input data).  
The greatest value for $\| {\bf x}^{(j)}_s\|$ that we encountered 
was approximately $3.8\times 10^{301} < 2^{1002}\,$
(this maximum was attained with $N=2^{21}$, $G=9$ and $j=J=5120$).
\par
In considering whether or not underflows might invalidate \eqref{tilde_y_errbound} 
we first recall the relevant discussion in Remarks~\mbox{C.2}~\mbox{(2)}, concerning \eqref{fast_hmmErrBound-2} and \eqref{Higham-2}. 
By the main conclusion reached there, we have only to check that 
$N^{-1} \| {\bf x}^{(j-1)}_s \| \geq 2^{-968}$ for $1\leq j\leq J$, $1\leq s\leq S$: subject to 
this being the case, underflow errors fail to make \eqref{tilde_y_errbound} invalid. 
This (it may be shown) is also a sufficient condition for 
\eqref{tilde_z_errbound} not to be invalidated by underflow errors: the empirical observation that 
$(\sum_{i=1}^{2M} \alpha_i^2)^{1/2} \in (\frac14 , \frac27)\,$ (when $12\leq n\leq 21$)  
is helpful in establishing this fact. 
It turns out that, for $12\leq n\leq 21$, 
we have $N^{-1}\| {\bf x}^{(j-1)}_s \| \geq 2^{-328}$ 
at all stages of the relevant computations.  
Therefore, in the cases that concern us, underflow errors do not invalidate 
\eqref{tilde_y_errbound} or \eqref{tilde_z_errbound},  
and so these errors (if they occur at all) must also fail to invalidate our conclusion 
that \eqref{Delta(A,g)Formula} implies \eqref{DieCast-2}.  
\par 
We come now (finally) to the cases where $n\in\{10, 11\}$: here the error analysis    
\eqref{tilde_y_errbound}--\eqref{Delta(A,g)Formula} does not apply. 
This has led us to modify how Steps~2,~3 and~4 of our algorithm are carried out in these cases. In particular, 
we use neither the function {\tt fast\_hmm()} nor any other implementation of Lemma~\mbox{C.1}. 
Instead we use the {\tt interval} package to compute an $N\times N$ matrix 
${\mathcal A}'$ of intervals, 
each of which contains the corresponding element of the matrix $gA$. 
Then, for Step~2 of the algorithm (where $1\leq j\leq J$),  
we use the {\tt interval} package  
and the data  ${\bf x}^{(j-1)}_1,\ldots , {\bf x}^{(j-1)}_S$ and ${\mathcal A}'$ to compute (for 
$1\leq s\leq S$) an $N\times 1$ matrix ${\mathcal X}^{(j)}_s$ of intervals, each containing the corresponding 
element of the $N\times 1$ matrix $gA{\bf x}^{(j-1)}_s$. 
We take  ${\bf x}^{(j)}_1,\ldots , {\bf x}^{(j)}_S\in{\mathbb R}^N$ to be the results 
returned when the ${\tt mid}()$ function (from the {\tt interval} package)  
is applied to  ${\mathcal X}^{(j)}_1,\ldots , {\mathcal X}^{(j)}_S$, respectively. 
This ensures that, for $1\leq s\leq S$, each element of the $N\times 1$ matrix ${\bf x}^{(j)}_s$ 
lies in the corresponding element of ${\mathcal X}^{(j)}_s$. 
We also introduce a refinement into Step~3, by replacing the constant factor $\Delta(A,g)$ with a variable one:  
$\Delta(A,g,j)\geq 
\max_{1\leq s\leq S} \| {\bf x}_s^{(j)} - gA{\bf x}_s^{(j-1)}\| / \| {\bf x}_s^{(j-1)}\|\,$ 
($1\leq j\leq J$). 
Using the {\tt interval} package, a satisfactory value for $\Delta(A,g,j)$ is computed 
from the  data ${\bf x}^{(j-1)}_s, {\mathcal X}^{(j)}_s\,$ 
($1\leq s\leq S$). Step~4 is correspondingly refined, so 
that we have 
$\Xi(\varrho_1) := \sum_{k=1}^J  \Delta(A,g,k)\xi_{k-1} \varrho_1^{-k}\,$ ($\varrho_1 > 0$), 
in place of \eqref{Def_Xi_rho1}. 
If we could have also taken this approach for each $n\in\{ 12,13,\ldots , 21\}$ 
then our treatment of those cases would have been greatly simplified 
(and put on a surer footing). This however was not practical, 
due to  limited computer memory and 
the great reduction in speed of computation that occurs with more intensive use of  
the {\tt interval} package (both these factors becoming more serious issues as $N$ increases).
\end{remarks}

\subsection{Probable upper bounds for the sizes of eigenvalues}

Our discussion in this section concerns the  
cases with $n\in\{ 10, 11,\ldots , 21\}$ and $N=2^n$. 
In each such case we have computed (as described in Section~\mbox{D.5}) 
a certain number $\tilde R$ that we call a `probable upper bound' for $\| A\|_2$. 
Lemma~\mbox{D.4} enables us to compute, 
for $1\leq m\leq M$ and $\pm\in\{ +, -\}$, a number $\tilde U^{\pm}_m$ 
such that 
\begin{equation}\label{Def-ProbUB}
\left| \varkappa^{\pm}_m (k_N')\right|\leq \tilde U^{\pm}_m\quad\text{if}\ \,\| A\|_2 \leq \tilde R\;. 
\end{equation} 
We call $\tilde U^{\pm}_m$ (as we do $\tilde R$) a 
`probable upper bound'. 
\par 
More specifically, in the case $n=21$ we obtain each probable upper bound $\tilde U^{\pm}_m$ 
by means of a computation almost identical to the computation (briefly described 
in  Section~\mbox{D.4}) by which we have obtained, also for $n=21$, 
both our bound $U^{+}_m$ in \eqref{Def-U^{+}_m} and 
our upper bound $U^{-}_m$ for $|\varkappa^{-}_m(k_N')|$: the only change being  
the substitution of the number $\tilde R$ in place of 
the upper bound  $R= 9.81212740819788 \times 10^{-3}$ from \eqref{SpecNormA21} 
(a substitution justified by virtue of 
the condition ``$\| A\|_2 \leq \tilde R$'' included in \eqref{Def-ProbUB}). 
Note that for $n=21$ we have $\tilde R < \frac14 R\,$ (see Table~\mbox{D-2} in the previous section). 
In the remaining cases (where $10\leq n \leq 20$) we compute 
$\tilde U^{\pm}_1,\ldots ,\tilde U^{\pm}_M$ similarly, using  
data (such as $\alpha^{+}_1,\alpha^{-}_1,\ldots ,\alpha^{+}_M,\alpha^{-}_M,$ and $\tilde R$)  
that relates to the case in question. 

\subsection{Our results in the cases where $10\leq n\leq 21$} 

Table~\mbox{D-3} (below) lists some of the lower and upper bounds ($L^{+}_m$ and $U^{+}_m$) for $\varkappa^{+}_m (k_N')$
obtained in respect of the case $n=21\,$ ($N=2^{21}$).  
The table includes, as well, the corresponding probable upper bounds,  $\tilde U^{+}_m$. 
Recall (from Sections~\mbox{D.4} and~\mbox{D.6}, and Table~\mbox{D-2}) that in computing $U^{+}_m$ we use a number 
$R \approx 9.812127408198\times 10^{-3}$ as our upper bound for $\| A\|_2$, 
whereas in computing $\tilde U^{+}_m$, for $N=2^{21}$, a hypothetical upper bound     
$\tilde R \approx 2.240542901683 \times 10^{-3}$ is used instead. 

\begin{remarks} 
\item{\it 1)} \quad Table~\mbox{D-3} contains just a small illustrative sample of the complete set of 
results that were obtained. 
For $n=21$ we obtained corresponding results for the indices $m\leq M=384$ not included in the table, 
and (for each positive integer $m\leq 384$) got numerical bounds, $L^{-}_m$ and $U^{-}_m$, for $|\varkappa^{-}_m (k_N')|$, 
as well as a probable upper bound $\tilde U^{-}_m$ satisfying \eqref{Def-ProbUB}. 
For $10\leq n\leq 20$, $1\leq m\leq 384$ and $\pm\in\{ +,-\}$, we got both 
a numerical bound $L^{\pm}_m$, satisfying $L^{\pm}_m \leq |\varkappa^{\pm}_m (k_N')|$, and a (complementary) 
probable upper bound $\tilde U^{\pm}_m$; though we did not 
(in these cases)  compute any upper bound $U^{\pm}_m\geq |\varkappa^{\pm}_m (k_N')|$,  
nor even (as a precursor to that) any upper bound $R$ for $\| A\|_2$, analogous to the bound \eqref{SpecNormA21}: 
see Remarks~\mbox{D.3}~\mbox{(2)} for some explanation of these omissions.   
\item{\it 2)}\quad By examining the relevant data and results (i.e. the numbers $\alpha^{\pm}_1,\ldots ,\alpha^{\pm}_M$ 
and numerical bounds 
$L^{\pm}_1,\ldots ,L^{\pm}_M$ and $U^{\pm}_1,\ldots ,U^{\pm}_M$) 
we have been able to ascertain, via \eqref{n=21intervals}, that for $n=21$ we have
$| \alpha^{\pm}_m - \varkappa^{\pm}_m(k_N')| < 3.2\times 10^{-11}$  whenever $1\leq m\leq 49$ 
and $\pm\in\{ + , -\}$; we find also that the corresponding $98$ relative errors 
(of form $| \alpha^{\pm}_m - \varkappa^{\pm}_m(k_N')| / | \alpha^{\pm}_m |$) 
are all less than $3.2\times 10^{-9}$. In contrast, for $n=21$, $\pm\in\{ + , -\}$ and $50\leq m\leq M=384$, 
we have $U^{\pm}_m - |\alpha^{\pm}_m| > 10^{-4}$ and 
$(U^{\pm}_m - |\alpha^{\pm}_m|)/ |\alpha^{\pm}_m| > 1.1 \times 10^{-2}$, and so 
are unable to rule out the possibility that the agreement between $\alpha^{\pm}_m$ 
and $\varkappa^{\pm}_m (k_N')$ might be less than $3$ significant figures.   
This dichotomous behaviour (of bounds on error terms) results from our use of the bound 
\eqref{SpecNormA21} for $\| A\|_2$: note that 
with $R := 9.81212740819788 \times 10^{-3}$ we have  
$R / \min\{ \alpha^{+}_m, |\alpha^{-}_m|\} < 0.995$ for $1\leq m\leq 49$, 
but $R / \max\{ \alpha^{+}_m, |\alpha^{-}_m|\} > 1.01$ for $50\leq m\leq M\,$
(when $n=21$). 
\item{\it 3)}\quad In light of what has been noted (in Section~\mbox{D.5}) regarding probabilities,  
it is reasonable for us to conjecture that $\| A\|_2\leq \tilde R$
in all cases of interest (i.e. for $10\leq n\leq 21$),  
so that, by \eqref{LBs_for_eigenmoduli} and \eqref{Def-ProbUB}, one has both 
$[ L^{+}_m , \tilde U^{+}_m] \ni \varkappa^{+}_m (k_N')$ and $[ -\tilde U^{-}_m , -L^{-}_m] \ni \varkappa^{-}_m (k_N')$ 
for $1\leq m\leq M=384$ and $10\leq n\leq 21$. By examination of the relevant data, we 
find that if our conjecture (just stated) is correct then, 
for $1\leq m\leq M=384$, $\pm\in\{ + , -\}$ and $10\leq n\leq 21$, one has both 
\begin{equation*} 
\left| \alpha^{\pm}_m - \varkappa^{\pm}_m(k_N')\right|  < 
\begin{cases} \min\left\{ 2^{17/2} {\tt u} N^{1/2} \,,\, 3.2\times 10^{-11}\right\} & \text{if $m\leq 382$} , \\
3.2\times 10^{-10} &\text{if $n=21$} , \\
2.0\times 10^{-6} &\text{otherwise} , 
\end{cases}
\end{equation*} 
and 
\begin{equation*} 
\frac{\left| \alpha^{\pm}_m - \varkappa^{\pm}_m(k_N')\right|}{\left| \alpha^{\pm}_m\right|}  < 
\begin{cases} \min\left\{ 2^{17} {\tt u} N^{1/2} \,,\, 1.5\times 10^{-8}\right\} & \text{if $m\leq 382$} , \\
1.5\times 10^{-7} &\text{if $n=21$} , \\
1.5\times 10^{-3} &\text{otherwise} .  
\end{cases}
\end{equation*} 
\end{remarks} 

\begin{table}[h!] 
\centering 
\begin{minipage}{170mm}
\resizebox{170mm}{!}{ 
\begin{tabular}[c]{|r|c|c|c|c|} 
\hline 
$m$ & $L^{+}_m$ & $\alpha^{+}_m$ & $\tilde U^{+}_m$ & $U^{+}_m$ \\ \hline 
$1$ & $7.96889109765\times 10^{-2}$ & $7.96889109771\times 10^{-2}$ & $7.96889109799\times 10^{-2}$ & $7.96889109801\times 10^{-2}$ \\ \hline
$2$ & $6.90322926751\times 10^{-2}$ & $6.90322926762\times 10^{-2}$ & $6.90322926785\times 10^{-2}$ & $6.90322926786\times 10^{-2}$ \\ \hline
$3$ & $5.52247170569\times 10^{-2}$ & $5.52247170587\times 10^{-2}$ & $5.52247170610\times 10^{-2}$ & $5.52247170612\times 10^{-2}$ \\ \hline
$4$ & $4.58979700002\times 10^{-2}$ & $4.58979700026\times 10^{-2}$ & $4.58979700051\times 10^{-2}$ & $4.58979700053\times 10^{-2}$ \\ \hline
$5$ & $4.40302596509\times 10^{-2}$ & $4.40302596531\times 10^{-2}$ & $4.40302596552\times 10^{-2}$ & $4.40302596555\times 10^{-2}$ \\ \hline
$6$ & $3.77564254334\times 10^{-2}$ & $3.77564254358\times 10^{-2}$ & $3.77564254382\times 10^{-2}$ & $3.77564254385\times 10^{-2}$ \\ \hline
$7$ & $3.37965989125\times 10^{-2}$ & $3.37965989148\times 10^{-2}$ & $3.37965989172\times 10^{-2}$ & $3.37965989175\times 10^{-2}$ \\ \hline
$8$ & $2.87683707857\times 10^{-2}$ & $2.87683707882\times 10^{-2}$ & $2.87683707908\times 10^{-2}$ & $2.87683707913\times 10^{-2}$ \\ \hline 
$9$ & $2.82035514184\times 10^{-2}$ & $2.82035514206\times 10^{-2}$ & $2.82035514228\times 10^{-2}$ & $2.82035514232\times 10^{-2}$ \\ \hline
$10$ & $2.79386995477\times 10^{-2}$ & $2.79386995498\times 10^{-2}$ & $2.79386995520\times 10^{-2}$ & $2.79386995524\times 10^{-2}$ \\ \hline
$11$ & $2.76253252619\times 10^{-2}$ & $2.76253252640\times 10^{-2}$ & $2.76253252662\times 10^{-2}$ & $2.76253252667\times 10^{-2}$ \\ \hline 
$12$ & $2.71188486835\times 10^{-2}$ & $2.71188486858\times 10^{-2}$ & $2.71188486883\times 10^{-2}$ & $2.71188486888\times 10^{-2}$ \\ \hline
$13$ & $2.47046510682\times 10^{-2}$ & $2.47046510704\times 10^{-2}$ & $2.47046510728\times 10^{-2}$ & $2.47046510733\times 10^{-2}$ \\ \hline
$14$ & $2.16319730642\times 10^{-2}$ & $2.16319730666\times 10^{-2}$ & $2.16319730691\times 10^{-2}$ & $2.16319730698\times 10^{-2}$ \\ \hline
$15$ & $2.12977886026\times 10^{-2}$ & $2.12977886048\times 10^{-2}$ & $2.12977886070\times 10^{-2}$ & $2.12977886077\times 10^{-2}$ \\ \hline
$16$ & $2.08797902564\times 10^{-2}$ & $2.08797902586\times 10^{-2}$ & $2.08797902609\times 10^{-2}$ & $2.08797902616\times 10^{-2}$ \\ \hline
$24$ & $1.56334007251\times 10^{-2}$ & $1.56334007272\times 10^{-2}$ & $1.56334007294\times 10^{-2}$ & $1.56334007306\times 10^{-2}$ \\ \hline
$32$ & $1.30254636639\times 10^{-2}$ & $1.30254636663\times 10^{-2}$ & $1.30254636690\times 10^{-2}$ & $1.30254636713\times 10^{-2}$ \\ \hline
$48$ & $9.92361438470\times 10^{-3}$ & $9.92361438695\times 10^{-3}$ & $9.92361438951\times 10^{-3}$ & $9.92361440820\times 10^{-3}$ \\ \hline
$49$ & $9.89056976310\times 10^{-3}$ & $9.89056976514\times 10^{-3}$ & $9.89056976747\times 10^{-3}$ & $9.89056978811\times 10^{-3}$ \\ \hline
$50$ & $9.70526482684\times 10^{-3}$ & $9.70526482904\times 10^{-3}$ & $9.70526483155\times 10^{-3}$ & $9.81212740820\times 10^{-3}$ \\ \hline
$64$ & $8.45564990163\times 10^{-3}$ & $8.45564990382\times 10^{-3}$ & $8.45564990637\times 10^{-3}$ & $9.81212740820\times 10^{-3}$ \\ \hline
$96$ & $6.37299745043\times 10^{-3}$ & $6.37299745273\times 10^{-3}$ & $6.37299745558\times 10^{-3}$ & $9.81212740820\times 10^{-3}$ \\ \hline
$128$ & $5.07486490043\times 10^{-3}$ & $5.07486490270\times 10^{-3}$ & $5.07486490575\times 10^{-3}$ & $9.81212740820\times 10^{-3}$ \\ \hline
$192$ & $3.77567706203\times 10^{-3}$ & $3.77567706447\times 10^{-3}$ & $3.77567706830\times 10^{-3}$ & $9.81212740820\times 10^{-3}$ \\ \hline
$256$ & $3.08683723978\times 10^{-3}$ & $3.08683724196\times 10^{-3}$ & $3.08683724612\times 10^{-3}$ & $9.81212740820\times 10^{-3}$ \\ \hline
$384$ & $2.24065583428\times 10^{-3}$ & $2.24065583651\times 10^{-3}$ & $2.24065615076\times 10^{-3}$ & $9.81212740820\times 10^{-3}$ \\ 
\hline 
\end{tabular}
}  
\vskip 3mm 
\ Table~D-3 (results for $N=2^{21}$, rounded to 12 significant digits) 
\end{minipage} 
\end{table} 

\printbibliography 

\end{refsection} 

\begin{refsection}[refsE.bib]

\section{Generating `random' numbers} 

Any satisfactory implementation of Step~1 of the algorithm described in Section~\mbox{D.5}  
requires some method by which vectors  ${\bf x}_1,\ldots ,{\bf x}_S$ 
can be chosen `at random' from the set $X_N$ defined in \eqref{DefSetXn}.   
We tried out three different methods of doing this, before deciding which we preferred. 
\par  
The first method that we tried out utilises Octave's {\tt rand()} function, 
which implements the `Mersenne Twister' pseudorandom number 
generator \cite{MN1998}.
When called with arguments $N$ and $S$, {\tt rand()}  
returns an $N\times S$ matrix $(u_{i,j})$ of  
pseudorandom double precision numbers that (in certain respects) simulate    
random samples from the uniform distribution on the interval $(0, 1)$. 
One can then use Octave's {\tt round()} 
function to compute the  $N\times S$ matrix $B$ with 
elements $b_{i,j} = 2\langle u_{i,j}\rangle - 1$, 
where $\langle u\rangle = 1$ if $u = \frac12$ and is otherwise the integer 
($0$ or $1$) nearest to $u$. One then has 
$b_{i,j} = 1$, if $u_{i,j}\geq \frac12 $; $b_{i,j} = -1$, otherwise. 
The elements of $B$ simulate random samples from the set $\{ 1, -1\}$. 
We take the columns of the matrix $N^{-1/2} B$ to be our chosen vectors ${\bf x}_1,\ldots , {\bf x}_S$.  
\par 
One deficiency of the approach just sketched is that the {\tt rand()} function 
simply cycles through its permitted `states' as it is used (producing one pseudorandom number 
per state encountered).  On first use of the {\tt rand()} function (in an Octave session)  
its state is initialized using `random' information, such as the current time of day.  
Each successive state of {\tt rand()} is then determined by the one before it. 
Thus the outcome of the above selection process (i.e. our choice of the matrix $B$) is determined by 
the initial state  of the {\tt rand()} function, prior to our use of it.
Therefore the number of possible outcomes cannot exceed the number of permitted states 
of the {\tt rand()} function. The latter number (the Mersenne Twister's `period') is $2^{19937} - 1$;  
one should compare it with the total number of subsets 
$\{ {\bf x}_1, -{\bf x}_1, \ldots , {\bf x}_S, -{\bf x}_S\}\subset X_N$ 
satisfying ${\bf x}_s \neq \pm {\bf x}_{s'}$ for $1\leq s < s' \leq S$, 
which is the binomial coefficient $C(2^{N-1}, S)$, and 
so is greater than or equal to $(2^{N-1} / S)^S$:  if $N$ is much greater than $\log_2(S) + 19937/S$, then
the vast majority of such subsets will be incapable of 
being produced by the application of the {\tt rand()} function described above, and so we must  
admit to a degree of failure in that approach to sampling `at random' from~$X_N$. 
Note in particular that when $S=32$ one has $\log_2(S) + 19937/S < 629 < 2^{10}$. 
On the other hand, since the Mersenne Twister has  
(see \cite[Table~II]{MN1998}) $19937$-dimensional equidistribution up to $1$-bit accuracy, 
we might hope that the above approach does simulate random sampling effectively when $NS < 19937$. 
\par
The problem (just noted) with our first approach to the selection of ${\bf x}_1, \ldots , {\bf x}_S$ 
has motivated us to also try out, as an alternative, a slight elaboration of that method. 
In doing so we utilise the fact that the 
{\tt rand()} function's state can be `reset' to a new `random' state (the extra source 
of randomness here might be derived from CPU time, wall clock time and current fraction of a second, 
though this is something that the documentation of the {\tt rand()} function does not fully clarify: 
we have, at least, been able to verify by experiment that the state after a reset is not determined 
by the state prior to the reset). We compute the $NS$ elements of the matrix $(u_{i,j})$ in batches of size 
$\nu := \min\{ N, 2^{14}\} < 19937$, with 
one call to {\tt rand()}, and one `random' reset (of the state), per batch. The relationships between the matrix 
$(u_{i,j})$, the matrix $B$ and the vectors  ${\bf x}_1,\ldots ,{\bf x}_S$ remain the same 
as they were for the first method.
\par 
The third method for selecting ${\bf x}_1,\ldots , {\bf x}_S$ that we tried out 
utilises the `RDRAND' hardware random number generator 
(available on many of the processors that Intel{\textregistered} has developed since 2012). 
We got the idea to try such an approach from the paper \cite{Ro2017}, which describes an application of 
RDRAND-based Monte Carlo simulation in the field of astrophysics.
\par 
A package (`randomgen') developed for the Python programming language  
gives one access to output from RDRAND. We wrote a short 
Python script to exploit this facility. 
Note that RDRAND, like the Mersenne Twister, is a pseudorandom number generator: its 
output is not truly random. Nevertheless, RDRAND is `reseeded' frequently with high-quality 
random numbers (`seeds') deriving from a non-determistic hardware source 
(see \cite[Sections~3.1 and~3.2]{In2018} for details). 
Our Python script produces a sequence of $64$-bit integers (i.e. integers in the interval $[0,2^{64}-1]$), 
each the result of one call to RDRAND. Between the production of one member of this sequence and the next 
an additional $1023$ calls are made to RDRAND, though the resulting ouputs (which are, again, $64$-bit integers) 
are completely disregarded: this many consecutive calls to RDRAND is sufficient to trigger a reseeding event 
(though, less happily, it also greatly increases the run-time of our script). 
We also introduced some artificial delays, each several microseconds in duration  (the aim, again, was to ensure reseeding). 
We believe that this forced reseeding (with seeds from a non-deterministic source)  
should make the output of our script more nearly `truly random' than it would otherwise be: 
this appears to be what is implied in \cite[Section~5.2.5]{In2018}. 
\par
Exactly $2^{21}$ $64$-bit integers were generated in one run of our script (taking just over $2$ hours). 
Their base-$2$ expansions, $2^{63} \beta_{i,1} + 2^{62} \beta_{i,2} + \ldots + 2^0 \beta_{i,64}\,$ ($1\leq i\leq 2^{21}$), were then computed. 
We had then, in $(\beta_{i,j})$, a $2^{21}\times 64$ matrix of random (or, at least, pseudorandom) 
samples from the set $\{0, 1\}$. 
For $N \in\{ 2^{10}, 2^{11}, ..., 2^{20}\}$  
we took,  
as our choice of vectors ${\bf x}_1,\ldots , {\bf x}_S\in X_N$, 
the columns of the matrix $(2^{-N/2}(2\beta_{i,j} - 1))_{N < i\leq 2N, j\leq 32}$; 
for $N = 2^{21}$ we took the remaining $32$ unused columns of the matrix  
$(2^{-N/2}(2\beta_{i,j} - 1))_{i\leq N, j\leq 64}\,$ (i.e. columns $33$ to $64$).  
Thus in every case we had $S=32$. 
\par 
The probable upper bounds $\tilde R(2^{10}),\ldots , \tilde R(2^{21})$ displayed in Table~\mbox{D-2} (in Section~\mbox{D.5})   
were obtained using the third of the above described methods for selecting 
${\bf x}_1,\ldots , {\bf x}_S$ `at random' from $X_N\,$ 
(i.e. the RDRAND-based method).  
Trials of the pair of {\tt rand()}-based methods 
yielded similar results:  in every instance, the alternative probable upper bound obtained for $\|A(N)\|_2$ was found to differ 
from the corresponding RDRAND-based probable upper bound by some factor $\phi$ satisfying $|\phi -1|<\frac{1}{10} J^{-1}\log Z$, 
where $J=J(N)$ and $Z=Z(N)$ are as in Table~\mbox{D-2}. 

\printbibliography 

\end{refsection}

\end{appendices} 

\end{document}